\documentclass[12pt,reqno]{amsart}
\usepackage{feynmp}
\usepackage{latexsym}
\usepackage{amsfonts}
\usepackage{amssymb}
\usepackage{amscd}
\usepackage{amsbsy}
\usepackage{amsmath}
\usepackage[dvips]{graphicx}
\usepackage{epsfig}    
\usepackage[all]{xy}
\xyoption{knot}
\xyoption{graph}

\def\eg{{\it e.g.\ }}
\def\cf{{\it cf.\ }}
\def\rhs{{\it r.h.s.\ }}
\def\lhs{{\it l.h.s.\ }}

\def\End{\mathop{{\rm End}}\nolimits}

\def\Hom{\mathop{{\rm Hom}}\nolimits}
\def\Ext{\mathop{{\rm Ext}}\nolimits}

\def\im{\mathop{{\rm im}}\nolimits}
\def\Tr{\mathop{{\rm Tr}}\nolimits}

\def\deg{ \mathop{{\rm deg}}\nolimits }
\def\p{^{\prime}}
\def\pp{^{\prime\prime}}

\def\rank{\mathop{{\rm rank}}\nolimits}

\def\del{ \partial }

\def\promod{\mathop{{\rm mod}}\nolimits}
\def\mod#1{\;(\promod #1)}

\def\End{\mathop{{\rm End}}\nolimits}

\def\pr#1#2{ \noindent{\em Proof of #1~\ref{#2}.} }
\def\proof{ \noindent{\em Proof.} }

\def\qed{ \hfill $\Box$ }

\def\lrbc#1{ \left( #1 \right) }
\def\lrbs#1{ \left[ #1 \right] }

\def\lrbl#1{ \left\langle #1 \right\rangle }

\def\lrbigc#1{ \big( #1 \big) }

\def\xG{ \mfG }

\parskip=\medskipamount                 
\thicklines                     
\newcommand{\bimn}[7]{\bibitem{#1}#2,
{\em #3}, { #4}\hspace{0.25em}{\bf
#5}\hspace{0.25em}(#6)\hspace{0.25em}{#7}.}

\newcommand{\bimna}[8]{\bibitem{#1}#2,
{\em #3}, { #4}\hspace{0.25em}{\bf
#5}\hspace{0.25em}(#6)\hspace{0.25em}{#7}\hspace{0.25em}(#8).}

\def\inbar{\vrule height1.5ex width.4pt depth0pt}
\def\IC{\relax\,\hbox{$\inbar\kern-.3em{\rm C}$}}
\def\IN{\relax{\rm I\kern-.18em N}}
\def\IQ{\relax\,\hbox{$\inbar\kern-.3em{\rm Q}$}}
\def\IR{\relax{\rm I\kern-.18em R}}
\def\ZZ{\relax{\sf Z\kern-.4em Z}}

  \def\cL{{\cal L}}

\newtheorem{theorem}{Theorem}[section]
\newtheorem{proposition}[theorem]{Proposition}
\newtheorem{corollary}[theorem]{Corollary}
\newtheorem{conjecture}[theorem]{Conjecture}
\newtheorem{lemma}[theorem]{Lemma}

\catcode`\@=11

\newif\if@fewtab\@fewtabtrue

\catcode`\@=11

\newif\if@fewtab\@fewtabtrue

{\count255=\time\divide\count255 by 60
\xdef\hourmin{\number\count255} \multiply\count255
by-60\advance\count255 by\time
\xdef\hourmin{\hourmin:\ifnum\count255<10 0\fi\the\count255}}
\def\ps@draft{\let\@mkboth\@gobbletwo
    \def\@oddhead{}
    \def\@oddfoot
      {\hbox to 7 cm{\footnotesize {\em Draft of \jobname:} \draftdate
       \hfil}\hskip -7cm\hfil\rm\thepage \hfil}
    \def\@evenhead{}\let\@evenfoot\@oddfoot}

\def\ceqno{\global\@fewtabfalse
    \ifcase\@eqcnt \def\@tempa{& & &}\or \def\@tempa{& &}
      \or \def\@tempa{&}
      \or\def\@tempa{}\fi\@tempa
{\rm(\theequation)}}

\def\aeqno#1{\global\@fewtabfalse
    \ifcase\@eqcnt \def\@tempa{& & &}\or \def\@tempa{& &}
      \or \def\@tempa{&}
      \or\def\@tempa{}\fi\@tempa
{\rm(\theequation,#1)}}

\def\label#1{\ifnum\draftcontrol=1
 \global\def\draftnote{$\scriptstyle #1$}\fi
 \@bsphack\if@filesw {\let\thepage\relax
   \def\protect{\noexpand\noexpand\noexpand}%
\xdef\@gtempa{\write\@auxout{\string
      \newlabel{#1}{{\@currentlabel}{\thepage}}}}}\@gtempa
   \if@nobreak \ifvmode\nobreak\fi\fi\fi
  \@esphack}

\def\alabel#1#2{\label{#1}\global\@fewtabfalse
    \ifcase\@eqcnt \def\@tempa{& & &}\or \def\@tempa{& &}
      \or \def\@tempa{&}
      \or\def\@tempa{}\fi\@tempa
{\hbox to 3cm{\phantom{\rm(\theequation,#2)} \draftnote
\hfil}\hskip -3cm {\rm(\theequation,#2)}}}

\def\clabel#1{\label{#1}\global\@fewtabfalse
    \ifcase\@eqcnt \def\@tempa{& & &}\or \def\@tempa{& &}
      \or \def\@tempa{&}
      \or\def\@tempa{}\fi\@tempa
{\hbox to 3cm{\phantom{\rm(\theequation)} \draftnote \hfil}\hskip
-3cm{\rm(\theequation)}}}

\def\eqnarray{\def\draftnote{{}}\global\@fewtabtrue
\stepcounter{equation}\let\@currentlabel=\theequation
\global\@eqnswtrue
\global\@eqcnt\z@\tabskip\@centering\let\\=\@eqncr
$$\halign to \displaywidth\bgroup\@eqnsel\hskip\@centering\@eqcnt\z@
  $\displaystyle\tabskip\z@{##}$&\global\@eqcnt\@ne
  \hskip 1\arraycolsep \hfil$\displaystyle{##}$\hfil
  &\global\@eqcnt\tw@ \hskip 1\arraycolsep
$\displaystyle\tabskip\z@{##}$ \hfil
\tabskip\@centering&\global\@eqcnt\thr@@\llap{##}\tabskip\z@ \cr}

\def\endeqnarray{\@@eqncr\egroup
      \global\advance\c@equation\m@ne$$\global\@ignoretrue}

\def\@eqnnum{\hbox to 3cm{\phantom{\rm(\theequation)} \draftnote
                         \hfil}\hskip -3cm {\rm(\theequation)}}

\def\@@eqncr{\let\@tempa\relax
    \ifcase\@eqcnt \def\@tempa{& & &}\or \def\@tempa{& &}
      \or \def\@tempa{&}
      \or\def\@tempa{}
\fi\@tempa \if@eqnsw \if@fewtab\@eqnnum\fi
\stepcounter{equation}\fi\global
\@eqnswtrue\global\@eqcnt\z@\global\@fewtabtrue\cr}

\def\draftcite#1{\ifnum\draftcontrol=1#1\else{}\fi}

\def\@lbibitem[#1]#2{\item{}\hskip -3cm \hbox to 2cm
{\hfil$\scriptstyle\draftcite{#2}$}\hskip
1cm[\@biblabel{#1}]\if@filesw
     {\def\protect##1{\string ##1\space}\immediate
      \write\@auxout{\string\bibcite{#2}{#1}}}\fi\ignorespaces}

\def\@bibitem#1{\item\hskip -3cm \hbox to 2cm
{\hfil $\scriptstyle\draftcite{#1}$}\hskip 1cm \if@filesw
\immediate\write\@auxout
       {\string\bibcite{#1}{\the\value{\@listctr}}}\fi\ignorespaces}

\def\draftdate{\number\month/\number\day/\number\year\ \ \ \hourmin }
 \global\def\draftcontrol{0}
\catcode`\@=12

\def\theequation{{\thesection.\arabic{equation}}}

\def\qq{\begin{eqnarray}}
\def\qqq{\end{eqnarray}}
\def\ee{\begin{eqnarray}}
\def\eee{\end{eqnarray}}
\def\rx#1{~(\ref{#1})}
\def\rxw#1{(\ref{#1})}
\def\ex#1{eq.\hspace*{-3pt}\rx{#1}}
\def\eex#1{eqs.\hspace*{-3pt}\rx{#1}}
\def\cx#1{~\cite{#1}}
\def\rw#1{~\ref{#1}}

\def\fg#1{Fig.~\ref{#1}}

\newlength{\shiftwidth}
\def\shift#1{&&\hbox to \shiftwidth{\hfill $\displaystyle#1$}}
\newlength{\sshiftwidth}
\def\sshift#1{\lefteqn{\hbox to
\sshiftwidth{\hfill$\displaystyle#1$}}}

\def\qbezier{\bezier{120}}
\def\DottedCircle{
\bezier{4}(0.966,-0.259)(1.04,0)(0.966,0.259)
\bezier{4}(0.966,0.259)(0.897,0.518)(0.707,0.707)
\bezier{4}(0.707,0.707)(0.518,0.897)(0.259,0.966)
\bezier{4}(0.259,0.966)(0,1.04)(-0.259,0.966)
\bezier{4}(-0.259,0.966)(-0.518,0.897)(-0.707,0.707)
\bezier{4}(-0.707,0.707)(-0.897,0.518)(-0.966,0.259)
\bezier{4}(-0.966,0.259)(-1.04,0)(-0.966,-0.259)
\bezier{4}(-0.966,-0.259)(-0.897,-0.518)(-0.707,-0.707)
\bezier{4}(-0.707,-0.707)(-0.518,-0.897)(-0.259,-0.966)
\bezier{4}(-0.259,-0.966)(0,-1.04)(0.259,-0.966)
\bezier{4}(0.259,-0.966)(0.518,-0.897)(0.707,-0.707)
\bezier{4}(0.707,-0.707)(0.897,-0.518)(0.966,-0.259) }
\def\Endpoint[#1]{
\ifcase#1 \put(1,0){\circle*{0.15}}
\or\put(0.866,0.5){\circle*{0.15}}
\or\put(0.5,0.866){\circle*{0.15}} \or\put(0,1){\circle*{0.15}}
\or\put(-0.5,0.866){\circle*{0.15}}
\or\put(-0.866,0.5){\circle*{0.15}} \or\put(-1,0){\circle*{0.15}}
\or\put(-0.866,-0.5){\circle*{0.15}}
\or\put(-0.5,-0.866){\circle*{0.15}} \or\put(0,-1){\circle*{0.15}}
\or\put(0.5,-0.866){\circle*{0.15}}
\or\put(0.866,-0.5){\circle*{0.15}} \fi}
\def\Arc[#1]{
\thicklines         
\ifcase#1 \bezier{25}(0.966,-0.259)(1.04,0)(0.966,0.259) \or
\bezier{25}(0.966,0.259)(0.897,0.518)(0.707,0.707) \or
\bezier{25}(0.707,0.707)(0.518,0.897)(0.259,0.966) \or
\bezier{25}(0.259,0.966)(0,1.04)(-0.259,0.966) \or
\bezier{25}(-0.259,0.966)(-0.518,0.897)(-0.707,0.707) \or
\bezier{25}(-0.707,0.707)(-0.897,0.518)(-0.966,0.259) \or
\bezier{25}(-0.966,0.259)(-1.04,0)(-0.966,-0.259) \or
\bezier{25}(-0.966,-0.259)(-0.897,-0.518)(-0.707,-0.707) \or
\bezier{25}(-0.707,-0.707)(-0.518,-0.897)(-0.259,-0.966) \or
\bezier{25}(-0.259,-0.966)(0,-1.04)(0.259,-0.966) \or
\bezier{25}(0.259,-0.966)(0.518,-0.897)(0.707,-0.707) \or
\bezier{25}(0.707,-0.707)(0.897,-0.518)(0.966,-0.259) \fi}
\def\DottedArc[#1]{
\ifcase#1 \bezier{4}(0.966,-0.259)(1.04,0)(0.966,0.259) \or
\bezier{4}(0.966,0.259)(0.897,0.518)(0.707,0.707) \or
\bezier{4}(0.707,0.707)(0.518,0.897)(0.259,0.966) \or
\bezier{4}(0.259,0.966)(0,1.04)(-0.259,0.966) \or
\bezier{4}(-0.259,0.966)(-0.518,0.897)(-0.707,0.707) \or
\bezier{4}(-0.707,0.707)(-0.897,0.518)(-0.966,0.259) \or
\bezier{4}(-0.966,0.259)(-1.04,0)(-0.966,-0.259) \or
\bezier{4}(-0.966,-0.259)(-0.897,-0.518)(-0.707,-0.707) \or
\bezier{4}(-0.707,-0.707)(-0.518,-0.897)(-0.259,-0.966) \or
\bezier{4}(-0.259,-0.966)(0,-1.04)(0.259,-0.966) \or
\bezier{4}(0.259,-0.966)(0.518,-0.897)(0.707,-0.707) \or
\bezier{4}(0.707,-0.707)(0.897,-0.518)(0.966,-0.259) \fi}
\def\Chord[#1,#2]{
\thinlines \ifnum#1>#2\Chord[#2,#1] \else\ifnum#1<#2 \ifcase#1
\ifcase#2 \or\qbezier(1,0)(0.516,0.138)(0.866,0.5)
\or\qbezier(1,0)(0.45,0.26)(0.5,0.866)
\or\qbezier(1,0)(0.327,0.327)(0,1)
\or\qbezier(1,0)(0.179,0.311)(-0.5,0.866)
\or\qbezier(1,0)(0.0536,0.2)(-0.866,0.5) \or\put(1, 0){\line(-2,
0){2}} \or\qbezier(1,0)(0.0536,-0.2)(-0.866,-0.5)
\or\qbezier(1,0)(0.179,-0.311)(-0.5,-0.866)
\or\qbezier(1,0)(0.327,-0.327)(0,-1)
\or\qbezier(1,0)(0.45,-0.26)(0.5,-0.866)
\or\qbezier(1,0)(0.516,-0.138)(0.866,-0.5) \fi \or\ifcase#2\or
\or\qbezier(0.866,0.5)(0.378,0.378)(0.5,0.866)
\or\qbezier(0.866,0.5)(0.26,0.45)(0,1)
\or\qbezier(0.866,0.5)(0.12,0.446)(-0.5,0.866)
\or\qbezier(0.866,0.5)(0,0.359)(-0.866,0.5)
\or\qbezier(0.866,0.5)(-0.0536,0.2)(-1,0) \or\put(0.866,
0.5){\line(-5, -3){1.73}}
\or\qbezier(0.866,0.5)(0.146,-0.146)(-0.5,-0.866)
\or\qbezier(0.866,0.5)(0.311,-0.179)(0,-1)
\or\qbezier(0.866,0.5)(0.446,-0.12)(0.5,-0.866)
\or\qbezier(0.866,0.5)(0.52,0)(0.866,-0.5) \fi \or\ifcase#2\or\or
\or\qbezier(0.5,0.866)(0.138,0.516)(0,1)
\or\qbezier(0.5,0.866)(0,0.52)(-0.5,0.866)
\or\qbezier(0.5,0.866)(-0.12,0.446)(-0.866,0.5)
\or\qbezier(0.5,0.866)(-0.179,0.311)(-1,0)
\or\qbezier(0.5,0.866)(-0.146,0.146)(-0.866,-0.5) \or\put(0.5,
0.866){\line(-3, -5){1}} \or\qbezier(0.5,0.866)(0.2,-0.0536)(0,-1)
\or\qbezier(0.5,0.866)(0.359,0)(0.5,-0.866)
\or\qbezier(0.5,0.866)(0.446,0.12)(0.866,-0.5) \fi
\or\ifcase#2\or\or\or \or\qbezier(0,1.)(-0.138,0.516)(-0.5,0.866)
\or\qbezier(0,1.)(-0.26,0.45)(-0.866,0.5)
\or\qbezier(0,1.)(-0.327,0.327)(-1,0)
\or\qbezier(0,1.)(-0.311,0.179)(-0.866,-0.5)
\or\qbezier(0,1.)(-0.2,0.0536)(-0.5,-0.866) \or\put(0, 1){\line(0,
-2){2}} \or\qbezier(0,1.)(0.2,0.0536)(0.5,-0.866)
\or\qbezier(0,1.)(0.311,0.179)(0.866,-0.5) \fi
\or\ifcase#2\or\or\or\or
\or\qbezier(-0.5,0.866)(-0.378,0.378)(-0.866,0.5)
\or\qbezier(-0.5,0.866)(-0.45,0.26)(-1,0)
\or\qbezier(-0.5,0.866)(-0.446,0.12)(-0.866,-0.5)
\or\qbezier(-0.5,0.866)(-0.359,0)(-0.5,-0.866)
\or\qbezier(-0.5,0.866)(-0.2,-0.0536)(0,-1) \or\put(-0.5,
0.866){\line(3, -5){1}}
\or\qbezier(-0.5,0.866)(0.146,0.146)(0.866,-0.5) \fi
\or\ifcase#2\or\or\or\or\or
\or\qbezier(-0.866,0.5)(-0.516,0.138)(-1,0)
\or\qbezier(-0.866,0.5)(-0.52,0)(-0.866,-0.5)
\or\qbezier(-0.866,0.5)(-0.446,-0.12)(-0.5,-0.866)
\or\qbezier(-0.866,0.5)(-0.311,-0.179)(0,-1)
\or\qbezier(-0.866,0.5)(-0.146,-0.146)(0.5,-0.866) \or\put(-0.866,
0.5){\line(5, -3){1.73}} \fi \or\ifcase#2\or\or\or\or\or\or
\or\qbezier(-1,0)(-0.516,-0.138)(-0.866,-0.5)
\or\qbezier(-1,0)(-0.45,-0.26)(-0.5,-0.866)
\or\qbezier(-1,0)(-0.327,-0.327)(0,-1)
\or\qbezier(-1,0)(-0.179,-0.311)(0.5,-0.866)
\or\qbezier(-1,0)(-0.0536,-0.2)(0.866,-0.5) \fi
\or\ifcase#2\or\or\or\or\or\or\or
\or\qbezier(-0.866,-0.5)(-0.378,-0.378)(-0.5,-0.866)
\or\qbezier(-0.866,-0.5)(-0.26,-0.45)(0,-1)
\or\qbezier(-0.866,-0.5)(-0.12,-0.446)(0.5,-0.866)
\or\qbezier(-0.866,-0.5)(0,-0.359)(0.866,-0.5) \fi
\or\ifcase#2\or\or\or\or\or\or\or\or
\or\qbezier(-0.5,-0.866)(-0.138,-0.516)(0,-1)
\or\qbezier(-0.5,-0.866)(0,-0.52)(0.5,-0.866)
\or\qbezier(-0.5,-0.866)(0.12,-0.446)(0.866,-0.5) \fi
\or\ifcase#2\or\or\or\or\or\or\or\or\or
\or\qbezier(0,-1.)(0.138,-0.516)(0.5,-0.866)
\or\qbezier(0,-1.)(0.26,-0.45)(0.866,-0.5) \fi
\or\ifcase#2\or\or\or\or\or\or\or\or\or\or
\or\qbezier(0.5,-0.866)(0.378,-0.378)(0.866,-0.5) \fi\fi\fi\fi}
\def\FullChord[#1,#2]{
\Endpoint[#1] \Endpoint[#2] \Arc[#1] \Arc[#2] \Chord[#1,#2] }
\def\EndChord[#1,#2]{
\Endpoint[#1] \Endpoint[#2] \Chord[#1,#2] }
\def\Picture#1{
\begin{picture}(2,1)(-1,-0.167)
#1
\end{picture}
}
\def\DottedChordDiagram[#1,#2]{
\Picture{\DottedCircle \FullChord[#1,#2]} }
\def\ZZ{ \mathbb{Z} }
\def\IQ{ \mathbb{Q} }
\def\IC{ \mathbb{C} }
\def\IR{ \mathbb{R} }

\def\bfx{ \mathbf{x} }
\def\bfy{ \mathbf{y} }

\def\cL{ \mathcal{L} }

\def\hlf{ {1\over 2} }

\def\xJ{ J }

\def\lg{local graph}
\def\lgs{\lg s}

\def\mf{matrix factorization}
\def\mfs{\mf s}

\def\xW{ W }

\def\xD{ D }
\def\xG{ G }

\def\xg{ g }

\def\xN{ N }

\def\xJ{ J }

\def\xi{ i }

\def\xWv#1{ \xW_{#1} }

\def\xDv#1{ \xD_{#1} }

\def\Zt{ \ZZ_2 }

\def\xDv#1{ \xD_{#1} }

\def\St{ S^2 }

\def\xDvv#1#2{ \xDv{#1,#2} }

\def\xgv#1{ \xg(#1) }

\def\xgm{ \gamma }

\def\Cnv#1{ \mathrm{C}#1 }

\def\Kmf#1#2{ (#1;#2) }
\def\yyn{ n }

\def\yym{ m }

\def\otR{ \otv{\xR} }
\def\Homv#1{ \Hom_{#1} }
\def\HomR{ \Homv{\xR} }

\def\yfv#1{ f_{#1} }

\def\Endv#1{ \End_{#1} }
\def\EndR{ \Endv{\xR} }

\def\homfly{HOMFLY-PT}
\def\homflyp{\homfly\ polynomial}

\def\tvar{2-variable}

\def\qgrd{$q$-grading}
\def\qgrdd{$q$-graded}
\def\tgrd{$t$-grading}
\def\qdgr{$q$-degree}

\def\clgr{closed graph}

\def\kngr{knotted graph}

\def\grtng{graph-tangle}
\def\grtngs{\grtng s}
\def\opgr{open graph}
\def\opgrs{\opgr s}

\def\mf{matrix factorization}
\def\mfs{\mf s}

\def\Kmf{Koszul matrix factorization}
\def\Kmfs{\Kmf s}

\def\Zt{ \ZZ_2 }

\def\Zgrdd{$\ZZ$-graded}

\def\Ztgrdd{$\Zt$-graded}

\def\xN{ N }

\def\Igr{$\zI$-graph}

\def\xL{ L }

\def\crs{\mathop{{\rm cr}}\nolimits}
\def\crsv#1{ \crs(#1) }
\def\crsv#1{ n }

\def\ncL{ \ncv{\xL} }

\def\dgq{ \dgv{q} }

\def\xr{ r }

\def\xrv#1{ \xr^{(#1)} }

\def\zx{ x }
\def\zy{ y }
\def\zxo{ \zx_1 }
\def\zxt{ \zx_2 }
\def\zyo{ \zy_1 }
\def\zyt{ \zy_2 }

\def\eltr#1#2#3{ [#3]_{#1#2} }

\def\tK{ \mathrm{K} }
\def\kcv#1{ \tK(#1) }

\def\kcyp{ \kcv{\yp} }

\def\kcbfp{ \kcv{\bfyp} }

\def\kmfv#1#2{ \tK(#1;#2) }
\def\kmfpq{ \kmfv{\yp}{\yq} }
\def\kmfbpq{ \kmfv{\bfyp}{\bfyq} }

\def\kmftwv#1#2#3#4{
\tK\begin{pmatrix}
#1, & #3
\\
#2, & #4
\end{pmatrix}
}

\def\MF{ \mathbf{MF} }
\def\MFv#1#2{ \MF(#1,#2) }

\def\ardpl{ \ar@{}[d]|{\oplus} }

\def\lxarr{\ar@/^2.5pc/@{-->}[rr]}
\def\lxadd{\ar@/^2.5pc/@{-->}[dd]}
\def\xarru{ \xdar@/^20pt/[rr] }
\def\xaddr{ \xdar@/^20pt/[dd] }

\def\xg{ g }

\def\yp{ p }
\def\bfyp{ \mathbf{\yp} }

\def\yq{ q }
\def\bfyq{ \mathbf{\yq} }

\def\yr{ r }

\def\xLo{ \xL_1 }
\def\xLt{ \xL_2 }

\def\qmqi{ q - \qi }

\def\qi{ q^{-1} }

\def\xLo{ \xL_1 }
\def\xLt{ \xL_2 }
\def\xLot{ \xLo \sqcup  \xLt }

\def\yG{ \hGm }
\def\yGp{ \hGm\p }
\def\yGv#1{ \yG_{#1} }

\def\id{ \mathrm{id} }

\def\xW{ W }
\def\xWv#1{ \xW_{#1} }

\def\yW{ W }

\def\xD{ D }

\def\arr#1{ \ar@<0.3ex>[r]^-{#1} }
\def\arl#1{ \ar@<0.3ex>[l]^-{#1} }

\def\bfxv#1{ \bfx_{#1} }

\def\Kf{Kauffman}
\def\Kfp{\Kf\ polynomial}
\def\Kfps{\Kfp s}
\def\Hpt{HOMFLY-PT}
\def\homflyp{\Hpt\ polynomial}

\def\ntrl{natural}
\def\Frb{Frobenius}
\def\Frba{\Frb\ algebra}
\def\mltp{multiplication}
\def\mltps{\mltp s}
\def\cmlt{co-multiplication}
\def\cmlts{\cmlt s}
\def\cnst{constituent}

\def\mf{matrix factorization}
\def\mfs{\mf s}
\def\Kmf{Koszul \mf}
\def\Kmfs{\Kmf s}
\def\Kcx{Koszul complex}
\def\Kcxs{\Kcx es}

\def\xxch{chain}
\def\xxcl{complex}
\def\ccomp{categorification complex}
\def\chcl{\xxch\ \xxcl}
\def\hteq{homotopy equivalence}
\def\hteqt{homotopy equivalent}
\def\htpc{homotopic}
\def\prmr{primary}
\def\scnd{secondary}

\def\scndh{\scnd\ homomorphism}
\def\scndhs{\scndh s}

\def\adj{adjacent}
\def\outr{outer}
\def\innr{inner}

\def\enfn{endo-functor}
\def\enfns{\enfn s}

\def\trns{translation}
\def\trnsf{\trns\ functor}

\def\mdl{vector space}
\def\mdls{\mdl s}
\def\grdmdl{\grdd\ \mdl}
\def\grdmdls{\grdd\ \mdls}
\def\qgrdmdl{\qgrdd\ \mdl}

\def\frlk{framed link}

\def\mphq{ \simeq }
\def\cchq{ \simeq }

\def\fvtx{4-vertex}
\def\fvtcs{4-vertices}

\def\Jr{Jacobi algebra}

\def\flggd{4-legged}

\def\vadj{virtually adjoint}

\def\rankv{ \dim }

\def\Igr{I-graph}

\def\ZtZ{$\ZZ_2\times\ZZ$}

\def\sotn{$\mathrm{SO}(2\xN+2)$}
\def\sotN{$\mathrm{SO}(2\xN)$}

\def\sunm{\mathrm{SU}(N)}
\def\sun{$\sunm$}

\def\msut{ \mathrm{SU}(2) }
\def\msutt{ \msut\times\msut }
\def\sutt{$\msutt$}
\def\sut{$\msut$}
\def\msof{ \mathrm{SO}(4) }
\def\sof{ $\msof$ }

\def\apmgv#1#2{ #1_{#2} }
\def\apmu#1{ \apmgv{#1}{\msut} }
\def\apmo#1{ \apmgv{#1}{\msof} }
\def\apmud#1{ \apmu{#1}^{\otimes 2} }
\def\apmusd#1{ \apmu{#1}^{\;\otimes 2} }

\def\apmsd#1{ #1^{\;\otimes 2} }

\def\ggu{ \apmu{\gg} }
\def\ggud{ \apmud{\gg} }

\def\yFsud{ \apmud{\yF} }
\def\yFso{ \apmo{\yF} }

\def\yFd{ \apmd{\yF} }
\def\yGd{ \apmd{\yG} }

\def\yGsud{ \apmud{\yG} }

\def\hlncrud{ \apmusd{\hlncr} }
\def\hlncro{ \apmo{\hlncr} }

\def\zarcxyooud{ \apmusd{\zarcxyoo} }
\def\zarcxyooo{ \apmo{\zarcxyoo} }

\def\hlparud{ \apmusd{\hlpar} }

\def\hlhorud{ \apmusd{\hlhor} }
\def\hlparo{ \apmo{\hlpar} }

\def\hlvirud{ \apmusd{\hlvir} }
\def\hlviro{ \apmo{\hlvir} }
\def\hlverudd{ \hlver_{\msutt} }

\def\hlncrsd{ \apmsd{\hlncr} }
\def\hlparsd{ \apmsd{\hlpar} }
\def\hlvirsd{ \apmsd{\hlvir} }
\def\hlhorsd{ \apmsd{\hlhor} }

\def\dgrsh{degree shift}

\def\tvar{2-variable}

\def\grdd{graded}

\def\dgrm{diagram}
\def\plgrdgs{\plgrdg s}

\def\ldgrm{link \dgrm}

\def\ntydeg{90^\circ}

\def\clgr{closed graph}
\def\clgrs{\clgr s}

\def\kngr{knotted graph}
\def\plgrdg{planar \gdgrm}
\def\gdgrm{graph \dgrm}

\def\eltr{elementary}

\def\opgr{open graph}
\def\opgrs{\opgr s}

\def\elopgr{\eltr\ \opgr}
\def\elopgrs{\elopgr s}

\def\tngl{tangle}
\def\tngls{\tngl s}

\def\eltngl{\eltr\ \tngl}
\def\eltngls{\eltngl s}

\def\grtngl{graph-\tngl}
\def\grtngls{\grtngl s}

\def\lg{leg}
\def\lgs{\lg s}

\def\mf{matrix factorization}
\def\mfs{\mf s}
\def\cntrb{contractible}

\def\tarc{2-arc}
\def\oarc{1-arc}

\def\crs{crossing}
\def\crss{\crs s}

\def\smvrt{semi-virtual}

\def\sdl{saddle}
\def\sdlmrp{\sdl\ morphism}
\def\sdlmrps{\sdlmrp s}
\def\prpr{proper}
\def\qprpr{`\prpr'}
\def\cmn{common}
\def\smcl{semi-closed}
\def\clsd{closed}

\def\knk{kink}

\def\cnvl{convolution}
\def\cnvls{\cnvl s}

\def\Ps{Postnikov system}
\def\Pss{\Ps s}

\def\sbs{subsystem}
\def\sbss{\sbs s}

\def\crnr{corner}
\def\csbs{\crnr\ \sbs}
\def\csbss{\csbs s}

\def\fcs{factorsystem}
\def\fcss{\fcs s}
\def\cfcs{corner \fcs}

\def\xct{cut}

\def\xctv#1{$#1$-\xct}
\def\Sctv{\xctv{\xSset}}
\def\xcnv{convex}
\def\xlcnv{\lngth-\xcnv}
\def\xcnvsst{\xcnv\ subset}
\def\xcnvssts{\xcnvsst s}
\def\xxct{\xcnv\ \xct}
\def\xxcts{\xxct s}

\def\lngth{length}
\def\lngthd{\lngth\ degree}
\def\tlngth{total length}
\def\nng{non-negative}

\def\lbndd{bounded}

\def\crnr{corner}
\def\smplfc{simplification}
\def\smplfd{simplified}

\def\excis{excision}
\def\cntrcn{contractible cone}
\def\cntrcns{\cntrcn s}
\def\cntrcnex{\cntrcn\ \excis}
\def\excse{excise}

\def\Ztgrdd{$\Zt$-graded}
\def\Ztgrdng{$\Zt$-grading}
\def\Ztdgr{$\Zt$-degree}

\def\Zdgr{$\ZZ$-degree}
\def\Zdgrs{\Zdgr s}
\def\Zgrdd{$\ZZ$-graded}
\def\Zgrdng{$\ZZ$-grading}

\def\qdgr{$q$-degree}

\def\MF{ \mathrm{MF} }
\def\PS{ \mathrm{PMF} }

\def\xL{ L }
\def\xLp{ \xL\p }
\def\Sth{ S^3 }

\def\xq{ q }
\def\xN{ N }

\def\catC{ \mathcal{C} }

\def\qgrd{$q$-grading}
\def\qgrdd{$q$-graded}
\def\tgrd{$t$-grading}

\def\qdgr{$q$-degree}

\def\dgq{ \deg_{\,q} }

\def\tdiff{twisted differential}
\def\tdiffs{\tdiff s}

\def\hmlg{ \mathrm{H} }

\def\hmlgMF{ \hmlg_{\mathrm{MF}} }
\def\hmlgMFv#1{ \hmlgMF^{#1} }
\def\hmlgMFb{ \hmlgMFv{\bullet} }
\def\hmlgMFi{ \hmlgMFv{i} }

\def\xcls{close}

\def\xg{ \gamma }

\def\xgv#1{ \xg_{#1} }
\def\xgo{ \xgv{1} }
\def\xgt{ \xgv{2} }

\def\xgi{ \xgv{i} }

\def\ZZn{ \ZZ^n }
\def\ZZnp{ \ZZ^{n\p} }

\def\xnvr{ m }
\def\xnpr{ n }

\def\tvm#1{ \tilde{#1} }
\def\nvm#1{ {#1} }

\def\nP{ P }

\def\nPvv#1#2{ \nP_{#1}(#2) }
\def\nPqv#1{ \nPvv{#1}{\xq} }
\def\nPLq{ \nPvv{\xL}{\xq} }

\def\nPGrq{ \nPqv{\xGr} }
\def\nPunq{ \nPqv{\xunkn} }

\def\nPgvv#1#2#3{ \nP_{#1;#2}(#3) }
\def\nPuvv#1#2{ \nPgvv{\msut}{#1}{#2} }
\def\nPovv#1#2{ \nPgvv{\msof}{#1}{#2} }
\def\nPuLq{ \nPuvv{\xL}{\xq} }
\def\nPoLq{ \nPovv{\xL}{\xq} }

\def\tCgu#1#2{ \tvm{C}_{#1}^{#2} }
\def\tCguv#1#2#3{ \tCgu{#1}{#2}\lrbc{#3} }
\def\tCuuv#1#2{ \tCguv{\msut}{#1}{#2} }
\def\tCubL{ \tCuuv{\bullet}{\xL} }

\def\tCgusv#1#2#3{ \tCgu{#1}{#2}\lrbs{#3} }
\def\tCuusv#1#2{ \tCgusv{\msut}{#1}{#2} }
\def\tCubs#1{ \tCuusv{\bullet}{\,\raise12pt\xybox{0;/r2pc/:#1}\,}}

\def\nC{ \nvm{C} }
\def\nCv#1{ \nC\lrbc{#1} }

\def\nCu#1{ \nC^{#1} }
\def\nCgu#1#2{ \nCu{#2}_{#1} }

\def\nCuv#1#2{ \nCu{#1}\lrbc{#2} }
\def\nCguv#1#2#3{ \nCgu{#1}{#2}\lrbc{#3} }
\def\nCuuv#1#2{ \nCguv{\msut}{#1}{#2} }
\def\nCouv#1#2{ \nCguv{\msof}{#1}{#2} }
\def\nCubL{ \nCuuv{\bullet}{\xL} }
\def\nCobL{ \nCouv{\bullet}{\xL} }
\def\nCusv#1#2{ \nCu{#1}\lrbs{\,\raise12pt\xybox{0;/r2pc/:#2}\,} }
\def\nCussv#1#2{ \nCu{#1}\lrbs{#2} }
\def\nCuLv#1{ \nCuv{#1}{\xL} }
\def\nCnL{ \nCuLv{n} }
\def\nCdv#1{ \nCuv{\bullet}{#1} }
\def\nCdsv#1{ \nCusv{\bullet}{#1} }
\def\nCdssv#1{ \nCussv{\bullet}{#1} }
\def\nCdL{ \nCdv{\xL} }
\def\nCdLp{ \nCdv{\xLp} }

\def\nCGr{ \nCv{\xGr} }

\def\xCugr#1#2#3{ C^{#1,#2}_{\qgrsv{#3}} }
\def\xCugrv#1#2#3#4{ \xCugr{#1}{#2}{#3}(#4) }
\def\xCugrL#1#2#3{ \xCugrv{#1}{#2}{#3}{\xL} }
\def\xCnijL{ \xCugrL{n}{j}{i} }

\def\yC{ C }
\def\yCvv#1#2{ \yC^{#1}_{#2} }
\def\yCNb{ \yCvv{\bullet}{\sunm} }
\def\yCNbL{ \yCNb(\xL) }
\def\yCNbLp{ \yCNb(\xLp) }

\def\dju{ \sqcup }
\def\jlg{ \# }
\def\jlgv#1#2{ \jlg_{#1,#2} }
\def\jlgij{ \jlgv{i}{j} }
\def\jlgino{ \jlgv{i}{n+1} }

\def\rIdv#1{ \mathord{ \yR(#1) } }

\def\vmvi#1{ #1 - #1^{-1} }

\def\qmqi{ \vmvi{\xq} }

\def\qtNpm{ q^{2\xN+1} - q^{-2\xN+1} }

\def\otRot{ \otimes_{\yR} }

\def\zgmf#1{ \hat{#1} }
\def\tgmf#1{ \tilde{#1} }

\def\gt{ \tau }
\def\gtp{ \gt\p }
\def\hgt{ \zgmf{\gt} }
\def\hgtp{ \zgmf{\gtp} }
\def\tgt{ \tgmf{\gt} }
\def\tgtin{ \tgt_{\xinit} }

\def\xg{ \gamma }
\def\xgp{ \xg\p}
\def\xgpp{ \xg^{\prime\prime} }

\def\xgm{ \xg_- }
\def\gg{ \zgmf{\xg} }
\def\ggp{ \gg\p }
\def\ggpp{ \gg^{\prime\prime}}
\def\ggm{ \gg_- }

\def\tg{ \tgmf{\xg} }

\def\tgv#1{ \tg_{#1} }
\def\tgin{ \tgv{\xinit} }
\def\tginv#1{ \tgv{\xinit,#1} }
\def\tginbr{ \tginv{\xbr} }
\def\tginbz{ \tginv{\bxz} }
\def\tgbr{ \tgv{\xbr} }
\def\tgbz{ \tgv{\bxz} }

\def\xgv#1{ \xg_{#1} }
\def\xgo{ \xgv{1} }
\def\xgt{ \xgv{2} }
\def\xgbk{ \xgv{\bxk} }
\def\xgbr{ \xgv{\xbr} }

\def\xinit{\mathrm{init}}

\def\ggv#1{ \gg_{#1} }
\def\ggo{ \ggv{1} }
\def\ggt{ \ggv{2} }
\def\ggi{ \ggv{i} }

\def\gginv#1{ \ggv{\xinit,#1} }
\def\gginbk{ \gginv{\bxk} }

\def\ggbk{ \ggv{\bxk} }

\def\ggbr{ \ggv{\xbr} }

\def\ggprp{ \ggv{\mprp} }

\def\ggtprp{ \ggv{\tprp} }

\def\otR{ \otimes_{\yR} }

\def\xG{ \Gamma }
\def\xGv#1{ \xG_{#1} }
\def\xGr{ \xGv{\xbr} }

\def\xGvv#1#2{ \xGv{#1;#2} }
\def\xGbrk{ \xGvv{\xbr}{\bxk} }

\def\hGvv#1#2{ \hGv{#1,#2} }
\def\hGbrk{ \hGvv{\xbr}{\bxk} }

\def\hG{ \zgmf{\xG} }

\def\tG{ \tgmf{\xG} }
\def\tGv#1{ \tG_{#1} }
\def\tGbr{ \tGv{\xbr} }

\def\hL{ \zgmf{\xL} }

\def\hGv#1{ \hG_{#1} }
\def\hGr{ \hGv{\xbr} }

\def\ncL{ \nc \xL }
\def\ncL{ n_{\xL} }

\def\an{ n }
\def\anv#1{ \an_{#1} }
\def\anbk{ \anv{\bxk} }

\def\chqd{ \chi_q }

\def\chqt{ \chi_{q,t} }

\def\dmv#1{ \dim_{#1} }
\def\dmq{ \dmv{q} }

\def\xr{ r }
\def\xrv#1{ \xr_{#1} }
\def\xri{ \xrv{i} }
\def\xbr{ \mathbf{\xr} }
\def\xbra{ |\xbr| }

\def\xLv#1{ \xL_{#1} }
\def\xLo{ \xLv{1} }
\def\xLt{ \xLv{2} }
\def\xLot{ \xLo\sqcup\xLt }

\def\zx{ x }
\def\zy{ y }
\def\zu{ u }
\def\zv{ v }

\def\zuv#1{ \zu_{#1} }
\def\zuo{ \zuv{1} }
\def\zut{ \zuv{2} }
\def\zuth{ \zuv{3} }
\def\zuf{ \zuv{4} }

\def\zxo{ \zx_1 }
\def\zxt{ \zx_2 }
\def\zxth{ \zx_3 }
\def\zxf{ \zx_4 }
\def\zxi{ \zx_i }

\def\zyo{ \zy_1 }
\def\zyt{ \zy_2 }
\def\zyth{ \zy_3 }
\def\zyf{ \zy_4 }
\def\zyi{ \zy_i }

\def\xone{ e }
\def\xones{ \xone^{\ast} }

\def\zxgv#1{ \zx_{(#1)} }
\def\zxgi{ \zxgv{i} }

\def\zxgtN{ \zxgv{2\xN} }

\def\zxsgv#1{ \zxs_{(#1)} }
\def\zxsgi{ \zxsgv{i} }

\def\zxuv#1{ \zx^{#1} }
\def\zxui{ \zxuv{i} }
\def\zxutN{ \zxuv{2\xN} }
\def\zxuis{ (\zxui)^\ast }
\def\zyus{ \zy^\ast }

\def\zxgv#1{ \hmx^{#1}(\xone) }
\def\zxgi{ \zxgv{i} }
\def\zxgtN{ \zxgv{2\xN} }

\def\zxsgv#1{ \lrbsc{\zxgv{#1}}^\ast }
\def\zxsgi{ \zxsgv{i} }

\def\zys{ \lrbsc{\zyg}^{\ast} }
\def\zyg{ \hmy(\xone) }

\def\hmxv#1{ \hmv{\zx}_{#1} }
\def\hmyv#1{ \hmv{\zy}_{#1} }

\def\hxo{ \hmxv{1} }
\def\hxt{ \hmxv{2} }
\def\hxth{ \hmxv{3} }
\def\hxf{ \hmxv{4} }

\def\hyo{ \hmyv{1} }
\def\hyt{ \hmyv{2} }
\def\hyth{ \hmyv{3} }
\def\hyf{ \hmyv{4} }

\def\hv{ \hat{v} }
\def\tla{ \tilde{a} }
\def\tlb{ \tilde{b} }

\def\IQbxy{ \IQ[\bfx,\bfy] }

\def\lrbsc#1{ \big( #1 \big) }

\def\lrbcsmo#1{ \lrbc{#1}\gszto }

\def\zxall{\zxo+\zxt+\zxth+\zxf}
\def\zyall{\zyo+\zyt+\zyth+\zyf}

\def\zbx{ \mathbf{x} }
\def\zby{ \mathbf{y} }
\def\zbu{ \mathbf{\zu} }

\def\zxfo{ \zxo+\zxf }
\def\zxft{ \zxt+\zxf }
\def\zxfth{ \zxth+\zxf }

\def\zyfo{ \zyo+\zyf }
\def\zyft{ \zyt+\zyf }
\def\zyfth{ \zy_3 + \zy_4 }

\def\zxtho{ \zxo+\zxth }
\def\zytho{ \zyo+\zyth }

\def\zA { A }
\def\zB { B }
\def\zC { C }

\def\ztC{ \tilde{C} }

\def\zAxy{ \zA(\zbx,\zby) }
\def\zBxy{ \zB(\zbx,\zby) }

\def\yp{ p }
\def\yq{ q }

\def\ypv#1{ \yp_{#1} }
\def\yqv#1{ \yq_{#1} }
\def\ypi{ \ypv{i} }
\def\yqi{ \yqv{i} }
\def\ypj{ \ypv{j} }
\def\yqj{ \yqv{j} }

\def\ypko{ \ypv{k+1} }
\def\yqko{ \yqv{k+1} }

\def\coltv#1#2{ \begin{pmatrix} #1\\#2 \end{pmatrix} }
\def\tcoltv#1#2#3#4{ \coltv{#1}{#3},\coltv{#2}{#4} }

\def\xlmb{ \lambda  }

\def\rtrv#1#2#3{ [#1,#2]_{#3} }
\def\rtrpv#1#2#3{ [#1,#2]_{#3}\p }
\def\rtrijl{ \rtrv{i}{j}{\xlmb} }
\def\rtrpijl{ \rtrpv{i}{j}{\xlmb} }
\def\rtrthoo{ \rtrv{3}{1}{1} }
\def\rtrfto{ \rtrv{4}{2}{1} }
\def\rtov#1#2{ [#1]_{#2} }
\def\rtoil{ \rtov{i}{\xlmb} }

\def\yconv{ \mathrm{conv} }

\def\xD{ D }
\def\xDp{ \xD\p }

\def\xDconv{ \xD_{\yconv} }

\def\xDv#1{ \xD_{#1} }

\def\xDpar{ \xDv{\ipar} }
\def\xDhor{ \xDv{\ihor} }
\def\xDvirt{ \xDv{\ivir} }

\def\xDvv#1#2{ \xDv{#1,#2} }
\def\xDAB{ \xDvv{\wA}{\wB} }
\def\xDAC{ \xDvv{\wA}{\wC} }
\def\xDMMp{ \xDvv{\yM}{\yMp} }

\def\xDS{ \xDv{\xSset} }

\def\xDbk{ \xDv{\bxk} }

\def\Zt{ \ZZ_2 }
\def\Ztgrd{$\Zt$-graded}

\def\xW{ W }
\def\xWp{ \xW\p }

\def\nW{ \nvm{\xW} }

\def\nWv#1#2{ \nW(#1,#2) }

\def\nWxy{ \nWv{\zx}{\zy} }
\def\nWxyi{ \nWv{\zxi}{\zyi} }
\def\nWhv#1{ \nW(#1) }
\def\nWhx{ \nWhv{\zx} }

\def\nWgv#1{ \nW_{#1} }
\def\nWu{ \nWgv{\msut} }
\def\nWo{ \nWgv{\msof} }
\def\nWoxy{ \nWo(\zx,\zy) }

\def\gWxy{ \nWxy }
\def\gWmxy{ \nWv{-\zx}{-\zy} }

\def\xWv#1{ \xW_{#1} }

\def\xWlg{ \xWv{\xmlg} }

\def\xsgm{ \sigma }
\def\sXv#1#2{ \xsgm_{#1#2} }

\def\sXot{ \sXv{1}{2} }

\def\sXoth{ \sXv{1}{3} }

\def\sXtth{ \sXv{2}{3} }

\def\hsgm{ \hat{\xsgm} }

\def\hXv#1#2{ \hsgm_{#1#2} }

\def\hXot{ \hXv{1}{2} }

\def\hXoth{ \hXv{1}{3} }

\def\hXtth{ \hXv{2}{3} }

\def\xmu{ \mu }
\def\xmui{ \xmu_i }

\def\abt#1#2{ #1#2^2 }
\def\atN#1{ #1^{2\xN+1} }

\def\xyt{ \abt{\zx}{\zy} }
\def\xtN{ \atN{\zx} }

\def\xw{ w }
\def\xwv#1{ \xw(#1) }
\def\xwx{ \xwv{\zx} }
\def\xwvv#1#2{ \xwv{#1,#2} }
\def\xwvvv#1#2#3{ \xwv{#1,#2,#3} }
\def\xwvf#1#2#3#4{ \xwv{#1,#2,#3,#4} }
\def\xwxot{ \xwvv{\zxo}{\zxt} }
\def\xwxott{ \xwvvv{\zxo}{\zxt}{\zxth} }

\def\yR{ \mathsf{R} }
\def\yRp{ \yR\p }
\def\yRuv#1{ \yR^{#1} }
\def\yRn{ \yRuv{n} }
\def\yRns{ (\yRn)^\ast }
\def\yRv#1{ \yR_{#1} }
\def\yRz{ \yRv{0} }
\def\yRo{ \yRv{1} }
\def\yRi{ \yRv{i} }

\def\yRzz{ \yRv{00} }
\def\yRzo{ \yRv{01} }
\def\yRoz{ \yRv{10} }
\def\yRoo{ \yRv{11} }
\def\yRij{ \yRv{ij} }

\def\yRp{ \yR\p }
\def\yRpv#1{ \yRp_{#1} }
\def\yRpz{ \yRpv{0} }
\def\yRpo{ \yRpv{1} }
\def\yRpj{ \yRpv{j} }

\def\yRt{ \yR^2 }
\def\yRtv#1{ \yRt_{#1} }
\def\yRtz{ \yRtv{0} }
\def\yRto{ \yRtv{1} }

\def\yRtvv#1#2{ \yRtv{#1,#2} }
\def\yRtzpar{ \yRtvv{\ipar}{0} }
\def\yRtopar{ \yRtvv{\ipar}{1} }
\def\yRtzvirt{ \yRtvv{\ivir}{0} }
\def\yRtovirt{ \yRtvv{\ivir}{1} }
\def\yRtzhor{ \yRtvv{\ihor}{0} }
\def\yRtohor{ \yRtvv{\ihor}{1} }

\def\yRs{ \yR^6 }
\def\yRsv#1{ \yRs_{#1} }
\def\yRsz{ \yRsv{0} }
\def\yRso{ \yRsv{1} }

\def\yRs{ \yR^6 }
\def\yRsv#1{ \yRs_{#1} }
\def\yRsz{ \yRsv{0} }
\def\yRso{ \yRsv{1} }

\def\yRf{ \yR^4 }
\def\yRfv#1{ \yRf_{#1} }
\def\yRfz{ \yRfv{0} }
\def\yRfo{ \yRfv{1} }

\def\yRprp{ \yR_{\mprp} }

\def\ybp{ \mathbf{\yp} }
\def\ybq{ \mathbf{\yq} }

\def\yM{ M }
\def\yMv#1{ \yM_{#1} }
\def\yMz{ \yMv{0} }
\def\yMo{ \yMv{1} }

\def\yMxp{ \yMv{\bfx\p} }
\def\yMyp{ \yMv{\bfy\p} }

\def\yMp{ \yM\p }
\def\yMs{ \yM^{\ast} }

\def\Homv#1{ \Hom_{#1} }
\def\HomR{ \Homv{\yR} }
\def\MFrm{ \mathrm{MF} }
\def\PMFrm{ \mathrm{PMF} }
\def\HomMF{ \Homv{\MFrm} }
\def\HomPMF{ \Homv{\PMFrm} }

\def\Endv#1{ \End_{#1} }
\def\EndMF{ \Endv{\MFrm} }

\def\MFv#1{ \MFrm_{#1} }
\def\MFvv#1#2{ \MFv{#1,#2} }
\def\MFW{ \MFv{\xW} }
\def\MFRW{ \MFvv{\yR}{\xW} }
\def\MFWf{ \MFv{\zWf} }
\def\MFWt{ \MFv{\zWt} }

\def\PMFvv#1#2{ \PMFrm_{#1}^{#2} }
\def\PMFWn{ \PMFvv{\xW}{\zan} }
\def\PMFWnp{ \PMFvv{\xW\p}{\zan\p} }
\def\PMFWfo{ \PMFvv{\zWf}{1} }

\def\ztcm#1#2{ [#1,#2]_{\mathrm{s} } }
\def\ztac#1#2{ \{#1,#2\} }

\def\Extmfv#1{ \Ext^{#1} }
\def\Extmfb{ \Extmfv{\bullet} }
\def\Extmfz{ \Extmfv{0} }

\def\Extmfzqv#1{ \Extmfz_{\qgrsv{#1}} }
\def\Extmfoqv#1{ \Extmfo_{\qgrsv{#1}} }
\def\ExtmfzqN{ \Extmfzqv{2\xN} }
\def\Extmfo{ \Extmfv{1} }
\def\Extmfi{ \Extmfv{i} }
\def\ExtmfoqNmo{ \Extmfoqv{\xN-1} }

\def\Extmfoqo{ \Extmfoqv{1} }

\def\ExtPmf{ \Ext_{\mathrm{P}} }
\def\ExtPmfv#1{ \ExtPmf^{#1} }
\def\ExtPmfb{ \ExtPmfv{\bullet} }
\def\ExtPmfz{ \ExtPmfv{0} }

\def\Endv#1{ \End_{#1} }
\def\EndR{ \Endv{\yR} }

\def\xal{ \mathrm{lng} }

\def\degv#1{ \deg_{#1} }
\def\degZt{ \degv{\Zt} }
\def\degZ{ \degv{\ZZ} }

\def\degl{ \deg_{\xal} }
\def\deglv#1{ \deg_{#1} }

\def\bdegl{ \mathbf{deg}_{\xal} }

\def\zp{ p }
\def\zq{ q }
\def\zr{ r }
\def\zxbp{ \mathbf{\zp} }
\def\zxbq{ \mathbf{\zq} }
\def\zxbr{ \mathbf{\zr} }
\def\zpv#1{ \zp_{#1} }
\def\zqv#1{ \zq_{#1} }
\def\zrv#1{ \zr_{#1} }
\def\zpo{ \zpv{1} }
\def\zpt{ \zpv{2} }
\def\zqo{ \zqv{1} }
\def\zqt{ \zqv{2} }
\def\zro{ \zrv{1} }
\def\zrt{ \zrv{2} }

\def\zP{ P }
\def\zPs{ \zP^{\ast} }
\def\zQ{ Q }
\def\zQs{ \zQ^{\ast} }

\def\ardpl{ \ar@{}[d]|{\oplus} }

\def\zPv#1{ \zP_{#1} }
\def\zQv#1{ \zQ_{#1} }
\def\zPcn{ \zPv{\yconv} }
\def\zQcn{ \zQv{\yconv} }

\def\zred{ \mathrm{r} }

\def\zPrcn{ \zPv{\zred,\yconv} }
\def\zQrcn{ \zQv{\zred,\yconv} }

\def\zPpv#1{ \zP\p_{#1} }
\def\zQpv#1{ \zQ\p_{#1} }

\def\zPprcn{ \zPpv{\zred,\yconv} }
\def\zQprcn{ \zQpv{\zred,\yconv} }

\def\zPi{ \zPv{i} }
\def\zQi{ \zQv{i} }

\def\zPpar{ \zPv{\ixpar} }
\def\zPhor{ \zPv{\ixhor} }
\def\zPvirt{ \zPv{\ixvir} }

\def\zQpar{ \zQv{\ixpar} }
\def\zQhor{ \zQv{\ixhor} }
\def\zQvirt{ \zQv{\ixvir} }

\def\zGv#1{ \xG_{#1} }
\def\zGpar{ \zGv{\ixpar} }
\def\zGhor{ \zGv{\ixhor} }
\def\zGV{ \zGv{\ixver} }

\def\znGpar{ \nCv{\zGpar} }
\def\znGhor{ \nCv{\zGhor} }
\def\znGV{ \nCv{\zGV} }

\def\znGr{ \nCv{\xGr} }

\def\Sf{ S_4 }
\def\Sth{ S_3 }

\def\St{ S_2 }

\def\zW{ W }
\def\zWf{ \zW_4 }
\def\zWfp{ \zW_{4,\,\mprp} }

\def\zWt{ \zW_2 }

\def\bfx{ \mathbf{\zzx} }

\def\hbfx{ \hmv{\bfx} }

\def\basW{ W }

\def\zzk{ \xnvr }
\def\zzn{ n }
\def\zzx{ x }

\def\yym{ m }
\def\yyn{ n }

\def\bfxv#1{ \bfx_{#1} }
\def\bfxlg{ \bfxv{\xmlg} }
\def\bfxi{ \bfxv{i} }
\def\bfxj{ \bfxv{j} }

\def\xmlg{ \mathfrak{l} }

\def\zWfxy{ \zWf(\bfx,\bfy) }
\def\zWfx{ \zWf(\bfx) }

\def\pcol#1#2{ \begin{pmatrix} #1_1 \\ \vdots \\ #1_{#2} \end{pmatrix} }

\def\hmv#1{ \hat{#1} }
\def\hmx{ \hmv{x} }
\def\hmy{ \hmv{y} }
\def\hmempt{ \hmv{{}} }

\def\hmr{ \hmv{r} }

\def\xhmfv#1{ \widehat{#1} }
\def\zhom{ h }
\def\zxhom{ h }
\def\zhomprp{ \zhom_{\mprp} }
\def\zhhom{ \xhmfv{\zhom} }
\def\zhomv#1#2{ \zhom_{#1,#2} }
\def\zhhomv#1#2{ \zhhom_{#1,#2} }
\def\zhomij{ \zhomv{i}{j} }
\def\zhhomij{ \zhhomv{i}{j} }

\def\zhhomn{ \zhhom_{\xnb} }

\def\zhhomtprp{ \prpv{\zhhom} }

\def\xgk{ k }

\def\xgkvv#1#2{ \xgk_{#1,#2} }
\def\xgkpq{ \xgkvv{\zp}{\zq} }

\def\grxhv#1{ \{#1 \} }

\def\grhm{ 2\xN-1 }

\def\grsv#1{ \left\{ #1 \right\} }

\def\grsk{ \grsv{\xgk} }

\def\grsph{ \grsv{\hlf} }
\def\grsmh{ \grsv{-\hlf} }

\def\grsev#1{ \left[#1\right] }
\def\grseo{ \grsev{1} }
\def\grsemo{ \grsev{-1} }
\def\grseph{ \grsev{\hlf} }
\def\grsemh{ \grsev{-\hlf} }

\def\hfpq{ \frac{1}{2}\,(\degZ\zp-\degZ\zq) }
\def\grshfpq{ \grsv{\xgkpq} }

\def\grkdl{ \grxhv{2\xN} }
\def\grkdlm{ \grxhv{-2\xN} }

\def\qgrsv#1{ \{#1\} }

\def\gsztv#1{ \left\langle #1 \right\rangle }
\def\gszto{ \gsztv{1} }

\def\gsztph{ \gsztv{\hlf} }
\def\gsztmh{ \gsztv{-\hlf} }

\def\tK{ \mathrm{K} }
\def\kcv#1{ \tK(#1) }

\def\kcyp{ \kcv{\yp} }

\def\kcbfp{ \kcv{\bfyp} }

\def\kmfv#1#2{ \tK(#1;#2) }
\def\kmfbv#1#2{ \tK\big(#1;#2\big) }
\def\kmfpq{ \kmfv{\yp}{\yq} }
\def\kmfpqi{ \kmfv{\yp_i}{\yq_i} }
\def\kmfbpq{ \kmfv{\bfyp}{\bfyq} }

\def\kmftwv#1#2#3#4{
\tK\begin{pmatrix}
#1, & #3
\\
#2, & #4
\end{pmatrix}
}

\def\kmffwv#1#2#3#4#5#6#7#8{
\tK\begin{pmatrix}
#1, & #5
\\
#2, & #6
\\
#3, & #7
\\
#4, & #8
\end{pmatrix}
}

\def\kctv#1#2{
\tK\begin{pmatrix}
#1
\\
#2
\end{pmatrix}
}

\def\mA{ A }

\def\mB{ B }

\def\mC{ C }

\def\mrf{ f }
\def\mrg{ g }

\def\yF{ F }
\def\yFp{ \yF\p }
\def\yFpp{ \yF\pp }
\def\yFprp{ \yF_{\mprp} }

\def\yFtprp{ \tprpv{\yF} }

\def\yG{ G }

\def\yGp{ \yG\p }
\def\yGpp{ \yG\pp }
\def\yGprp{ \yG_{\mprp} }

\def\yH{ H }
\def\yFv#1{ \yF_{#1} }
\def\yFz{ \yFv{0} }
\def\yFo{ \yFv{1} }
\def\yFi{ \yFv{i} }
\def\yFj{ \yFv{j} }
\def\yFbvv#1#2{ \yFv{#1,#2} }
\def\yFbki{ \yFbvv{\bxk}{i} }
\def\yGv#1{ \yG_{#1} }
\def\yGz{ \yGv{0} }
\def\yGo{ \yGv{1} }
\def\yHv#1{ \yH_{#1} }
\def\yHz{ \yHv{0} }
\def\yHo{ \yHv{1} }
\def\yHp{ \yH\p }
\def\yHpp{ \yH\pp }
\def\yHpp{ \yH^{\prime\prime} }
\def\yHprp{ \yH_{\mprp} }

\def\yf{ f }
\def\yfv#1{ \yf_{#1} }

\def\yfvv#1#2{ \yfv{#1,#2} }
\def\yfbij{ \yfvv{\bxi}{\bxj} }

\def\yfS{ \yfv{\xSset} }
\def\yfbS{ \yfv{\bSset} }

\def\yg{ g }
\def\ygv#1{ \yg_{#1} }

\def\ygSp{ \ygv{\xSsetp} }

\def\sxhv#1{ \mathcal{F}_{#1} }
\def\sxhxS{ \sxhv{\xSset} }

\def\idv#1{ \id_{#1} }
\def\idS{ \idv{\xSset} }

\def\yFi{ \yFv{i} }

\def\xeclv#1{ [#1] }
\def\xeclXk{ \xeclv{\yXk} }
\def\xeclXlpl{ \xeclv{\yXblpl} }

\def\yX{ X }
\def\yXprp{ \prpv{\yX} }
\def\yXv#1{ \yX_{#1} }
\def\yXz{ \yXv{0} }
\def\yXo{ \yXv{1} }

\def\yXk{ \yXv{k} }
\def\yXm{ \yXv{m} }
\def\yXl{ \yXv{l} }
\def\yXlp{ \yXl\p }
\def\yXkz{ \yXv{k,0} }
\def\yXkmo{ \yXv{k-1} }
\def\yXl{ \yXv{l} }
\def\yXp{ \yX\p }
\def\yXpp{ \yX^{\prime\prime} }

\def\ycAv#1{ \mathcal{A}(#1) }
\def\ycAXk{ \ycAv{\yXk} }
\def\ycAXpk{ \ycAv{\yXpk} }

\def\yY{ Y }
\def\yYinv{ \yY^{-1} }
\def\yYp{ \yY\p }

\def\yYv#1{ \yY_{#1} }
\def\yYk{ \yYv{k} }

\def\yZ{ Z }
\def\yZp{ \yZ\p }

\def\wX{ X }
\def\wXp{ \wX\p }

\def\yFu{ \yF_{\yup} }
\def\yFd{ \yF_{\ydown} }
\def\yGu{ \yG_{\yup} }
\def\yGd{ \yG_{\ydown} }

\def\ymrpv#1{ {#1}^{\prime} }
\def\yFp{ \ymrpv{\yF} }
\def\yGp{ \ymrpv{\yG} }

\def\xmlt{ \mathrm{m} }
\def\xcmlt{ \Delta }
\def\ymrmv#1{ #1_{\xmlt} }
\def\ymrDv#1{ #1_{\xcmlt} }
\def\yFm{ \ymrmv{\yF} }

\def\yHD{ \ymrDv{\yH} }

\def\xfis{ f_{\cchq} }
\def\xfisv#1{ f_{\cchq,#1} }
\def\xfisz{ \xfisv{0} }
\def\xfiso{ \xfisv{1} }
\def\xfist{ \xfisv{2} }

\def\xgisv#1{ g_{\cong,#1} }
\def\xgisz{ \xgisv{0} }
\def\xgiso{ \xgisv{1} }

\def\Kmpl{ \mathbf{K} }
\def\Kmplv#1{ \Kmpl(#1) }

\def\Obv#1{ \mathord{\mathrm{Ob}\lrbc{#1}} }

\def\zK{ K }
\def\zKcmn{ \zK_{\xxcmn} }

\def\Cn{ \mathrm{Cone} }
\def\Cnv#1{ \Cn\lrbc{#1} }

\def\Cnvv#1#2{ \Cn_{#1}\lrbc{#2} }

\def\KmplMF{ \Kmplv{\MF} }
\def\KmplMFWf{ \Kmplv{\MFWf} }

\def\CnMFv#1{ \Cnvv{\MF}{#1} }
\def\CnPSv#1{ \Cnvv{\PS}{#1} }

\def\Cvl{ \mathrm{Conv} }
\def\Cvlvv#1#2{ \Cvl_{#1}\lrbc{#2} }
\def\CvlMFv#1{ \Cvlvv{\MF}{#1} }

\def\xfrbth{ \save [].[rr]!C="g1" *[F]\frm{} \restore }

\def\xCthwv#1#2#3#4#5{ {#1} \ar[r]^-{#2}="a" & {#3} \ar[r]^-{#4}="b" &
{#5} \save {{[].[ll]}."a"}."b"*[F]\frm{} \restore }

\def\biCnMFthv#1#2#3#4#5{ \xfrbth {#1} \ar[r]^{#2} & {#3} \ar[r]^{#4} & {#5} }

\def\bicnmfFGsh{ \biCnMFthv{\hlparsh}{\yF}{\hlvirsh}{\yG}{\hlhorsh} }
\def\bCnMFthv#1#2#3#4#5{ \xymatrix{ \biCnMFthv{#1}{#2}{#3}{#4}{#5} } }

\def\cnmfFGsh{ \bCnMFthv{\hlparsh}{\yF}{\hlvirsh}{\yG}{\hlhorsh} }

\def\ttNo{ 2(2\xN+1) }
\def\spxot{ \siztN \hmx_1^i\,\hmx_2^{2\xN-i} }

\def\xlst#1#2{ (#1_{1},\dotsc, #1_{#2}) }
\def\ylst#1#2{ #1_{1},\dotsc, #1_{#2} }

\def\xJ{ J }
\def\xJv#1{ \xJ_{#1} }
\def\JrW{ \xJv{\nW} }
\def\JrWp{ \JrW\p }
\def\JrWpdbl{ \JrWp\otimes\JrWp }
\def\JrWqv#1{ \xJv{\nW,\qgrsv{#1} } }

\def\TrF{ \Tr_{\mathrm{F}} }

\def\botoncL{ \bigotimes_{i=1}^{\ncL} }

\def\bopiztN{ \bigoplus_{i=0}^{2\xN} }

\def\bopiztNmo{ \bigoplus_{i=0}^{2\xN-1} }

\def\bopiZjt{ \bigoplus_{\begin{subarray}{l}
i\in\ZZ \\ j\in\Zt
\end{subarray}} }

\def\bopbkzn{ \bigoplus_{\bxk\in\ZZ^\zan} }

\def\snjZit{ \sum_{\begin{subarray}{r}
i,n\in\ZZ \\ j\in\Zt
\end{subarray}} }

\def\sxbr{ \sum_{\xbr} }

\def\siztN{ \sum_{i=0}^{2\xN} }

\def\ychi{ \chi }
\def\ychiv#1{ \ychi_{#1} }
\def\yin{ \mathrm{in} }
\def\yout{ \mathrm{out} }
\def\yup{ \mathrm{u} }
\def\ydown{ \mathrm{d} }
\def\ychin{ \ychiv{\yin} }
\def\ychout{ \ychiv{\yout} }
\def\ychiz{ \ychiv{0} }
\def\ychio{ \ychiv{1} }

\def\ytchi{ \tilde{\ychi} }
\def\ytchin{ \ytchi_{\yin} }
\def\ytchout{ \ytchi_{\yout} }

\def\ychinv#1{ \ychiv{\yin,#1} }
\def\ychoutv#1{ \ychiv{\yout,#1} }

\def\ychinprp{ \ychinv{\mprp} }
\def\ychinprpp{ \ychinprp\p}
\def\ychoutprp{ \ychoutv{\mprp} }

\def\ychvin#1{ \ychiv{#1,\yin} }
\def\ychvout#1{ \ychiv{#1,\yout} }

\def\ychsSin{ \ychvin{\xSset} }
\def\ychsSpin{ \ychvin{\xSset\p} }
\def\ychsSSpin{ \ychvin{\xSSpset} }
\def\ychbSout{ \ychvout{\bSset} }
\def\ychbSpout{ \ychvout{\bSset\p} }
\def\ychsSSpout{ \ychvout{\bSSpset} }

\def\ychsv#1{ \ychiv{\ast,#1} }
\def\ychsu{ \ychsv{\yup} }
\def\ychsd{ \ychsv{\ydown} }

\def\ychinu{ \ychinv{\yup} }
\def\ychind{ \ychinv{\ydown} }
\def\ychini{ \ychinv{i} }
\def\ychoutu{ \ychoutv{\yup} }
\def\ychoutd{ \ychoutv{\ydown} }
\def\ychouti{ \ychoutv{i} }

\def\yFu{ \yF_{\yup} }
\def\yFd{ \yF_{\ydown} }
\def\yGu{ \yG_{\yup} }
\def\yGd{ \yG_{\ydown} }

\def\ychpin{ \ychpiv{\yin} }
\def\ychpout{ \ychpiv{\yout} }

\def\ychip{ \ychi\p }
\def\ychipv#1{ \ychip_{#1} }
\def\ychpin{ \ychipv{\yin} }
\def\ychpout{ \ychipv{\yout} }

\def\grso{ \grsv{1} }
\def\grsmo{ \grsv{-1} }
\def\grst{ \grsv{2} }
\def\grsmt{ \grsv{-2} }
\def\grsmN{ \grsv{-\xN} }
\def\grsN{ \grsv{\xN} }
\def\grstmN{ \grsv{2-\xN} }

\def\grsmNo{ \grsv{1-\xN} }
\def\yRomNo{ \yRo\grsmNo }

\def\yRzmo{ \yRz\grsmo }
\def\yRzso{ \yRz\grso }
\def\yRoso{ \yRo\grso }
\def\yRomo{ \yRo\grsmo }
\def\yRomN{ \yRo\grsmN }
\def\yRotmN{ \yRo\grstmN }

\def\hlparsh{ \hlpar\grsmo }
\def\hlvirsh{ \hlvir\gszto }

\def\hlhorsh{ \hlhor\grso }

\def\hlhorsho{ \hlhor\grsmo }
\def\hlparsho{ \hlpar\grso }

\def\xxa{ a }
\def\xxb{ b }

\def\smiztN{ \sum_{i=0}^{2\xN} }
\def\smion{ \sum_{i=1}^n }
\def\smiok{ \sum_{i=1}^k }

\def\ssetZ{ \mathcal{S} }

\def\ssetZm{ \ssetZ_- }
\def\xaAm{ \xaA_- }
\def\crv#1{ |#1| }
\def\crssetZ{ \crv{\ssetZ} }

\def\mprp{ \mathrm{p} }
\def\tprp{ \widetilde{\mprp} }
\def\prpv#1{ #1_{\mprp} }
\def\tprpv#1{ #1_{\tprp} }

\def\arcxyoo{ \lxoarvv{1}{2} }
\def\zarcxyoo{ \hloarvv{1}{2} }

\def\zan{ n }

\def\HomRngl{ \Hom_{\yR;\geq 0} }
\def\EndRngl{ \End_{\yR;\geq 0} }

\def\ZZinf{ \ZZ_{\infty} }

\def\xsbsv#1#2{ #1_{#2} }
\def\Asbsv#1{ \xsbsv{\xaA}{#1} }
\def\Bsbsv#1{ \xsbsv{\xaB}{#1} }
\def\ABsbsv#1{ \xsbsv{(\xaA\otimes\xaB)}{#1} }
\def\AsbsS{ \Asbsv{\xSset} }
\def\AsbxS{ \Asbsv{\bSset} }
\def\BsbsS{ \Bsbsv{\xSset} }
\def\BsbsSp{ \Bsbsv{\xSsetp} }
\def\BsbxS{ \Bsbsv{\bSset} }
\def\ABsbsS{ \ABsbsv{\xSset\times\xSset\p} }

\def\xSset{ \mathcal{S} }
\def\xSsetp{ \xSset\p }
\def\bSset{ \bar{\xSset} }
\def\bSsetp{ \bSset\p }
\def\xSSpset{ \xSset\times\xSsetp }
\def\bSSpset{ \bSset\times\bSsetp }

\def\sbslv#1{ (\leq #1) }
\def\sbsgv#1{ (\geq #1) }
\def\sbslbk{ \sbslv{\bxk} }
\def\sbsgbk{ \sbsgv{\bxk} }
\def\Asbslv#1{ \xaA_{\sbslv{#1}} }
\def\Asbsgv#1{ \xaA_{\sbsgv{#1}} }

\def\Asbslbk{ \Asbslv{\bxk} }
\def\Asbsgbk{ \Asbsgv{\bxk} }

\def\bxi{ \mathbf{i} }
\def\bxj{ \mathbf{j} }
\def\bxk{ \mathbf{k} }
\def\bxv{ \mathbf{v} }
\def\bxw{ \mathbf{w} }
\def\bxz{ \mathbf{0} }
\def\bxl{ \mathbf{l} }
\def\bxlp{ \bxl\p }
\def\bxe{ \mathbf{e} }
\def\bxm{ \mathbf{m} }

\def\xaAv#1{ \xaA_{#1} }
\def\xaAbi{ \xaAv{\bxi} }
\def\xaAbk{ \xaAv{\bxk} }
\def\xaAbl{ \xaAv{\bxl} }
\def\xaAblp{ \xaAv{\bxlp} }

\def\xaAp{ \xaA\p }
\def\xaApv#1{ \xaAp_{#1} }
\def\xaApbk{ \xaApv{\bxk} }

\def\yXbk{ \yXv{\bxk} }
\def\yXbvv#1#2{ \yXv{#2,#1} }
\def\yXbkl{ \yXbvv{\bxk}{\bxl} }
\def\yXbkk{ \yXbvv{\bxk}{\bxk} }
\def\yXblpl{ \yXbvv{\bxl}{\bxlp} }
\def\yXblplz{ \yXbvv{\bxl;0}{\bxlp} }

\def\yXblk{ \yXbvv{\bxl}{\bxk} }

\def\yXpv#1{ \yXp_{#1} }
\def\yXpk{ \yXpv{k} }
\def\yXpm{ \yXpv{m} }

\def\DyXk{ \Delta\yXk }

\def\xtlv#1{ |#1| }

\def\xbt{ \beta }
\def\xnb{ n }

\def\wrth{ \mathrm{w} }
\def\wrthL{ \wrth(\xL) }
\def\ncmpv#1{ \# #1 }
\def\xnp{ n_{+} }
\def\xnn{ n_{-} }

\def\wA{ A }
\def\wB{ B }
\def\wC{ C }
\def\wF{ F }
\def\wFp{ \wF\p }
\def\wG{ G }
\def\wGp{ \wG\p }

\def\xDA{ \xDv{\wA} }
\def\xDAp{ \xDA\p }
\def\xDB{ \xDv{\wB} }
\def\xDC{ \xDv{\wC} }

\def\xdz{ d_0 }

\def\Cnthbgv#1#2#3#4#5#6{ \Cnv{
\xymatrix{{#1} \ar[r]^-{#2} \ar@/_1pc/@{-->}[rr]_{#6} & {#3} \ar[r]^-{#4} & {#5} } }
}

\def\ythbgdv#1#2#3#4#5#6{
\xymatrix{{#1} \ar[r]^-{#2} \ar@/_2pc/@{-->}[rr]_{#6} & {#3} \ar[r]^-{#4} & {#5} }
}

\def\ythbguv#1#2#3#4#5#6{
\xymatrix{ {#1} \ar[r]^-{#2} \ar@/^2pc/@{-->}[rr]^{#6} & {#3} \ar[r]^-{#4} & {#5} }
}

\def\smpm#1{ \lrbc{\! \begin{smallmatrix} #1 \end{smallmatrix} \!}
}

\def\Ctwwv#1#2#3{
{#1} \ar[r]^-{#2}="f" & {#3}
\save {[].[l]}."f"
*[F]\frm{} \restore
}

\def\Cbdv#1#2#3{
{#1} \ar[r]^-{#2}="f" & {#3}
\save {[].[l]}."f"
*[F]\frm{} \restore
}

\def\Pbdv#1#2#3{
{#1} \ar[r]^-{#2}="f" & {#3}
\save {[].[l]}."f"
*[F.]\frm{} \restore
}

\def\Cbtwv#1#2#3#4#5#6
{
{#1} \ar[r]^-{#2}="f" \ar@/_2pc/@{-->}[rr]_{#6}="x" & {#3} \ar[r]^-{#4} & {#5}
\save {{[].[ll]}."f"}."x"
*[F]\frm{} \restore
}

\def\Pbtwv#1#2#3#4#5#6
{
{#1} \ar[r]^-{#2}="f" \ar@/_2pc/@{-->}[rr]_{#6}="x" & {#3} \ar[r]^-{#4} & {#5}
\save {{[].[ll]}."f"}."x"
*[F.]\frm{} \restore
}

\def\Pbtwshv#1#2#3#4#5#6
{
{#1\grsmo} \ar[r]^-{#2}="f" \ar@/_1.5pc/@{-->}[rr]_{#6}="x" & {#3\gszto} \ar[r]^-{#4} &
{#5\grso}
\save {{[].[ll]}."f"}."x"
*[F.]\frm{} \restore
}

\def\Cbtwshv#1#2#3#4#5#6
{
{#1\grsmo} \ar[r]^-{#2}="f" \ar@/_1.5pc/@{-->}[rr]_{#6}="x" & {#3\gszto} \ar[r]^-{#4} &
{#5\grso}
\save {{[].[ll]}."f"}."x"
*[F]\frm{} \restore
}

\def\Cbtwshuv#1#2#3#4#5#6
{
{#1\grsmo} \ar[r]^-{#2}="f" \ar@/_2pc/@{-->}[rr]_{#6}="x" & {#3\gszto} \ar[r]^-{#4} &
{#5\grso}
\save {{[].[ll]}."f"}."x"
*[F]\frm{} \restore
}

\def\cxy#1{\vcenter{ \xymatrix{#1} } }

\def\Pssdl{ \Pbtwshv{\hlpar}{\wF}{\hlvir}{\wG}{\hlhor}{\wX} }
\def\Pssdlh{ \Pbtwshv{\hlhor}{\wFp}{\hlvir}{\wGp}{\hlpar}{\wXp} }

\def\Psgmp{ \Pbtwshv{\gg}{\yF}{\ggp}{\yG}{\ggpp}{\wX} }
\def\Cnsdl{ \Cbtwshv{\hlpar}{\wF}{\hlvir}{\wG}{\hlhor}{\wX} }
\def\Cnsdr{ \Cbtwshv{\hlhor}{\wFp}{\hlvir}{\wGp}{\hlpar}{\wXp} }
\def\Cnsdrprp{ \Cbtwshv{\hlhorp}{\yFprp}{\hlvirp}{\yGprp}{\hlparp}{\yXprp} }
\def\Pssdrprp{ \Pbtwshv{\hlhorp}{\yFprp}{\hlvirp}{\yGprp}{\hlparp}{\yXprp} }

\def\Cbtv#1#2#3#4#5
{
{#1} \ar[r]^-{#2}="f" & {#3} \ar[r]^-{#4} & {#5}
\save {[].[ll]}."f"
*[F]\frm{} \restore
}

\def\Cbttv#1#2#3#4
{
{#1} \ar[r]^-{#2}="f" & {#3}  & {#4}
\save {[].[ll]}."f"
*[F]\frm{} \restore
}

\def\Cbtov#1#2#3#4
{
{#1}  & {#2} \ar[r]^-{#3}="f" & {#4}
\save {[].[ll]}."f"
*[F]\frm{} \restore
}

\def\Cbtvsht#1#2#3#4{ \Cbttv{#1\grsmo}{#2}{#3\gszto}{#4\grso} }

\def\CbtwABCX{ \Cbtwv{\mA}{\wF}{\mB\gszto}{\wG}{\mC}{\wX} }
\def\PbtwABCX{ \Pbtwv{\mA}{\wF}{\mB\gszto}{\wG}{\mC}{\wX} }

\def\Cbtwvsh#1#2#3#4#5#6{ \Cbtwv{#1\grsmo}{#2}{#3\gszto}{#4}{#5\grso}{#6} }
\def\Cbtwgny{
\Cbtwvsh{\hlpary}{\yF}{\hlviry}{\yG}{\hlhory}{\yX} }

\def\Cbtany{
\Cbtvsht{\hlhcpy\otimes\JrWp}{\yFm}{\hlhcpy}{\hlhcpy}
}

\def\smmpver{ \lrbc{\! \begin{smallmatrix} 0\\0\\ \id
\end{smallmatrix} \! } }

\def\smmphor{ \lrbc{\! \begin{smallmatrix} \id& 0& 0
\end{smallmatrix} \! } }

\def\stphor{ \lrbc{\! \begin{smallmatrix} \id& 0
\end{smallmatrix} \! } }

\def\vcirc{ \circ }

\def\aPsk#1{ \mathord{\nP\lrbs{\,\raise12pt\xybox{0;/r2pc/:#1}\,}
} }

\def\akinksk{
\xygraph{
!{0;/r1pc/:}
[d(0.5)]
!{\hover}
!{\hcap}
[ld]!{\xcapv@(0)}
[uuu]!{\xcapv@(0)}
}
}

\def\akinkskn{
\xygraph{
!{0;/r1pc/:}
[d(0.5)]
!{\hunder}
!{\hcap}
[ld]!{\xcapv@(0)}
[uuu]!{\xcapv@(0)}
}
}

\def\avlinsk{
\xygraph{
!{0;/r1pc/:}
!{\xcapv[2]@(0)}
}
}

\def\bvertsk{
\xygraph{
!{0;/r1pc/:}
[d(1)l(1.8)]
{\zbendv[2]@(0)}
[u]
{\sbendv[2]@(0)}
[d][r(0.2)]
*{\bullet}
}
}

\def\ipar{
\xygraph{
!{0;/r0.6pc/:}
[u(0.3)l(1.3)]
{\xunoverv}
}
}

\def\ihor{
\xygraph{
!{0;/r0.6pc/:}
[u(0.3)l(1.3)]
{\xunoverh}
}
}

\def\iver{
\xygraph{
!{0;/r0.3pc/:}
[d(0.5)l(1.8)]
{\zbendv[2]@(0)}
[u]
{\sbendv[2]@(0)}
[d(0.6)][r(0.2)]
*{\bullet}
}
}

\def\ivir{
\xygraph{
!{0;/r0.3pc/:}
[d(0.5)l(1.8)]
{\zbendv[2]@(0)}
[u]
{\sbendv[2]@(0)}
[d(0.6)][r(0.2)]
}
}

\def\ivir{
\xygraph{
!{0;/r0.3pc/:}
[d(0.5)l(1.8)]
{\zbendv[2]@(0)}
[u]
{\sbendv[2]@(0)}
[d(0.6)][r(0.2)]
*{\vcirc}
}
}

\def\ixv#1{ \!\!\!\!\xybox{0;/r1pc/:#1}\,\, }
\def\ixpar{ \ixv{\ipar} }
\def\ixver{ \ixv{\iver} }
\def\ixhor{ \ixv{\ihor} }
\def\ixvir{ \ixv{\ivir} }

\def\lpcr{
\xygraph{
!{0;/r1pc/:}
[l(1.2)u(0.25)]
{\xoverv}
}
}

\def\lncr{
\xygraph{
!{0;/r1pc/:}
[l(1.2)u(0.25)]
{\xunderv}
}
}

\def\lhor{
\xygraph{
!{0;/r1pc/:}
[l(1.2)u(0.25)]
{\xunoverh}
}
}

\def\lverx{
\xygraph{
!{0;/r0.5pc/:}
[l(1.4)d(0.5)]
{\zbendv[2]@(0)}
[u]
{\sbendv[2]@(0)}
[d(0.6)][r(0.65)]
*{\bullet}
}
}

\def\lvrbx{
\xygraph{
!{0;/r0.5pc/:}
[l(1.4)d(0.5)]
{\zbendv[2]@(0)}
[u]
{\sbendv[2]@(0)}
[d(0.6)][r(0.65)]
*{\bullet}
[u(1)]
*{.}
[u(-2)]
*{.}
}
}

\def\lvbhx{
\xygraph{
!{0;/r0.5pc/:}
[l(1.4)d(0.5)]
{\zbendv[2]@(0)}
[u]
{\sbendv[2]@(0)}
[d(0.6)][r(0.65)]
*{\bullet}
[l(1)]
*{.}
[r(2)]
*{.}
}
}

\def\lvirx{
\xygraph{
!{0;/r0.5pc/:}
[l(1.4)d(0.5)]
{\zbendv[2]@(0)}
[u]
{\sbendv[2]@(0)}
[d(0.6)][r(0.65)]
}
}

\def\lvirx{
\xygraph{
!{0;/r0.5pc/:}
[l(1.4)d(0.5)]
{\zbendv[2]@(0)}
[u]
{\sbendv[2]@(0)}
[d(0.56)][r(0.61)]
*{\vcirc}
}
}

\def\loarx{
\xygraph{
!{0;/r0.5pc/:}
[l(1.4)d(0.5)]
[u(0.5)]{\xcaph[2]@(0)}
[d(0.6)][r(0.65)]
}
}

\def\loarxvv#1#2{
\xygraph{
!{0;/r0.5pc/:}
[l(1.4)d(0.5)]
[u(0.5)]{\xcaph[2]@(0)<{#1}>{#2}}
[d(0.6)][r(0.65)]
}
}

\def\loarxvv#1#2{ {\scriptstyle #1}\lxoar {\scriptstyle #2} }

\def\lpar{
\xygraph{
!{0;/r1pc/:}
[l(1.2)u(0.25)]
{\xunoverv}
}
}

\def\lxv#1{\xybox{0;/r1pc/:(0,0)*{#1},(-0.7,0)*{},(0.9,0)*{},(0,0.7)*{},(0,-1)*{}}
}

\def\lxpar{ \lxv{\lpar} }
\def\lxhor{ \lxv{\lhor} }
\def\lxver{ \lxv{\lverx} }
\def\lxvrb{ \lxv{\lvrbx} }
\def\lxvbh{ \lxv{\lvbhx} }
\def\lxvir{ \lxv{\lvirx} }
\def\lxpcr{ \lxv{\lpcr} }
\def\lxncr{ \lxv{\lncr} }
\def\lxoar{ \lxv{\loarx} }

\def\lxoarvv#1#2{ {\scriptstyle #1} \lxoar {\scriptstyle #2} }
\def\hloarvv#1#2{ \widehat{\lxoarvv{#1}{#2}} }

\def\loarxins{
\xygraph{
!{0;/r0.5pc/:}
[l(1.4)d(0.5)]
[u(0.5)]{\xcaph[2]@(0)<>><}
[d(0.6)][r(0.65)]
}
}
\def\lxoarins{ \lxv{\loarxins} }
\def\hloarins{ \mfgrv{\lxoarins} }
\def\lxoarinsvv#1#2{ {\scriptstyle #1} \lxoarins {\scriptstyle #2}
}
\def\hloarinsvv#1#2{ \widehat{\lxoarinsvv{#1}{#2}} }

\def\lzcir{\xybox{0;/r1pc/:(0,0)*{\lcirz},(-1.8,0.8)*{}
}
}

\def\lzcir{ \bigcirc }

\def\hlzcir{ \mfgrv{\lzcir} }
\def\nCcir{ \nCv{\lzcir} }
\def\nCcirgv#1{ \nC_{\xunkn,\grxhv{#1} } }
\def\nCcirgmN{ \nCcirgv{-2\xN} }
\def\nCscir{ \nCcir^\ast }

\def\xunkn{ \mathrm{unkn} }
\def\nCcir{ \nC_{\xunkn} }
\def\hlhrcir{ \hlhor\otimes\JrWp }

\def\arnCcir{ \hlhcpy\otimes\JrWp }

\def\mfgrv#1{ \widehat{#1} }
\def\hlpcr{ \mfgrv{\lxpcr} }
\def\hlncr{ \mfgrv{\lxncr} }
\def\hlhor{ \mfgrv{\lxhor} }
\def\hlver{ \mfgrv{\lxver} }
\def\hlvir{ \mfgrv{\lxvir} }
\def\hlpar{ \mfgrv{\lxpar} }
\def\hloar{ \mfgrv{\lxoar} }

\def\psgrv#1{ \widetilde{#1} }
\def\tlvrb{ \psgrv{\lxvrb} }
\def\tlvbh{ \psgrv{\lxvbh} }
\def\tlncr{ \psgrv{\lxncr} }
\def\tlpcr{ \psgrv{\lxpcr} }

\def\hlvirsh{ \hlvir\gszto }

\def\hlpv#1{ #1_{\mprp} }
\def\hlparp{ \hlpv{\hlpar} }
\def\hlvirp{ \hlpv{\hlvir} }
\def\hlhorp{ \hlpv{\hlhor} }

\def\hlpartp{ \tprpv{\hlpar} }

\def\tlvrbp{ \prpv{\tlvrb} }
\def\tlvbhp{ \prpv{\tlvbh} }

\def\lyv#1{\xybox{0;/r1pc/:(0,0)*{#1},(-0.7,0)*{},(1.3,0)*{},(0,0.7)*{},(0,-1)*{}}
}

\def\lhory{
\xygraph{
!{0;/r1pc/:}
[l(1.2)u(0.25)]
{\huntwist}
[r(0.5)]
{\hcap}
}
}

\def\lpary{
\xygraph{
!{0;/r1pc/:}
[l(1.2)u(0.25)]
{\huncross}
[r(0.5)]
{\hcap}
}
}

\def\lvery{
\xygraph{
!{0;/r0.5pc/:}
[l(1.5)d(0.5)]
[u]
{\zbendh[2]}
[d]
{\sbendh[2]}
[r(1)][u]
{\hcap[2]}
[d(0.5)][l(0.5)]
*{\bullet}
}
}

\def\lviry{
\xygraph{
!{0;/r0.5pc/:}
[l(1.5)d(0.5)]
[u]
{\zbendh[2]}
[d]
{\sbendh[2]}
[r(1)][u]
{\hcap[2]}
[d(0.5)][l(0.5)]
}
}

\def\lviry{
\xygraph{
!{0;/r0.5pc/:}
[l(1.5)d(0.5)]
[u]
{\zbendh[2]}
[d]
{\sbendh[2]}
[r(1)][u]
{\hcap[2]}
[d(0.5)][l(0.5)]
*{\vcirc}
}
}

\def\ltwpy{
\xygraph{
!{0;/r0.5pc/:}
[l(1.5)d(0.5)]
[u]
{\htwist[2]}
[d]
[r(1)][u]
{\hcap[2]}
[d(0.5)][l(0.5)]
}
}

\def\lxhory{ \lyv{\lhory} }
\def\hlhory{ \mfgrv{\lxhory} }

\def\lxpary{ \lyv{\lpary} }
\def\hlpary{ \mfgrv{\lxpary} }

\def\lxvery{ \lyv{\lvery} }
\def\hlvery{ \mfgrv{\lxvery} }

\def\lxviry{ \lyv{\lviry} }
\def\hlviry{ \mfgrv{\lxviry} }

\def\lxtwpy{ \lyv{\ltwpy} }
\def\hltwpy{ \mfgrv{\lxtwpy} }

\def\hlhorz{ \mfgrv{\lxhorz} }
\def\hlparz{ \mfgrv{\lxparz} }
\def\hlvirz{ \mfgrv{\lxvirz} }

\def\lparz{
\xygraph{
!{0;/r1pc/:}
[l(1.1)u(0.25)]
{\huncross}
[r(0.5)]
{\hcap}
[l(0.5)]
{\hcap-}
}
}

\def\lvirz{
\xygraph{
!{0;/r0.5pc/:}
[l(1.2)d(0.5)]
[u]
{\zbendh[2]}
[d]
{\sbendh[2]}
[r(1)][u]
{\hcap[2]}
[l]
{\hcap[-2]}
[d(0.5)][l(0.5)]
}
}

\def\lvirz{
\xygraph{
!{0;/r0.5pc/:}
[l(1.2)d(0.5)]
[u]
{\zbendh[2]}
[d]
{\sbendh[2]}
[r(1)][u]
{\hcap[2]}
[l]
{\hcap[-2]}
[d(0.5)][r(0.8)]
*{\vcirc}
}
}

\def\lhorz{
\xygraph{
!{0;/r1pc/:}
[l(1.1)u(0.25)]
{\huntwist}
[r(0.5)]
{\hcap}
[l(0.5)]
{\hcap-}
}
}

\def\lzv#1{\xybox{0;/r1pc/:(0,0)*{#1},(-0.7,0)*{},(1.3,0)*{},(0,0.7)*{},(0,-1)*{}}
}

\def\lxparz{ \lzv{\lparz} }
\def\lxhorz{ \lzv{\lhorz} }
\def\lxvirz{ \lzv{\lvirz} }

\def\lycv#1{\xybox{0;/r1pc/:(0,0)*{#1},(-0.3,0)*{},(0.9,0)*{},(0,0.7)*{},(0,-1)*{}}
}

\def\lhcpy{
\xygraph{
!{0;/r1pc/:}
[l(0.95)u(0.25)]
{\hcap}
}
}

\def\lxhcpy{ \lycv{\lhcpy} }
\def\hlhcpy{ \mfgrv{\lxhcpy} }

\def\lycnv#1#2#3{\xybox{0;/r1pc/:(0,0)*{#1},(-0.6,0)*{},(0.9,0)*{},(0,0.7)*{},(0,-1)*{},
(-0.2,-0.6)*{\scriptstyle #2},
(-0.2,0.5)*{\scriptstyle #3}
}
}

\def\lhcpny{
\xygraph{
!{0;/r1pc/:}
[l(0.95)u(0.25)]
{\hcap}
}
}

\def\lxhcpothy{ \lycnv{\lhcpny}{1}{3} }
\def\hlhcpothy{ \mfgrv{\lxhcpothy} }

\def\rttw{
\xygraph{
!{0;/r2.0pc/:}
!{\vtwist}
!{\vtwistneg}
}
}

\def\rtutw{
\xygraph{
!{0;/r0.5pc/:}
[rrr]
!{\xcapv[2]@(0)}
[d]
!{\xbendd-}
[l]!{\xbendd}
[l]!{\xbendu}
[ld]!{\xbendu-}
[ld]!{\xcapv[2]@(0)}
[rrrruuuuuuu]!{\xcapv[2]@(0)}
[dl]!{\xbendu}
[dl]!{\xbendu-}
[dl]!{\xbendd-}
[l]!{\xbendd}
!{\xcapv[2]@(0)}
}
}

\def\bgrtbx#1{ \widehat{
\xybox{0;/r1pc/:
(-4,0)*{#1},
(-2,-2)*{},
(2,0)*{}
}
}
}

\def\bgrtbx#1{
\mspace{-20mu}
\widehat{
\raise22pt\xybox{0;/r1pc/:
(-5,0)*{#1},
(0,-4.5)*{},
(1,0.5)*{}
}
}\mspace{-20mu}
}

\def\bgrtbwx#1{
\mspace{-20mu}
\raise24pt\xybox{0;/r1pc/:
(-5,0)*{#1},
(0,-4.5)*{},
(1,0.2)*{}
}
\mspace{-20mu}
}

\def\hbrttw{ \bgrtbx{\rttw} }
\def\hbrtutw{ \bgrtbx{\rtutw} }

\def\hbrttww{ \bgrtbwx{\rttw} }
\def\hbrtutww{ \bgrtbwx{\rtutw} }

\def\lrtho{
\xygraph{
!{0;/r2pc/:}
[l][u(1.5)]
!{\vcross}
[rru]!{\xcapv@(0)}
[l]!{\vcross}
[lu]!{\xcapv@(0)}
!{\vcross}
[rru]!{\xcapv@(0)}
}
}

\def\lvrtho{
\xygraph{
!{0;/r2pc/:}
[l][u(1.5)]
!{\sbendv}
[l]!{\zbendv}[ld]
[rru]!{\xcapv@(0)}
[l]!{\sbendv}
[l]!{\zbendv}[ld]
[lu]!{\xcapv@(0)}
!{\vcross}
[rru]!{\xcapv@(0)}
[u(1.5)][l(0.5)]*{\vcirc}
[u(1)][l]*{\vcirc}
}
}

\def\lrtht{
\xygraph{
!{0;/r2pc/:}
[u(1.5)]
!{\vcross}
[r]!{\xcapv@(0)}
[l]!{\vcross}
[uul]!{\vcross}
!{\xcapv@(0)}
[uuu]!{\xcapv@(0)}
}
}

\def\lvrtht{
\xygraph{
!{0;/r2pc/:}
[u(1.5)]
!{\vcross}
[r]!{\xcapv@(0)}
[l]
!{\sbendv}
[l]!{\zbendv}[ld]
[uul]
!{\sbendv}
[l]!{\zbendv}[ld]
!{\xcapv@(0)}
[uuu]!{\xcapv@(0)}
[d(0.5)][r(0.5)]*{\vcirc}
[d][r]*{\vcirc}
}
}

\def\xxrv#1{
\xybox{
0;/r1pc/:
(-0.5,0)*{#1},
(-1,-3)*{},
(4,2)*{}
}
}

\def\xvxrv#1{
\xybox{
0;/r1pc/:
(-0.5,0)*{#1},
(-2,-4)*{},
(5,3)*{}
}
}

\def\xrtho{ \xxrv{\lrtho} }
\def\xrtht{ \xxrv{\lrtht} }

\def\xvrtho{ \xvxrv{\lvrtho} }
\def\xvrtht{ \xvxrv{\lvrtht} }

\def\smgrthbv#1{
\raise11pt\xybox{0;/r1pc/:
(0,0)*{#1},
(-0.7,-2)*{},
(0.7,0.2)*{}
}
}

\def\tgrtbv#1{
\widetilde{
\raise11pt\xybox{0;/r1pc/:
(0,0)*{#1},
(-0.7,-2)*{},
(0.7,0.2)*{}
}
}
}

\def\bgrtbv#1{
\widehat{
\raise11pt\xybox{0;/r1pc/:
(0,0)*{#1},
(-0.7,-2)*{},
(0.7,0.2)*{}
}
}
}

\def\tbgrtbv#1{
\widetilde{
\raise11pt\xybox{0;/r1pc/:
(0,0)*{#1},
(-0.7,-2)*{},
(0.7,0.2)*{}
}
}
}

\def\abgrtbv#1{
{
\raise11pt\xybox{0;/r1pc/:
(0,0)*{#1},
(-0.7,-2)*{},
(0.7,0.2)*{}
}
}
}

\def\lrttw{
\xygraph{
!{0;/r1.0pc/:}
[l(1.5)]
!{\vtwist}
!{\vtwistneg}
}
}

\def\parhor{
\xygraph{
!{0;/r1.0pc/:}
[l(1.5)]
!{\vuntwist}
!{\vuncross}
}
}

\def\parvir{
\xygraph{
!{0;/r1.0pc/:}
[l(1.5)]
!{\vuntwist}
!{\sbendv}
[l]!{\zbendv}
}
}

\def\parvir{
\xygraph{
!{0;/r1.0pc/:}
[l(1.5)]
!{\vuntwist}
!{\sbendv}
[l]!{\zbendv}
[d(0.5)l(0.5)]{\vcirc}
}
}

\def\parpar{
\xygraph{
!{0;/r1.0pc/:}
[l(1.5)]
!{\vuntwist}
!{\vuntwist}
}
}

\def\horhor{
\xygraph{
!{0;/r1.0pc/:}
[l(1.5)]
!{\vuncross}
!{\vuncross}
}
}

\def\horvir{
\xygraph{
!{0;/r1.0pc/:}
[l(1.5)]
!{\vuncross}
!{\sbendv}
[l]!{\zbendv}
}
}

\def\horvir{
\xygraph{
!{0;/r1.0pc/:}
[l(1.5)]
!{\vuncross}
!{\sbendv}
[l]!{\zbendv}
[d(0.5)l(0.5)]{\vcirc}
}
}

\def\horpar{
\xygraph{
!{0;/r1.0pc/:}
[l(1.5)]
!{\vuncross}
!{\vuntwist}
}
}

\def\virhor{
\xygraph{
!{0;/r1.0pc/:}
[l(1.5)]
!{\sbendv}
[l]!{\zbendv}
[dl]
!{\vuncross}
}
}

\def\virhor{
\xygraph{
!{0;/r1.0pc/:}
[l(1.5)]
!{\sbendv}
[l]!{\zbendv}
[d(0.5)l(0.5)]{\vcirc}[u(0.5)r(0.5)]
[dl]
!{\vuncross}
}
}

\def\virpar{
\xygraph{
!{0;/r1.0pc/:}
[l(1.5)]
!{\sbendv}
[l]!{\zbendv}
[dl]
!{\vuntwist}
}
}

\def\virpar{
\xygraph{
!{0;/r1.0pc/:}
[l(1.5)]
!{\sbendv}
[l]!{\zbendv}
[d(0.5)l(0.5)]{\vcirc}[u(0.5)r(0.5)]
[dl]
!{\vuntwist}
}
}

\def\virvir{
\xygraph{
!{0;/r1.0pc/:}
[l(1.5)]
!{\sbendv}
[l]!{\zbendv}
[dl]
!{\sbendv}
[l]!{\zbendv}
}
}

\def\virvir{
\xygraph{
!{0;/r1.0pc/:}
[l(1.5)]
!{\sbendv}
[l]!{\zbendv}
[d(0.5)l(0.5)]{\vcirc}[u(0.5)r(0.5)]
[dl]
!{\sbendv}
[l]!{\zbendv}
[d(0.5)l(0.5)]{\vcirc}
}
}

\def\verver{
\xygraph{
!{0;/r1.0pc/:}
[l(1.5)]
!{\sbendv}
[l]!{\zbendv}
[d(0.5)l(0.5)]{\bullet}[u(0.5)r(0.5)]
[dl]
!{\sbendv}
[l]!{\zbendv}
[d(0.5)l(0.5)]{\bullet}
}
}

\def\hparhor{ \bgrtbv{\parhor} }

\def\hparvir{ \bgrtbv{\parvir} }
\def\hparpar{ \bgrtbv{\parpar} }

\def\hhorhor{ \bgrtbv{\horhor} }

\def\hhorvir{ \bgrtbv{\horvir} }
\def\hhorpar{ \bgrtbv{\horpar} }

\def\hvirhor{ \bgrtbv{\virhor} }

\def\hvirvir{ \bgrtbv{\virvir} }
\def\hvirpar{ \bgrtbv{\virpar} }

\def\hverver{ \bgrtbv{\verver} }

\def\hrttw{ \bgrtbv{\lrttw} }

\def\arttw{ \abgrtbv{\lrttw} }

\def\xsmgrtbv#1{
\raise11pt\xybox{0;/r0.3pc/:
(0,0)*{#1},
(-1.7,-2.8)*{},
(1.7,0.5)*{}
}
}

\def\xgrtbv#1{
\raise11pt\xybox{0;/r1pc/:
(0,0)*{#1},
(-1.7,-2.8)*{},
(1.7,0.5)*{}
}
}

\def\lhpfw{
\xygraph{
!{0;/r1.0pc/:}
[l(2.5)]
[u(0.1)]
!{\vcap}
[rr]!{\vcap}
[ll]!{\xcapv[2]@(0)}
[rrru]!{\xcapv[2]@(0)}
[ld]!{\vcap[-1]}
[ll]!{\vcap[-1]}
[ruu]!{\vtwistneg}
!{\vtwistneg}
}
}

\def\cparhor{
\xygraph{
!{0;/r1.0pc/:}
[l(2.5)]
[u(0.1)]
!{\vcap}
[rr]!{\vcap}
[ll]!{\xcapv[2]@(0)}
[rrru]!{\xcapv[2]@(0)}
[ld]!{\vcap[-1]}
[ll]!{\vcap[-1]}
[u(2)r(2.5)]
[l(1.5)]
!{\vuntwist}
!{\vuncross}
}
}

\def\cparver{
\xygraph{
!{0;/r1.0pc/:}
[l(2.5)]
[u(0.1)]
!{\vcap}
[rr]!{\vcap}
[ll]!{\xcapv[2]@(0)}
[rrru]!{\xcapv[2]@(0)}
[ld]!{\vcap[-1]}
[ll]!{\vcap[-1]}
[u(2)r(2.5)]
[l(1.5)]
!{\vuntwist}
!{\sbendv}
[l]!{\zbendv}
[d(0.5)l(0.5)]{\bullet}
}
}

\def\cparvir{
\xygraph{
!{0;/r1.0pc/:}
[l(2.5)]
[u(0.1)]
!{\vcap}
[rr]!{\vcap}
[ll]!{\xcapv[2]@(0)}
[rrru]!{\xcapv[2]@(0)}
[ld]!{\vcap[-1]}
[ll]!{\vcap[-1]}
[u(2)r(2.5)]
[l(1.5)]
!{\vuntwist}
!{\sbendv}
[l]!{\zbendv}
}
}

\def\cparvir{
\xygraph{
!{0;/r1.0pc/:}
[l(2.5)]
[u(0.1)]
!{\vcap}
[rr]!{\vcap}
[ll]!{\xcapv[2]@(0)}
[rrru]!{\xcapv[2]@(0)}
[ld]!{\vcap[-1]}
[ll]!{\vcap[-1]}
[u(2)r(2.5)]
[l(1.5)]
!{\vuntwist}
!{\sbendv}
[l]!{\zbendv}
[d(0.52)l(0.5)]{\vcirc}
}
}

\def\cparpar{
\xygraph{
!{0;/r1.0pc/:}
[l(2.5)]
[u(0.1)]
!{\vcap}
[rr]!{\vcap}
[ll]!{\xcapv[2]@(0)}
[rrru]!{\xcapv[2]@(0)}
[ld]!{\vcap[-1]}
[ll]!{\vcap[-1]}
[u(2)r(2.5)]
[l(1.5)]
!{\vuntwist}
!{\vuntwist}
}
}

\def\chorhor{
\xygraph{
!{0;/r1.0pc/:}
[l(2.5)]
[u(0.1)]
!{\vcap}
[rr]!{\vcap}
[ll]!{\xcapv[2]@(0)}
[rrru]!{\xcapv[2]@(0)}
[ld]!{\vcap[-1]}
[ll]!{\vcap[-1]}
[u(2)r(2.5)]
[l(1.5)]
!{\vuncross}
!{\vuncross}
}
}

\def\chorvir{
\xygraph{
!{0;/r1.0pc/:}
[l(2.5)]
[u(0.1)]
!{\vcap}
[rr]!{\vcap}
[ll]!{\xcapv[2]@(0)}
[rrru]!{\xcapv[2]@(0)}
[ld]!{\vcap[-1]}
[ll]!{\vcap[-1]}
[u(2)r(2.5)]
[l(1.5)]
!{\vuncross}
!{\sbendv}
[l]!{\zbendv}
}
}

\def\chorvir{
\xygraph{
!{0;/r1.0pc/:}
[l(2.5)]
[u(0.1)]
!{\vcap}
[rr]!{\vcap}
[ll]!{\xcapv[2]@(0)}
[rrru]!{\xcapv[2]@(0)}
[ld]!{\vcap[-1]}
[ll]!{\vcap[-1]}
[u(2)r(2.5)]
[l(1.5)]
!{\vuncross}
!{\sbendv}
[l]!{\zbendv}
[d(0.52)l(0.5)]{\vcirc}
}
}

\def\chorpar{
\xygraph{
!{0;/r1.0pc/:}
[l(2.5)]
[u(0.1)]
!{\vcap}
[rr]!{\vcap}
[ll]!{\xcapv[2]@(0)}
[rrru]!{\xcapv[2]@(0)}
[ld]!{\vcap[-1]}
[ll]!{\vcap[-1]}
[u(2)r(2.5)]
[l(1.5)]
!{\vuncross}
!{\vuntwist}
}
}

\def\cvirhor{
\xygraph{
!{0;/r1.0pc/:}
[l(2.5)]
[u(0.1)]
!{\vcap}
[rr]!{\vcap}
[ll]!{\xcapv[2]@(0)}
[rrru]!{\xcapv[2]@(0)}
[ld]!{\vcap[-1]}
[ll]!{\vcap[-1]}
[u(2)r(2.5)]
[l(1.5)]
!{\sbendv}
[l]!{\zbendv}
[dl]
!{\vuncross}
}
}

\def\cvirhor{
\xygraph{
!{0;/r1.0pc/:}
[l(2.5)]
[u(0.1)]
!{\vcap}
[rr]!{\vcap}
[ll]!{\xcapv[2]@(0)}
[rrru]!{\xcapv[2]@(0)}
[ld]!{\vcap[-1]}
[ll]!{\vcap[-1]}
[u(2)r(2.5)]
[l(1.5)]
!{\sbendv}
[l]!{\zbendv}
[d(0.52)l(0.5)]{\vcirc}[u(0.52)r(0.5)]
[dl]
!{\vuncross}
}
}

\def\cvirpar{
\xygraph{
!{0;/r1.0pc/:}
[l(2.5)]
[u(0.1)]
!{\vcap}
[rr]!{\vcap}
[ll]!{\xcapv[2]@(0)}
[rrru]!{\xcapv[2]@(0)}
[ld]!{\vcap[-1]}
[ll]!{\vcap[-1]}
[u(2)r(2.5)]
[l(1.5)]
!{\sbendv}
[l]!{\zbendv}
[dl]
!{\vuntwist}
}
}

\def\cvirpar{
\xygraph{
!{0;/r1.0pc/:}
[l(2.5)]
[u(0.1)]
!{\vcap}
[rr]!{\vcap}
[ll]!{\xcapv[2]@(0)}
[rrru]!{\xcapv[2]@(0)}
[ld]!{\vcap[-1]}
[ll]!{\vcap[-1]}
[u(2)r(2.5)]
[l(1.5)]
!{\sbendv}
[l]!{\zbendv}
[d(0.52)l(0.5)]{\vcirc}[u(0.52)r(0.5)]
[dl]
!{\vuntwist}
}
}

\def\cvirvir{
\xygraph{
!{0;/r1.0pc/:}
[l(2.5)]
[u(0.1)]
!{\vcap}
[rr]!{\vcap}
[ll]!{\xcapv[2]@(0)}
[rrru]!{\xcapv[2]@(0)}
[ld]!{\vcap[-1]}
[ll]!{\vcap[-1]}
[u(2)r(2.5)]
[l(1.5)]
!{\sbendv}
[l]!{\zbendv}
[dl]
!{\sbendv}
[l]!{\zbendv}
}
}

\def\cvirvir{
\xygraph{
!{0;/r1.0pc/:}
[l(2.5)]
[u(0.1)]
!{\vcap}
[rr]!{\vcap}
[ll]!{\xcapv[2]@(0)}
[rrru]!{\xcapv[2]@(0)}
[ld]!{\vcap[-1]}
[ll]!{\vcap[-1]}
[u(2)r(2.5)]
[l(1.5)]
!{\sbendv}
[l]!{\zbendv}
[d(0.52)l(0.5)]{\vcirc}[u(0.52)r(0.5)]
[dl]
!{\sbendv}
[l]!{\zbendv}
[d(0.52)l(0.5)]{\vcirc}
}
}

\def\cverhor{
\xygraph{
!{0;/r1.0pc/:}
[l(2.5)]
[u(0.1)]
!{\vcap}
[rr]!{\vcap}
[ll]!{\xcapv[2]@(0)}
[rrru]!{\xcapv[2]@(0)}
[ld]!{\vcap[-1]}
[ll]!{\vcap[-1]}
[u(2)r(2.5)]
[l(1.5)]
!{\sbendv}
[l]!{\zbendv}
[d(0.5)l(0.5)]{\bullet}[u(0.5)r(0.5)]
[dl]
!{\vuncross}
}
}

\def\cverpar{
\xygraph{
!{0;/r1.0pc/:}
[l(2.5)]
[u(0.1)]
!{\vcap}
[rr]!{\vcap}
[ll]!{\xcapv[2]@(0)}
[rrru]!{\xcapv[2]@(0)}
[ld]!{\vcap[-1]}
[ll]!{\vcap[-1]}
[u(2)r(2.5)]
[l(1.5)]
!{\sbendv}
[l]!{\zbendv}
[d(0.5)l(0.5)]{\bullet}[u(0.5)r(0.5)]
[dl]
!{\vuntwist}
}
}

\def\cverver{
\xygraph{
!{0;/r1.0pc/:}
[l(2.5)]
[u(0.1)]
!{\vcap}
[rr]!{\vcap}
[ll]!{\xcapv[2]@(0)}
[rrru]!{\xcapv[2]@(0)}
[ld]!{\vcap[-1]}
[ll]!{\vcap[-1]}
[u(2)r(2.5)]
[l(1.5)]
!{\sbendv}
[l]!{\zbendv}
[d(0.5)l(0.5)]{\bullet}[u(0.5)r(0.5)]
[dl]
!{\sbendv}
[l]!{\zbendv}
[d(0.5)l(0.5)]{\bullet}
}
}

\def\cverhor{
\xygraph{
!{0;/r1.0pc/:}
[l(2.5)]
[u(0.1)]
!{\vcap}
[rr]!{\vcap}
[ll]!{\xcapv[2]@(0)}
[rrru]!{\xcapv[2]@(0)}
[ld]!{\vcap[-1]}
[ll]!{\vcap[-1]}
[u(2)r(2.5)]
[l(1.5)]
!{\sbendv}
[l]!{\zbendv}
[d(0.5)l(0.5)]{\bullet}[u(0.5)r(0.5)]
[dl]
!{\vuncross}
}
}

\def\xparhor{ \xgrtbv{\cparhor} }
\def\xparver{ \xgrtbv{\cparver} }
\def\xparvir{ \xgrtbv{\cparvir} }
\def\xparpar{ \xgrtbv{\cparpar} }

\def\xhorhor{ \xgrtbv{\chorhor} }

\def\xhorvir{ \xgrtbv{\chorvir} }
\def\xhorpar{ \xgrtbv{\chorpar} }

\def\xvirhor{ \xgrtbv{\cvirhor} }

\def\xvirvir{ \xgrtbv{\cvirvir} }
\def\xvirpar{ \xgrtbv{\cvirpar} }

\def\xverhor{ \xgrtbv{\cverhor} }
\def\xverver{ \xgrtbv{\cverver} }

\def\xverpar{ \xgrtbv{\cverpar} }

\def\Cparhor{ \nCv{\xparhor} }
\def\Cparver{ \nCv{\xparver} }

\def\Cparpar{ \nCv{\xparpar} }

\def\Chorhor{ \nCv{\xhorhor} }

\def\Chorpar{ \nCv{\xhorpar} }

\def\Cverhor{ \nCv{\xverhor} }
\def\Cverver{ \nCv{\xverver} }

\def\Cverpar{ \nCv{\xverpar} }

\def\wxgrtbv#1{
\widehat{
\raise11pt\xybox{0;/r1pc/:
(0,-1.5)*{#1},
(-1.7,-2.8)*{},
(1.7,0.5)*{}
}
}
}

\def\Chparhor{ \wxgrtbv{\xparhor} }

\def\Chparvir{ \wxgrtbv{\xparvir} }
\def\Chparpar{ \wxgrtbv{\xparpar} }

\def\Chhorhor{ \wxgrtbv{\xhorhor} }

\def\Chhorvir{ \wxgrtbv{\xhorvir} }
\def\Chhorpar{ \wxgrtbv{\xhorpar} }

\def\Chvirhor{ \wxgrtbv{\xvirhor} }

\def\Chvirvir{ \wxgrtbv{\xvirvir} }
\def\Chvirpar{ \wxgrtbv{\xvirpar} }

\def\xhpfw{ \xgrtbv{\lhpfw} }

\def\tgrtbv#1{
\widetilde{
\raise11pt\xybox{0;/r1pc/:
(0,0)*{#1},
(-1.3,-4)*{},
(1.3,1.2)*{}
}
}
}

\def\bverver{
\xygraph{
!{0;/r2.0pc/:}
[l(1.5)]
[u(0.55)]
!{\sbendv}
[l]!{\zbendv}
[d(0.5)l(0.5)]{\bullet}
[l(0.4)]*{.}
[r(0.8)]*{.}
[r(0.4)][l(0.8)]
[u(0.5)r(0.5)]
[dl]
!{\sbendv}
[l]!{\zbendv}
[d(0.5)l(0.5)]{\bullet}
[u(0.3)]*{.}
[d(0.6)]*{.}
[d(0.3)][u(0.6)]
}
}

\def\xtgrtbv#1{
\widetilde{
\raise11pt\xybox{0;/r1pc/:
(0,0)*{#1},
(-1.1,-3.8)*{},
(1.1,1)*{}
}
}
}

\def\xbverver{
\xygraph{
!{0;/r1.5pc/:}
[l(1.5)]
[u(0.55)]
!{\sbendv}
[l]!{\zbendv}
[d(0.5)l(0.5)]{\bullet}
[l(0.4)]*{.}
[r(0.8)]*{.}
[r(0.4)][l(0.8)]
[u(0.5)r(0.5)]
[dl]
!{\sbendv}
[l]!{\zbendv}
[d(0.5)l(0.5)]{\bullet}
[u(0.3)]*{.}
[d(0.6)]*{.}
[d(0.3)][u(0.6)]
}
}

\def\xbparver{
\xygraph{
!{0;/r1.5pc/:}
[l(1.5)]
[u(0.55)]
!{\vuntwist}
[ur]
[d(0.5)l(0.5)]
[l(0.4)]
[r(0.8)]
[r(0.4)][l(0.8)]
[u(0.5)r(0.5)]
[dl]
!{\sbendv}
[l]!{\zbendv}
[d(0.5)l(0.5)]{\bullet}
[u(0.3)]*{.}
[d(0.6)]*{.}
[d(0.3)][u(0.6)]
}
}

\def\xbhorver{
\xygraph{
!{0;/r1.5pc/:}
[l(1.5)]
[u(0.55)]
!{\vuncross}
[ur]
[d(0.5)l(0.5)]
[l(0.4)]
[r(0.8)]
[r(0.4)][l(0.8)]
[u(0.5)r(0.5)]
[dl]
!{\sbendv}
[l]!{\zbendv}
[d(0.5)l(0.5)]{\bullet}
[u(0.3)]*{.}
[d(0.6)]*{.}
[d(0.3)][u(0.6)]
}
}

\def\xbverpar{
\xygraph{
!{0;/r1.5pc/:}
[l(1.5)]
[u(0.55)]
!{\sbendv}
[l]!{\zbendv}
[d(0.5)l(0.5)]{\bullet}
[l(0.4)]*{.}
[r(0.8)]*{.}
[r(0.4)][l(0.8)]
[u(0.5)r(0.5)]
[dl]
!{\vuntwist}
[d(0.5)l(0.5)]
[u(0.3)]
[d(0.6)]
[d(0.3)][u(0.6)]
}
}

\def\xbverhor{
\xygraph{
!{0;/r1.5pc/:}
[l(1.5)]
[u(0.55)]
!{\sbendv}
[l]!{\zbendv}
[d(0.5)l(0.5)]{\bullet}
[l(0.4)]*{.}
[r(0.8)]*{.}
[r(0.4)][l(0.8)]
[u(0.5)r(0.5)]
[dl]
!{\vuncross}
[d(0.5)l(0.5)]
[u(0.3)]
[d(0.6)]
[d(0.3)][u(0.6)]
}
}

\def\xbparpar{
\xygraph{
!{0;/r1.5pc/:}
[l(1.5)]
[u(0.55)]
!{\vuntwist}
[ur]
[d(0.5)l(0.5)]
[l(0.4)]
[r(0.8)]
[r(0.4)][l(0.8)]
[u(0.5)r(0.5)]
[dl]
!{\vuntwist}
[d(0.5)l(0.5)]
[u(0.3)]
[d(0.6)]
[d(0.3)][u(0.6)]
}
}

\def\xbparhor{
\xygraph{
!{0;/r1.5pc/:}
[l(1.5)]
[u(0.55)]
!{\vuntwist}
[ur]
[d(0.5)l(0.5)]
[l(0.4)]
[r(0.8)]
[r(0.4)][l(0.8)]
[u(0.5)r(0.5)]
[dl]
!{\vuncross}
[d(0.5)l(0.5)]
[u(0.3)]
[d(0.6)]
[d(0.3)][u(0.6)]
}
}

\def\xbhorpar{
\xygraph{
!{0;/r1.5pc/:}
[l(1.5)]
[u(0.55)]
!{\vuncross}
[ur]
[d(0.5)l(0.5)]
[l(0.4)]
[r(0.8)]
[r(0.4)][l(0.8)]
[u(0.5)r(0.5)]
[dl]
!{\vuntwist}
[d(0.5)l(0.5)]
[u(0.3)]
[d(0.6)]
[d(0.3)][u(0.6)]
}
}

\def\xbhorhor{
\xygraph{
!{0;/r1.5pc/:}
[l(1.5)]
[u(0.55)]
!{\vuncross}
[ur]
[d(0.5)l(0.5)]
[l(0.4)]
[r(0.8)]
[r(0.4)][l(0.8)]
[u(0.5)r(0.5)]
[dl]
!{\vuncross}
}
}

\def\xbnegpos{
\xygraph{
!{0;/r1.5pc/:}
[l(1.5)]
[u(0.55)]
!{\vtwist}
[ur]
[d(0.5)l(0.5)]
[l(0.4)]
[r(0.8)]
[r(0.4)][l(0.8)]
[u(0.5)r(0.5)]
[dl]
!{\vtwistneg}
}
}

\def\xtgrtbv#1{
\widetilde{
\raise11pt\xybox{0;/r1pc/:
(0,0)*{#1},
(-1.1,-2.5)*{},
(1.1,1)*{}
}
}
}

\def\xtverver{ \xtgrtbv{\xbverver} }
\def\xtparver{ \xtgrtbv{\xbparver} }
\def\xthorver{ \xtgrtbv{\xbhorver} }
\def\xtverpar{ \xtgrtbv{\xbverpar} }
\def\xtverhor{ \xtgrtbv{\xbverhor} }

\def\xtnegpos{ \xtgrtbv{\xbnegpos} }

\def\xhgrtbv#1{
\widehat{
\raise11pt\xybox{0;/r1pc/:
(0,0)*{#1},
(-1.1,-2.5)*{},
(1.1,1)*{}
}
}
}

\def\xhparpar{ \xhgrtbv{\xbparpar} }
\def\xhparhor{ \xhgrtbv{\xbparhor} }
\def\xhhorpar{ \xhgrtbv{\xbhorpar} }
\def\xhhorhor{ \xhgrtbv{\xbhorhor} }

\def\xabx{ \mathbf{x} }
\def\xay{ y }
\def\xayv#1{ \xay_{#1} }
\def\xayo{ \xayv{1} }
\def\xayt{ \xayv{2} }

\def\xaQbx{ \IQ[\xabx] }
\def\xaQbxy{ \IQ[\xabx,\xay] }

\def\xaQxyot{ \IQ[\xabx,\xayo,\xayt] }

\def\xan{ \xnpr }
\def\xanv#1{ \xan_{#1} }
\def\xano{ \xanv{1} }
\def\xant{ \xanv{2} }

\def\xaW{ \xW }

\def\xaWb{ \bar{\xaW} }

\def\xaM{ M }
\def\xaMb{ \bar{M} }
\def\xaMbp{ \xaMb\p }
\def\xaMp{ \xaM\p }

\def\sbf{subfactorization}

\def\xap{ p }
\def\xaq{ q }
\def\xapv#1{ \xap_{#1} }
\def\xaqv#1{ \xaq_{#1} }
\def\xapo{ \xapv{1} }
\def\xapt{ \xapv{2} }
\def\xaqo{ \xaqv{1} }
\def\xaqt{ \xaqv{2} }

\def\ynp{ \xay^{\xnpr} - \xap }
\def\ynpq{ (\ynp)\,\xaq }

\def\ynpv#1{ \xayv{#1}^{\xanv{#1}}-\xapv{#1} }
\def\ynpo{ \ynpv{1} }
\def\ynpt{ \ynpv{2} }
\def\ynpqv#1{ (\ynpv{#1})\,\xaqv{#1} }
\def\ynpqo{ \ynpqv{1} }
\def\ynpqt{ \ynpqv{2} }

\def\xaotxy{ \otimes_{\xaQbxy} }
\def\xaotxyot{ \otimes_{\xaQxyot} }
\def\kmpynpq{ \kmfv{\ynp}{\xaq} }
\def\kmpynqpv#1{ \kmfv{\ynpv{#1}}{\xaqv{#1}} }
\def\kmpynqpo{ \kmpynqpv{1} }
\def\kmpynqpt{ \kmpynqpv{2} }

\def\xadgy{ \deg_{\xay} }

\def\xscr{ \scriptstyle }

\def\xspmttv#1#2#3#4
{\begin{pmatrix}
\xscr #1 & \xscr #2
\\
\xscr #3 & \xscr #4
\end{pmatrix}
}

\def\xspmthtv#1#2#3#4#5#6
{\begin{pmatrix}
\xscr #1 & \xscr #2 & \xscr #3
\\
\xscr #4 & \xscr #5 & \xscr #6
\end{pmatrix}
}

\def\xspmththv#1#2#3#4#5#6#7#8#9
{\begin{pmatrix}
\xscr #1 & \xscr #2 & \xscr #3
\\
\xscr #4 & \xscr #5 & \xscr #6
\\
\xscr #7 & \xscr #8 & \xscr #9
\end{pmatrix}
}

\def\xxcmn{ \mathrm{cmn} }

\def\id{ \mathrm{id} }
\def\idcmn{ \id_{\xxcmn} }

\def\xDp{ \xD\p }

\def\xaarl{ \save
[]+<-0.4cm,0cm>*i{BAC}
\ar@(ul,dl)_-{-\xD}
\restore }

\def\xaarr{
\save []+<0.1cm,0cm>*i{CAB}
\ar@(ur,dr)^-{\xD}
\restore }

\def\xaf{ f }
\def\xaC{ C }
\def\xaCp{ \xaC\p }
\def\xaA{ A }
\def\xaAp{ \xaA\p }
\def\xaAz{ \bar{\xaA} }
\def\xaB{ B }

\def\xmaAv#1{ \xaA_{#1} }

\def\xag{ g }
\def\xagp{ \xag\p }
\def\xah{ h }
\def\xahp{ \xah\p }

\def\xagv#1{ \xag_{#1} }
\def\xagA{ \xagv{\xaA} }
\def\xagB{ \xagv{\xaB} }
\def\xagpv#1{ \xagp_{#1} }
\def\xagpA{ \xagpv{\xaA} }
\def\xagpB{ \xagpv{\xaB} }
\def\xahv#1{ \xah_{#1} }

\def\xahB{ \xahv{\xaB} }
\def\xahpv#1{ \xahp_{#1} }
\def\xahpA{ \xahpv{\xaA} }

\def\zzd{ d }
\def\zzdp{ \zzd\p }

\def\cnmff{ \Cnvv{\MF}{\xaf} }
\def\cgcnf{ \xymatrix@1{ {\xaC} \ar[r]^-{\xag} & {\CnMFv{\xaf}}} }
\def\cgcnfp{ \xymatrix@1{ {\CnMFv{\xaf}} \ar[r]^-{\xagp} & {\xaCp}} }

\def\JrWpmz{ \JrWpmv{0} }
\def\JrWpmtN{ \JrWpmv{2\xN} }
\def\JrWpmztN{ \JrWpmv{0,2\xN} }

\def\JrWpv#1{ \xJv{\nW;#1}\p}
\def\JrWpmz{ \JrWpv{\hat{0} } }
\def\JrWpmtN{ \JrWpv{\widehat{2\xN}} }
\def\JrWpmztN{ \JrWpv{\hat{0},\widehat{2\xN}} }

\def\JrWpmz{ \JrWpv{-} }
\def\JrWpmtN{ \JrWpv{+} }
\def\JrWpmztN{ \JrWpv{\pm} }

\def\xdar{ \ar@{-->} }

\def\xunk{ \ast }

\def\nWitv#1#2{ \nW(#1,#2) }
\def\nWitxot{ \nWitv{\zxo}{\zxt} }
\def\kmfbv#1#2{ \tK\big(#1;#2\big) }

\def\nWithv#1#2#3{ \nW(#1,#2,#3) }

\def\zzxfm{ \zxth+\zxf-\zxo-\zxt }

\def\loarh{
\xygraph{
!{0;/r1pc/:}
[l(1.4)d(0.5)]
[u(0.5)]
!{\xcaph[1.5]@(0)=<}
!{\xcaph[1.5]@(0)}
[d(0.6)][r(0.65)]
}
}

\def\lxhv#1{\xybox{0;/r1pc/:(0,0)*{#1},(-0.7,0)*{},(2.7,0)*{},(0,0.7)*{},(0,-1)*{}}
}

\def\lxoarh{ \lxhv{\loarh} }

\def\loarhvv#1#2{ {\scriptstyle #1} \lxoarh {\scriptstyle #2} }

\def\bpar{
\xygraph{
!{0;/r2pc/:}
[l(1.2)u(0.25)]
!{\xunoverv=<}
}
}

\def\bpru{
\xygraph{
!{0;/r2pc/:}
[l(1.2)u(0.25)]
!{\xunoverv}
}
}

\def\bcrn{
\xygraph{
!{0;/r2pc/:}
[l(1.2)u(0.25)]
!{\xoverv=<}
}
}

\def\bcnw{
\xygraph{
!{0;/r2pc/:}
[l(1.2)u(0.25)]
!{\xoverv}
}
}

\def\bcrp{
\xygraph{
!{0;/r2pc/:}
[l(1.2)u(0.25)]
!{\xunderv=<}
}
}

\def\bcpw{
\xygraph{
!{0;/r2pc/:}
[l(1.2)u(0.25)]
!{\xunderv}
}
}

\def\bpro{
\xygraph{
!{0;/r2pc/:}
[l(1.2)u(0.25)]
!{\xcapv=<}
[ur]
!{\xcapv-=<}
}
}

\def\bhor{
\xygraph{
!{0;/r2pc/:}
[l(1.2)u(0.25)]
[u(0.25)]
!{\xbendr[0.5]=<}
[u][l(0.5)]
!{\xbendl[0.5]=<}
[l(0.5)][d(0.25)]
!{\xcaph}
}
}

\def\bhru{
\xygraph{
!{0;/r2pc/:}
[l(1.2)u(0.25)]
[u(0.25)]
!{\xbendr[0.5]}
[u][l(0.5)]
!{\xbendl[0.5]}
[l(0.5)][d(0.25)]
!{\xcaph}
}
}

\def\bhro{
\xygraph{
!{0;/r2pc/:}
[l(1.2)u(0.25)]
[u(0.25)]
!{\xbendr[0.5]}
[u][l(0.5)]
!{\xbendl[0.5]=<}
[l(0.5)][d(0.25)]
!{\xcaph=>}
}
}

\def\bhru{
\xygraph{
!{0;/r2pc/:}
[l(1.2)u(0.25)]
[u(0.25)]
!{\xbendr[0.5]}
[u][l(0.5)]
!{\xbendl[0.5]}
[l(0.5)][d(0.25)]
!{\xcaph}
}
}

\def\bvir{
\xygraph{
!{0;/r2pc/:}
[l(1.2)d(0.75)]
!{\zbendv[1]@(0)}
[l]
!{\sbendv[1]@(0)}
[l(0.5)][u(0.5)]{\vcirc}[r(0.5)][d(0.5)]
[d(0.6)][r(0.65)]
}
}

\def\lbv#1{\xybox{0;/r1pc/:(0,0.5)*{#1},(-0.7,0)*{},(2,0)*{},(0,1.5)*{},(0,-1.5)*{}}
}

\def\lbpar{
{\UseComputerModernTips
\lbv{\bpar} }
}

\def\lbpru{
{\UseComputerModernTips
\lbv{\bpru} }
}

\def\lbcrn{
{\UseComputerModernTips
\lbv{\bcrn} }
}

\def\lbcnw{
\lbv{\bcnw}
}

\def\lbcpw{
\lbv{\bcpw}
}

\def\lbcrp{
{\UseComputerModernTips
\lbv{\bcrp} }
}

\def\lbpro{
{\UseComputerModernTips
\lbv{\bpro} }
}

\def\lbhro{
{\UseComputerModernTips
\lbv{\bhro} }
}

\def\lbhru{
{\lbv{\bhru} }
}

\def\lbviu{
{\lbv{\bvir} }
}

\def\lbvir{
{\UseComputerModernTips
\xybox{0;/r1pc/:
(1.7,-1)*{}="a";
(-0.2,1)*{}="b";
{\ar"a";"b"},
(-0.2,-1)*{}="c";
(1.7,1)*{}="d";
{\ar"c";"d"},
(0.75,0)*{\vcirc},
(-0.7,0)*{},(2,0)*{},(0,1.5)*{},(0,-1.5)*{}}
}
}

\def\lbver{
{\UseComputerModernTips
\xybox{0;/r1pc/:
(0,0.5)*{\bhor}
,(-0.7,0)*{},(2,0)*{},(0,1.5)*{},(0,-1.5)*{},
(0.6,0.45)*{\bullet},
(0.6,-0.55)*{\bullet},
(0.6,0.45)*\cir<1pt>{} = "a";
(0.6,-0.55)*\cir<1pt>{} = "b";
**\dir{=}
}
}
}

\def\lbvru{
{\UseComputerModernTips
\xybox{0;/r1pc/:
(0,0.5)*{\bhru}
,(-0.7,0)*{},(2,0)*{},(0,1.5)*{},(0,-1.5)*{},
(0.6,0.45)*{\bullet},
(0.6,-0.55)*{\bullet},
(0.6,0.45)*\cir<1pt>{} = "a";
(0.6,-0.55)*\cir<1pt>{} = "b";
**\dir{=}
}
}
}

\def\hbpar{ \mfgrv{\lbpar} }
\def\hbpru{ \mfgrv{\lbpru} }
\def\hbpro{ \mfgrv{\lbpro} }
\def\hbvir{ \mfgrv{\lbvir} }
\def\hbver{ \mfgrv{\lbver} }
\def\hbvru{ \mfgrv{\lbvru} }
\def\hbhro{ \mfgrv{\lbhro} }
\def\hbhru{ \mfgrv{\lbhru} }
\def\hbviu{ \mfgrv{\lbviu} }

\def\hbcrn{ \mfgrv{\lbcrn} }
\def\hbcrp{ \mfgrv{\lbcrp} }
\def\hbcnw{ \mfgrv{\lbcnw} }
\def\hbcpw{ \mfgrv{\lbcpw} }

\def\hbparprp{ \hbpar_{\mprp} }
\def\hbvirprp{ \hbvir_{\mprp} }
\def\hbverprp{ \hbver_{\mprp} }

\def\hbpruprp{ \hbpru_{\mprp} }
\def\hbviuprp{ \hbviu_{\mprp} }
\def\hbhruprp{ \hbhru_{\mprp} }
\def\hbvruprp{ \hbvru_{\mprp} }

\def\sbpar{
{\UseComputerModernTips
\xygraph{
!{0;/r1pc/:}
[l(1.2)]
!{\xunoverv=<}
}
}
}

\def\sbcrp{
{\UseComputerModernTips
\xygraph{
!{0;/r1pc/:}
[l(1.2)]
!{\xoverv=<}
}
}
}

\def\sbcrn{
{\UseComputerModernTips
\xygraph{
!{0;/r1pc/:}
[l(1.2)]
!{\xunderv=<}
}
}
}

\def\sbpro{
{\UseComputerModernTips
\xygraph{
!{0;/r1pc/:}
[l(1.2)]
!{\xcapv=<}
[ur]
!{\xcapv-=<}
}
}
}

\def\sbhor{
{\UseComputerModernTips
\xygraph{
!{0;/r1pc/:}
[l(1.2)]
[u(0.25)]
!{\xbendr[0.5]=<}
[u][l(0.5)]
!{\xbendl[0.5]=<}
[l(0.5)][d(0.25)]
!{\xcaph}
}
}
}

\def\sbhru{
{\UseComputerModernTips
\xygraph{
!{0;/r1pc/:}
[l(1.2)]
[u(0.25)]
!{\xbendr[0.5]}
[u][l(0.5)]
!{\xbendl[0.5]}
[l(0.5)][d(0.25)]
!{\xcaph}
}
}
}

\def\sbhro{
{\UseComputerModernTips
\xygraph{
!{0;/r1pc/:}
[l(1.2)]
[u(0.25)]
!{\xbendr[0.5]}
[u][l(0.5)]
!{\xbendl[0.5]=<}
[l(0.5)][d(0.25)]
!{\xcaph=>}
}
}
}

\def\lspar{\xybox{0;/r1pc/:(0,0.5)*{\sbpar},(-0.4,0)*{},(1,0)*{},(0,0.6)*{},(0,-0.5)*{}}}
\def\lspro{\xybox{0;/r1pc/:(0,0.5)*{\sbpro},(-0.4,0)*{},(1,0)*{},(0,0.6)*{},(0,-0.5)*{}}}
\def\lshro{\xybox{0;/r1pc/:(0,0.5)*{\sbhro},(-0.4,0)*{},(1,0)*{},(0,0.6)*{},(0,-0.5)*{}}}
\def\lshru{\xybox{0;/r1pc/:(0,0.5)*{\sbhru},(-0.4,0)*{},(1,0)*{},(0,0.6)*{},(0,-0.5)*{}}}
\def\lscrp{\xybox{0;/r1pc/:(0,0.5)*{\sbcrp},(-0.4,0)*{},(1,0)*{},(0,0.6)*{},(0,-0.5)*{}}}
\def\lscrn{\xybox{0;/r1pc/:(0,0.5)*{\sbcrn},(-0.4,0)*{},(1,0)*{},(0,0.6)*{},(0,-0.5)*{}}}

\def\lsvir{
{\UseComputerModernTips
\xybox{0;/r0.5pc/:
(1.7,-1)*{}="a";
(-0.2,1)*{}="b";
{\ar"a";"b"},
(-0.2,-1)*{}="c";
(1.7,1)*{}="d";
{\ar"c";"d"},
(0.75,0)*{\vcirc},
(-0.5,0)*{},(2,0)*{},(0,1.3)*{},(0,-1.3)*{}
}
}
}

\def\lsver{
\xybox{0;/r0.5pc/:
(0,1)*{\sbhor}
,(-0.7,0)*{},(2,0)*{},(0,1.3)*{},(0,-1.3)*{},
(0.6,0.45)*{\bullet},
(0.6,-0.55)*{\bullet},
(0.6,0.45)*\cir<1pt>{} = "a";
(0.6,-0.55)*\cir<1pt>{} = "b";
**\dir{=}
}
}

\def\lsvru{
\xybox{0;/r0.5pc/:
(0,1)*{\sbhru}
,(-0.7,0)*{},(2,0)*{},(0,1.3)*{},(0,-1.3)*{},
(0.6,0.45)*{\bullet},
(0.6,-0.55)*{\bullet},
(0.6,0.45)*\cir<1pt>{} = "a";
(0.6,-0.55)*\cir<1pt>{} = "b";
**\dir{=}
}
}

\def\hspar{ \mfgrv{\lspar} }
\def\hspro{ \mfgrv{\lspro} }
\def\hshro{ \mfgrv{\lshro} }
\def\hsvir{ \mfgrv{\lsvir} }
\def\hsver{ \mfgrv{\lsver} }
\def\hsvru{ \mfgrv{\lsvru} }
\def\hshru{ \mfgrv{\lshru} }

\def\hsparprp{ \hspar_{\mprp} }
\def\hsvirprp{ \hsvir_{\mprp} }
\def\hsverprp{ \hsver_{\mprp} }
\def\hsvruprp{ \hsvru_{\mprp} }
\def\hshruprp{ \hshru_{\mprp} }

\def\hlf{ \frac{1}{2} }
\def\spexlpl{ \Extmfo(\xaAbl,\xaAblp) }

\numberwithin{equation}{section}

\CompileMatrices

\begin{document}

\setlength{\unitlength}{3947sp}
\setlength{\unitlength}{1mm}

\begin{titlepage}
\vfill
\begin{center}
{\large \bf
Virtual crossings, convolutions
and a categorification of the
\sotN\ \Kfp
}\\

\bigskip

\bigskip
\centerline{\textsc{M. Khovanov\footnotemark[1]}}

\centerline{\em Department of Mathematics, Columbia University}
\centerline{\em 2990 Broadway}
\centerline{\em 509 Mathematics Building}
\centerline{\em New York, NY 10027}
\centerline{{\em E-mail address:} {\tt khovanov@math.columbia.edu}}

\bigskip
\centerline{\textsc{L. Rozansky}\footnote[1]{ This work was
supported by NSF Grants DMS-0407784 and DMS-0509793.
}
}

\centerline{\em Department of Mathematics, University of North
Carolina} \centerline{\em CB \#3250, Phillips Hall}
\centerline{\em Chapel Hill, NC 27599} \centerline{{\em E-mail
address:} {\tt rozansky@math.unc.edu}}

\vfill
{\bf Abstract}

\end{center}
\begin{quotation}

We suggest a categorification procedure for the \sotN\ one-variable
specialization of the two-variable \Kfp.
The construction has many similarities with the \Hpt\
categorification: a planar graph formula for the polynomial is
converted into a complex of graded vector spaces, each of them being
the
homology of a \Ztgrdd\ differential vector space associated to a graph
and constructed using matrix factorizations.
This time, however, the elementary \mfs\ are not Koszul; instead, they are
\cnvls\ of chain complexes of \Kmfs.

We prove that the homotopy class of the resulting complex
associated to a \dgrm\ of a link is invariant under the first two
Reidemeister moves and conjecture its invariance
under the third move.

\end{quotation} \vfill \end{titlepage}

\pagebreak
\tableofcontents

\section{Introduction}
\label{Sect1}
\subsection{\sotn\ \Kfp\ and graded complexes}
The \sotn\ \Kfp\cx{Ka} $$\nPLq\in\ZZ[q^{\pm 1}]$$ is an invariant of an
unoriented \frlk\ $\xL\in S^3$, which satisfies the skein relation
\ee
\label{1.0a1}
\aPsk{\xunderh} - \aPsk{\xoverh}=
(q-q^{-1})\lrbc{\aPsk{\xunoverv} -
\aPsk{\xunoverh}},
\eee
the change of framing relation
\ee
\label{1.0a2}
\aPsk{\akinksk} = q^{2\xN+1}\aPsk{\avlinsk}
\eee
(throughout the paper we assume that links are endowed with the
blackboard framing),
and the multiplicativity property
\ee
\label{1.0a3}
\nPqv{\xLot} = \nPqv{\xLo} \nPqv{\xLt},
\eee
where $\xLot$ is the disjoint union of the links $\xLo$ and
$\xLt$.
The conditions\rx{1.0a1}--\rxw{1.0a3} determine $\nPLq$ uniquely\cx{Ka},\cx{KV}. In
particular, the \Kfp\ of the unknot is equal to
\ee
\label{1.2}
\nPunq = \frac{\qtNpm}{\qmqi}+1.
\eee

We will describe a conjectural categorification of the \sotn\ \Kfp.
Namely, to a \ldgrm\ $\xL$ we associate a chain
complex
$\nCdL$ of $\ZZ_2\times\ZZ$-graded $\IQ$-vector spaces
%
\ee
\label{1.3a}
\cdots  \rightarrow\nCnL \rightarrow\nCuLv{n+1}
\rightarrow\cdots,\qquad
\nCnL  = \bopiZjt \xCnijL.
\eee
%
We prove that, as an object
in the homotopy
category of complexes of \ZtZ-\grdmdls,
the complex $\nCdL$
behaves nicely under the first two Reidemeister moves:
\begin{theorem}
\label{th.lr1}
The complex $\nCdL$ is homotopy invariant up to degree shifts under the first
Reidemeister move:
\ee
\label{1.3a1}
\nCdsv{\akinkskn} \cchq \nCdsv{\avlinsk}\grsv{-2\xN-1}\gszto\grsemo,
\eee
Here $\grsv{\cdot}$, $\gsztv{\cdot}$, and $\grsv{\cdot}$ denote the
shifts in the homological degree of the complex\rx{1.3a}, in the
\Ztdgr\ and in the \Zdgr\ respectively.
\end{theorem}


\begin{theorem}
\label{th.lr2}
The complex $\nCdL$ is homotopy invariant under the second
Reidemeister move:
\ee
\label{1.3a2}
\nCdssv{ \hbrttww} \cchq \nCdssv{\hbrtutww}.
\eee
\end{theorem}

Further, we conjecture that
\begin{conjecture}
\label{cj.lr3}
The complex $\nCdL$ is homotopy invariant under the third
Reidemeister move:
\ee
\label{1.3a3}
\nCdssv{\xrtho} \cchq \nCdssv{\xrtht}.
\eee
\end{conjecture}
\begin{conjecture}
\label{cj.dgzt}
The whole complex $\nCdL$ has a homogeneous \Ztdgr
\ee
\label{1.3a4}
\degZt \nCdL = \ncL\;\;\mod{2},
\eee
where $\ncL$ is the number of crossings in the diagram
$\xL$.
\end{conjecture}

We will show that Conjecture\rw{cj.lr3} implies that the
\grdd\ Euler characteristic of the complex $\nCdL$
\ee
\label{1.6}
\chqd\lrbigc{\nCdL } = \snjZit(-1)^{j+n} q^i\rankv\xCnijL
\eee
equals the \Kfp:
\ee
\label{1.7}
\chqd\lrbigc{\nCdL }=\nPLq.
\eee

Our approach to the construction of the categorification complex
$\nCdL$ is similar to that of\cx{KR1}: it follows the alternating
sum formula for the polynomial $\nPLq$, which expresses it in
terms of polynomial invariants of planar graphs. The
categorification of elementary \opgrs\ is provided by \mfs,
and this time the basic polynomial is $\nWxy = \xyt+ \xtN$, as suggested
by Gukov and Walcher\cx{GV}.

\subsection{Closed graphs and the alternating sum formula}
%

Kauffman and Vogel\cx{KV} extended the \Kfp\ from links to
knotted 4-valent graphs. They defined the
invariant of graphs by presenting the 4-vertex as a linear
combination of crossings and their resolutions:
\ee
\label{1.7a1}
\begin{split}
\aPsk{\bvertsk} &= - \aPsk{\xoverv} + q\aPsk{\xunoverv} +
\qi\aPsk{\xunoverh}
\\
&= - \aPsk{\xunderv} + q\aPsk{\xunoverh} + \qi\aPsk{\xunoverv}.
\end{split}
\eee
The second equality follows from the
skein relation\rx{1.0a1} and expresses the $\ntydeg$ rotational
symmetry of the 4-vertex.

The formula\rx{1.7a1} presents the \Kfp\
of a \kngr\ as a linear combination of \Kfps\ of links constructed
by resolving each 4-vertex in three possible ways, thus reducing
the computation of the \Kfp\ of a graph to that of links.

The relation\rx{1.7a1} can
also be played backwards, that is, a crossing can be presented as
a linear combination
\ee
\label{1.7a2}
\aPsk{\xunderv} = q\aPsk{\xunoverh} -\aPsk{\bvertsk}
+\qi\aPsk{\xunoverv}.
\eee
In this form, it allows us to express the \Kfp\ of a link in terms of the
polynomials of \plgrdgs. Namely, let $\xL$ be a \ldgrm\ with
$\ncL$ enumerated crossings. To a multi-index $\xbr=\xlst{\xr}{\ncL}$,
$\xr_i\in\{-1,0,1\}$, we associate a \plgrdg\ $\xGr$ constructed
by resolving all crossings of $\xL$: the $i$-th crossing is
resolved in $\xr_i$-th way. Then, according to \ex{1.7a2},
\ee
\label{e2.3}
\nPLq =  (-1)^{\ncL}
\sxbr (-1)^{\xbra} q^{-\xbra} \nPGrq,
\eee
where $\xbra=\sum_{i=1}^{\ncL} \xr_i$.
We call 4-valent graphs $\xGr$ \emph{closed}.

\subsection{The outline of the categorification}
\label{ss.outcat}
Similar to the \homflyp\ case, the categorification of the \Kfp\
is based on turning the alternating sign sum\rx{e2.3} into a
complex of \grdmdls. To each \clgr\ $\xGr$ we
associate a \ZtZ-\grdmdl\ $\nCGr$ (in fact, we conjecture that
$\degZt\nCGr=0$; the appearance of the $\Zt$ degree $\ncL$ in the
formula\rx{1.3a4} is due to the $\Zt$ degree shift in
\ex{e2.3x2}).
We assemble all spaces $\nCGr$
into an $\ncL$-complex by placing them according to their
multi-indices $\xbr$ (interpreted as $\ncL$-dimensional coordinates)
at the nodes of a $\ZZ^{\ncL}$-lattice
and defining appropriate differential maps
between adjacent spaces. Namely, if three graphs
$\zGhor,\zGV,\zGpar$ result from the same resolution of all
crossings of $\xL$, except for an $i$-th crossing $\lncr$, where they
differ according to their indices, then we construct the linear
maps forming the chain complex
\ee
\label{e2.3x1}
\xymatrix{
{\znGhor\grso} \ar[r]^-{\ychini} &
{\znGV}
\ar[r]^-{\ychouti} &{\znGpar\grsmo}}. 
\eee
%
The complexes\rx{e2.3x1} represent the action of $\ncL$ mutually
anti-commuting
differentials, each differential being related to a resolution of a particular
crossing of $\xL$.
The categorification complex
$\nCdL$ results from collapsing the $\ncL$-complex into a
single complex by taking the sum of individual differentials
and shifting the \Ztdgr\ of the whole complex by $\ncL$ units.
The homological degree of each space $\nCGr$ is set to be equal to
$\xbra$.

For example, the complex of the Hopf link is presented in
\fg{fg.e2.3x1x}.
\begin{figure}
\begin{multline}
\nonumber
\nCdv{\xhpfw} =
\\
\nonumber
\lrbs{\vcenter{
\xymatrix{
{\Chorhor\grst}
\ar[r]^-{\ychinu}
\ar[d]^-{-\ychind}
&
{\Cverhor\grso}
\ar[r]^-{\ychoutu}
\ar[d]^-{\ychind}
&
{\Cparhor}
\ar[d]^-{-\ychind}
\\
{\Cparver\grso}
\ar[r]^-{\ychinu}
\ar[d]^-{-\ychoutd}
&
{\Cverver}
\ar[r]^-{\ychoutu}
\ar[d]^-{\ychoutd}
&
{\Cparver\grsmo}
\ar[d]^-{-\ychoutd}
\\
{\Chorpar}
\ar[r]^-{\ychinu}
&
{\Cverpar\grsmo}
\ar[r]^-{\ychoutu}
&
{\Cparpar\grsmt}
}
}
},
\end{multline}
\caption{The categorification complex of the Hopf link}
\label{fg.e2.3x1x}
\end{figure}
There the maps $\ychsu$ and $\ychsd$ are related to the upper and
lower crossing of the Hopf link diagram.

In order to construct the spaces $\znGr$ and the morphisms between them, we break the
\ldgrm\ $\xL$ and the corresponding \clgrs\ $\xGr$ into
simple pieces by cutting across each edge of the \dgrm\ $\xL$. The
\dgrm\ $\xL$ splits into \eltngls\ $\lncr$, while the
\clgrs\ $\xGr$ split into \eltr\ \emph{\opgrs} $\lxhor$, $\lxver$ and
$\lxpar$. To each \elopgr\ $\xg$ we associate a \mf\ $\gg$,
and to an \eltngl\ $\lxncr$ we associate a complex of \mfs
\ee
\label{e2.3x2}
\hlncr = \lrbcsmo{
\xymatrix{
{\hlhor\grso} \ar[r]^-{\ychin} & {\hlver} \ar[r]^-{\ychout} &
{\hlpar\grsmo}
}
}
\eee
with a special choice of morphisms $\ychin$ and $\ychout$, which
are local versions of the morphisms\rx{e2.3x1}.
The homological degree of the middle \mf\ in this complex is set
to zero.
If a \clgr\ $\xGr$ consists of $\ncL$ \elopgrs\ $\xgi$, then we
set
\begin{align}
\label{e2.3x3}
\nCGr& = \hmlgMF(\hGr),&
\hGr &= \botoncL\ggi,
\\
\label{e2.3x4}
\nCdL& = \hmlgMF(\hL),&
\hL &= \botoncL\hlncr_i,
\end{align}
and $\hmlgMF$ is the \mf\ homology
(its definition is given in subsection\rw{ss.cmfs}; the precise definition
of the tensor product appearing in the
formulas for $\hGr$ and $\hL$ will be given in
subsection\rw{ss.gtamf}). As a tensor product of
complexes\rx{e2.3x2}, $\nCdL$ will have an expected structure of
the total complex of the $\ncL$-dimensional cube-like complex.

The first challenge of the categorification of the \Kfp\ is to choose the
right \mf\ $\hlver$ and morphisms $\ychin,\ychout$. It turns out that in
contrast to the \sun\ and \Hpt\ cases we can not present $\hlver$
as a \Kmf\ (its rank is not a power of 2). Instead, $\hlver$ is realized as a
\cnvl\ of a complex of \Kmfs.

First of all, we introduce a new type of 4-vertex for graphs: a
\emph{virtual crossing} $\lxvir$. Its \mf\ $\hlvir$ is just a
tensor product of two \oarc\ \mfs\ (similarly to $\hlpar$ and $\hlhor$),
as if the segments of $\lxvir$
did not cross. Then we define two special morphisms
\ee
\label{e2.x5}
\xymatrix{ {\hlpar} \ar[r]^-{\yF} &{\hlvir} \ar[r]^-{\yG} &{\hlhor}},
\eee
which we call \emph{saddle} morphisms, and we prove that
\ee
\label{e2.x6}
\yG\yF\mphq 0.
\eee
%

The category of \mfs\ is triangulated, which means that to a
morphism $\mA\!\xrightarrow{\mrf}\!\mB$ between two \mfs\ one can
associate a \mf\ $\CnMFv{\mrf}$ called the cone of $\mrf$. The cone
comes with a pair of natural morphisms
\ee
\xaB\rightarrow\CnMFv{\mrf}\rightarrow\xaA\gszto,
\eee
which form two sides
of the exact triangle of $\mrf$.

A similar construction can be applied to a chain of two
morphisms $\xaA\xrightarrow{\mrf}\xaB\xrightarrow{\mrg}\xaC$, such
that
$\mrg\mrf\mphq 0$.
The resulting \mf\ is called a \emph{\cnvl}
and we denote it by a frame box around the chain:
\fbox{$\xaA\xrightarrow{\mrf}\xaB\xrightarrow{\mrg}\xaC$}
(note that here we omitted the \scnd\ homomorphism in order to
simplify the diagrams; the role of \scnd\ morphisms will be
explained in Section\rw{s.pst}). The module of the convolution is a sum of
modules $\xaA\oplus\xaB\oplus\xaC$. The
convolution comes with two natural morphisms, which form a chain complex of \mfs:

\ee
\label{e.nm}
\xymatrix@C=1.5cm{
{\xaC} \ar[r]^{\smpm{0 \\ 0 \\ \id}} &
\xCthwv{\xaA}{\mrf}{\xaB}{\mrg}{\xaC} \ar[r]^{\smpm{\id & 0 & 0}} &
{\xaA}}.
\eee
This allows us to define $\hlver$ as a convolution of the chain\rx{e2.x5}
\ee
\label{e2.x9}
\hlver = \cnmfFGsh
\eee
and use the natural morphisms of\rx{e.nm} as $\ychin$ and $\ychout$,
so that the definition\rx{e2.3x2} of the categorification of a
crossing becomes
%
%
%
\begin{multline}
\label{ee.crch}
\hlncr =
\\
\lrbcsmo{
\xymatrix@C=1.2cm{
{\hlhorsh} \ar[r]^-{\smpm{0 \\ 0 \\ \id}} & \bicnmfFGsh \ar[r]^-{\smpm{\id & 0 & 0}} &
{\hlparsh}}
}.
\end{multline}
In view of the explicit form\rx{e.nm} of the morphisms $\ychin$
and $\ychout$, we may depict the diagram\rx{ee.crch} informally as
%
%
\ee
\label{eh.x1}
\hlncr =\lrbcsmo{
\vcenter{
\xymatrix@C=1.5cm{
&
{\hlparsh}
\ar[d]^-{\yF}
\save
[]+<1cm,0cm>*{}="a1",
\ar"a1";"1,3"^-{\id}
\restore
&
{\hlparsh}
\\
&
{\hlvirsh}
\ar[d]^-{\yG}
\\
{\hlhorsh}
&
{\hlhorsh}
\save
{[].[u].[uu]}*[F]\frm{}
\restore
\save
"3,2"+<-1cm,0cm>*{}="a2",
\ar"3,1";"a2"^-{\id}
\restore
}
}
}
\eee

If we substitute the definition\rx{e2.x9} into the diagram of \fg{fg.e2.3x1x},
then the complex for the Hopf link takes the form
depicted in \fg{fg.hpfcm}.
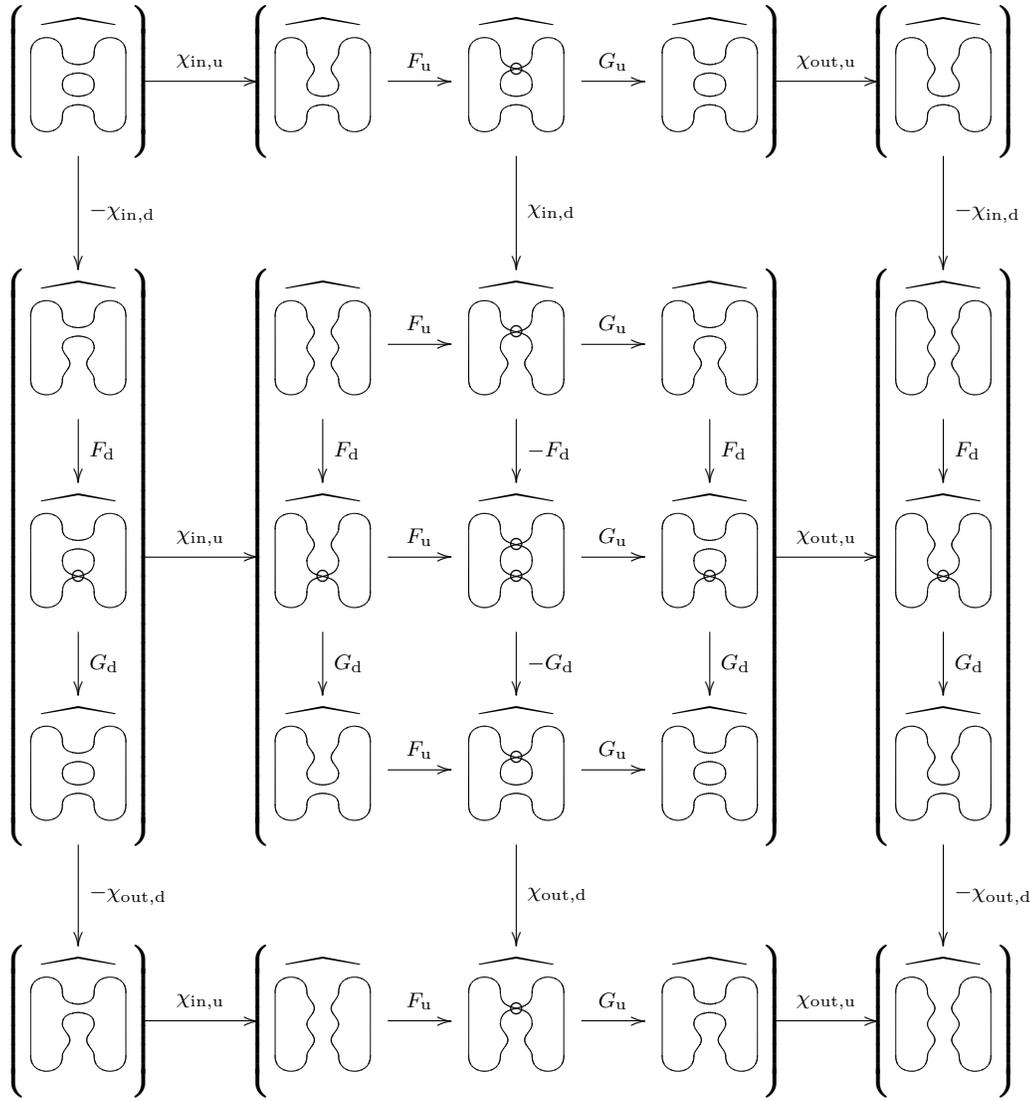
\begin{figure}
\ee
\nonumber
\quad
\vcenter{
\xy
\xymatrix"aa"
{{\Chhorhor}}
\POS*\frm{(}\POS*\frm{)}
\POS(32.5,0)
\xymatrix"ba"
{
{\Chparhor} \ar["aa"];[]^-{\ychinu}
\ar[r]^-{\yFu}
&
{\Chvirhor}
\ar[r]^-{\yGu}
& {\Chhorhor}
}
\POS*\frm{(}\POS*\frm{)}
\POS(115,0)
\xymatrix"ca"{
{\Chparhor}\ar["ba"rr];[]^-{\ychoutu} }
\POS*\frm{(}\POS*\frm{)}
\POS(0,-35)
\xymatrix"ab"{
{\Chhorpar} \ar["aa"];[]^-{-\ychind}
\ar[d]^-{\yFd}
\\
{\Chhorvir}
\ar[d]^-{\yGd}
\\
{\Chhorhor}
}
\POS*\frm{(}\POS*\frm{)}
\POS(32.5,-35)
\xymatrix"bb"{
{\Chparpar}
\ar[r]^-{\yFu}
\ar[d]^-{\yFd}
&
{\Chvirpar}
\ar[r]^-{\yGu}
\ar[d]^-{-\yFd}
\ar["ba"];[]^-{\ychind}
&
{\Chhorpar}
\ar[d]^-{\yFd}
\\
{\Chparvir}
\ar[r]^-{\yFu}
\ar[d]^-{\yGd}
\ar["ab"];[]^-{\ychinu}
&
{\Chvirvir}
\ar[r]^-{\yGu}
\ar[d]^-{-\yGd}
&
{\Chhorvir}
\ar[d]^-{\yGd}
\\
{\Chparhor}
\ar[r]^-{\yFu}
&
{\Chvirhor}
\ar[r]^-{\yGu}
&
{\Chhorhor}
}
\POS*\frm{(}\POS*\frm{)}
\POS(115,-35)
\xymatrix"cb"{
{\Chparpar}
\ar[d]^-{\yFd}
\ar["ca"];[]^{-\ychind}
\\
{\Chparvir}
\ar[d]^-{\yGd}
\ar["bb"rr];[]^{\ychoutu}
\\
{\Chparhor}
}
\POS*\frm{(}\POS*\frm{)}
\POS(0,-125)
\xymatrix"ac"{
{\Chhorpar} \ar["ab"dd];[]^-{-\ychoutd}
}
\POS*\frm{(}\POS*\frm{)}
\POS(32.5,-125)
\xymatrix"bc"{
{\Chparpar}
\ar[r]^-{\yFu}
\ar["ac"];[]^-{\ychinu}
&
{\Chvirpar}
\ar[r]^-{\yGu}
\ar["bb"dd];[]^-{\ychoutd}
&
{\Chhorpar}
}
\POS*\frm{(}\POS*\frm{)}
\POS(115,-125)
\xymatrix"cc"{
{\Chparpar}
\ar["cb"dd];[]^-{-\ychoutd}
\ar["bc"rr];[]^-{\ychoutu}
}
\POS*\frm{(}\POS*\frm{)}\endxy
}
\quad
\eee
\caption{Detailed structure of the Hopf link categorification complex}
\label{fg.hpfcm}
\end{figure}
We omitted there degree shifts and \scnd\ homomorphisms and used an abbreviated notation
\ee
\label{ee.crch1}
\Big(\ast \Big) = \hmlgMF\lrbc{\boxed{\ast}},
\eee
where $\ast$ is a chain of \mf\ morphisms.

\subsection{A categorification complex of a virtual link}
Since we have introduced the \mf\ $\hlvir$ for a virtual crossing,
it is natural to extend the definition\rx{e2.3x4} of the categorification
complex of a link to virtual links.

Virtual links were invented by L.~Kauffman\cx{Kfv} and their topological meaning was clarified
by G.~Kuperberg\cx{GKp} (for the details see V.~Manturov's book\cx{Mntr}). From the
combinatorial point of view, they are equivalence classes of
virtual link diagrams. A virtual link diagram is a
4-valent planar graph with two types of 4-valent vertices: an
ordinary crossing $\lncr$ and a virtual crossing $\lxvir$. Two
diagrams are equivalent (and thus present the same virtual link),
if one can be transformed into another by a sequence of special
moves: ordinary Reidemeister moves, virtual Reidemeister moves and
\smvrt\ Reidemeister moves.

The virtual Reidemeister moves are the ordinary
Reidemeister moves, in which all crossings are replaced by virtual
crossings. A \smvrt\ Reidemeister move is a third Reidemeister
move, in which two crossings are replaced by virtual crossings:
\ee
\label{ee.vrm}
\xvrtho \cchq \xvrtht\;.
\eee

Let $\xL$ be a virtual link diagram with $\ncL$ 4-vertices. We break it into
pieces
$\xgi$ ($1\leq i\leq\ncL$)
by cutting across its edges. We call these pieces \eltngls. The \eltngls\ $\xgi$ coming from
the breakup of $\xL$ are either
crossings $\lncr$ or virtual crossings $\lxvir$. We define the
categorification complex $\nCdL$ of a virtual link diagram $\xL$
similar to \ex{e2.3x4}:
\ee
\label{ee.vct}
\nCdL = \hmlgMF(\hL),\qquad
\hL &= \botoncL\ggi,
\eee
\begin{theorem}
\label{th.vadj}
If two virtual link diagrams $\xL$ and $\xLp$ are related by a virtual or
a \smvrt\ Reidemeister move, then their complexes $\nCdL$ and
$\nCdLp$ are \hteqt.
\end{theorem}
We will prove this theorem in subsection\rw{ss.vadj}.
\subsection{Acknowledgements}
We want to thank Yasuyoshi Yonezawa for helpful suggestion regarding the
manuscript and Vassily Manturov for a discussion of properties of
virtual links. L.R. is indebted to his wife Ruth for her invaluable support.
This work was
supported by NSF Grants DMS-0407784 and DMS-0509793.

\section{Matrix factorizations and graphs}

\subsection{A category of \mfs}
\label{ss.cmfs}
Let us recall the basic definitions and notations related to \mfs.
For a polynomial $\xW$ in a polynomial ring $\yR$, we
define a (homotopy) category of \mfs\ $\MFRW$ (or simply $\MFW$). Its objects
are \Ztgrdd\ free
$\yR$-modules
$\yM=\yMz\oplus\yMo$ equipped with the \emph{\tdiff}
$\xD\in\EndR(\yM)$ such that
\ee
\label{e2.d0}
\degZt\xD=1,\qquad\xD^2=\xW\,\id.
\eee
%
Thus
$\xD$ splits into a sum of two homomorphisms $\xD=\zP+\zQ$
\ee
\label{e2.d1}
\xymatrixnocompile{ \yMo \ar[r]^-{\zP} & \yMz \ar[r]^-{\zQ} & \yMo
},\qquad \zP\zQ = \zQ\zP = \xW\,\id.
\eee
If $\xW\neq 0$, then the conditions\rx{e2.d0} imply
\ee
\label{e2.d1x}
\rank \yMz = \rank\yMo.
\eee
If $\xW=0$, however, then the relation\rx{e2.d1x} no longer holds.

Sometimes, when $\yR=\IQ[\bfx]$, $\bfx=\ylst{\zzx}{n}$, we may use
a notation $\yMxp$ (for a subset $\{\bfx\p\}\subset\{\bfx\}$) instead of
simply $\yM$. The advantage is that we can then denote by $\yMyp$
the \mf\ obtained from $\yMxp=\yM$ by renaming some of the variables $\bfx$, namely $\bfx\p$,
into $\bfy\p$.

The
\tdiff\ $\xD$ generates a differential $d$ ($d^2=0$) acting on $\yR$-homomorphisms
$\yF\in\HomR(\yM,\yMp)$ between the underlying $\yR$-modules $\yM$ and $\yMp$ by the formula
\ee
\label{e2.d2}
d\yF = \ztcm{\xDMMp}{\yF},
\eee
where we use the notation
\ee
\label{e2.d2x}
\xDMMp=\xDv{\yM}+\xDv{\yMp},
\eee
while $\ztcm{\cdot}{\cdot}$ is the \Ztgrdd\ commutator.
In our conventions a composition of two homomorphisms is zero,
unless the target of the first one matches the domain of the
second one. Then we
define the \Ztgrd\ extensions as the homology of $d$:
\ee
\label{e2.d3}
\Extmfb(\yM,\yMp) = \ker d/\im d,\qquad \Extmfb(\yM,\yMp) =
\Extmfz(\yM,\yMp)\oplus\Extmfo(\yM,\yMp),
\eee
and the morphisms between $\yM$ and $\yMp$ in the category $\MFRW$ are
defined as
\ee
\label{e2.d3a1}
\HomMF(\yM,\yMp) = \Extmfz(\yM,\yMp).
\eee

If $\zP=\id$, $\zQ=\xW\,\id$ or $\zP=\xW\,\id$, $\zQ=\id$, then the
identity endomorphism in $\End(\yM)$ becomes a trivial element in
$\Extmfz(\yM,\yM)$, so the corresponding \mf\ is isomorphic to the zero object in
$\MFW$. We call the latter \mfs\ \emph{\cntrb}.

If $\yW=0$, then $\xD^2=0$ and we can define the homology of the
\mf\ $\yM$ as
\ee
\label{e2.d3b1}
\hmlgMFb(\yM) = \ker\xD/\im \xD.
\eee
If in addition to that $\yR=\IQ$, then
\ee
\label{e2.d3b2}
\yM\cchq \hmlgMFb(\yM)
\eee
as objects in the category $\MFvv{\IQ}{0}$ of \Ztgrdd\ differential \mdls.

The \trnsf\ $\gszto$ turns the \mf\rx{e2.d1} into
\ee
\label{e2.d1xa}
\yM\gszto = \lrbc{
\xymatrixnocompile@1{ \yMo \ar[r]^-{-\zQ} & \yMz \ar[r]^-{-\zP} & \yMo
}
},
\eee
the double \trns\ $\gsztv{2}$ acts as the identity and
\ee
\label{e2.d1x1}
\HomMF(\yM,\yMp\gsztv{1}) = \Extmfv{1}(\yM,\yMp).
\eee

For $\yM\in\Obv{\MFRW}$ and $\yMp\in\Obv{\MFvv{\yRp}{\xWp}}$ we define a
tensor product
\ee
\label{e2.d1z1}
\yM\otimes\yMp\in\Obv{\MFvv{\yR\otimes\yRp}{\xW+\xWp}}
\eee
as a \mf, whose module is the tensor product of modules
$\yM\otimes\yMp$ and whose \tdiff\ is the graded sum of \tdiffs
\ee
\label{e2.d1x2}
\xD\otimes \id + (-1)^{\degZt}\,\id\otimes\xDp.
\eee
Similarly, for $\yM\in\Obv{\MFRW}$, $\yMp\in\Obv{\MFvv{\yR}{\xWp}}$ we define
a tensor product
\ee
\label{e2.d1z2}
\yM\otimes_{\yR} \yMp\in\Obv{\MFvv{\yR}{\xW+\xWp}}
\eee
as a \mf, whose module is the tensor product of modules
$\yM\otimes_{\yR}\yMp$ and whose \tdiff\ is\rx{e2.d1x2}.

The conjugate of a \mf\rx{e2.d1} is the \mf\ $\yMs$ of $-\xW$:
\ee
\label{e2.d4}
\xymatrixnocompile{ \yMs_1 \ar[r]^-{\zQs} & \yMs_0 \ar[r]^-{-\zPs} & \yMs_1
}
\eee
It satisfies the property that for any \mf\ $\yMp$ of $\xW$,
\ee
\label{e2.d5}
\Extmfb(\yM,\yMp) \cong \hmlgMFb(\yMp\otimes_{\yR}\yMs)
\eee
(note that according to\rx{e2.d1z2}
$\yMp\otimes_{\yR}\yMs\in\Obv{\MFvv{\yR}{0}}$, so the \rhs is well-defined).

For a homomorphism $\yR\xrightarrow{\zhom}\yRp$ between two
polynomial algebras there is a functor
\ee
\label{e2.d3a2}
\xymatrixnocompile{ \MFRW \ar[r]^-{\zhhom} & \MFvv{\yRp}{\zhom(\xW)}
}
\eee
which takes a \mf\ $\yM$ to $\yM\otimes_{\yR}\yRp$.

There is a homomorphism
\ee
\label{e2.d3a3}
\xymatrixnocompile{\yR \ar[r]^-{ \hmempt } & \EndMF(\yM)},
\eee
which turns an element $r\in\yR$ into a multiplication by $r$:
\ee
\label{e2.d3a4}
\hmr(v)=rv,\qquad v\in\yM.
\eee
If $\yR= \IQ[\bfx]$, $\bfx=\xlst{x}{\xnvr}$, then the homomorphism\rx{e2.d3a3}
factors through the \Jr\
\ee
\label{e2.d3a5}
\JrW = \IQ[\bfx]/(\del_{x_1}\xW,\dotsc,\del_{x_\xnvr}\xW),
\eee
that is, the homomorphism\rx{e2.d3a3} is a composition of
homomorphisms
\ee
\label{e2.d3a6}
\xymatrixnocompile{\IQ[\bfx] \ar[r]
\ar@/^2pc/[rr]^{\hmempt}
 & \JrW\ar[r] & \EndMF(\yM)}.
\eee

Finally, if $\yR$ is a \Zgrdd\ ring and $\yW$ has an even homogeneous \Zdgr\ $\degZ\yW=2N$, then
one can define a category of \Zgrdd\ \mfs\ of $\yW$. Its objects
are \mfs\ with $(\ZZ_2\times\ZZ)$-graded modules $\yM$, whose \tdiffs\
have a homogeneous \Zdgr\ $\degZ\xD=N$.


\subsection{\Kmfs}
Many of our constructions are based upon \Kmfs. For two polynomials $\yp,\yq\in\yR$ the
\emph{\Kmf} $\kmfpq$ of $\xW=\yp\yq$
is a free  \Ztgrd\ $\yR$-module $\yRz\oplus\yRo$ of rank $(1,1)$
with the \tdiff\
$\xD$, whose action is
\ee
\label{ea2.12a}
\xymatrixnocompile{
\yRo\ar[r]^-{\yp} &\yRz \ar[r]^-{\yq} &\yRo.
}
\eee
It is easy to verify that
\ee
\label{ea2.a12a}
\hmv{\yp},\hmv{\yq}\mphq 0,
\eee
that is, both $\hmv{\yp}$ and $\hmv{\yq}$ are homotopic to 0.


If both polynomials $\yp,\yq$ have homogeneous (although, maybe,
unequal) \Zdgrs, then the \Kmf\rx{ea2.12a} can be made \Zgrdd\
by shifting the degree of $\yRo$ by
\ee
\label{ea2.12x01}
\xgkpq = \hfpq.
\eee
Hence when dealing with \grdd\
\Kmfs\ we use the notation
\ee
\label{ea2.12x1}
\kmfpq = (
\xymatrixnocompile{
\yRo\grshfpq\ar[r]^-{\yp} &\yRz \ar[r]^-{\yq} &\yRo\grshfpq
}
),
\eee
where $\grsv{*}$ denotes the \Zgrdng\ shift.
$\kmfpq\grsk$ denotes the \Kmf\ in which the \Zgrdng\
of both modules $\yRz,\yRo$ is shifted by extra $\xgk$ units
relative to that of\rx{ea2.12x1}.
Obviously,
\ee
\label{ea2.21ax}
\kmfpq\grsk\gszto \cong \kmfv{-\yq}{-\yp}\grsv{\xgkpq+\xgk},\qquad
\big(\kmfpq\grsk\big)^\ast \cong \kmfv{\yq}{-\yp}\grsv{-\xgk}.
\eee

For two columns
\ee
\label{ea2.12a3}
\ybp =\pcol{\yp}{n},\qquad \ybq =
\pcol{\yq}{n},\qquad
\qquad\yp_i,\yq_i\in\yR,
\quad i=1,\ldots,n,
\eee
we define the \Kmf\ $\kmfbpq$ of the polynomial
$\sum_{i=1}^n\yp_i\yq_i$ as the tensor product of the elementary
ones:
\ee
\label{ea2.12a4}
\kmfbpq = \bigotimes_{i=1}^n \kmfpqi.
\eee

The case $n=2$ plays an important role in our computations, so in
order to set the basis notations, we present it explicitly. Thus
consider a \Kmf
\ee
\label{ea2.12y1}
\kmftwv{\zpo}{\zpt}{\zqo}{\zqt} =
(\yRo\xrightarrow{\zpo}\yRz\xrightarrow{\zqo}\yRo)\otimes
(\yRpo\xrightarrow{\zpt}\yRpz\xrightarrow{\zqt}\yRpo)
\eee
We choose the generators of rank 2 submodules of $\Zt$-degree 0
and 1 according to the splitting
\ee
\label{ea2.12y2}
\yRtz = \yRzz\oplus\yRoo,\quad
\yRto = \yRzo\oplus\yRoz,\qquad\text{where}\quad
\yRij = \yRi\otimes\yRpj.
\eee
Then the \Kmf\rx{ea2.12y1} has the form
\ee
\label{ea2.12y3}
\xymatrixnocompile{
\yRto \ar[r]^-{\zP} & \yRtz \ar[r]^-{\zQ} & \yRto
}
\eee
where
\ee
\label{ea2.12y4}
\zP = \begin{pmatrix}\zpt & \zpo\\ \zqo & - \zqt\end{pmatrix},
\quad
\zQ = \begin{pmatrix}\zqt & \zpo\\ \zqo & -\zpt \end{pmatrix}.
\eee

In fact, the columns of ring elements\rx{ea2.12a3} should be
considered as elements of $\yRn$ and $\yRns$:
\ee
\label{ea2.12z1}
\ybp \in\yRn,\qquad\ \ybq\in\yRns,\qquad \xW =  \ybq\neg\ybp.
\eee
The implication is that a particular choice of polynomials $\ypi$
and $\yqi$ representing $\ybp$ and $\ybq$ depends on the choice of
free generators of $\yRn$. A change of generators modifies
$\ypi$'s and $\yqi$'s, but the \Kmf\ $\kmfbpq$ remains the same. In
particular, it is invariant under the following row
transformations, acting on $i$-th and $j$-th rows of the
columns\rx{ea2.12a3}:
\ee
\label{ea2.12z2}
\xymatrixnocompile@C=1.5cm{
{\tcoltv{\ypi}{\yqi}{\ypj}{\yqj}} 
\ar@{|->}
[r]^-{\rtrijl} &
{\tcoltv{\ypi}{\yqi-\xlmb\yqj}{\ypj + \xlmb\ypi}{\yqj}
},
}
\qquad\xlmb\in\yR.
\eee
A conjugation of this transformation by the \trnsf\ $\gszto$
(which may act by switching $\ypj$ and $\yqj$) yields another
useful equivalence:
\ee
\label{ea2.12z3}
\xymatrixnocompile@C=1.5cm{
{\tcoltv{\ypi}{\yqi}{\ypj}{\yqj}} 
\ar@{|->}
[r]^-{\rtrpijl} &
{\tcoltv{\ypi}{\yqi-\xlmb\ypj}{\ypj}{\yqj+\xlmb\ypi}
},
}
\qquad\xlmb\in\yR.
\eee
For an invertible element $\xlmb\in\yR$ there is another
equivalence transformation acting on the $i$-th row of $\ybp$ and $\ybq$:
\ee
\label{ea2.12z3x}
\xymatrixnocompile{(\ypi,\yqi) \ar@{|->}[r]^-{\rtoil} &
(\xlmb\ypi,\xlmb^{-1}\yqi)
}
\eee

\subsection{Two theorems}
The following two theorems present us with convenient technical
tools for establishing equivalences between \Kmfs.

\begin{theorem}
\label{th.equiv}
If $\ylst{\yp}{n}\in\yR$ form a regular sequence and
\ee
\label{ea2.12z3a}
\smion \ypi\,\yqi = \smion\ypi\,\yqi\p,
\eee
then the following two \Kmfs\ are isomorphic:
\ee
\label{ea2.12z4}
\kmfv{\ybp}{\ybq} \cong \kmfv{\ybp}{\ybq\p}
\eee
\end{theorem}
\proof
We prove the theorem by induction over $n$. If $n=1$, then the
claim is obvious. Suppose that the theorem holds for $n=k$ and
consider the case $n=k+1$. In view of the condition\rx{ea2.12z3a},
\ee
\label{ea2.12z5}
\ypko(\yqko\p-\yqko) = - \smiok \ypi(\yqi\p - \yqi).
\eee
The definition of a regular sequence implies that $\ypko$ is not a
zero divisor in the quotient $\yR/\xlst{\yp}{k}$, hence it follows
from \ex{ea2.12z5} that
\ee
\label{ea2.12z6}
\yqko\p-\yqko \in \rIdv{\ylst{\yp}{k}},
\eee
that is, there exist $\ylst{\xlmb}{k}\in\yR$ such that
\ee
\label{ea2.12z7}
\yqko\p-\yqko = \smiok \xlmb_i\ypi.
\eee
Therefore, a composition of transformations $\rtrpv{i}{k+1}{\xlmb_i}$ for
$1\leq i\leq k$ turns $\kmfv{\ybp}{\ybq}$ into a \Kmf
\ee
\label{ea2.12xx1}
\tK\begin{pmatrix}
\yp_1, & \yq_1-\xlmb_1\yp_{k+1}
\\
\vdots & \vdots
\\
\yp_k, & \yq_k - \xlmb_k\yp_{k+1}
\\
\yp_k, &\yq\p_{k+1}
\end{pmatrix}
\cong
\tK\begin{pmatrix}
\yp_1, & \yq_1-\xlmb_1\yp_{k+1}
\\
\vdots & \vdots
\\
\yp_k, & \yq_k - \xlmb_k\yp_{k+1}
\end{pmatrix}
\otR
\kmfv{\yp_{k+1}}{\yq_{k+1}\p}
\eee
Since
\ee
\label{ea2.12xx2}
\sum_{i=1}^k \yp_i (\yq_i-\xlmb_i\yp_{k+1}) = \sum_{i=1}^k
\yp_i\yq_i\p,
\eee
then by the inductive assumption
\ee
\label{ea2.12xx3}
\tK\begin{pmatrix}
\yp_1, & \yq_1-\xlmb_1\yp_{k+1}
\\
\vdots & \vdots
\\
\yp_k, & \yq_k - \xlmb_k\yp_{k+1}
\end{pmatrix}
\cong
\tK\begin{pmatrix}
\yp_1, & \yq_1\p
\\
\vdots & \vdots
\\
\yp_k, & \yq_k\p
\end{pmatrix}
\eee
and the tensor product in the \rhs of \ex{ea2.12xx1} is isomorphic
to $\kmfv{\ybp}{\ybq\p}$, which proves the theorem.\qed


Consider three polynomials: $\xaWb\in\xaQbx$ and
$\xap,\xaq\in\xaQbxy$. Let us fix a positive integer $\xnpr$ such that
\ee
\label{xat3}
\xadgy \xap < \xnpr.
\eee
Let $\xaM$ be a \mf\ of the polynomial
\ee
\label{xat2}
\xaW = \xaWb - (\ynp)\,\xaq.
\eee
In view of the condition\rx{xat3}, the quotient
\ee
\label{xat5}
\xaMp=\xaM/(\ynp)\xaM
\eee
is a free module over $\xaQbx$. Hence it represents a \mf\ of the
polynomial $\xaWb$ over that algebra.

A tensor product of $\xaM$ with a \Kmf\ of a special form
\ee
\label{xat1}
\xaMb = \kmpynpq\,\xaotxy\,\xaM
\eee
is also a \mf\ of $\xaWb$.
\begin{theorem}
\label{txat}
$\xaMp$ and $\xaMb$ are \hteqt\ as \mfs\ over $\xaQbx$ via
a natural \hteq\ morphism
\ee
\label{zat1}
\xymatrixnocompile{\xaMb\ar[r]^-{\xfis} &\xaMp}.
\eee
\end{theorem}
\proof
The module of the tensor product $\xaMb = \kmpynpq\,\xaotxy\,\xaM$
is the sum of modules $\xaM\oplus\xaM\gszto$, and its
\tdiff\ is the sum of \tdiffs\ of the modules and the homomorphisms in the diagram
\ee
\label{yat5}
\xymatrixnocompile@C=2cm{
\xaM\gszto
\ar@<0.5ex>[r]^-{\ynp}
\xaarl
&
\xaM
\ar@<0.5ex>[l]^-{\yq}
\xaarr
},
\eee
where $\xD$ is the \tdiff\ of $\xaM$.

The \mf\ $\xaMb$ contains a \sbf\
\ee
\label{zat2}
\xaMbp = \kmfbv{1}{\ynpq}\xaotxy\xaM,
\eee
as demonstrated by the following diagram
\ee
\label{zat3}
\vcenter{\xymatrixnocompile@R=1cm@C=2cm{
\xaMbp
\ar@{^{(}->}[d]
&
\xaM\gszto
\ar[d]_-{1}
\ar@<0.5ex>[r]^-{1}
&
\xaM
\ar[d]^-{\ynp}
\ar@<0.5ex>[l]^-{\ynpq}
\\
\xaMb
&
\xaM\gszto
\ar@<0.5ex>[r]^-{\ynp}
&
\xaM
\ar@<0.5ex>[l]^-{\yq}
}
}
\eee
Its top line represents $\xaMbp$, its bottom line represents
$\xaMb$ and the vertical homomorphisms represent the injection.

It is obvious from the diagram\rx{zat3} that
$\xaMb/\xaMbp=\xaM/(\ynp)\xaM=\xaMp$.
The \Kmf\ $\kmfbv{1}{\ynpq}$ is contractible, hence $\xaMbp$ is
also contractible and $\xaMb\cong\xaMb/\xaMbp$. Together with the previous isomorphism,
this proves the \hteq\ $\xaMb\cong\xaMp$. The equivalence is
established by the quotient map
\ee
\label{zat4}
\xymatrixnocompile{\xaMb\ar[r]^-{\xfis} & \xaMb/\xaMbp=\xaMp},
\eee
which is natural in $\xaM$.\qed

Different \hteq\ maps $\xfis$ commute with each other. Namely,
consider five polynomials: $\xaWb\in\xaQbx$, $\xapo\in\IQ[\xabx,\xayo]$ and
$\xapt,\xaqo,\xaqt\in\xaQxyot$ with the conditions
\ee
\label{zat5}
\deg_{\xayo}\xapo<\xano,\qquad \deg_{\xayt}\xapt<\xant.
\eee

Let $\xaM$ be a \mf\ of the polynomial
\ee
\label{zat6}
\xaW = \xaWb - \ynpqo - \ynpqt.
\eee
Its tensor product with two \Kmfs
\ee
\label{zat7}
\xaMb = \kmpynqpo\xaotxyot\kmpynqpt\xaotxyot\xaM
\eee
is a \mf\ of $\xaWb$ and we consider it over the algebra $\xaQbx$.
A double application of Theorem\rw{txat} leads to the \hteq\
\ee
\label{zat9}
\xaMb\cchq\xaMp,
\eee
where
\ee
\label{zat10}
\xaMp = \xaM/\big((\ynpo)\xaM + (\ynpt)\xaM\big)
\eee
is a \mf\ of $\xaWb$ over $\xaQbx$. However, one could establish the \hteq\rx{zat9}
by either applying Theorem\rw{txat} first to $\xayo,\xapo,\xaqo$
and then to $\xayt,\xapt,\xaqt$ or in the opposite order, thus
obtaining two \hteq\ homomorphisms
\ee
\label{zat11}
\xymatrixnocompile@C=2.5cm{\xaMb \ar@/^10pt/[r]^-{\xfist\xfiso}
\ar@/_10pt/[r]_-{\xfiso\xfist}
& \xaMp}
\eee
\begin{theorem}
\label{th.heq}
The \hteq\ homomorphisms\rx{zat11} are equal:
\ee
\label{zat12}
\xfist\xfiso=\xfiso\xfist.
\eee
\end{theorem}
\proof
The \mf\ $\xaMb$ in\rx{zat7} has the structure
\ee
\label{zat13}
\xaMb = \lrbc{\vcenter{
\xymatrixnocompile@C=2cm@R=2cm{
\xaM
\ar@<0.5ex>[r]^-{\ynpo}
\ar@<-0.5ex>[d]_-{-(\ynpt)}
&
\xaM\gszto
\ar@<0.5ex>[l]^-{\xaqo}
\ar@<-0.5ex>[d]_-{\ynpt}
\\
\xaM\gszto
\ar@<0.5ex>[r]^-{\ynpo}
\ar@<-0.5ex>[u]_-{-\xaqt}
&
\xaM
\ar@<0.5ex>[l]^-{\xaqo}
\ar@<-0.5ex>[u]_-{\xaqt}
}}}
\eee
This \mf\ has a \sbf\ $\xaMbp$:
\ee
\label{zat14}
\xaMbp = \lrbc{\vcenter{
\xymatrixnocompile@C=2cm@R=2cm{
\xaM
\ar@<0.5ex>[r]^-{\ynpo}
\ar@<-0.5ex>[d]_-{-(\ynpt)}
&
\xaM\gszto
\ar@<0.5ex>[l]^-{\xaqo}
\ar@<-0.5ex>[d]_-{\ynpt}
\\
\xaM\gszto
\ar@<0.5ex>[r]^-{\ynpo}
\ar@<-0.5ex>[u]_-{-\xaqt}
&
(\ynpo)\xaM+(\ynpt)\xaM
\ar@<0.5ex>[l]^-{\xaqo}
\ar@<-0.5ex>[u]_-{\xaqt}
}}}
\eee
A double application of the proof of Theorem\rw{txat} to $\xaMb$
indicates that $\xaMp$  is contractible and both compositions
$\xfist\xfiso$ and $\xfiso\xfist$ are canonical homomorphisms
between a module and its quotient:
\ee
\label{zat15}
\xymatrixnocompile@C=2.5cm{\xaMb \ar@/^10pt/[r]^-{\xfist\xfiso}
\ar@/_10pt/[r]_-{\xfiso\xfist}
& \xaMb/\xaMbp=\xaMp}.
\eee
\qed

\subsection{Graphs, \tngls\ and associated \mfs}
\label{ss.gtamf}

In our categorification constructions, an \emph{\opgr}
is a \plgrdg\ with special univalent vertices called \emph{\lgs}.
$\lxhor$, $\lxver$ and $\lxvir$ are examples of \opgrs.
The \lgs\ are enumerated and oriented as `in' or `out', but sometimes we drop these
decorations on the pictures.

If an
\opgr\ has no \lgs, then it is called \emph{closed}. The graphs
$\xGr$ resulting from the resolution of the crossings of $\cL$ are
examples of \clgrs.

More generally, a \grtngl\ is a \plgrdg\ with \lgs\ and also with
special 4-valent vertices called \emph{\crss}. These vertices are
decorated with the under-over selection for the incident edges. A
link diagram $\cL$ and an \eltngl\ $\lxncr$ are examples of
\grtngls, the first one being, in fact, closed.

There are two operations on the sets of \clgrs\ (and similarly on
\grtngls). The first operation is a disjoint union
\ee
\label{ex.1}
\xymatrixnocompile{
(\xgo,\xgt) \ar@{|->}[r]^-{\dju} & \xgo\dju\xgt
}
\eee

We call two legs of an \opgr\ or a \grtngl\ $\gg$ \xcls\ if they can
be connected by a line in the complement of $\gg$. The second
operation joins two \xcls\ \lgs\ $i$ and $j$ with opposite orientations:
\ee
\label{ex.2}
\xymatrixnocompile{
\xg \ar@{|->}[r]^-{\jlgij}
&\jlgij(\xg)
}.
\eee

Our categorification is a map
from \opgrs\
to \mfs\ and from \grtngls\ to complexes of \mfs
\ee
\label{ex.2a}
\xg\mapsto\gg,\qquad\gt\mapsto\hgt.
\eee
The categorification map
should convert the operations\rx{ex.1}
and\rx{ex.2} into appropriate algebraic counterparts.

First of all, we choose a basic set of variables
$\bfx=\xlst{\zzx}{\zzk}$ and a basic polynomial
$\basW(\bfx)\in\IQ[\bfx]$. To a set $\xmlg$ of $\zzn$ enumerated oriented \lgs\
(of an \opgr\ or a \grtngl) we associate a polynomial of
$\zzn\zzk$ variables
\ee
\label{ex.3}
\bfxlg = \ylst{\bfx}{\zzn},\qquad \bfxi =
\zzx_{1,i},\dotsc,\zzx_{\zzk,i},
\eee
according to the formula
\ee
\label{ex.4}
\xWlg(\bfxlg) = \smion\xmui\,\basW(\bfxi),
\eee
where $\xmui$ is the orientation of the $i$-th leg: $\xmui=1$ if
it is oriented outwards, and $\xmui=-1$ if it is oriented inwards.

To an \opgr\ $\xg$ we associate a \mf\
\ee
\label{ex.5}
\gg\in\Obv{\MFv{\xWv{\xmlg(\xg)}}},
\eee
where $\xmlg(\xg)$ is a set of oriented legs of $\xg$. Similarly,
to a \grtngl\ $\gt$ we associate a complex of \mfs\ up to
homotopy:
\ee
\label{ex.6}
\hgt\in\Obv{\Kmplv{\MFv{\xWv{\xmlg(\gt)}}}}.
\eee

If $\xG$ is a \clgr\ then $\xmlg(\xG)=\emptyset$, so
$\xWv{\xmlg(\xG)}=0$, and as a `\mf', $\hG$ is isomorphic to its own homology
\ee
\label{ex.7}
\hG\cchq\nCv{\xG},
\eee
where by definition
\ee
\nCv{\xG}=\hmlgMF(\hG)
\eee
(\cf \ex{e2.3x3}).

Finally, we specify the functoriality rules for the
maps\rx{ex.2a}.
To a disjoint union of \opgrs\ $\xgo\sqcup\xgt$ we associate the
tensor product of \mfs
\ee
\label{ex.8}
\widehat{\xgo\sqcup\xgt} = \ggo\otimes\ggt,
\eee
and to the joining of legs we associate the identification of
the variables $\bfxi$ and $\bfxj$:
\ee
\label{ex.9}
\widehat{\jlgij(\xg)} = \zhhomij(\gg),
\eee
where $\zhhomij$ is the functor\rx{e2.d3a2} corresponding to the
algebra homomorphism
\ee
\label{ex.10a}
\xymatrixnocompile{
\IQ[\bfx_{\xmlg(\xg)}]\ar[r]^-{\zhomij}
& \IQ[\bfx_{\xmlg(\xg)}]/
(\zzx_{1,i}-
\zzx_{1,j},\dotsc, \zzx_{\zzk,i}-\zzx_{\zzk,j})
}.
\eee
It is important to note that the \rhs of \ex{ex.9} should be
considered as a \mf\ over the ring
$\IQ[\bfx_{\xmlg(\jlgij(\xg))}]$, that is, the multiplication by
$\bfxi$ and $\bfxj$ is not a part of the module structure, since
these variables are no longer related to legs. Nevertheless,
$\hmv{\bfx}_i$ and $\hmv{\bfx}_j$ remain collections of endomorphisms of
$\widehat{\jlgij(\xg)}$, and
\ee
\label{ex.11}
\hmv{\bfx}_i=\hmv{\bfx}_j
\eee
in view of the quotient\rx{ex.10a}.

The functorial properties of the categorification map allow us to
define it first for the \eltr\ \opgrs\ and \grtngls, which are
just single vertices together with their legs, and for the \oarc\
graph $\lxoar$, and then derive the categorification of more complicated
objects by cutting them into \eltr\ pieces and assembling the
corresponding \mfs\ with the help of \eex{ex.8} and\rx{ex.10a}. The
\mf\ $\hloar$
is determined by the
principal property of $\lxoar$:
gluing it to a leg of another graph $\xg$ would produce the same
graph. In other words, if $\xg$ has $n$ \lgs\ enumerated by
numbers $1\dotsc n$, then
\ee
\label{ex.12}
\jlgino\lrbc{\xg\; \dju\;\;
\loarxvv{n+1}{n+2}
} = \xgp,\qquad 1\leq i\leq n,
\eee
where $\xgp$ is the same graph as $\xg$, except the $i$-th leg of
$\xg$ is labeled as $(n+2)$-th.
Then the formulas\rx{ex.8} and\rx{ex.10a} require that
\ee
\label{ex.13}
\gg_{\bfxi}\otimes_{\IQ[\bfxi]}\lrbc{\hloarvv{i}{n+2}} \cchq
\gg_{\bfxv{n+2}}.
\eee
%

Note that the property\rx{ex.13} of the \oarc\ \mf\ allows us to
present the \rhs of \ex{ex.9} as a tensor product with \oarc\ \mf,
so that \ex{ex.9} can be written as
\ee
\label{ex.14}
\widehat{\jlgij(\xg)} \cchq \gg\otimes_{\IQ[\bfxi,\bfxj]} \lrbc{
\hloarvv{i}{j}
}.
\eee
This formula is similar to definition of the Hochschild
homology of a bimodule, $\lrbc{
\hloarvv{i}{j}}$ playing the role of the algebra resolution.

\section{Arc \mfs}
\subsection{Basic polynomial}
We begin the categorification of the \sotn\ \Kfp\ by choosing the
basic algebra $\IQ[x,y]$ and the basic polynomial
\ee
\label{ea2.4}
\nWxy = \xyt+ \xwx,\qquad \xwx = \xtN,
\eee
which was previously considered by Gukov and Walcher\cx{GV}.
This polynomial has a homogeneous \qdgr
\ee
\label{ea2.4x1}
\dgq\nWxy = 4 \xN + 2,
\eee
if we set
\ee
\label{ea2.6}
\dgq x = 2,\qquad \dgq y = 2\xN,
\eee
and this allows us to use \grdd\ \mfs\ for categorification. In
particular, the \tdiff\ should have a homogeneous \qdgr
\ee
\label{ea2.13a2}
\dgq\xD = 2\xN+1.
\eee

The basic polynomial\rx{ea2.4} is an odd function of its combined argument:
\ee
\label{ea2.15}
\gWmxy = -\gWxy.
\eee
This allows us to impose an additional constraint on the
categorification map\rx{ex.2a}:
if $\xgp$ is obtained from
$\xg$ by reversing the orientation of the $i$-th leg, then $\ggp$
is obtained from $\gg$ by changing the signs of the variables
$\zxi,\zyi$, that is,
\ee
\label{ea2.15z}
\ggp_{\zxi,\zyi} = \gg_{-\zxi,-\zyi}.
\eee
Therefore, from now on we assume that all legs of graphs and tangles are oriented
outwards, unless specified otherwise, and whenever we need to join two legs, we change the
orientation of one of them with the help of \ex{ea2.15z}. In other words, the homomorphism $\zhomij$ in the joining
formula\rx{ex.9} should be modified from \ex{ex.10a} to
\ee
\label{ex.10}
\xymatrixnocompile{
\IQ[\bfx_{\xmlg(\xg)}]\ar[r]^-{\zhomij}
& \IQ[\bfx_{\xmlg(\xg)}]/
(\zzx_{1,i}+
\zzx_{1,j},\dotsc, \zzx_{\zzk,i}+\zzx_{\zzk,j})
},
\eee
%
because
$i$-th and $j$-th legs are both oriented outwards.


We introduce a convenient notation for iterated
differences of $\xwx$, which we define recursively as
\ee
\label{ea.defdif}
\xw(\zx_1,\ldots,\zx_{n-1},\zx_n) =
\frac{\xw(\zx_1,\ldots,\zx_{n-2},\zx_{n-1})-\xw(\zx_1,\ldots,\zx_{n-2},\zx_{n})}{\zx_{n-1}-\zx_{n}
}.
\eee
For example,
\ee
\label{ea2.13a}
\xwxot = \frac{ \xwv{\zxo}-\xwv{\zxt} }{\zxo-\zxt},\qquad
\xwxott = \frac{ \xwvv{\zxo}{\zxt} -
\xwvv{\zxo}{\zxth}}{\zxt-\zxth}.
\eee
Note that all these polynomials are symmetric functions of their
arguments.

\subsection{\Jr}
The \Jr\ of $\nW$ is defined as the quotient
\ee
\label{ea2.6a}
\JrW = \IQ[x,y]/(\del_x\nW,\del_y\nW) = \IQ[x,y]/(y^2 + (2\xN+1)x^{2\xN},xy).
\eee
It is \grdd\ with a \qgrdd\ basis
\ee
\label{ea2.6ya}
1,x\dotsc x^{2\xN},y,\qquad \dgq x^i=2i,\quad\dgq y =2\xN.
\eee
%
Hence the dimensions of \qgrdd\ subspaces of $\JrW$ are
\ee
\label{ea2.6xa2}
\dim\JrWqv{2i} =
\begin{cases}
1& \text{if $0\leq i\leq 2\xN$ and $i\neq \xN$,}
\\
2& \text{if $i=\xN$},
\end{cases}
\eee
and the
\grdd\ dimension of $\JrW$
is equal to the \Kfp\ of the unknot\rx{1.2} up to a \dgrsh:
\ee
\label{ea2.6a1}
\dmq \JrW = \frac{q^{4\xN+2} - 1}{q^2-1} + q^{2\xN} =
q^{2\xN}\nPunq.
\eee
This agreement indicates that \ex{ea2.4} is a suitable
choice for the basic polynomial.

Let
\ee
\label{ee.fr1}
1^*,x^*\dotsc (x^{2\xN})^*,y^*\in\JrW^*
\eee
denote the dual basis of\rx{ea2.6ya}, and for $v\in\JrW$ let
$\hv:\;\JrW\xrightarrow{v}\JrW$ denote the linear operator of
multiplication by $v$. Then the algebra multiplication map
\ee
\label{ee.fr2}
\xymatrixnocompile{\JrW\otimes\JrW \ar[r]^-{\xmlt} &\JrW}
\eee
can be presented as
\ee
\label{ee.fr4}
\xmlt = \siztN \hmx^i\otimes (x^i)^* + \hmy\otimes y^*.
\eee

The \Jr\ $\JrW$ acquires a \Frba\ structure, if we
choose the \Frb\ trace to be
\ee
\label{ee.fr5}
\TrF = - \frac{1}{\ttNo}\, (x^{2\xN})^*
\eee
(the origin of the normalization factor will become clear later).
Then the \cmlt\ map
\ee
\label{ee.fr6}
\xymatrixnocompile@1{\JrW  \ar[r]^-{\xcmlt} & \JrW\otimes\JrW}
\eee
can be presented as
\ee
\label{ee.fr7}
\xcmlt = 2\lrbc{\hmy\otimes y - (2\xN+1)\siztN \hmx^i\otimes
x^{2\xN-i}
}.
\eee

\subsection{The 1-arc \Kmf\ and the unknot space}

To a \oarc\ \opgr\ we associate a \Kmf
\ee
\label{ea2.16c}
\zarcxyoo=\kmftwv{\zyt+\zyo}{\zxt+\zxo}{\zxt\,(\zyt-\zyo)}{\zyo^2 + \xwvv{-\zxo}{\zxt}}.
\eee
of the polynomial
\ee
\label{ey.c1}
\zWt = \nW(\zxo,\zyo) + \nW(\zxt,\zyt).
\eee
The left column of this \Kmf\ contains sums rather than differences,
because in our conventions both legs of the \oarc\ graph
$\lxoarvv{1}{2}$ are oriented outwards.
Theorem\rw{txat} indicates that the \mf\rx{ea2.16c} satisfies the
property\rx{ex.13}. Also note that in view of \ex{ea2.a12a},
\ee
\label{ea2.a16}
\hxt\mphq - \hxo,\qquad \hyt\mphq-\hyo.
\eee

According to \ex{ea2.21ax}, the dual to the \oarc\ \mf\ is the
\mf\ of the \oarc\ graph in which both legs are oriented inwards:
\ee
\label{ea2.16cx}
\lrbc{\zarcxyoo}^\ast \cong \hloarinsvv{1}{2}\grkdl.
\eee
%


The unknot diagram $\lzcir$ can be constructed by joining legs 1 and 2 of the arc
$\arcxyoo$, hence, in accordance with \eex{e2.3x3} and\rx{ex.9}, the corresponding
space $\nCcir$
is the homology of the differential \mdl\ $\hlzcir$
constructed by setting
\ee
\label{ea2.a16cx}
-\zxt=\zxo=\zx,\qquad-\zyt=\zyo=\zy
\eee
in the \Kmf\rx{ea2.16c}:
\ee
\label{ea2.16c2}
%
\xymatrixnocompile{
\yRzz \ardpl\ar[rr]^-{\zy^2+\xw\p(\zx)}
\ar[rrd]^{2\zx\zy}
&& \yRzo\grxhv{1-2\xN} \ardpl
\ar[rrd]^-{2\zx\zy}
&& \yRzz \ardpl
\\
\yRoo\grkdlm &&
\yRoz\grxhv{-1}
\ar[rr]^{-(\zy^2 +\xw\p(\zx))}
&& \yRoo\grkdlm
}
\eee
This is a Koszul complex $\kctv{2\zx\zy}{\zy^2+\xw\p(\zx)}$. Since
the polynomials $\zx\zy,\zy^2+\xw\p(\zx)$ form a regular sequence,
the homology of the complex\rx{ea2.16c2} appears only in the lower
right corner and, in fact, as a \qgrdmdl\ it is isomorphic to the \Jr\rx{ea2.6a} up to
a \qdgr\ shift:
\ee
\label{ea2.16ac3}
\hmlgMFi\lrbc{\hlzcir} \cong
\begin{cases}
\JrWp,
& \text{if $i=0$},\\
0, & \text{if $i=1$},
\end{cases}
\qquad \text{where $\JrWp=\JrW\grkdlm$.}
\eee
Thus
\ee
\label{ea2.16c3}
\nCcir \cong
\JrWp,
\qquad\degZt\nCcir = 0.
\eee
It follows from this equation and from \ex{ea2.6a1} that the
relation\rx{1.7} holds for the crossingless diagram of the unknot.

$\nCcir$ is isomorphic to $\JrWp$ not only as a \grdmdl\ but
also as a module over $\JrW$ whose action on $\nCcir$ is
generated by
\ee
\label{ey.b2}
\hmx=\hxo=-\hxt,\qquad\hmy=\hyo=-\hyt.
\eee
This module structure allows us to introduce a convenient basis of
$\nCcir$.
Let the first basis element $\xone$ be a non-zero element of the 1-dimensional lowest
\qdgr\ subspace
\ee
\label{ey.c2}
\nCcirgmN\subset\nCcir,
\eee
then the rest of the basis is generated by the action of the elements of $\JrW$ on $\xone$:
\ee
\label{ey.b1}
\xone,\hmx(\xone),\zxgv{2}\dotsc\zxgtN,\zyg\in\nCcir.
\eee
Later we will endow $\nCcir$ with an algebra structure, and this
will turn\rx{ea2.16c3} into an algebra isomorphism. Then $\xone$
will be chosen to be the unit of $\JrW$.

Finally, we describe the endomorphism algebra of the \oarc\
\mf\rx{ea2.16c}. The partial derivatives
\ee
\label{ea2.16d1}
\del_x\xW = y^2 + (2\xN+1)x^{2\xN},\qquad
\del_y\xW = 2xy
\eee
form a regular sequence, hence
\ee
\label{ea2.16d2}
\EndMF(\zarcxyoo)=
\Extmfz(\zarcxyoo,\zarcxyoo)
\cong \JrW,\qquad \Extmfo(\zarcxyoo,\zarcxyoo)=0.
\eee
This isomorphism follows directly
from\rx{ea2.16ac3}. Indeed, according to the relations\rx{e2.d5},
\rxw{ea2.16cx} and to the joining legs formula\rx{ex.14},
\ee
\label{ea2.16d5}
\begin{split}
\Extmfi(\zarcxyoo,\zarcxyoo)
&= \hmlgMFi(\zarcxyoo \otRot (\zarcxyoo)^\ast) =
\hmlgMFi(\zarcxyoo\otRot\hloarins)\grkdl
\\
&=
\hmlgMFi\lrbc{\hlzcir}\grkdl =
\begin{cases}
\JrW & \text{if $i=0$},\\
0 & \text{if $i=1$},
\end{cases}
\end{split}
\eee
where $\yR = \IQ[x_1,y_1,x_2,y_2]$.
The isomorphism $\JrW\xrightarrow{\cong} \EndMF(\zarcxyoo)$
can be generated by the homomorphism $\;\hmempt\;$ in
two different ways:
\ee
\label{ea2.a16d2}
(x,y)\mapsto (\hxo,\hyo)\qquad\text{or}\quad
(x,y)\mapsto (\hxt,\hyt).
\eee
Both versions differ by a sign in view of \ex{ea2.a16}.

The calculations\rx{ea2.16d5} leading to \ex{ea2.16d2} can be
generalized to find the morphisms between multi-arc graphs.

\begin{theorem}
\label{th.arcgl}
Consider a system of $2m$ distinct 1-vertices (\lgs). Let $\xg$ and
$\xgp$ be two \opgrs\ which are the unions of $m$ \oarc\ graphs
connecting these $2m$ legs. Then
\ee
\label{ee.mrarc}
\Extmfz(\xg,\xgp) = (\JrW)^{\otimes m_0} \grxhv{2\xN(m-m_0)},\qquad
\Extmfo(\xg,\xgp) = 0,
\eee
where $m_0$ is the number of circles formed by joining $\xg$ and
$\xgp$ together through their common \lgs.
\end{theorem}
\proof
Let $\xgm$ be the graph $\xg$, in which we reversed the
orientations of all \lgs, so that they are now oriented inwards.
Then
\ee
\label{ee.mrarc1}
\Extmfi(\gg,\ggp) = \hmlgMFi(\ggp\otRot(\gg)^*) =
\hmlgMFi(\ggp\otRot\ggm) \grxhv{ 2\xN m},
\eee
and the formula\rx{ee.mrarc} follows from the basic property of
the \oarc\ \mf\ and from \ex{ea2.16c3}.\qed

\subsection{\tarc\ \Kmfs}
\subsubsection{Category, \enfns\ and \sdlmrps}
\label{sss.endf}
Four \lgs\ $1,2,3,4$ can be connected by two arcs in three different
ways:
\ee
\label{eb2.a16}
\xybox{0;/r4pc/:
\xunoverh,
(0.8,0.2)*{3},
(2.2,0.2)*{4},
(0.8,-1.2)*{1},
(2.2,-1.2)*{2}
}
\qquad\qquad
\xybox{0;/r4pc/:
\xunoverv,
(0.8,0.2)*{3},
(2.2,0.2)*{4},
(0.8,-1.2)*{1},
(2.2,-1.2)*{2}
}
\qquad\qquad
\xybox{0;/r4pc/:
(0.5,-0.5)*{\zbendv@(0)},
(0.5,0)*{\sbendv@(0)},
(1.5,-0.503)*{\vcirc},
(0.8,0.2)*{3},
(2.2,0.2)*{4},
(0.8,-1.2)*{1},
(2.2,-1.2)*{2}
}
\eee
The third connection requires the arcs to intersect in the plane of the diagram, however when
we assign the \mf\ to the \opgr\ $\lxvir$, we treat this intersection as
non-existent (that is, virtual), while considering $\lxvir$ to be a
disjoint union of two arcs similar to $\lxpar$ and $\lxhor$.

The \mfs\ 
$\hlpar,\hlvir,\hlhor$ belong to the same category of \mfs\ of the
polynomial
\ee
\label{eb2.16}
\zWfxy = \sum_{i=1}^4 \nWxyi
\eee
over the algebra
\ee
\label{eb2.16p}
\yR=\IQbxy,\qquad \bfx=\ylst{x}{4},\quad\bfy=\ylst{y}{4}.
\eee
The structure of \tarc\ \mfs\ is simple: in view of \ex{ex.8} they are tensor
products of \mfs\ corresponding to constituent arcs. In
particular, the virtual nature of the crossing $\lxvir$ implies
that
\ee
\label{eb2.16o}
\hlvir = \hloarvv{1}{4}\otimes\hloarvv{2}{3}.
\eee

The symmetric group $\Sf$ acts on the set of graphs\rx{eb2.a16} by
permuting their \lgs, \eg
\ee
\label{eb2.16a}
\xymatrixnocompile{
\lxpar \ar@{|->}[r]^{\sXot}  & \lxvir 
}.
\eee
The polynomial\rx{eb2.16} is invariant under the simultaneous
permutations
of the variables $\bfx$ and $\bfy$, hence the elements $\xsgm\in\Sf$
act as \enfns\ on the category $\MFWf$.

Theorem\rw{th.arcgl} describes the spaces $\Extmfi(\gg,\ggp)$ for
\tarc\ graphs $\xg$ and $\xgp$:
\ee
\label{eb2.16a2}
\Extmfz(\gg,\ggp) =
\begin{cases}
\JrW\otimes\JrW & \text{if $\xg=\xgp$,}
\\
\JrW\grkdl & \text{if $\xg\neq \xgp$,}
\end{cases}
\qquad\qquad
\Extmfo(\gg,\ggp) = 0.
\eee
Suppose that $\xg\neq\xgp$. According to \ex{eb2.16a2},
\ee
\label{eb2.16xa4}
\dim\ExtmfzqN\lrbc{\gg,\ggp} = 1.
\eee
Therefore, up to a constant factor, there exists a unique non-trivial morphism
\ee
\label{eb2.16xa5}
\xymatrixnocompile@1{\gg\ar[r]^-{\yF}&\ggp},\qquad \dgq\yF=2\xN.
\eee
We call it \emph{the \sdlmrp} and we denote it either as $\yF$, or
$\yG$, or $\yH$
depending on the choice of $\xg$ and $\xgp$.

The uniqueness of the \sdlmrp\ (up to a constant factor)
for a given pair of graphs
$\xg\neq\xgp$ implies that the elements of $\Sf$
transform \sdlmrps\ into \sdlmrps\
(these functors are images of group elements, hence they
are invertible and
can not transform a \sdlmrp\ into a trivial morphism).

\subsubsection{Proper 2-arc \Kmfs}

We have to know an explicit presentation of a saddle morphism
to perform computations with it.
According to \eex{ea2.16c} and\rx{eb2.16o}, the \tarc\ \mfs\ $\gg$
are presented as 4-row \Kmfs, so their modules have rank (8,8),
which makes explicit formulas for morphisms rather unwieldy.
However, as we shall see, all \tarc\ \mfs\ can be presented as
tensor products
\ee
\label{ea2.23a}
\gg = \zKcmn\otR\ggtprp,
\eee
where $\zKcmn$ is a `\cmn' 2-row \Kmf\ (the same for all
$\gg$), while $\ggtprp$ are `\prpr' 2-row \Kmfs. The \sdlmrps\ act
as identity on $\zKcmn$, that is
\ee
\label{ea2.33x2}
\vcenter{\xymatrixnocompile{
\gg \ar[rrr]^-{\yF} \ar@{=}[d] &&& \ggp \ar@{=}[d]
\\
\zKcmn\otR\ggtprp  \ar[rrr]^-{\id\otR\yFtprp}&&&
\zKcmn\otR\ggtprp\p
}}
\eee
for a suitable map $\xymatrixnocompile@1{\ggtprp \ar[r]^-{\yFtprp} &
\ggtprp\p}$.

In order to describe this construction explicitly, we introduce a
\qprpr\ algebra
\ee
\label{ey.2}
\yRprp = \IQ[\zxbp,\zxbq,\zxbr,\zC],\qquad\text{where}\quad
\zxbp=(\zpo,\zpt),\quad\zxbq=(\zqo,\zqt),\quad\zxbr=(\zro,\zrt),
\eee
and a \qprpr\ polynomial
\ee
\label{ea2.28}
\zWfp = \zpo\zqt\zrt+\zqo\zpt\zrt+\zro\zpt\zqt +\zpo\zqo\zro\zC.
\eee
We also define a homomorphism
\ee
\label{ey.3}
\xymatrixnocompile{\yRprp\ar[r]^-{\zhomprp} & \yR}
\eee
by the formulas
\begin{align}
\label{ey.3a}
\nonumber
\zhomprp(\zpo) & = \zxft, & \zhomprp(\zqo)&=\zxfth, & \zhomprp(\zro) &= \zxfo,
\\
\zhomprp(\zpt) & = \zyft, & \zhomprp(\zqt)&=\zyfth, & \zhomprp(\zrt) &= \zyfo.
\end{align}
and
\ee
\label{ea2.22}
\zhomprp(\zC) =\ztC(\bfx),
\eee
where
\ee
\label{ea2.22x}
\ztC(\bfx) = \xwvf{\zxo}{\zxt}{-\zxth}{\zxf} +
\xwvf{\zxo}{-\zxo}{\zxt}{\zxf} + \xwvf{\zxo}{-\zxo}{-\zxt}{\zxf}
\eee
and \ex{ea.defdif} defines the polynomial $w$.

Let $\Sth\subset\Sf$ be the subgroup permuting the \lgs\
$(1,2,3)$.
Consider its permutation action on the triplet of
variables $(\zxbr,\zxbp,\zxbq)$. Note that the formulas\rx{ey.3a} are
invariant with respect to the simultaneous action of $\Sth$ on
$(\zxo,\zxt,\zxth)$, $(\zyo,\zyt,\zyth)$
and $(\zxbr,\zxbp,\zxbq)$.

The polynomial $\zWfp$ is invariant under the action of $\Sth$ on $(\zxbr,\zxbp,\zxbq)$,
hence the elements $\xsgm\in\Sth$ act as \enfns\ $\hsgm$ on the category
of \mfs\ $\MFvv{\yRprp}{\zWfp}$. In particular, they permute the
\mfs\
\begin{align}
\label{ea2.25}
\hlparp &=
\kmftwv{\zpo}{\zpt}{\zqt\zrt+\zqo\zro\zC}{\zqt\zro+\zqo\zrt},
\\
\nonumber
\\
\label{ea2.26}
\hlvirp &=
\kmftwv{\zro}{\zrt}{\zpt\zqt+\zpo\zqo\zC}{\zpo\zqt+\zpt\zqo},
\\
\nonumber
\\
\label{ea2.27}
\hlhorp &=
\kmftwv{\zqo}{\zqt}{\zpt\zrt+\zpo\zro\zC}{\zpo\zrt+\zpt\zro},
\end{align}
equivariantly with respect to their action on the \tarc\ graphs:
\ee
\label{ey.4}
\hsgm (\ggprp) = \widehat{\lrbsc{\xsgm(\xg)  }}_{\mprp}.
\eee

Define
\ee
\label{ey.5}
\ggtprp = \zhhomtprp(\ggprp),
\eee
where $\zhhomtprp$ is the functor\rx{e2.d3a2} induced by the
homomorphism\rx{ey.3}.

\begin{proposition}
\label{pr.factprp}
Three \tarc\ \mfs\ factor as
\ee
\label{ea2.23a1}
\gg = \zKcmn\otR\ggtprp,
\eee
where
\ee
\label{ea2.24}
\zKcmn = \kmftwv{\zxall}{\zyall}{\zAxy}{\zBxy}
\eee
and
\begin{align}
\label{ea2.20}
\zAxy & = -(\zyfth)
(\zyo+\zyt+\zyf)-\zyo\zyt + \xwvv{-\zxo}{\zxth}
+(\zxft)\xwvvv{\zxo}{\zxt}{-\zxth}
\\
\nonumber
&\hspace{0.8in}-(\zxft)(\zxfo)\big(\xwvf{\zxo}{-\zxo}{\zxt}{\zxf}+\xwvf{\zxo}{-\zxo}{-\zxt}{\zxf}
\big)
\\
\label{ea2.21}
\zBxy & = \zxo\zyo + \zxt\zyt + \zxth\zyth - \zxf\zyf.
\end{align}
\end{proposition}
\proof
We will prove this proposition for $\xg=\lxpar$, the proofs for
other graphs are similar. According to our definitions,
\begin{multline}
\label{ey.6}
\zKcmn\otR\hlpartp =
\\
\kmffwv{\zxall}{\zyall}{\zxft}{\zyft}
{\zAxy}{\zBxy}{(\zyfo)(
\zyfth
) +
(\zxfo)(\zxfth)\,\ztC(\bfx)}{(\zxfo)(
\zyfth
)+(\zyfo)(\zxfth)}
\end{multline}
A straightforward computation shows that this is a \mf\ of the
polynomial $\zWfxy$. On the other hand, the \mf\ $\hlpar$ is a
4-row \Kmf\ of the same polynomial, which is transformed by the
combination $\rtrthoo\circ\rtrfto$ of the
transformations\rx{ea2.12z2} into the following form
\ee
\label{ey.7}
\begin{split}
\hlpar &=
\kmffwv{\zxtho}{\zytho}{\zxft}{\zyft}{\zyo^2+\xwvv{-\zxo}{\zxth}}
{\zxth\,(\zyth-\zyo)}{\zyt^2+\xwvv{-\zxt}{\zxf}}{\zxf\,(\zyf-\zyt)}
\\
&=
\kmffwv{\zxall}{\zyall}{\zxft}{\zyft}{\zyo^2+\xwvv{-\zxo}{\zxth}}
{\zxth\,(\zyth-\zyo)}{\zyt^2-\zyo^2+\xwvv{-\zxt}{\zxf}-\xwvv{-\zxo}{\zxth}}{\zxf\,(\zyf-\zyt)
-\zxth\,(\zyth-\zyo)}.
\end{split}
\eee
The latter \Kmf\ shares the same first column with the
\mf\rx{ey.6}. Since the polynomials in that column form a regular
sequence, the claim of the proposition follows from
Theorem\rw{txat}.\qed

\section{Saddle morphisms}

\subsection{A \prpr\ \sdlmrp}
According to the formulas\rx{ea2.12y3} and\rx{ea2.12y4}, the \prpr\ \mf\ $\ggprp$
has an explicit form
\ee
\label{ea2.28z1}
\xymatrixnocompile{\yRto \ar[r]^-{\zPi} & \yRtz \ar[r]^-{\zQi} & \yRto,
}
\eee
where the homomorphisms $\zPi$ and $\zQi$ are presented by
matrices
\begin{align}
\label{ea2.35}
\zPpar & = \begin{pmatrix}
\zpt & \zpo \\
\zqt\zrt + \zqo\zro\zC & -(\zqt\zro+\zqo\zrt)
\end{pmatrix}
, &
\zQpar & = \begin{pmatrix}
\zqo\zrt + \zqt\zro& \zpo\\
\zqt\zrt+\zqo\zro\zC & -\zpt
\end{pmatrix},
\\
\label{ea2.36}
\zPvirt & = \begin{pmatrix}
\zrt & \zro \\
\zpt\zqt + \zpo\zqo\zC & - (\zpo\zqt+\zpt\zqo)
\end{pmatrix}
, &
\zQvirt & = \begin{pmatrix}
\zpo\zqt+\zpt\zqo & \zro \\
\zpt\zqt + \zpo\zqo\zC & -\zrt
\end{pmatrix},
\\
\label{ea2.37}
\zPhor & = \begin{pmatrix}
\zqt & \zqo\\
\zpt\zrt+\zpo\zro\zC & -(\zpo\zrt+\zpt\zro)
\end{pmatrix}
, &
\zQhor & = \begin{pmatrix}
\zpo\zrt+\zpt\zro & \zqo \\
\zpt\zrt+\zpo\zro\zC & -\zqt
\end{pmatrix}.
\end{align}

Consider the following homomorphism between two \mf\ modules
\ee
\label{ea2.33x1}
\xymatrix{
{\hlparp} \ar
@<-0.55ex>
[d]^-{\yFprp}
&
{\yRto} \ar[rr]^-{\zPpar}  \ar[d]^-{\yFo} && {\yRtz} \ar[rr]^-{\zQpar} \ar[d]^-{\yFz} &&
{\yRto} \ar[d]^-{\yFo}
\\
{\hlvirp}
&
{\yRto} \ar[rr]^-{\zPvirt}  && {\yRtz} \ar[rr]^-{\zQvirt}
&& {\yRto}
},
\eee
where
\ee
\label{ea2.38c1}
\yFz  =
\begin{pmatrix}
0 & 1\\ \zqt^2 - \zqo^2\zC & 0
\end{pmatrix}
,\qquad
\yFo =
\begin{pmatrix}
\zqt & -\zqo \\ \zqo\zC & -\zqt
\end{pmatrix}.
\eee
By using expressions\rx{ea2.35} and\rx{ea2.36}, it is easy to verify that the
diagram\rx{ea2.33x1} is commutative, that is, $[\xD,\yFprp]=0$. This means that $\yFprp$
defines a \mf\ morphism. The isomorphisms\rx{ea2.23a} allow us to
extend $\yFprp$ to a morphism between full \tarc\ \mfs:
\ee
\label{ey.8}
\xymatrix@1{
{\hlpar} \ar[r]^-{\yF} & {\hlvir}
},\qquad\text{where}\quad
\yF = \id\otR\yFtprp,\quad\yFtprp = \zhhomtprp(\yFprp)
\eee
(\cf diagram\rx{ea2.33x2}).
\begin{proposition}
\label{p.prsd}
The morphism\rx{ey.8} is a \sdlmrp.
\end{proposition}
\proof
It is easy to verify that
\ee
\label{ey.9}
\degZt\yF=\degZt\yFprp=0,\qquad \dgq\yF=\dgq\yFprp=2\xN,
\eee
so it remains to verify that $\yF\not\sim 0$. Indeed, the matrix
$\yFz$ has a unit entry, so the matrix of $\yF$ also has a unit
entry. However, all entries of the \tdiff\ matrices\rx{ea2.35}
and\rx{ea2.36}
belong to the ideal generated by the variables
$x_1,\ldots,x_4,y_1,\ldots,y_4$,
hence there does not exist a homomorphism
$\yX\in\HomR\lrbc{\hlpar,\hlvir}$, such that $\yF=[\xD,\yX]$.\qed

The \sdlmrps\ for other pairs of \tarc\ graphs $\xg\neq\xgp$ can
be obtained from \eex{ea2.38c1} and\rx{ey.8} by the
simultaneous
permutation action of $\Sth$ on the \lgs\ of the graphs and on the
variables $(\zxbr,\zxbp,\zxbq)$.

\subsection{Compositions of \sdlmrps}
Let us find a composition of two \sdlmrps
\ee
\label{ea2.40x2}
\xymatrixnocompile@C=1.5cm{\gg \ar[r]^-{\yF} & \ggp \ar[r]^-{\yG} & \ggpp},\qquad
\xgp\neq\xg,\xgpp.
\eee
We should consider two different cases: $\xg=\xgpp$ and $\xg\neq \xgpp$. In
both cases the symmetric group $\Sf$ acts transitively on such triplets
of \tarc\ graphs $(\xg,\xgp,\xgpp)$, so it is sufficient to
perform the calculation just for one triplet, and the answer for
others could be deduced by permutation of \lgs\ and variables.

\begin{proposition}
\label{p.sdlcomult}
The composition of two \sdlmrps\
\ee
\label{ea2.40y1}
\xymatrix@1{ {\hlpar} \ar[r]^-{\yF} & {\hlvir} \ar[r]^-{\yH} & {\hlpar} }
\eee
is equal to
\ee
\label{ea2.40y2}
\yH\yF \mphq 2\hmy_1\hmy_2 - 2(2\xN+1)\,
\spxot.
\eee
\end{proposition}
\proof
Consider the composition of \prpr\ \sdlmrps
\ee
\label{ea2.40x4}
\xymatrix{
{\hlparp} \ar@<-0.55ex>[d]^-{\yFprp}
&
{\yRto} \ar[rr]^-{\zPpar}  \ar[d]^-{\yFo} &&
{\yRtz} \ar[rr]^-{\zQpar} \ar[d]^-{\yFz} &&
{\yRto}
\ar[d]^-{\yFo}
\\
{\hlvirp} \ar@<-0.55ex>[d]^-{\yHprp}
&
{\yRto} \ar[rr]^-{\zPvirt} \ar[d]^-{\yHo} && {\yRtz} \ar[rr]^-{\zQvirt} \ar[d]^-{\yHz}
&& {\yRto}
\ar[d]^-{\yHo}
\\
{\hlparp}
&
{\yRto} \ar[rr]^-{\zPpar} && {\yRtz} \ar[rr]^-{\zQpar} &&{\yRto}
}
\eee
Since $\yHprp = \hsgm_{12}(\yFprp)$, and $\xsgm_{12}$ switches
$\zxbr$ and $\zxbp$ while leaving $\zxbq$ and $\zC$ intact, we
conclude from \ex{ea2.38c1} that
\ee
\label{ea2.40x6}
\yHz=\yFz,\qquad\yHo=\yFo,
\eee
and matrix multiplication shows that
\ee
\label{ea2.40x5}
\yHprp\yFprp = (\zqt^2 - \zqo^2\zC)\,\id.
\eee
Hence
\ee
\label{ea2.40y1a}
\yH\yF = (\hyth+\hyf)^2
-(\hxth+\hxf)^2
\ztC(\hbfx),
\eee
where $\ztC$ is defined by \ex{ea2.22x}.
Now it is an elementary algebra exercise to transform the \rhs of
this formula into the \rhs of \ex{ea2.40y2} with the help of
\eex{ea2.22x}, \rx{ea2.4} and the
following relations between the endomorphisms of $\hlpar$:
\begin{gather}
\label{ea2.40y2a}
\hmx_3\mphq -\hmx_1,\qquad\hmy_3\mphq -\hmy_1,\qquad
\hmx_4\mphq -\hmx_3,\qquad\hmy_4\mphq-\hmy_3,
\\
\label{ea2.40y3a}
\hmy_1 \mphq - (2\xN+1)\,\hmx_1^{2\xN},\qquad
\hmy_2 \mphq - (2\xN+1)\,\hmx_2^{2\xN}.
\end{gather}
\qed

\begin{proposition}
\label{pr.nh}
The composition of two \sdlmrps\
\ee
\label{ea2.40y3}
\xymatrix@1{
{\hlpar} \ar[r]^-{\yF} & {\hlvir} \ar[r]^-{\yG} & {\hlhor}
}
\eee
is null-homotopic:
\ee
\label{ey.13}
\yG\yF\mphq 0.
\eee
\end{proposition}
\proof
By definition, \ex{ey.13} means that
there exists an $\yR$-module homomorphism
\ee
\xymatrix@1{ {\hlpar}\ar[r]^-{\yX}&{\hlhor}},
\eee
such that
%
\ee
\label{ea2.40a1}
\degZt\yX=1,\qquad
\dgq\yX = \grhm
\eee
and
\ee
\label{ea2.41}
\yG\yF = -\{\xD,\yX\}
\eee
%
(we put a minus sign in for future
convenience).

Consider the \prpr\ \sdlmrps \ee
\label{ea2.40x7}
\xymatrix{
{\hlparp} \ar[d]^-{\yFprp}
&
{\yRto} \ar[rr]^-{\zPpar}  \ar[d]^-{\yFo} && {\yRtz} \ar[rr]^-{\zQpar} \ar[d]^-{\yFz} &&
{\yRto}
\ar[d]^-{\yFo}
\\
{\hlvirp} \ar[d]^-{\yGprp}
&
{\yRto} \ar[rr]^-{\zPvirt} \ar[d]^-{\yGo} && {\yRtz} \ar[rr]^-{\zQvirt} \ar[d]^-{\yGz}
&& {\yRto}
\ar[d]^-{\yGo}
\\
{\hlhorp}
&
{\yRto} \ar[rr]^-{\zPhor} && {\yRtz} \ar[rr]^-{\zQhor} &&{\yRto}
}
\eee
where the second morphism is
presented by the matrices
\ee
\label{ea2.39}
\yGz  =
\begin{pmatrix}
0 & 1 \\ \zpt^2 - \zpo^2 \zC & 0
\end{pmatrix}
,\qquad
\yGo =
\begin{pmatrix}
\zpt & -\zpo \\ \zpo\zC & -\zpt
\end{pmatrix},
\eee
in accordance with the relation $\yGprp=\hsgm_{13}\hsgm_{12}(\yFprp)$. On the
diagram
\ee
\label{ea2.40}
\xymatrix{
{\yRto} \ar[dd]_-{\yGo\yFo}\ar[rr]^-{\zPpar}  &&
{\yRtz} \ar[dd]_-{\yGz\yFz} \ar[rr]^-{\zQpar} \ar[lldd]_-{\yXz} &&
{\yRto} \ar[dd]^-{\yGo\yFo} \ar[lldd]_-{\yXo}
\\
\\
{\yRto} \ar[rr]^-{\zPhor} &&
{\yRtz} \ar[rr]^-{\zQhor} &&
{\yRto}
}
\eee
representing the composition of \sdlmrps\ we choose an $\yR$-module homomorphism
\ee
\label{ey.11}
\yXprp\in\HomR(\hlparp,\hlhorp),\qquad
\yXprp=\yXz+\yXo
\eee
presented by the matrices
\ee
\label{ea2.43}
\yXz = \begin{pmatrix} -\zqt & 0 \\ \zqo\zC & 0\end{pmatrix},
\qquad
\yXo = \begin{pmatrix} 0 & 0 \\ \zpo\zC &  \zpt\end{pmatrix}.
\eee
A direct computation shows that these matrices satisfy the
relations
\ee
\label{ea2.42}
\yGz\yFz = -(\zPhor \yXz + \yXo\zQpar),\qquad
\yGo\yFo = -(\zQhor \yXo + \yXz\zPpar),
\eee
which means that $\yXprp$ satisfies the properties
\ee
\label{ea2.40b1}
\degZt\yXprp=1,\qquad\dgq\yXprp = \grhm,\qquad
\yGprp\yFprp = -\{\xD,\yXprp\},
\eee
and if we choose
\ee
\label{ey.13x}
\yX= \id\otR\yXprp,
\eee
then $\yX$ would
satisfy\rx{ea2.40a1} and\rx{ea2.41}.\qed

\subsection{Semi-closed and \clsd\ \sdlmrps}
In order to exhibit the structure of the categorification
complex  $\nCdL$ and to prove its Reidemeister move invariance
we have to find what happens to the \sdlmrp\ when
we join together one or two pairs of legs of the \tarc\ \opgrs. As
a result of this joining, a \tarc\ graph becomes either an \oarc\ graph, or a \oarc\ graph
with a disjoint circle, or a circle, or a pair of circles. Hence the
initial \tarc\ \mfs\ can be simplified by \hteq\ transformations,
and we want to know how this simplification affects the \sdlmrps.

\subsubsection{Semi-closed \sdlmrps}
Let us join one pair of legs in the \sdlmrps.
The endomorphism
symmetry $\Sf$ allows us to consider only three cases, namely, the
saddle morphisms
\ee
\label{ea2.74}
\xymatrix{
{\hlpar} \ar@/^/[r]^-{\yF} & {\hlvir} \ar@/^/[l]^-{\yH} \ar[r]^-{\yG} &
{\hlhor}
},
\eee
in which we join together legs 2 and 4. After that the
diagram becomes
\ee
\label{ea2.74a}
\xymatrix{
{\hlpary} \ar@/^/[r]^-{\yF} & {\hlviry} \ar@/^/[l]^-{\yH} \ar[r]^-{\yG} &
{\hlhory}
},
\eee
its homomorphisms resulting from taking the quotient of
the \sdl\ homomorphisms\rx{ea2.74} by the ideal $(\zxt+\zxf,\zyt+\zyf)$
according to \ex{ex.10}.
The second and the third graphs in the diagram\rx{ea2.74a} are just \oarc s, so
we pass to \hteqt\ \mfs\ and present this diagram as
%
\ee
\label{ey.29}
\xymatrix{
{\hlhcpy\otimes\nCcir} \ar@/^/[r]^-{\yFp} & {\hlhcpy} \ar@/^/[l]^-{\yHp} \ar[r]^-{\yGp} &
{\hlhcpy}
}.
\eee
The first \mf\ in this diagram splits into a direct sum of \oarc\ \mfs
\ee
\label{ey.31}
\hlhcpy\otimes\nCcir = \bopiztN \lrbc{\hlhcpy\otimes \zxgi} \oplus
\lrbc{\hlhcpy\otimes\zy},
\eee
where $\zxgi$ and $\zy$ form the basis\rx{ey.b1} of the space
$\nCcir$.
Hence we will express the \smcl\ \sdlmrps\ $\yFp$ and $\yGp$
in terms of the endomorphisms of the \oarc\ \mf\ $\hlhcpothy$.
According to \ex{ea2.16d5},
\ee
\label{ey.32}
\EndMF\lrbc{\hlhcpothy} = \JrW,
\eee
where $\JrW$ is generated by
\ee
\label{ey.33}
\hmx=\hxo=-\hxth,\qquad\hmy=\hyo=-\hyth.
\eee

Recall that the space $\nCcir$ has a special basis\rx{ey.b1}.
Let us introduce the dual basis for the dual space $\nCscir$:
\ee
\label{ey.30}
\xones,(\hmx(\xone))^*,\zxsgv{2}\dotsc\zxsgv{2\xN},\zys\in\nCscir.
\eee

\begin{proposition}
\label{pr.sdlc}
With the appropriate choice of the basis element $\xone$,
the homomorphisms $\yFp$ and $\yGp$ of the diagram\rx{ey.29} are
\htpc\ to the following homomorphisms
\begin{align}
\label{ea2.82}
\yFp & \mphq\yFm= 
\lrbc{\smiztN \hmx^i\otimes\zxsgi +
\hmy\otimes\zys
},
\\
\label{ea2.83}
\yHp &\mphq \yHD= 2\lrbc{\hmy\otimes\zyg - (2\xN+1)\smiztN \hmx^i\otimes\zxgv{2\xN-i} }.
\end{align}
%
The homomorphism $\yGp$ is null-homotopic:
\ee
\label{ea2.86}
\yGp\mphq 0.
\eee
\end{proposition}
The indices ${\mathrm m}$ and $\Delta$ indicate that the
homomorphisms $\yFm$ and $\yHD$ are related to the
multiplication\rx{ee.fr4}
and comultiplication\rx{ee.fr7} in $\JrW$, as we will see shortly.

\proof
We derive the expressions for the \smcl\ \sdlmrps\ not by a direct computation, but rather
from three properties of the \sdlmrps\rx{ea2.74}. First, a \sdlmrp\ has homogeneous
\qdgr\ $2\xN$, so
\ee
\label{ea2.76a1}
\dgq\yFp=\dgq\yHp=2\xN.
\eee
and these morphisms must have a form
\begin{align}
\label{ea2.81}
\yFp & \mphq \smiztN a_i\,\hmx^i\otimes\zxsgi + \tla\, \hmy\otimes\zys
+ a\p_1\,
\hmy\otimes\zxsgv{\xN} + a\p_2 \hmx^{\xN}\otimes\zys,
\\
\label{ea2.81x}
\yHp &\mphq \smiztN b_i\,\hmx^i\otimes\zxsgv{2\xN-i} + \tlb\,\hmy\otimes\zys +
b\p_1\,\hmy\otimes\zxsgv{\xN} + b\p_2\, \hmx^{\xN}\otimes\zys,
\end{align}
where the coefficients $a$'s and $b$'s are rational numbers.

The second property is the composition formula\rx{ea2.40y2}, which
we apply to the composition $\yHp\yFp$:
\ee
\label{ee.comp}
\yHp\yFp\mphq
2\hmy_1\hmy_2 - 2(2\xN+1)\,
\frac{\hmx_1^{2\xN+1}-\hmx_2^{2\xN+1}}{\hmx_1-\hmx_2}.
\eee
It imposes relations among the coefficients of expressions\rx{ea2.81}, which imply that
\ee
\label{eh.1}
\yFp\mphq\xxa\yFm,\qquad\yHp\mphq\xxb\yHD,\qquad\xxa,\xxb\in\IQ,\qquad \xxa\xxb=2.
\eee
Indeed, if we apply \ex{ee.comp} to $\hlhcpy\otimes\xone$, then we
find that $\yHp\mphq\xxb\yHD$ and $a_0\xxb=2$. The other
coefficients
of\rx{ea2.81} are determined by applying the
relation\rx{ee.comp} to $\hlhcpy\otimes\zxgi$ and to
$\hlhcpy\otimes\zyg$.

Since $\xxa\xxb=2$, then $\xxa\neq 0$ and we can rescale the
basis element $\xone$ so that $\xxa=1$, $\xxb=2$. Thus we
obtain
\eex{ea2.82} and\rx{ea2.83}.

Finally, Proposition\rw{pr.nh} says that
\ee
\label{ea2.76a6}
\yGp\yFp\mphq 0.
\eee
At the same time, according to \ex{ea2.82}, the morphism $\yFp$ acts on the submodule
$\hlhcpothy\otimes \xone
\subset\hlhcpothy\otimes\nCcir$ as $\id$, so $\yGp$ must be null-homotopic. \qed

\subsubsection{Closed \sdlmrps}
Let us consider the result of joining legs 2 and 4, as well as
legs 1 and 3 in the graphs of
the diagram\rx{ea2.74}:
\ee
\label{eh.3}
\xymatrix{
{\hlparz} \ar@/^/[r]^-{\yF} & {\hlvirz} \ar@/^/[l]^-{\yH} \ar[r]^-{\yG} &
{\hlhorz}
}.
\eee
The first graph in this diagram is a pair of disjoint circles,
while the second and the third graphs are circles. Therefore, we
can use the \mf\ \hteq\ in order to present the diagram\rx{eh.3} as
\ee
\label{eh.4}
\xymatrix{
{\nCcir\otimes\nCcir}
\ar@/^/[r]^-{\yFpp} & {\nCcir} \ar@/^/[l]^-{\yHpp} \ar[r]^-{\yGpp} &
{\nCcir}
}.
\eee
The morphisms on this diagram can be obtained from the morphisms
of the diagram\rx{ey.29} by imposing the joining relation
$\hxo=-\hxth$ in the expressions\rx{ea2.82},\rx{ea2.83}
and\rx{ea2.86}. In fact, since $\hxth$ never appears there, the
expressions remain the same. First of all, we find that
\ee
\yGpp = 0.
\eee
Second, we observe that the operators $\hmx$ and $\hmy$ appearing
in the formulas for $\yFm$ and $\yHD$ are the same as the
operators $\hmx$ and $\hmy$ of \eex{ee.fr4} and\rx{ee.fr7}. Thus,
if we establish a canonical isomorphism
\ee
\label{eh.5}
\nCcir \cong \JrWp=\JrW\grkdlm
\eee
by the basis elements correspondence
\ee
\label{eh.6}
\zxgi \leftrightarrow x^i, \qquad\zyg \leftrightarrow y,
\eee
then the diagram\rx{eh.4} becomes
\ee
\label{eh.7}
\xymatrixnocompile{
\JrWp\otimes\JrWp
\ar@/^/[r]^-{\xmlt} & \JrWp \ar@/^/[l]^-{\xcmlt} \ar[r]^-{0} &
\JrWp
},
\eee
while the homomorphisms\rx{ea2.82} and\rx{ea2.83} take the form
\begin{align}
\label{ea2.82x}
\yFp & \mphq\yFm= 
\smiztN \hmx^i\otimes\zxuis +
\hmy\otimes\zyus,
\\
\label{ea2.83x}
\yHp &\mphq \yHD= 2\lrbc{\hmy\otimes\zy - (2\xN+1)\smiztN \hmx^i\otimes\zxuv{2\xN-i} }.
\end{align}

\section{Multi-dimensional \Pss\ and their \cnvls}
\label{s.pst}
The category of \mfs\ is triangulated.
Relation\rx{ey.13} means that the diagram\rx{ea2.40y3} is a \chcl\
of \mf\ morphisms. This allows us to define the \mf\ $\hlver$ as
its \cnvl. But first let us review the basic facts about \Pss\ and
their \cnvls. The general theory of \Pss\ within triangulated
categories is described in the book\cx{GM} (Chapter IV, exercises). We adapt its approach
to the specific case of \mfs.

Some proofs in this section are omitted, because they are
standard exercises in homological algebra.
%
\subsection{A multi-dimensional \Ps}
Let us introduce multi-index notations: for a positive integer $n$
and $v_i\in\ZZ$ we denote $\bxv = \xlst{v}{\zan}$ and $\xtlv{\bxv}
= \smion v_i$. Also $\bxv\geq \bxw$ means that $v_i\geq w_i$ for
all $i=1\dotsc \zan$, and $\bxv>\bxw$ means that $\bxv\geq \bxw$
and $\bxv\neq\bxw$. The zero vector is denoted $\bxz=(0,\dotsc,0)$,
and the vectors $\ylst{\bxe}{n}$ form the standard basis in the
lattice $\ZZn$.

We introduce a new $\ZZn$ grading, which we call \emph{\lngth},
with multi-degree
%
\ee
\label{ef.1}
\bdegl = (\deglv{1},\dotsc,\deglv{\zan}),
\eee
and assume that $\bdegl x=\bxz$ for all $x\in\yR$.
We also introduce the \tlngth\ degree:
\ee
\label{ef.2}
\degl = \xtlv{\bdegl}.
\eee
%
A \lngth-graded $\yR$-module $\xaA$
is called \emph{\lbndd}
if its grading
expansion
\ee
\label{ef.2a}
\xaA = \bopbkzn \xaAbk, \quad\bdegl\xaAbk = \bxk
\eee
has only finitely many
non-trivial modules $\xaAbk$.
A homomorphism
$\yX\in\HomR(\xaA,\xaB)$ between two \lngth-graded $\yR$-modules
is called \nng\ if it is a sum of homomorphisms of \nng\ \lngthd s:
\ee
\label{ef.3}
\yX = \sum_{\bxk\geq \bxz} \yXbk,\qquad\bdegl\yXbk = \bxk.
\eee
Let $\HomRngl(\xaA,\xaB)$ denote the space of all such
homomorphisms.

For a fixed polynomial $\xW\in\yR$ we define the category $\PMFWn$ of
$n$-dimensional \Pss. Its objects are \lngth-graded \lbndd\ \mfs\ $\xaA$ of
$\xW$ with \nng\ \tdiffs\ $\xDA$:
\ee
\xaA = \bopbkzn
\xaAbk,\qquad\xDA:\;\xaAbk\rightarrow\bigoplus_{\bxl\geq\bxk}\xaAbl,\qquad
\xDA^2 = \xW\,\id_{\xaA}.
\eee
For two such \mfs\ $\xaA$ and $\xaB$, we
define the differential $d$ acting on the space $\HomRngl(\xaA,\xaB)$ by the
formula\rx{e2.d2}:
\ee
\label{ef.4}
d = \ztcm{\xDAB}{\cdot}.
\eee
Then we define the $\Zt$-graded space
\ee
\label{ef.5}
\ExtPmfb(\xaA,\xaB) = \ker d/\im d,
\eee
and the space of $\PMFWn$-morphisms between $\xaA$ and
$\xaB$:
\ee
\label{ef.6}
\HomPMF(\xaA,\xaB) = \ExtPmfz(\xaA,\xaB).
\eee

The category $\PMFWn$ is triangulated. For a morphism
$\xymatrixnocompile@1{\xaA\gszto \ar[r]^-{f}& \xaB}$
we define the
cone $\CnPSv{f}$ as a \Ps\ whose module is the sum $\xaA\oplus
\xaB$ and whose \tdiff\ $\xD$ is the sum of individual \tdiffs\ and $f$:
$\xD=\xDA+\xDB+f$.

The tensor product\rx{e2.d1z1} creates a bifunctor
\ee
\label{ef.23a}
\xymatrixnocompile@1{
\PMFvv{\yR_1,\xW_1}{\zan_1}\times
\PMFvv{\yR_2,\xW_2}{\zan_2}
\ar[r]^-{\otimes} &
\PMFvv{\yR_1\otimes\yR_2,\xW_1+\xW_2}{\zan_1+\zan_2}
},
\eee

A \emph{\cnvl} is a functor
\ee
\label{ef.7}
\xymatrixnocompile@C=2.5cm{\PMFWn \ar[r]^-{\CvlMFv{\cdot}} & \MFW}
\eee
which turns a \Ps\ $\xaA$ into an ordinary \mf\
$\CvlMFv{\xaA}$
by `forgetting' about its \lngth\ grading. It extends to a functor
between the corresponding homotopy categories of complexes by
acting on individual chain modules:
\ee
\label{ef.7a}
\xymatrixnocompile@C=2.5cm{\Kmplv{\PMFWn} \ar[r]^-{\CvlMFv{\cdot}} & \Kmplv{\MFW}
}.
\eee

In order to reveal the inner structure of a \Ps\ and its relation
to an $n$-\xxcl\ of \mfs, we split the module and the \tdiff\ of a \Ps\
$\xaA$ according to the \lngth-grading:
\ee
\label{ef.8}
\xaA = \bopbkzn \xaAbk,\quad\bdegl\xaAbk = \bxk,\qquad
\xDA = \sum_{k\geq 0}
\yXk,\quad
\degl \yXk = k.
\eee
We further split the \tdiff\ $\xDA$ according to the domain and range
with respect to the splitting of $\xaA$:
\ee
\label{ef.9}
\xDA = \sum_{\bxk\leq\bxl} \yXbkl,\qquad\yXbkl\in
\HomR(\xaAbk,\xaAbl),
\eee
so that
\ee
\label{ef.9x}
\yXk = \sum_{\substack{\bxl,\bxlp\in\ZZ^n\\ \xtlv{\bxlp-\bxl} = k }}\yXblpl.
\eee

Each $\xaAbk$, equipped with a \tdiff\ $\yXbkk$,
is a \mf. The $\yR$-module map $\yXbkl$ has the
internal $\Zt$-degree 1.

If we split both sides of the basic relation
\ee
\label{ef.10}
\xDA^2 = \xW\,\id_{\xaA}
\eee
according to the \tlngth, then the \lngth-0 part says that
$\yXz^2=\xW\,\id_{\xaA}$, or equivalently
\ee
\label{ef.11}
\yXbvv{\bxk}{\bxk}^2 = \xW\,\id_{\xaAbk}.
\eee
This means that the modules $\xaAbk$ together with the \tdiffs\
$\xDbk=\yXbvv{\bxk}{\bxk}$ form \mfs\ of $\xW$. We call them
\emph{\cnst} \mfs\ of the \Ps\ $\xaA$.

Let us introduce a differential
\ee
\label{ef.12}
\xdz = \ztcm{\yXz}{\cdot},
\eee
acting on $\EndRngl(\xaA)$. Then the \lngth-$k$ ($k\geq 1$) part of \ex{ef.10}
can be put in a form
\ee
\label{ef.13}
\xdz\yXk = - \sum_{l=1}^{k-1} \yXv{k-l}\yXl.
\eee
In particular, the \lngth-1 part says that $\xdz\yXo=0$, which
means that the homomorphism
$\yFv{\bxk,i}=\yXbvv{\bxk}{\bxk+\bxe_i}$ is a morphism
between the \adj\ \mfs\ $\xaAbk\gszto$ and $\xaAv{\bxk+\bxe_i}$.
Moreover, the \lngth-2 part of \eex{ef.13} implies that
\ee
\label{ef.14}
\ztac{\yFi}{\yFj}\mphq 0,\qquad\text{where}\quad\yFi =
\sum_{\bxk}\yFbki,
\eee
which means that the morphisms $\yFi$ can be interpreted as
differentials of an $n$-\xxcl\ over the homotopy category of \mfs. This complex is formed by placing
the \mfs\ $\xaAbk$ at the
nodes of an $n$-dimensional lattice $\ZZ^n$ in accordance with
their multi-indices $\bxk$, so that each pair of adjacent \mfs\ $\xaAbk$ and $\xaAv{\bxk+\bxe_i}$
is connected by a differential $\yFbki$ (see the diagram inside
the dotted box of \ex{ef.15}).

We refer to
$\yFi$ as \emph{\prmr} homomorphisms (or differentials), since they determine
the $n$-\xxcl\ of \mfs, and we call
$\yXk$ ($k\geq 2$) \emph{\scnd} homomorphisms, because their role
is in part to `correct' the consequences of the distinction between the module chain
differentials for which the relation $\ztac{\yFi}{\yFj}=0$
holds, and the \mf\ chain differentials satisfying \ex{ef.14}.

\begin{proposition}
If two \Pss\ $\xaA$ and $\xaAp$ are isomorphic in the category
$\PMFWn$, then there exist homotopy equivalences
\ee
\label{ef.14a}
\xaAbk \cchq\xaAbk\p
\qquad
\text{in $\MFW$}
\eee
making the following diagrams commutative:
\ee
\label{ef.14ax}
\vcenter{
\xymatrixnocompile{
\xaAbk \ar[r]^-{\yFbki}
\ar[d]^-{\cchq}
&
\xaAv{\bxk+\bxe_i}
\ar[d]^-{\cchq}
\\
\xaAbk\p \ar[r]^-{\yFbki\p}
&
\xaAv{\bxk+\bxe_i}\p
}
}
\eee
\end{proposition}
\proof
If a \hteq\ between $\xaA$ and $\xaAp$ within $\PMFWn$ is
established by the homomorphisms
\ee
\xymatrixnocompile@1{\xaA \ar@<0.3ex>[r]^{f} & \xaB \ar@<0.3ex>[l]^{f\p} },
\eee
then the \hteq s\rx{ef.14a} are established by their \lngth-0
parts $f_0$ and $f\p_0$.\qed

We denote a \Ps\ built upon an $n$-\xxcl\ by encircling it with a
dotted frame:
\def\xFF{[F.]}
\ee
\label{ef.15}
\vcenter{
\xymatrixnocompile{
& \vdots \ar[d] & \vdots \ar[d]
\\
\cdots \ar[r] &
\xaAbk \ar[r]^-{\yFbki} \ar[d]
&
\xaAv{\bxk+\bxe_i}
\ar[r] \ar[d]
& \cdots
\\
& \vdots  & \vdots &
\save {[u].[ulll].[l].[uul]}
*[F.]\frm{} \restore
}
}
\eee
Since the \Ps\ might not be determined by this complex uniquely, we
may also add the arrows carrying the \scnd\ morphisms to the
picture. The \cnvl\ of a \Ps\ is denoted by a solid frame:
\ee
\label{ef.15a}
\vcenter{\xymatrixnocompile{
& \vdots \ar[d] & \vdots \ar[d]
\\
\cdots \ar[r] &
\xaAbk \ar[r]^-{\yFbki} \ar[d]
&
\xaAv{\bxk+\bxe_i}
\ar[r] \ar[d]
& \cdots
\\
& \vdots  & \vdots &
\save {[u].[ulll].[l].[uul]}*[F]\frm{} \restore
}
}
\eee

\subsection{Building a \Ps\ from the bottom up}
A \Ps\ $\xaA$ can be constructed from the bottom up in two stages.
At first we choose its \cnst\ \mfs\ $\xaAbk$. At the
second stage we solve the equations\rx{ef.13} inductively on the
value of $k$ starting at $k=1$. The inductive procedure is
possible, since the \rhs of \ex{ef.13} contain only the
homomorphisms of \lngth\ up to $k-1$.

Let us find out how various choices involved in this build-up
affect the emerging \Ps.
At the first stage we choose  $\xaAbk$, $\bxk\in\ssetZ$, where $\ssetZ\subset\ZZn$ is a finite
subset.
%
\begin{theorem}
\label{th.hps}
Let $\xaA$ be a \Ps\ presented by \mfs\ $\xaAbk$, $\bxk\in\ssetZ$
and homomorphisms $\yXm$, $m\geq 1$.
Given a family of \hteqt\ \mfs 
%
\ee
\label{ef.16}
\xaApbk\cchq\xaAbk,\qquad\bxk\in\ssetZ,
\eee
there exists a system of homomorphisms $\yXpm$, $m\geq 1$ such that the
\Ps\ $\xaAp$ presented by $\xaApbk$ and $\yXpm$ is \hteqt\ to $\xaA$.
\end{theorem}
In other words, the set of equivalence classes of all \Pss\ based on
\mfs\ $\xaAbk$ is determined by \hteq\ classes of $\xaAbk$.

\proof
Let $\crssetZ$ denote the number of elements in $\ssetZ$. We will
prove the theorem by induction over $\crssetZ$. If $\crssetZ=1$,
then the claim is obvious. Suppose that the claim is true for sets
of $k$ elements. Consider a \Ps\ $\xaA$ consisting of $k+1$
constituent modules. Since the set $\ssetZ$ is finite,
it contains an element $\bxk$ such that
$\bxl\not\geq \bxk$ for all $\bxl\in\ssetZm=\ssetZ\setminus\{\bxk\}$.
Non-negativity of the \tdiff\ $\xD$ of $\xaA$ implies that
\ee
\label{ef.16a}
\yXbkl = 0,\qquad \bxl\in\ssetZm.
\eee
Let $\xaAm$ be the \Ps\ presented by \mfs\ $\xaAbl$,
$\bxl\in\ssetZm$ and homomorphisms $\yXblpl$,
$\bxl,\bxlp\in\ssetZm$. The sum
$f = \sum_{\bxl\in\ssetZm} \yXblk$ defines a morphism between
$\xaAm$ and $\xaAbk$. In view of \ex{ef.16a}, the \Ps\ cone of $f$
is isomorphic to $\xaA$:
\ee
\label{ef.16a1}
\CnPSv{f}\cong\xaA.
\eee
By the induction
assumption there exists a \Ps\ $\xaAm\p\cchq\xaAm$, which is based
upon the \mfs\ $\xaAbl\p$, $\bxl\in\ssetZm$.
Fix an equivalence homomorphism $g:\;\xaAm\p\rightarrow\xaAm$
and set $f\p=fg$.
The cone of $f\p$ is a \Ps\ based upon \mfs\ $\xaAbl\p$, $\bxl\in\ssetZ$,
so the claim of the theorem follows from the \hteq\ of the cones
$\CnPSv{f}\cchq\CnPSv{f\p}$ and from the isomorphism\rx{ef.16a1}.
\qed

At the second stage we choose the solutions of \eex{ef.13}
inductively over $k$ starting with $k=1$.
Let us fix the \mfs\ $\xaAbk$ and a sequence of homomorphisms $\yXo,\dotsc,\yXkmo$
which solve the equations\rx{ef.13} up to  the \lngth\ $k-1$. Consider
the equation\rx{ef.13} at the \lngth\ $k$. Note that its \rhs is $\xdz$-closed,
however it is not necessarily $\xdz$-exact, so not all choices of
homomorphisms up to \lngth\ $k-1$ may lead to a \Ps.
For a solution $\yXk$ of the equation\rx{ef.13}, let $\ycAXk$
denote the set of \hteq\ classes of \Ps, which can be built upon
the sequence of homomorphisms $\yXo,\dotsc,\yXkmo,\yXk$.

Let $\yXk$ and $\yXpk$ be two solution of \ex{ef.13}.
Since the \rhs of this equation is fixed, then
$\xdz\yXk=\xdz\yXbk$, so
the
difference $\DyXk=\yXpk-\yXk$ is $\xdz$-closed:
\ee
\label{ef.17}
\xdz\,\DyXk = 0.
\eee
\begin{theorem}
\label{th.hl}
If $\DyXk$ is  $\xdz$-exact (that is, if there exists a
homomorphism $\yYk$ ($\degZt\yYk=0$, $\degl\yYk=k$) such that $\DyXk=\xdz\yYk$),
then $\ycAXk=\ycAXpk$.
\end{theorem}
\proof
We will prove the inclusion $\ycAXk\subset\ycAXpk$, the other
inclusion $\ycAXpk\subset\ycAXk$ is established similarly.

Suppose that a \Ps\ $\xaA\in\ycAXk$ is presented by a \tdiff\
$\xDA = \sum_{l\geq 0}\yXl$. A \nng\ homomorphism $\yY = \id +
\yYk$
is invertible:
\ee
\label{ee.yinv}
\yYinv = \id + \sum_{i=1}^{\infty}(-1)^i \yYk^i.
\eee
Hence it establishes an isomorphism between the \Pss\ with the
\tdiff\ $\xDA$ and with the conjugated \tdiff\
$\xDAp = \yYinv\xDA\yY$.
It is easy to see that since $\yY$ and $\yYinv$ are equal to
identity up to the terms of \lngth\ $k$, then the low \lngth\
components of $\xDAp=\sum_{l\geq 0}\yXlp$ are
\ee
\label{ee.llcm}
\yXp_0 = \yX_0,\quad\dotsc\quad,\yXp_{k-1}=\yX_{k-1},\quad\yXpk = \yXp +
\xdz\yYk.
\eee
Hence the \Ps\ $\xaA$ is also an element of $\ycAXpk$.\qed

Let us choose a particular solution $\yXkz$ of \ex{ef.13}. Then, in view of \ex{ef.17}, for any solution
$\yXk$ of that equation we can define a class
\ee
\label{ecr.2}
[\yXk]=\yXk-\yXkz\in\Extmfo_{k}(\xaAz,\xaAz),
\eee
where
$\Extmfo_{k}$ is the \tlngth\ $k$ summand of $\Extmfo$ and
$\xaAz$ is a \lngth-\grdd\ \mf, which is a direct sum of
\mfs\ $\xaAbk$:
\ee
\label{ecr.1}
\xaAz = \bopbkzn \xaAbk
\eee
(in other words, its module is $\xaA$ and
its \tdiff\ is just $\yXz$, which is the \lngth-0 part of $\xD$).
Theorem\rw{th.hl} says that
the set $\ycAXk$
depends only on the
class $\xeclXk$.

Let us split a solution $\yXk$ into components according to
\ex{ef.9x}.
%
%
An individual component $\yXblpl$ satisfies
the equation
\ee
\label{ef.20}
\xdz\yXblpl = - \sum_{\bxl<\bxm<\bxlp}\yXbvv{\bxm}{\bxlp}
\yXbvv{\bxl}{\bxm}
\eee
and defines a class
\ee
\label{ecr.3}
\xeclXlpl = \yXblpl -\yXblplz \in \Extmfo(\xaAbl,\xaAblp).
\eee
According to Theorem\rx{th.hl}, the set $\ycAXk$ is determined by all classes $\xeclXlpl$.
\begin{corollary}
\label{cor.ext}
If for a pair of multi-indexes $\bxl,\bxlp\in\ZZ^n$, $\xtlv{\bxlp-\bxl} =
k$ we have
\ee
\label{ef.22}
\Extmfo(\xaAbl,\xaAblp) = 0,
\eee
then the set $\ycAXk$ does not depend on the choice of a solution
to \ex{ef.20}.
\end{corollary}

\subsection{Corner \sbss\ and \fcss}

Fix a subset $\xSset\in\ZZn$.
For a
\lngth-\grdd\ module $\xaA$, an \emph{\Sctv\ submodule} is defined
as
\ee
\label{ek.1}
\xaA\supset\AsbsS = \bigoplus_{\bxi\in\xSset}\xaAbi.
\eee
For a homomorphism $\yf\in\HomR(\xaA,\xaB)$ between two
\lngth-\grdd\ modules, an \emph{\Sctv\ homomorphism}
$\yfS\in\HomR(\AsbsS,\BsbsS)$ is defined as
\ee
\label{ek.2}
\yfS = \sum_{\bxi,\bxj\in\xSset} \yfbij.
\eee

A subset $\xSset\in\ZZn$ is called \emph{\xlcnv} (or simply \xcnv)
if for any triplet of multi-indices $\bxi,\bxj\in\xSset$, $\bxk\in\ZZn$ such that
$\bxi<\bxk<\bxj$, it turns out that $\bxk\in\xSset$. If $\xaA$ is
\Ps\ with a \tdiff\ $\xD$, then it is easy to verify that $\xDS^2 = \xW
\idS$, so that $\AsbsS$ is a \Ps\ with a \tdiff\ $\xDS$.

\begin{theorem}
Fix a \xcnv\ subset $\xSset\in\ZZn$. The map $\sxhxS$,
which associates $\AsbsS$
to a \Ps\ $\xaA$ and the morphism $\yfS$ to a \Ps\
morphism $\yf$, is an \enfn
\ee
\label{ek.3}
\xymatrixnocompile{
\PMFWn \ar[r]^-{\sxhxS} & \PMFWn
}.
\eee
\end{theorem}
We call $\sxhxS$ an \emph{\Sctv\ functor}. It behaves properly
under the tensor product functor\rx{ef.23a}.
\begin{theorem}Fix two \xcnvssts\ $\xSset\in\ZZn$ and
$\xSsetp\in\ZZnp$. Let $\xaA$ and $f$ be an object and a morphism
of $\PMFWn$, and let $\xaB$ and $g$ be an object and a morphism of
$\PMFWnp$. Then
%
\ee
\label{ef.33}
\AsbsS\otimes\BsbsSp =
\ABsbsS,\qquad
\yfS\otimes\ygSp = (\yf\otimes\yg)_{\xSSpset}.
\eee
\end{theorem}

For a subset $\xSset\subset\ZZn$ denote the complementary subset
$\bSset=\ZZn\setminus\xSset$. Suppose that a pair of complementary
subsets $\xSset,\bSset$ satisfies the property
\ee
\label{ek.x3}
\bxi
\not<
\bxj\qquad\text{for all $\bxi\in\xSset$, $\bxj\in\bSset$.}
\eee
Then both $\xSset$ and $\bSset$ are \xcnv.
Moreover, $\AsbsS$ is a \sbs\ of
$\xaA$, and $\AsbxS$ is a \fcs: $\AsbxS\cong\xaA/\AsbsS$.
Denote by $\ychsSin$ and $\ychbSout$ the injection and surjection
homomorphisms
\ee
\label{ef.28}
\xymatrixnocompile{\AsbsS\ar[r]^-{\ychsSin} & \xaA},\qquad
\xymatrixnocompile{\xaA\ar[r]^-{\ychbSout} &
\AsbxS=\xaA/\AsbsS}.
\eee
\begin{theorem}
\hfill
\begin{enumerate}
\item[a.]
The homomorphisms $\ychsSin$ and $\ychbSout$ are \Ps\ morphisms.
If $\xaA\cchq\xaAp$ in $\PMFWn$ then
$\ychsSin\mphq\ychsSin\p$ and $\ychbSout\mphq\ychbSout\p$.
\item[b.]
The morphisms $\ychin$ and $\ychout$ are natural, that is, for any
\Ps\
morphism $\xymatrixnocompile@1{\xaA\ar[r]^-{\yf}&\xaB}$ the following diagrams are commutative:
\ee
\label{ef.31a}
\vcenter{\xymatrixnocompile{\AsbsS \ar[r]^-{\ychsSin}
\ar[d]_-{\yfS}
& \xaA \ar[d]^-{\yf}
\\
\BsbsS \ar[r]^-{\ychsSin} &\xaB
}
}\qquad\qquad
\vcenter{\xymatrixnocompile{\xaA \ar[r]^-{\ychbSout} \ar[d]_-{\yf}
&
\AsbxS \ar[d]^-{\yfbS}
\\
\xaB \ar[r]^-{\ychbSout} & \BsbxS
}
}
\eee
\item[c.]
If $\xSset\p,\bSset\p\subset\ZZn$ is another pair of complementary
subsets satisfying the condition\rx{ek.x3} and
$\xSset\p\subset\xSset$, then
\ee
\label{ef.32}
\ychbSout\ychsSpin = 0.
\eee
\end{enumerate}
\end{theorem}

The natural morphisms\rx{ef.28}
behave properly under the tensor product functor\rx{ef.23a}
\begin{theorem}For two pairs of
complementary subsets $\xSset,\bSset\in\ZZn$, $\;\;\xSsetp,\bSsetp\in\ZZnp$
\ee
\label{ef.34}
\ychsSin\otimes\ychsSpin = \ychsSSpin,\qquad
\ychbSout\otimes\ychbSpout = \ychsSSpout.
\eee
%
\end{theorem}

Two types of pairs of complementary subsets
$\xSset,\bSset\subset\ZZn$ satisfying the condition\rx{ek.x3}
are particularly important.
For $\bxk\in(\ZZinf)^{\times n}$, where $\ZZinf=\ZZ\cup\{-\infty,+\infty\}$
is an ordered set, define two subsets of $\ZZn$:
\ee
\label{ek.4}
\sbslbk = \{\bxj\in\ZZn\;|\; \bxj\leq\bxk\},
\qquad \sbsgbk =
\{\bxi\in\ZZn\;|\;\bxi\geq\bxk\}.
\eee
The condition\rx{ek.x3} is satisfied by the pair
$\xSset=\sbsgbk$, $\bSset=\ZZn\setminus\sbsgbk$ and by the pair
$\xSset=\ZZn\setminus\sbslbk$, $\bSset=\sbslbk$. We call
$\Asbsgbk$ \emph{a \csbs} and we call $\Asbslbk$ \emph{a
\cfcs}.

\subsection{Graded \Pss}
A special feature of the \mfs\ appearing in the categorification of
the \sotn\ \Kfp\ is their \qgrd, for which
\ee
\label{ef.35}
\dgq\xW = 4\xN + 2,\qquad \dgq\xD = 2\xN + 1.
\eee
The \prmr\ homomorphisms appearing in convolutions producing the
\mfs\ $\gg$ and $\hG$ (associated with an \opgr\ $\xg$ and a
\clgr\ $\xG$)
will be \sdlmrps, so according to
\ex{eb2.16xa5},
\ee
\label{ef.36}
\dgq\yXv{\bxe_i} = 2\xN.
\eee
To guarantee that
$\dgq\xDA = 2\xN+1$, we introduce a (relative) \qdgr\ shift
for the constituent modules $\xaAbk\grxhv{\xtlv{\bxk}}$
and impose a condition on the \qdgr\ of the \scnd\ homomorphisms $\yXv{\bxk}$:
\ee
\label{ef.37}
\dgq\yXk = 2\xN-k+1.
\eee
As a result, the space that describes the
dependence of the set of \Pss\ on
the choice of
solutions
$\yXblpl$
to \ex{ef.20} is the subspace of $\Extmfo(\xaAbl,\xaAblp)$ of \qdgr\ $2\xN-\xtlv{\bxlp-\bxl}+1$:
\ee
\label{ef.38}
\Extmfo_{\qgrsv{2\xN-\xtlv{\bxlp-\bxl}+1}}(\xaAbl,\xaAblp).
\eee

\subsection{Examples of \Pss\ and \cnvls}
\label{ss.expss}
We will consider three simple examples of \Pss\ of \lngth\ 0, 1 and
2 with $n=1$.

A \Ps\ $\xaA$ of \lngth\ 0 is just an ordinary \mf, and its \cnvl\
is equal to that \mf:
\ee
\label{ef.39}
\boxed{\xaA} = \xaA.
\eee

Consider a morphism
\ee
\label{ef.40}
\xymatrixnocompile@1{\xaA\ar[r]^{\yF} & \xaB},\qquad\degZt\yF=0.
\eee
In order to construct a \Ps\ out of it, we shift the \Ztdgr\ of
$\xaA$:
\ee
\label{ef.41}
\xymatrixnocompile{\Pbdv{\mA\gszto}{\wF}{\mB}}\;\;.
\eee
Its \cnvl\ is isomorphic to the cone of the
morphism\rx{ef.40}:
\ee
\label{ed.42}
\xymatrixnocompile{\Cbdv{\mA\gszto}{\wF}{\mB}} = \CnMFv{\wF},
\eee
and the \cnvls\ of natural morphisms
\ee
\label{ed.24a4}
\xymatrixnocompile{
\mB \ar[r]^-{\ychin}
&
\Cbdv{\mA\gszto}{\wF}{\mB}
\ar[r]^-{\ychout}
&
\mA\gszto
},\qquad
\ychin=\begin{pmatrix}0\\ \id_{\xaB} \end{pmatrix},
\quad
\ychin=\begin{pmatrix} \id_{\xaA\gszto} & 0 \end{pmatrix}
\eee
form two sides of the exact triangle related to\rx{ef.40}.

Now consider a chain of two morphisms
\ee
\label{ea2.24a5}
\xymatrixnocompile{
\wA \ar[r]^-{\wF} & \wB \ar[r]^-{\wG} & \wC
},\qquad\degZt\wF=\degZt\wG=0,\quad \wG\wF\sim 0.
\eee
The latter condition implies that there exists a homomorphism
$\wX\in\HomR(\wA,\wC)$ such that
\ee
\label{ea2.24a6}
\degZt\wX=1,
\qquad
\wG\wF = - \ztac{\xDAC}{\wX},
\eee
where $\xDAC= \xDA + \xDC$.
If we shift the \Ztdgr\ of $\wB$, then we can form a \Ps\ based on
the chain\rx{ea2.24a5}
\ee
\label{ea2.24a7}
\xymatrixnocompile{ \PbtwABCX }\;\;.
\eee
Its module is a sum of modules
\ee
\label{ea2.24a8}
\wA\oplus\wB\gszto\oplus\wC
\eee
and its \tdiff\ is a sum of the \mf\ \tdiffs\ and of the
homomorphisms $\wF$, $\wG$ and $\wX$:
\ee
\label{ea2.24a9}
\xD = \xDA - \xDB + \xDC + \wF + \wG + \wX.
\eee
The effect of the choice of $\wX$ satisfying \ex{ea2.24a6} on the
class of this \Ps\ is (relatively) parametrized by the elements of
$\Extmfo(\wA,\wC)$, so if
\ee
\label{ea2.24a9b}
\dim\Extmfo(\wA,\wC)=0,
\eee
then the \Ps\rx{ea2.24a7} is determined by the
complex\rx{ea2.24a5} uniquely.
%

The convolutions of two natural morphisms
\ee
\label{ea2.24a10}
\xymatrixnocompile{
\wC \ar[r]^-{\ychin} & \CbtwABCX \ar[r]^-{\ychout} & \wA
},\qquad
\ychin = \begin{pmatrix} 0 \\ 0 \\ \id_{\wC} \end{pmatrix},\quad
\ychout = \begin{pmatrix} \id_{\wA} & 0 & 0 \end{pmatrix}
\eee
form a chain complex:
\ee
\label{ea2.24a11}
\ychout\ychin=0.
\eee

\section{Categorifiction of the \fvtx, graphs and tangles}
\subsection{A \cnvl\ of the chain of \sdlmrps}
We are going to build a \Ps\ upon
the \chcl\rx{ea2.40y3}, following the steps outlined in subsection\rw{ss.expss}.
According to \ex{eb2.16a2}, the subspace
$\Extmfo(\hlpar,\hlhor)$ is trivial, hence the \Ps\ is determined
uniquely. The \scnd\ homomorphism $\yX$ of\rx{ea2.24a6} is
provided by \eex{ey.13x},\rx{ea2.43}. We denote the resulting \Ps\
as
\ee
\label{eg.1}
\tlvrb = \cxy{\Pssdl}\;\;.
\eee
The \qgrd\ shifts in the first and the third modules
are required in order to ensure that the \tdiff
\ee
\label{ea2.43b2}
\xDconv = \xDpar - \xDvirt + \xDhor + \yF + \yG +
\yX
\eee
has a homogeneous degree in accordance with \ex{ea2.13a2}. For the
same reason the homomorphism $\yX$ must have a degree prescribed
by \ex{ef.37}, which is that of \ex{ea2.40a1}.
Two markers $\cdot$ in the notation $\tlvrb$ are needed, because
this \Ps\ is not
invariant under the $\ntydeg$ rotation.

The uniqueness of the \sdlmrp\ up to a constant factor and the
fact that \sdlmrps\ determine the \Ps\rx{eg.1} uniquely allows us to
define similar \Pss\ for any
triplet of distinct \tarc\ graphs:
\ee
\label{ey.22}
\cxy{\Psgmp}\;\;,\qquad \xg\neq \xgp\neq \xgpp.
\eee
In particular,
\ee
\label{eg.2}
\tlvbh = \cxy{\Pssdlh}\;\;,
\eee
where $\yFp$ and $\yGp$ are corresponding \sdlmrps.

\begin{theorem}
\label{th.rotinv}
The \cnvl\ of the \Ps\rx{eg.1} is invariant under the $\ntydeg$
rotation, that is
\ee
\label{ea2.51}
\vcenter{\xymatrix{\Cnsdl}}\;\;\; \cchq\;\;\;
\vcenter{\xymatrix{\Cnsdr}}\;.
\eee
\end{theorem}
\proof
The \Ps\rx{eg.1} factors similarly to
\mfs\rx{ea2.23a}:
\ee
\label{ea2.43b4}
\tlvrb = \zKcmn\otR\tlvrbp,\qquad\!\!\!\!\text{where}\quad
\tlvrbp = \cxy{\Pssdrprp}\;.
\eee
Therefore, it suffices to prove the rotation invariance for the \prpr\
\cnvls
\ee
\label{ea2.51a1}
\CvlMFv{\tlvrbp} \cchq \CvlMFv{\tlvbhp}.
\eee
%

The \cnvl
\ee
\label{eg.2x}
\cxy{\Cnsdrprp}
\eee
in the \lhs of \ex{ea2.51a1} has an
explicit presentation
\ee
\label{ea2.43b2x1}
\xymatrixnocompile{
\yRso \ar[r]^-{\zPcn} & \yRsz \ar[r]^-{\zQcn} & \yRso
},
\eee
where the homomorphisms $\zPcn$ and $\zQcn$ are presented by the
matrices
\ee
\label{ea2.43b2x2}
\zPcn =
\begin{pmatrix}
\zPpar & 0 & 0
\\
\yFo & -\zQvirt & 0
\\
\yXo &\yGz &  \zPhor
\end{pmatrix}
,\qquad\zQcn =
\begin{pmatrix}
\zQpar & 0 & 0
\\
\yFz & -\zPvirt & 0
\\
\yXz & \yGo &  \zQhor
\end{pmatrix}
\eee
in the basis corresponding to the splitting
\ee
\label{ea2.43b2x3}
\yRsz = \yRtzpar\oplus\yRtovirt\oplus\yRtzhor
,\qquad \yRso =
\yRtopar\oplus\yRtzvirt\oplus\yRtohor
.
\eee
The diagonal entries of the matrices\rx{ea2.43b2x2}
are given by formulas\rx{ea2.35}--\rxw{ea2.37}.
We simplify this presentation\rx{ea2.43b2x2} by using the fact
that the matrices $\yFz$ and $\yGz$ have unit entries. Consider
the isomorphism of two \mfs\ (the bottom one being\rx{ea2.43b2x1})
established by the commutative diagram
\ee
\label{ea2.52}
\xymatrix{
{\yRo\oplus\yRo\oplus\yRfo} \ar[rrr]^-{\zWfp\oplus 1\oplus\zPrcn} \ar[dd]^-{\xfiso} &&&
{\yRz\oplus\yRz\oplus\yRfz} \ar[rrr]^-{1\oplus\zWfp\oplus\zQrcn} \ar[dd]^-{\xfisz} &&&
{\yRo\oplus\yRo\oplus\yRfo}            \ar[dd]^-{\xfiso}
\\
\\
{\yRso} \ar[rrr]^-{\zPcn} &&& \yRsz \ar[rrr]^-{\zQcn} &&& {\yRso}
}
\eee
where
\begin{align}
\label{ea2.53}
\zPrcn & =
\begin{pmatrix}
\zpt & \zpo & 0 & 0
\\
\zqt & -\zqo & \zqt\zro & \zqo\zro
\\
\zqo\zC & -\zqt & -\zqt\zrt & -\zqo\zrt
\\
\zpo\zC & \zpt & \zpt\zrt + \zpo\zro\zC & -(\zpo\zrt + \zpt\zro)
\end{pmatrix}
,
\\
\label{ea2.53a}
\zQrcn & =
\begin{pmatrix}
\zqo\zrt+\zqt\zro & \zpo\zrt & \zpo\zro & 0
\\
\zqt\zrt+\zqo\zro\zC & -\zpt\zrt & -\zpt\zro & 0
\\
-\zqt & \zpt & -\zpo & \zqo
\\
\zqo\zC & \zpo\zC & -\zpt & -\zqt
\end{pmatrix}
\end{align}
and
\ee
\label{ea2.54}
\xfisz =
\begin{pmatrix}
0 & 0 & 1 & 0 & 0 & 0
\\
1 & 0 & 0 & \zrt & \zro & 0
\\
0 & -\zro & 0 & 1 & 0 & 0
\\
0 & \zrt & 0 & 0 & 1 & 0
\\
0 & 1 & 0 & 0 & 0 & 0
\\
0 & 0 & 0 & 0 & 0 & 1
\end{pmatrix}
,\qquad
\xfiso =
\begin{pmatrix}
\zpo & 0 & 1 & 0 & 0 & 0
\\
-\zpt & 0 & 0 & 1 & 0 & 0
\\
1 & 0 & 0 & 0 & 0 & 0
\\
0 & 1 & 0 & 0 & -\zqt & -\zqo
\\
0 & 0 & 0 & 0 & 1 & 0
\\
0 & 0 & 0 & 0 & 0 & 1
\end{pmatrix}
\eee
It is easy to verify that the latter matrices are invertible and hence define an
isomorphism $\xfis$ of two \mfs. The top \mf\ decomposes into a
sum of three \mfs, the first two being
\ee
\label{ea2.55}
\xymatrixnocompile{
\yRo \ar[r]^-{\zWfp} & \yRz \ar[r]^-{1} & \yRo
},\qquad
\xymatrixnocompile{
\yRo \ar[r]^-{1} & \yRz \ar[r]^-{\zWfp} & \yRo
}
\eee
and the third (reduced) one being a rank-$(4,4)$ \mf
\ee
\label{ea2.56}
\xymatrixnocompile{
\yRfo \ar[r]^-{\zPrcn} & \yRfz \ar[r]^-{\zQrcn} & \yRfo.
}
\eee
The \mfs\rx{ea2.55} are \cntrb, hence the \lhs\ of \ex{ea2.51a1}
is \hteqt\ to\rx{ea2.56}.

We construct the presentation of the \rhs \cnvl\ of
\ex{ea2.51a1} by using the \enfn\ action of the symmetric group
$\Sth$.
Namely, $\sXtth\in\Sth$ transforms the triplet of graphs
\ee
\label{ea2.57}
\sXtth
\lrbc{\lxpar,\lxvir,\lxhor} = \lrbc{\lxhor,\lxvir,\lxpar}.
\eee
Therefore
\ee
\label{ea2.58}
\hXtth
\lrbc{\tlvrbp}=\tlvbhp\;\;.
\eee
The functor $\hXtth$ acts by switching the variables $\zxbp$ and $\zxbq$, while leaving
$\zxbr$ intact, so
the \rhs of \ex{ea2.51a1} can be presented similarly to\rx{ea2.56}
as
\ee
\label{ea2.60}
\xymatrixnocompile{
\yRfo \ar[rr]^-{\zPprcn} && \yRfz \ar[rr]^-{\zQprcn} && \yRfo,
}
\eee
where
\begin{align}
\label{ea2.61}
\zPprcn & =
\begin{pmatrix}
\zqt & \zqo & 0 & 0
\\
\zpt & -\zpo & \zpt\zro & \zpo\zro
\\
\zpo\zC & -\zpt & -\zpt\zrt & -\zpo\zrt
\\
\zqo\zC & \zqt & \zqt\zrt + \zqo\zro\zC & -\zqo\zrt - \zqt\zro
\end{pmatrix},
\\
\label{ea2.62}
\zQprcn & =
\begin{pmatrix}
\zpo\zrt+\zpt\zro & \zqo\zrt & \zqo\zro & 0
\\
\zpt\zrt+\zqo\zro\zC & -\zqt\zrt & -\zqt\zro & 0
\\
-\zpt & \zqt & -\zqo & \zpo
\\
\zpo\zC & \zqo\zC & -\zqt & -\zpt
\end{pmatrix}.
\end{align}
Now the isomorphism\rx{ea2.51a1} is established by the commutative
diagram
\ee
\label{ea2.63}
\xymatrix{
{\yRfo} \ar[rr]^-{\zPrcn} \ar[dd]^-{\xgiso}
&& {\yRfz} \ar[rr]^-{\zQrcn} \ar[dd]^-{\xgisz}
&& {\yRfo} \ar[dd]^-{\xgiso}
\\
\\
{\yRfo} \ar[rr]^-{\zPprcn} && {\yRfz} \ar[rr]^-{\zQprcn} &&
{\yRfo}
}
\eee
in which the isomorphism matrices $\xgisz$, $\xgiso$ are
\ee
\label{ea2.64}
\xgisz =
\begin{pmatrix}
0 & 1 & 0 & 0
\\
1 & 0 & 0 & 0
\\
0 & 0 & 0 & 1
\\
0 & 0 & 1 & 0
\end{pmatrix}
,
\qquad
\xgiso =
\begin{pmatrix}
1 & 0 &\zro & 0
\\
0 & -1 & 0 & \zro
\\
0 & 0 & -1 & 0
\\
0 & 0 & 0 & 1
\end{pmatrix}.
\eee
\qed

\subsection{The \fvtx\ \mf\ and the crossing complex}

Now that we have proved the isomorphism\rx{ea2.51}, we choose the
\mf\ associated with the \fvtx\ graph to be the \cnvl\rx{ea2.51}:
%
\begin{align}
\label{ea2.65}
\hlver & = \CvlMFv{\tlvrb}=
\vcenter{\xymatrix{\Cnsdl}}
\\
\intertext{or, equivalently,}
\label{ea2.65x}
\hlver & =
\CvlMFv{\tlvbh}
=
\vcenter{\xymatrix{\Cnsdr}}\;\;.
\end{align}
%
This \mf\ has a
$90^\circ$-rotation symmetry as the picture of the graph $\lxver$ suggests.

The \Ps\rx{eg.1} has \mfs\ $\hlparsh$ and $\hlhorsh$ as its left
and right \csbss, hence they are connected to it by \ntrl\
morphisms of \Pss
%
\begin{gather}
\label{ey.23}
\xymatrix{{\hlhorsh} \ar[r]^-{\ychin} &
\Pssdl
\ar[r]^-{\ychout} &
{\hlparsh}},
\\
\label{ey.24}
\ychin = \begin{pmatrix} 0\\ 0\\ \id \end{pmatrix},\qquad
\ychout = \begin{pmatrix} \id & 0 & 0 \end{pmatrix},
\\
\label{ey.25}
\degZt\ychin=\degZt\ychout=0,\qquad
\dgq\ychin=\dgq\ychout =0,\qquad
\ychout\ychin=0.
\end{gather}
The latter relation allows us to
interpret the diagram\rx{ey.23} as a \chcl\
in the homotopy category of \Ps\ complexes $\Kmplv{\PMFWfo}$ and
define
\ee
\label{eg.4}
\begin{split}
\tlncr
& =
\lrbc{
\xymatrix{{\hlhorsh} \ar[r]^-{\ychin} &
{\tlvrb} \ar[r]^-{\ychout} &
{\hlparsh}
}
}\gszto
\\
&= \lrbc{
\xymatrix{ {\hlhorsh} \ar[r]^-{\ychin} &
\Pssdl
\ar[r]^-{\ychout} &
{\hlparsh} }
}\gszto.
\end{split}
\eee
Now
we define the complex of \mfs\ $\hlncr$ for the \eltr\ \crs\ as the \cnvl\
functor\rx{ef.7a} applied to the \rhs of \ex{eg.4}:
%
\ee
\label{ey.26}
\begin{split}
\hlncr & =
\CvlMFv{\tlncr} =
\CvlMFv{
\cxy{ {\hlhorsh} \ar[r]^-{\ychin}
&
{\tlvrb}\ar[r]^-{\ychout}
&
{\hlparsh}
}
}\gszto
\\
& =
\lrbc{
\xymatrix{
{\hlhorsh} \ar[r]^-{
\ychin
}
&
\Cnsdl
\ar[r]^-{
\ychout
}
&
{\hlparsh}
}
}\gszto.
\end{split}
\eee
%
This definition coincides with \ex{e2.3x2} and the resulting
complex should be considered as an object in the homotopy category
of \mf\ complexes $\Kmplv{\MFWf}$.

The rotational invariance of $\hlver$ allows us to
present a $\ntydeg$ rotation version of \ex{ey.26} as
\ee
\label{ey.27}
\begin{split}
\hlpcr & =
\CvlMFv{\tlpcr} =
\CvlMFv{
\cxy{ {\hlparsho} \ar[r]^-{\ychin}
&
{\tlvbh} \ar[r]^-{\ychout}
&
{\hlhorsho}
}
}
\\
& =
\lrbc{
\xymatrix{
{\hlparsho} \ar[r]^-{
\ychin
}
&
\Cnsdr
\ar[r]^-{
\ychout
}
&
{\hlhorsho}
}
}\gszto.
\end{split}
\eee

\subsection{A categorification complex}
\label{ss.smplf}
\subsubsection{A \mf\ of an \opgr\ as a \cnvl}

Let $\xg$ be an \opgr\ with $\yym$ \fvtcs. Its \mf\ $\gg$ is
constructed according to the standard procedure outlined in
subsection\rw{ss.outcat}: we cut all edges of $\xg$ connecting \fvtcs, so
that $\xg$ splits into \eltr\ pieces: circles, arcs and \eltr\ \fvtx\
graphs. To each piece we associate its \mf\ according to
\eex{eh.5},\rx{ea2.16c} and either \ex{ea2.65} or \ex{ea2.65x}.
Then we join the pieces back together with the help of the
operations\rx{ex.8} and\rx{ex.9}.

The structure of $\gg$ can be described more precisely if this
\mf\ is constructed in a
slightly different way. The \cnvl\
functor commutes with the gluing operations\rx{ex.8} and\rx{ex.9}.
Hence, instead of assembling the \cnvls\rx{ea2.65} and\rx{ea2.65x} we
can
assemble the original \Pss\rx{eg.1} and\rx{eg.2}. The assembly of \eltr\ pieces into
the graph $\xg$ results in an $\yym$-dimensional \Ps\ $\tgin$,
which will be simplified by \hteq\ transformations within the
category of \Pss, thus yielding $\tg$. Finally we will apply the \cnvl\ functor:
\ee
\label{el.1}
\gg = \CvlMFv{\tg}.
\eee
This procedure can be presented schematically by the diagram:
\ee
\xymatrixnocompile@C=2cm{
*+[F]\txt{graph \\ $\xg$}
\ar[r]^-{\txt{\tiny assembly\\ \tiny of elementary\\ \tiny Postnikov \\ \tiny systems}}
&
*+[F]\txt{initial\\Postnikov\\system\\ $\tgin$}
\ar[r]^-{\txt{\tiny simplification}}
&
*+[F]\txt{simplified\\Postnikov\\system\\ $\tg$}
\ar[r]^-{\txt{\tiny convolution}}
&
*+[F]\txt{matrix\\factorization\\ $\gg$}
}
\eee

For $\bxk\in\ZZ^{\yym}$ ($-1\leq k_i\leq 1$)
let $\xgbk$
be an \opgr\ constructed
from $\xg$ by
replacing each \fvtx\ $\lxver_i$ ($1\leq i\leq \yym$) with either $\lxpar$, or $\lxvir$, or
$\lxhor$ in accordance with the value of $k_i$ ($-1$, or $0$, or $1$) and the
choice\rx{ea2.65} or\rx{ea2.65x} for the presentation of that
\fvtx.
Initially, the \cnst\ \mfs\ of the
$\yym$-dimensional \Ps\ $\tgin$
are \mfs\ $\gginbk$, constructed by applying the assembly
procedure to the \tarc\ \mfs\ $\hlpar$, $\hlvir$ and $\hlhor$,
which combines their graphs together into $\xgbk$. The homomorphisms between
$\gginbk$
within $\tgin$ are exactly the homomorphisms of the \Pss\rx{eg.1} and\rx{ey.22}
(up to appropriate negative signs related to the \Ztgrdng) applied
to the corresponding factors of the tensor products which
form $\gginbk$. In particular, the \prmr\ morphisms of $\tgin$ are
\sdlmrps\ reconnecting the \tarc\ subgraphs sitting inside the
graphs $\xgbk$.

The graphs $\xgbk$ have a simple structure: they are disjoint
unions of circles and arcs. Therefore each \cnst\ \mf\ $\gginbk$ is
\hteqt\ to the \mf\ $\ggbk$ which, according to our rules, is just
a tensor product of circle spaces $\nCcir=\JrWp$ and \oarc\
\mfs\ of the type\rx{ea2.16c}. If we replace $\gginbk$ with
$\ggbk$, then the \prmr\ homomorphisms between $\ggbk$ will be \hteqt\
to \sdlmrps, \smcl\ \sdlmrps\rx{ey.29} and \clsd\
\sdlmrps\rx{eh.7}, applied to the \mfs, corresponding to circle and \oarc\ graphs,
which are the
connected components of the graphs $\xgbk$. Hence, according to
theorems\rw{th.hps} and\rw{th.hl}, the system $\tgin$ is \hteqt\ in the category of
\Pss\ to
a \emph{\smplfd} system $\tg$, whose \cnst\ \mfs\ are $\ggbk$ and whose
\prmr\ homomorphisms are \sdlmrps\ as well as their \smcl\ and
\clsd\ versions, acting on the \mfs\ of connected components of
the graphs $\xgbk$.

The \scnd\ homomorphisms of the \Ps\ $\tgin$ originate from the
homomorphisms $\yX$ of\rx{eg.1} and\rx{ey.22}, hence they all have
\lngth\ 2. However, the \hteq\ transformation from $\tgin$ to
$\tg$ may generate \scnd\ homomorphisms of higher \lngth. The only
restriction that we can impose on their structure so far follows
from Theorem\rw{th.arcgl}, according to which the spaces
$\spexlpl$
are trivial for \scnd\ homomorphisms of even \lngth, so the
\hteq\ class of $\tg$ does not depend on the choice of a
particular solution of \ex{ef.20}.


Thus we construct the graph \mf\ $\gg$ by first reducing the \Ps\
$\tgin$ to $\tg$ and then applying the \cnvl\ functor $\gg =
\CvlMFv{\tg}$, which just turns the dotted frame boxes into solid
ones.

As an example of the simplification procedure, consider the initial \Ps\ for the graph
$\smgrthbv{\verver}$ as depicted in \fg{fg.em.1} (this graph will appear in the proof of the second
Reidemeister move invariance).
\def\lxarr{\ar@/^2.5pc/@{-->}[rr]}
\def\lxadd{\ar@/^2.5pc/@{-->}[dd]}
\begin{figure}
\ee
\nonumber
\tgrtbv{\bverver}_{\xinit} =
\vcenter{
\xymatrix@C=2cm@R=2cm{
{\hhorpar}
\lxarr^-{}="t"
\lxadd
\ar[r] \ar[d]
&
{\hvirpar}
\lxadd
\ar[r] \ar[d]
&
{\hparpar}
\lxadd^-{}="r"
\ar[d]
\\
{\hhorvir} \lxarr
\ar[r] \ar[d]
&
{\hvirvir} \ar[r] \ar[d]
&
{\hparvir} \ar[d]
\\
{\hhorhor}
\lxarr
\ar[r]
&
{\hvirhor} \ar[r]
&
{\hparhor}
\save
[].[uull]."t"."r"*[F.]\frm{}
\restore
}
}
%
\eee
\caption{An example of an initial \Ps}
\label{fg.em.1}
\end{figure}
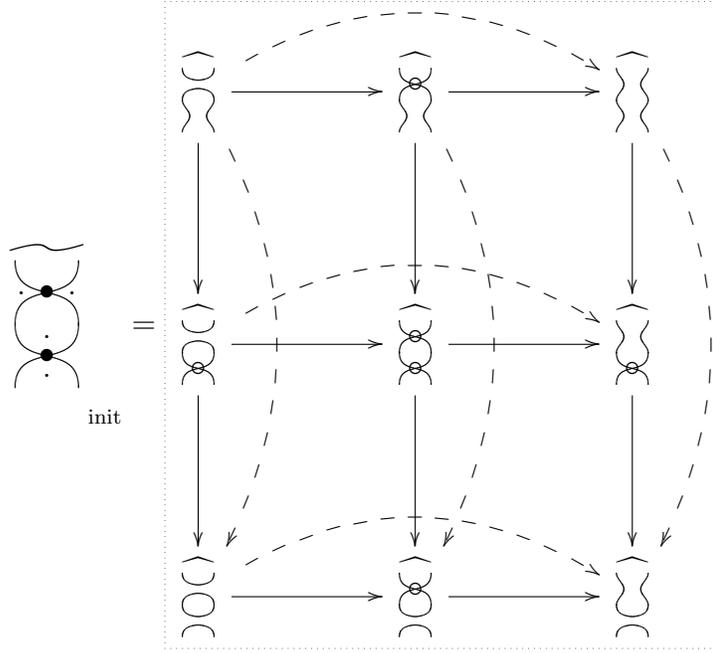
There the solid arrows denote the \sdlmrps, while the dashed arrows
denote the \scnd\ morphisms $\yX$ of\rx{eg.1}. The horizontal
arrow morphisms act on \mfs\ coming from the top \fvtx, while vertical
ones act on the bottom \mfs.
Note that we omitted the degree shifts
$\gsztv{\xtlv{\bxk}}\grxhv{\xtlv{\bxk}}$, which should
accompany the \cnst\ \mfs\ $\ggbk$ in the \Ps\ of \fg{fg.em.1}.

The simplification of the \cnst\ \mfs\ and \prmr\
(that is, \lngth\ 1)
morphisms leads to the \hteqt\ \Ps\ in \fg{fg.ey.56}.
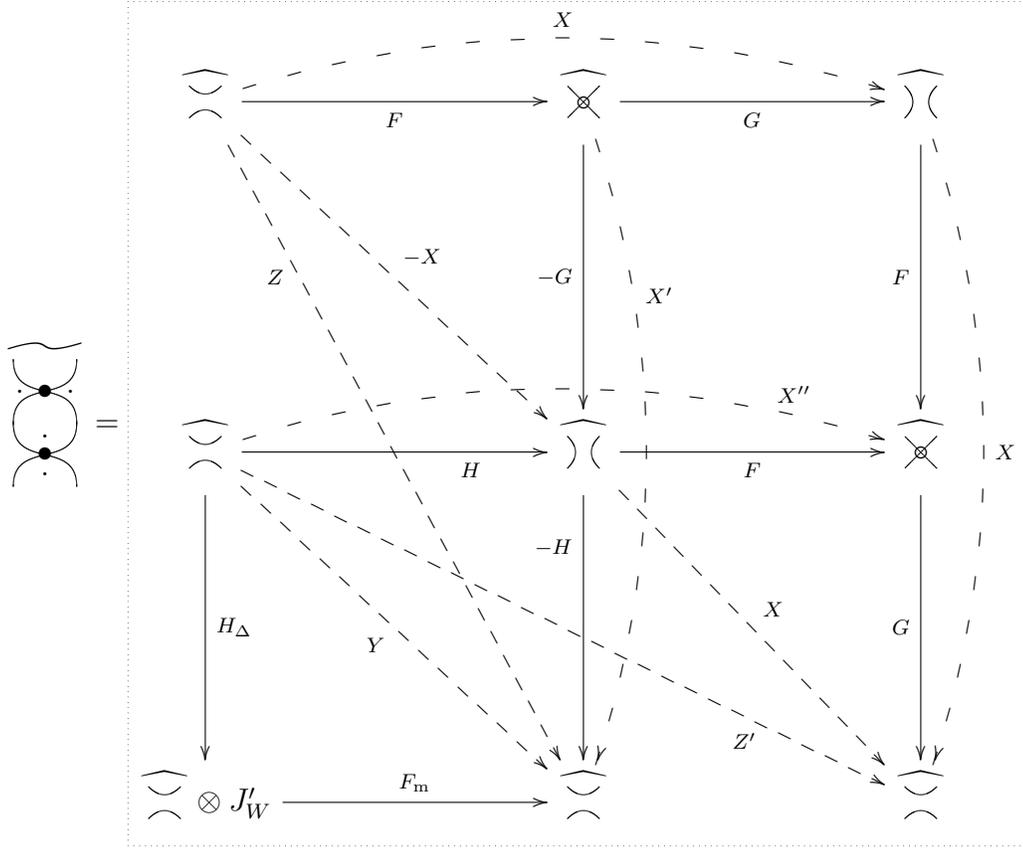
\begin{figure}
\ee
\nonumber
\tgrtbv{\bverver} =
\vcenter{
\xymatrix@C=100pt@R=100pt{
{\hlhor}
\ar[r]_-{\yF}
\ar@{-->}[ddr]_-(0.2){\yZ}
\ar@/^2pc/@{-->}[rr]^-{\yX}
="xarr"
\ar@{-->}[dr]^-{-\yX}
&
{\hlvir}
\ar[r]_-{\yG}
\ar[d]_-{-\yG}
\ar@/^2pc/@{-->}[dd]^-(0.25){\yXp}
&
{\hlpar} \ar[d]_-{\yF}
\ar@/^2pc/@{-->}[dd]^{\yX}="yarr"
\\
{\hlhor}
\ar[r]_-(0.75){\yH}
\ar[d]^-{\yHD}
\ar@{-->}[rd]_-{\yY}
\ar@{-->}[rrd]_-(0.8){\yZp}
\ar@/^2pc/@{-->}[rr]^-(0.85){\yXpp}
&
{\hlpar}
\ar[r]_-{\yF}
\ar[d]_-(0.2){-\yH}
\ar@{-->}[dr]^-{\yX}
&
{\hlvir}
\ar[d]_-{\yG}
\\
{\hlhrcir}
\ar[r]^-{\yFm}
&
{\hlhor}
&
{\hlhor}
\save
{[].[ll]}.[uu]!C="cbb"
\restore
\save
"cbb"."xarr"!C="cbbv"
\restore
\save
"cbb"."yarr"!C="cbbw"
\restore
\save
"cbbv"."yarr"!C="cbbt"*[F.]\frm{}
\restore
}
}
\eee
\caption{An example of a simplified \Ps}
\label{fg.ey.56}
\end{figure}
In that diagram we omitted two
\smcl\
saddle homomorphisms which are
equal to zero (see \ex{ea2.86}).
The homomorphisms $\yF$, $\yG$ and $\yH$ are
\sdlmrps\ related by the action of the symmetry group
$\Sf$ permuting
the legs. $\yX$ is the \scndh\rx{ea2.43} and the homomorphisms $\yXp$, $\yXpp$ are
related to it
by the leg permutation symmetry action. $\yY$, $\yYp$ and $\yZ$
are other \scndhs\ (we will not attempt to determine them).

The origin of the \smplfd\ \cnst\ \mfs\ and \prmr\ homomorphisms
in the diagram of \fg{fg.ey.56} is obvious.
Let us explain the structure of the \scnd\ homomorphisms. We
follow the step-by-step procedure outlined after the
Corollary\rw{cor.ext}. Namely, for any given $k\geq 2$ we first
choose the homomorphisms of \lngth\ $k$ which satisfy the
conditions\rx{ef.20} and then add to them the
representatives of the spaces\rx{ef.38}.
\begin{itemize}
\item[($k=2$)]
We included all homomorphisms that might be needed to satisfy the conditions\rx{ef.20}.
The spaces $\spexlpl$ are trivial, so we do not have to add any
extra homomorphisms. In particular, we set the homomorphism
$\xymatrix@1{{\hlvir} \ar[r] &{\hlvir}}$ from $\xmaAv{0,-1}$ to
$\xmaAv{1,0}$ equal to zero.

\item[($k=3$)]
We marked by $\yZ$ and $\yZp$ the homomorphisms needed to satisfy the conditions\rx{ef.20}.
As for the spaces\rx{ef.38}, it follows easily from \ex{ee.mrarc}
that they are non-trivial only for two homomorphisms that have already
been marked, hence we set homomorphisms
$\xmaAv{-1,-1}\rightarrow \xmaAv{1,0}$ and
$\xmaAv{0,-1}\rightarrow\xmaAv{1,1}$ equal to zero.
%

\item[($k=4$)]
We set the only homomorphism to zero. This choice satisfies \ex{ef.20},
and the corresponding space\rx{ef.38} is trivial
\end{itemize}
\subsubsection{A complex of a \tngl}
Let $\gt$ be a \tngl\ with $\yyn$ crossings. Its \ccomp\ $\hgt$ is
a complex of \mfs\
constructed by assembling the \eltr\ crossing complexes\rx{eg.4}.
For $\xbr\in\ZZ^{\yyn}$ ($-1\leq \xr_i\leq 1$) let $\xgbr$ be an
\opgr\ constructed from $\gt$ by replacing each crossing
$\lxncr_i$ ($1\leq i\leq \yyn$) with either $\lxhor$, or $\lxver$,
or $\lxpar$ in accordance with the value of $\xr_i$ ($-1$, or $0$, or
$1$). Then the complex $\hgt$ is formed by the \mfs\ $\ggbr$, each
placed at the homological degree $\xbra$, the differential being the
sum of \ntrl\ morphisms $\ychin$ and $\ychout$ of\rx{ey.26}.

We are going to construct a simplified \hteqt\ version of the
complex $\hgt$ by first constructing the complex of \Pss, then
simplifying them as in the previous subsection and finally applying
the \cnvl\ functor to each \Ps\ forming the complex. Thus we start
by applying the assembly process not to the \mfs\ $\hlncr$ but to
their underlying \Ps\ complexes $\tlncr$ of \ex{eg.4}. The result
is a complex $\tgtin$ consisting of \Pss\ $\tginbr$,
its differential being the sum of \ntrl\ morphisms $\ychin$ and
$\ychout$.

Since the left and the right \Pss\ (consisting of a single \cnst\
\mf) in the complex\rx{eg.4} are \xxcts\
of the middle \Ps\ (the left one is a \cfcs\ and
the right one is a \csbs), then all \Pss\ $\tginbr$ are \xxcts\ of
the middle system $\tginbz$. Let $\tgbz$ be a \hteqt\ \smplfc\
of $\tginbz$. Then its \xxcts\ $\tgbr$ are the \smplfc\ of $\tginbr$,
and all \ntrl\ morphisms between them remain intact. Thus we
obtain a complex of \Pss\ $\tgt$, which is a \hteqt\ \smplfd\ version of
$\tgtin$: its \Pss\ are $\tgbr$ and its differential is still the
sum of \ntrl\ morphisms $\ychin$ and $\ychout$. Finally, we apply
the \cnvl\ functor: $\hgt = \CvlMFv{\tgt}$.

As an example, let us consider the \tngl\ $\arttw$ appearing in the second
Reidemeister move. The initial form of its \Pss\ complex is
is depicted in \fg{fg.ey.53y}.
\def\lxarr{\ar@/^2.2pc/@{-->}[rr]}
\def\lxadd{\ar@/^2pc/@{-->}[dd]}
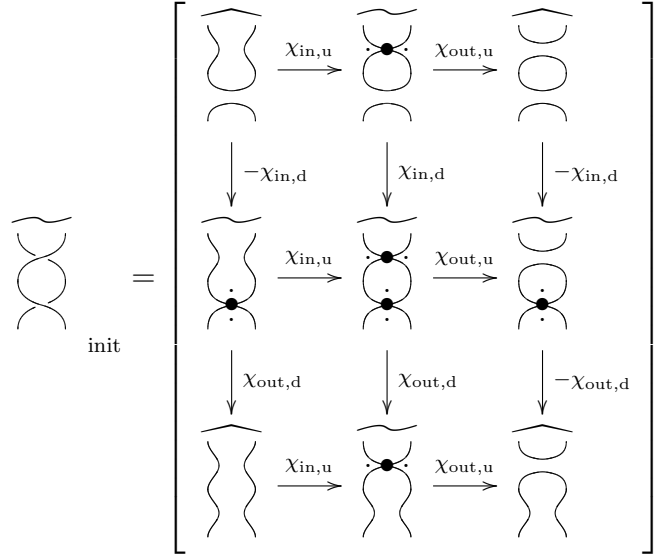
\begin{figure}
\begin{equation}
\nonumber
\xymatrixnocompile{\xtnegpos}_{\xinit}
=
\lrbs{\vcenter{
\xymatrix{
{\xhparhor}
\ar[r]^-{\ychinu} \ar[d]^-{-\ychind}
&
{\xtverhor}
\ar[r]^-{\ychoutu} \ar[d]^-{\ychind}
&
{\xhhorhor}
\ar[d]^-{-\ychind}
\\
{\xtparver}
\ar[r]^-{\ychinu} \ar[d]^-{\ychoutd}
&
{\xtverver}
\ar[r]^-{\ychoutu} \ar[d]^-{\ychoutd}
&
{\xthorver}
\ar[d]^-{-\ychoutd}
\\
{\xhparpar}
\ar[r]^-{\ychinu}
&
{\xtverpar}
\ar[r]^-{\ychoutu}
&
{\xhhorpar}
}
}
}
\end{equation}
\caption{Initial complex for the tangle diagram in the second Reidemeister move}
\label{fg.ey.53y}
\end{figure}
A more detailed version of the same complex, in which we substituted
the \Pss\rx{eg.1} for the \fvtcs, is presented in \fg{fg.ey.53}.
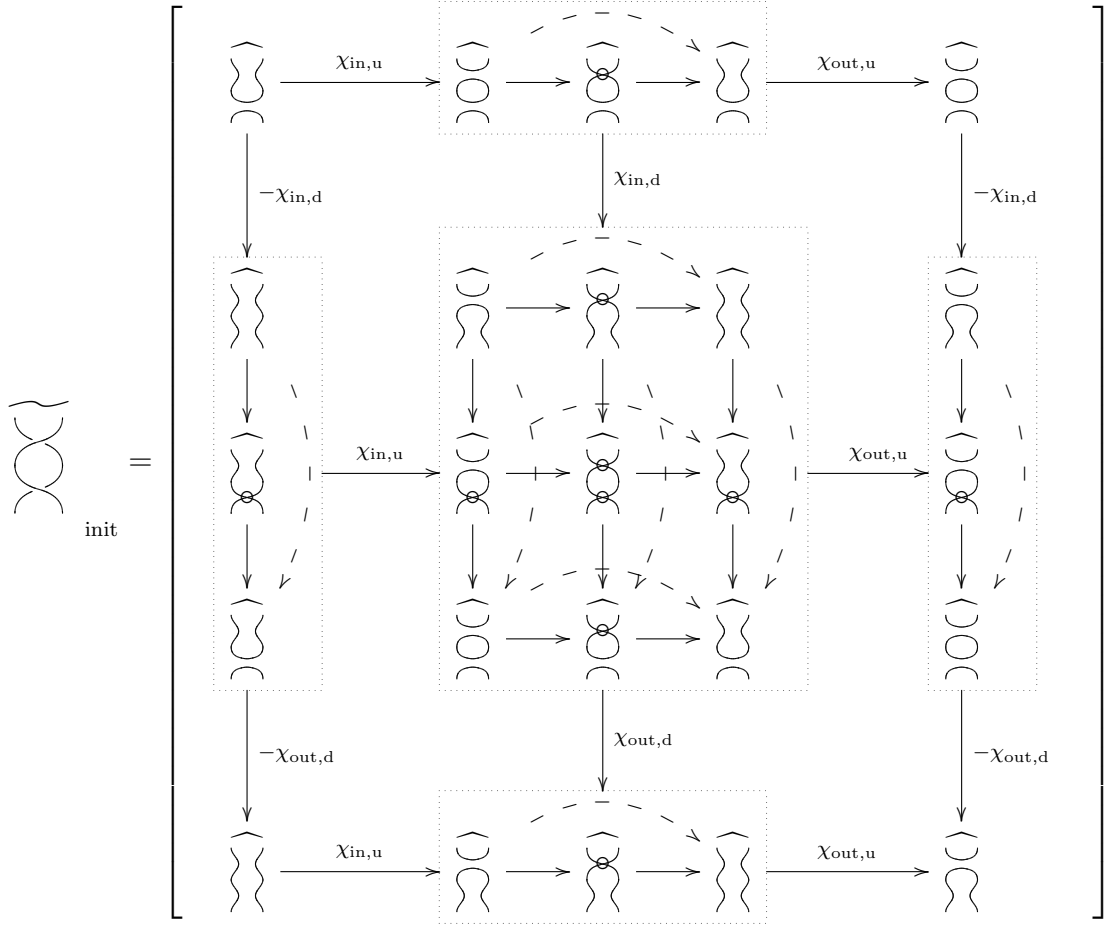
\begin{figure}
\ee
\nonumber
\xymatrixnocompile{\xtnegpos}_{\xinit}
=
\lrbs{
\quad
\vcenter{
\xy
\xymatrix"aa"
{{\hparhor}}
\POS(30,0)
\xymatrix"ba"
{
{\hhorhor}\lxarr^-{}="mu" \ar["aa"];[]^-{\ychinu} \ar[r] & {\hvirhor} \ar[r] &
{\hparhor}
\save [].[ll]."mu"*[F.]\frm{} \restore
}
%
\POS(95,0)
\xymatrix"ca"{ {\hhorhor}\ar["ba"rr];[]^-{\ychoutu} }
\POS(0,-30)
\xymatrix"ab"{
{\hparpar}
\lxadd^-{}="ml"
\ar["aa"];[]^-{-\ychind}
\ar[d]
\\
{\hparvir} \ar[d]
\\
{\hparhor}
\save [].[uu]."ml"*[F.]\frm{}
\restore
}
%
\POS(30,-30)
\xymatrix"bb"{
{\hhorpar}
\lxarr^-{}="mmo"
\lxadd
\ar[r] \ar[d]
&
{\hvirpar}
\lxadd
\ar[r] \ar[d]
\ar["ba"];{[]."mmo"}^-{\ychind}
&
{\hparpar}
\lxadd^-{}="mmt"
\ar[d]
\\
{\hhorvir}
\lxarr
\ar[r] \ar[d]
\ar
"ml";[]^-{\ychinu}
&
{\hvirvir} \ar[r] \ar[d]
&
{\hparvir} \ar[d]
\\
{\hhorhor}
\lxarr
\ar[r]
&
{\hvirhor} \ar[r]
&
{\hparhor}
\save
[].[uull]."mmo"."mmt"!C="tst"*[F.]\frm{}
\restore
}
%
\POS(95,-30)
\xymatrix"cb"{
{\hhorpar}
\lxadd^-{}="mr"
\ar[d]
\ar["ca"];[]^{-\ychind}
\\
{\hhorvir} \ar[d]
\ar
"mmt";[]^{\ychoutu}
\\
{\hhorhor}
\save
[].[uu]."mr"*[F.]\frm{}
\restore
}
%
\POS(0,-105)
\xymatrix"ac"{
{\hparpar} \ar["ab"dd];[]^-{-\ychoutd}
}
\POS(30,-105)
\xymatrix"bc"{
{\hhorpar}
\lxarr^-{}="dm"
\ar[r]
\ar["ac"];[]^-{\ychinu}
&
{\hvirpar} \ar[r]
\ar["bb"dd];"dm"^-{\ychoutd}
&
{\hparpar}
\save
[].[ll]."dm"*[F.]\frm{}
\restore
}
%
\POS(95,-105)
\xymatrix"cc"{
{\hhorpar}
\ar["cb"dd];[]^-{-\ychoutd}
\ar["bc"rr];[]^-{\ychoutu}
}
\endxy
}
\quad
}
\eee
\caption{Detailed structure of the initial complex for the tangle diagram in the second Reidemeister move}
\label{fg.ey.53}
\end{figure}
There the horizontal morphisms $\ychin$ and $\ychout$ are related
to the top crossing, while the vertical ones are related to the
bottom crossing.

In order to simplify the complex of \fg{fg.ey.53}, first, we simplify the
central \Ps\ to the form of \fg{fg.ey.56}. Then the other \Pss\ are
simplified as its sub- and \fcss. The result is the complex of \fg{fg.el.2}
\begin{figure}
\begin{multline}
\nonumber
\xymatrixnocompile{\xtnegpos} =
\\
\nonumber
\lrbs{
\vcenter{
\xy
\xymatrix"aa"
{{\hlhor}}
\POS(35,0)
\xymatrix@C=40pt"ba"
{
{\hlhrcir}\ar["aa"];[]^-{\ychinu} \ar[r]^-{\yFm}
& {\hlhor} 
& {\hlhor}
\save [].[ll]!C="cba"  \restore}
\POS*\frm{.}
\POS(125,0)
\xymatrix"ca"{ {\hlhrcir}\ar["ba"rr];[]^-{\ychoutu} }
\POS(0,-30)
\xymatrix"ab"{
{\hlpar} \ar@/^2pc/@{-->}[dd]^-{\yX}="xab"
\ar["aa"];[]^-{\ychind}
\ar[d]_-{\yF}
\\
{\hlvir} \ar[d]_-{\yG} 
\\
{\hlhor}
\save
{[].[uu]}."xab"!C*[F.]\frm{}
\restore
}
%
%
%
\POS(35,-30)
\xymatrix@C=40pt"bb"{
{\hlhor} \ar[r]
\ar@{-->}[ddr]
\ar@/^2pc/@{-->}[rr]
^{}="xarr"
\ar@{-->}[dr]
&
{\hlvir} \ar[r] \ar[d]
\ar@/^2pc/@{-->}[dd]
&
{\hlpar} \ar[d]
\ar@/^2pc/@{-->}[dd]^{}="yarr"
\\
{\hlhor}
\ar[r] \ar[d]
\ar@{-->}[rd]
\ar@{-->}[rrd]
\ar@/^2pc/@{-->}[rr]
&
{\hlpar} \ar[r] \ar[d]
\ar@{-->}[dr]
&
{\hlvir} \ar[d]
\\
{\hlhrcir} \ar[r]
&
{\hlhor}
&
{\hlhor}
\save
{[].[ll]}.[uu]!C="cbb"
\restore
\save
"cbb"."xarr"!C="cbbv"
\restore
\save
"cbb"."yarr"!C="cbbw"
\restore
\save
"cbbv"."yarr"!C="cbbt"*[F.]\frm{}
\restore
\ar"xab";"cbb"^-{\ychinu}
\ar"cba";"cbbv"^-{\ychind}
}
%
%
%
\POS(125,-30)
\xymatrix"cb"{
{\hlhor}
\ar["ca"];[]^{\ychind}
\\
{\hlhor} \ar[d]^-{\yHD}
\\
{\hlhrcir}
\save
[].[uu]!C="ccb"*[F.]\frm{}
\restore
\ar"cbbw";"ccb"^-{\ychoutu}
}
%
\POS(0,-105)
\xymatrix"ac"{
{\hlpar} \ar["ab"dd];[]^-{\ychoutd}
}
\POS(35,-105)
\xymatrix"bc"{
{\hlhor} \ar[r]^-{\yF}
\ar@/^2pc/@{-->}[rr]^{\yX}="hbc"
\ar["ac"];[]^-{\ychinu}
&
{\hlvir} \ar[r]^-{\yG}
&
{\hlpar}
\save {[].[ll]!C="cbc"}."hbc"!C="cbct"*[F.]\frm{} \restore
\ar"cbb";"cbct"^-{\ychoutd}
}
%
\POS(125,-105)
\xymatrix"cc"{
{\hlhor}
\ar["cb"dd];[]^-{\ychoutd}
\ar["bc"rr];[]^-{\ychoutu}
}
\endxy
}
}
\end{multline}
\caption{Simplified complex of the tangle diagram in the second
Reidemeister move}
\label{fg.el.2}
\end{figure}
in which the differential is still a sum of \ntrl\ morphisms relating
a \Ps\ to its \csbss\ and \fcss.

\subsubsection{A complex of a link}
A link $\cL$ is a \tngl\ without \lgs.
Let $\xGr$ be the \clgrs\ constructed by replacing the
crossings of $\cL$ with \tarc\ or \fvtx\ \eltr\ graphs according
to the values of $\xri$. Further, let $\xGbrk$ be the \clgrs,
constructed by replacing the \fvtcs\ of $\xGr$ with \tarc\ graphs
according to the values of $k_i$. Since a graph $\xGbrk$ is
closed, it is a disjoint union of $\anbk$ circles. Hence the spaces
$\hGbrk$ are tensor powers of the unknot spaces $(\JrWp)^{\otimes\anbk}$, and the
\prmr\ morphisms of the \smplfd\ \Ps\ $\tGbr$ are \clsd\
\sdlmrps\ of\rx{eh.7}, which are \mltps\ $\xmlt$, if two circles
coalesce into one, \cmlts\ $\xcmlt $, if one circle splits into two, or
zero if one circle reconnects into another circle.

Since the \tdiffs\ in the \cnst\ \mfs\ $(\JrWp)^{\otimes\anbk}$ are zero,
particular solutions to the equations\rx{ef.20} may be chosen to
be zero. Therefore
the presence of \scndhs\ in a \smplfd\ \Ps\ $\tGbr$ is
due exclusively to non-trivial spaces $\spexlpl$. This means
that $\tGbr$ contains \scnd\ morphisms of only odd \lngth.

As an example of a link complex, let us consider the \smplfd\
version of the Hopf link complex in \fg{fg.hpfcm} (see Section\rw{Sect1}; note that we have
omitted the \scndhs\ in that diagram). All \cnst\ graphs of its
\cnvls\ are either single or double circles, so the corresponding
\smplfd\ \mfs\ are either $\JrWp$ or $\JrWpdbl$. Now we choose the
\prmr\ and \scndhs\ of the \Ps\ underlying the middle \cnvl, while the differentials of other \cnvls\
are determined by the fact that they are \cnvls\ of the \xxcts\ of
the middle \Ps. The result is the complex depicted in
\fg{fg.ee.hpff}, in which
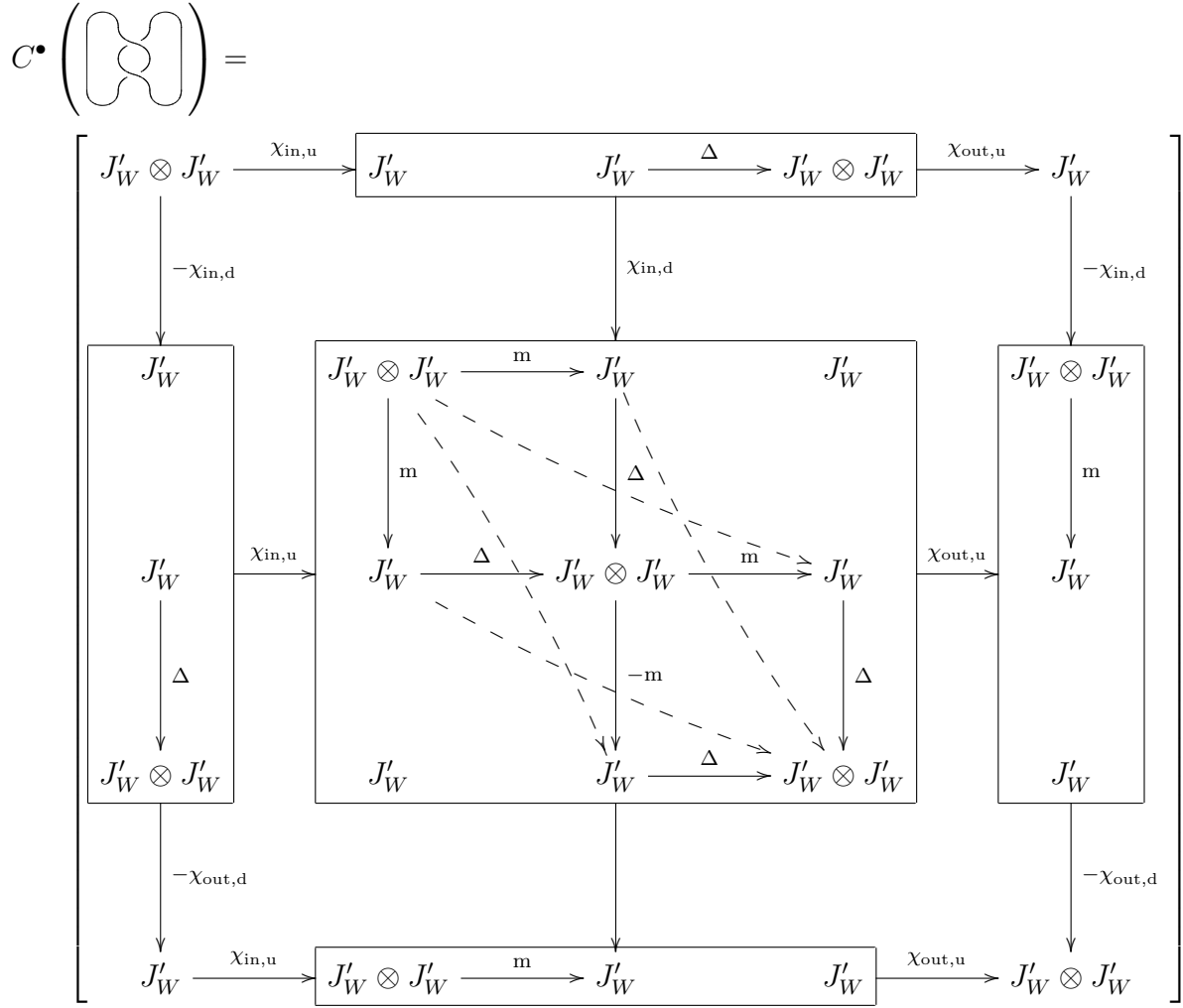
\begin{figure}
\begin{multline}
\nonumber
\nCdv{\xhpfw} =
\\
\nonumber
\lrbs{
\vcenter{
\xymatrix@C=1.1cm@R=2cm{
{\JrWpdbl}
\save
[]="lu"
\restore
&
{\JrWp}
&
{\JrWp}
\ar[r]^-{\xcmlt}="mua"
&
{\JrWpdbl}
\save
{[].[ll]}!C="mu"
\ar "lu";"mu"^-{\ychinu}
\restore
\save
{"mu"."mua"}!C*[F]\frm{}
\restore
&
{\JrWp}
\ar"mu";[]^-{\ychoutu}
\\
{\JrWp}
&
{\JrWpdbl}
\ar[r]^-{\xmlt}="mmao"
\ar[d]^{\xmlt}
&
{\JrWp}
\ar[d]^-{\xcmlt}
&
{\JrWp}
&
{\JrWpdbl}
\ar[d]^-{\xmlt}
\\
{\JrWp}
\ar[d]^-{\xcmlt}
&
{\JrWp}
\ar[r]^-{\xcmlt}
&
{\JrWpdbl}
\ar[r]^-{\xmlt}
\ar[d]^-{-\xmlt}
&
{\JrWp}
\ar[d]^-{\xcmlt}="mmat"
&
{\JrWp}
\\
{\JrWpdbl}
&
{\JrWp}
&
{\JrWp}\ar[r]^-{\xcmlt}
&
{\JrWpdbl}
&
{\JrWp}
\save
{[].[uu]}!C="rm"
\ar "1,5";"rm"^-{-\ychind}
\restore
\save
{"2,5"."4,5"}!C*[F]\frm{}
\restore
\save
{"2,2"."4,4"}!C="mm"
\ar "mm";"rm"^-{\ychoutu}
\restore
\save
{"mm"."mmat"."mmao"}!C*[F]\frm{}
\ar "1,3";{"mm"."mmao"}^-{\ychind}
\restore
\save
{"2,1"."4,1"}!C="lm"
\ar "lm";"mm"^-{\ychinu}
\ar "1,1";"lm"^-{-\ychind}
\restore
\save
{"2,1"."4,1"}!C*[F]\frm{}
\restore
\\
{\JrWp}
\ar"lm";[]^-{-\ychoutd}
&
{\JrWpdbl}
\ar[r]^-{\xmlt}="mda"
&
{\JrWp}
&
{\JrWp}
&
{\JrWpdbl}
\save
{"5,2"."5,4"."mda"}!C*[F]\frm{}
\ar"4,3";"5,3"
\ar"5,1";"5,2"^-{\ychinu}
\ar"5,4";"5,5"^-{\ychoutu}
\ar"4,5";"5,5"^-{-\ychoutd}
\ar@{-->}@/^6pt/"2,2";"4,3"
\ar@{-->}@/_6pt/"3,2";"4,4"
\ar@{-->}@/_6pt/"2,2";"3,4"
\ar@{-->}@/_6pt/"2,3";"4,4"
\restore
}
}
}
\end{multline}
\caption{The simplified complex for the Hopf link}
\label{fg.ee.hpff}
\end{figure}
%
we omitted the grading shifts as well as the \prmr\ morphisms which
are zero.

\section{The Reidemeister moves}

\subsection{Excision of a \cntrcn}
\label{ss.exc}
The proof of Reidemeister move invariance will require a
simplification of \mfs\ and their complexes, which goes beyond the
simplification of the \Pss\ described in Subsection\rw{ss.smplf}.
In fact, we will need just one elementary trick, which we call
\emph{an \excis\ of a \cntrcn}.
Let $\catC$ be a triangulated
category (we have two examples in mind: the category of \mfs\
$\MF$ and the homotopy category of their complexes $\Kmplv{\MF}$).
Let $\grseo$ be the translation functor of $\catC$.
\begin{lemma}
\label{l.excc}
Consider two objects $\xaA,\xaB\in\Obv{\catC}$ and a morphism
$\xaA
\xrightarrow{\xaf}\xaB$.
Suppose that either $\xaA$ or $\xaB$ is a zero-object in $\catC$
(that is, it is contractible). Then
\ee
\label{el.3x}
\Cnvv{\catC}{\xaf}\cchq
\begin{cases}
\xaA\grseo,& \text{if $\xaB\cchq 0$,}
\\
\xaB,&\text{if $\xaA\cchq 0$.}
\end{cases}
\eee
\end{lemma}

In our future examples contractible objects will be the cones of
the identity morphism $\Cnvv{\catC}{\id}\cchq 0$.

Now let us assume that $\catC$ is the category of \mfs\ $\MF$ and consider two \mfs\
$\xaC,\xaCp\in\Obv{\MF}$ and two morphisms, relating them to the
cone $\CnMFv{\xaf}$:
\ee
\label{el.3}
\xymatrixnocompile{\xaC \ar[r]^-{\xag} & \CnMFv{\xaf}},\qquad
\xymatrixnocompile{\CnMFv{\xaf} \ar[r]^-{\xagp} & \xaCp}.
\eee
As an $\yR$-module, the cone is a sum
\ee
\label{el.4}
\CnMFv{\xaf} = \xaA\oplus\xaB,
\eee
and we present the homomorphisms $\xag$ and $\xagp$ as a column and a row according
to this decomposition:
\ee
\label{el.5}
\xag = \begin{pmatrix} \xagA \\ \xagB \end{pmatrix},\qquad
\xagp = \begin{pmatrix} \xagpA & \xagpB \end{pmatrix}.
\eee
Since $\xag$ and $\xagp$ are \mf\ morphisms, they commute with the \tdiff, which means that
\begin{align}
\label{el.6}
\zzd \xagA & = 0, & \zzd \xagB & = -\xaf\xagA,
\\
\label{el.7}
\zzdp \xagpB & = 0, & \zzdp\xagpA & = (-1)^{\degZt\xagpB} \xagpB\xaf,
\end{align}
where
\ee
\zzd = \ztcm{\xDA+\xDB + \xDC}{\cdot},\qquad \zzdp=
\ztcm{\xDA+\xDB+\xDv{\xaCp}}{\cdot}.
\eee
The first conditions indicate that $\xagA$ and $\xagB$ are \mf\
morphisms, while the second conditions imply that $\xagB$ and
$\xagpA$ are morphisms iff
\ee
\label{el.8}
\xaf\xagA=0,\qquad\xagpB\xaf = 0.
\eee
The following lemma describes what happens to the
morphisms\rx{el.3} if one of the \mfs\ $\xaA$ or $\xaB$ is
contractible, and therefore, according to \ex{el.3x}, the cone of $\xaf$ is \hteqt\ to
either $\xaA$ or $\xaB$.
\begin{lemma}
\label{l.css}
\begin{enumerate}
\item If the \mf\ $\xaA$ is contractible, then $\cnmff\cchq\xaB$
and
\begin{itemize}
\item[a.]
the morphisms $\cgcnfp$ and $\xymatrixnocompile@1{\xaB\ar[r]^-{\xagpB} &
\xaCp}$ are equivalent in the category $\MF$;
\item[b.] there exists a morphism $\xahB\in\HomMF(\xaC,\xaB)$ such
that the morphisms\linebreak $\cgcnf$ and $\xymatrixnocompile@1{\xaC\ar[r]^-{\xahB} &
\xaB}$ are equivalent in $\MF$;
\item[c.]furthermore, if $\xagB$ is a morphism (that is,
$\xaf\xagA=0$), then the morphisms $\cgcnf$ and
$\xymatrixnocompile@1{\xaC\ar[r]^-{\xagB} &
\xaB}$ are equivalent in $\MF$.
\end{itemize}

\item If the \mf\ $\xaB$ is contractible, then $\cnmff\cchq\xaA$
and
\begin{itemize}
\item[a.] the morphisms $\cgcnf$ and
$\xymatrixnocompile@1{\xaC\ar[r]^{\xagA} &\xaA}$ are equivalent in $\MF$;
\item[b.] there exists a morphism
$\xahpA\in\HomMF(\xaA,\xaCp)$
such that the morphisms\linebreak $\cgcnfp$ and $\xymatrixnocompile@1{\xaA\ar[r]^-{\xahpA} & \xaCp}$
are equivalent in $\MF$.
\item[c.] furthermore, if $\xagpA$ is a morphism (that is, if
$\xagpB\xaf=0$), then the morphisms $\cgcnfp$ and
$\xymatrixnocompile@1{\xaA\ar[r]^-{\xagpA} & \xaCp}$ are equivalent in $\MF$.

\end{itemize}

\end{enumerate}

\end{lemma}

It will be convenient to use the following abbreviated notations
for the sums of graded subspaces of $\JrWp$:
\ee
\label{el.31}
\JrWpmz = \zy\IQ\oplus\bigoplus_{i=1}^{2\xN} \zx^i\IQ,\qquad
\JrWpmtN = \zy\IQ\oplus\bigoplus_{i=0}^{2\xN-1} \zx^i\IQ,\qquad
\JrWpmztN = \zy\IQ\oplus\bigoplus_{i=1}^{2\xN-1} \zx^i\IQ.
\eee
%


\subsection{First Reidemeister move}
We are going to prove Theorem\rw{th.lr1} by establishing the
\hteq\ (up to a degree shift) between the complex of a \knk\ \tngl\ and
the \mf\ of the \oarc\ graph.

\begin{lemma}
\label{l.r1}
The following complexes are \hteqt\ in the category
$\Kmplv{\MFWt}$:
%
%
\ee
\label{ea2.86a}
\hltwpy \cchq \hlhcpy\grsv{-2\xN-1}\gszto\grsemo,
\eee
\end{lemma}

Equation\rx{1.3a1} follows if we tensor-multiply both sides of
\ex{ea2.86a} by the complex associated with the
\tngl\ representing the rest of the link.

\proof
The \mf\ $\hltwpy$ is presented initially by the complex of \mfs\rx{ey.26}, in which
legs 2 and 4 are connected:
\begin{multline}
\label{ey.34}
\hltwpy_{\xinit} =
\lrbc{
\xymatrix{
{\hlhory\grso} \ar[r]^-{\ychpin} &
{\hlvery}
\ar[r]^-{\ychpout} & {\hlpary\grsmo}
}
}\gszto
\\
=\lrbc{
\xymatrix@C=.6cm{
{\hlhory\grso} \ar[rr]^-{\smmpver}
&&
\Cbtwgny \ar[rr]^-{\smmphor}
&&
{\hlpary\grsmo}
}
}\gszto.
\end{multline}

After we simplify the \cnst\ \mfs\ and \prmr\ homomorphisms of the underlying complex of \Pss\ as
described in Subsection\rw{ss.smplf}, the complex\rx{ey.34}
becomes
\ee
\label{ey.38}
\xymatrix{ {\hlhcpy\grso} \ar[r]^-{\smmpver} &
\Cbtany \ar[rr]^-{\smmphor} && {\arnCcir\grsmo}
}.
\eee
We can set the \scndhs\ equal to zero, because this
choice satisfies the condition\rx{ef.20} and
the corresponding space\rx{ef.22} for \lngth-2 homomorphisms is trivial.

Now it remains to apply the \cntrcnex\ procedure in order to
reduce the complex\rx{ey.38} to the \rhs of \ex{ea2.86a}.
The middle \cnvl\ in the complex\rx{ey.38}, which represents the
\smcl\ 4-vertex \mf\ $\hlvery$, splits:
\ee
\label{ey.38a}
\hlvery =
\xymatrix{
\Ctwwv{\arnCcir\grsmo}{\yFm}{\hlhcpy\gszto}
} \;\;\oplus \hlhcpy\grso,
\eee
and the
complex\rx{ey.38} also splits: it is a direct sum of two complexes, the first
being a contractible complex
\ee
\label{ey.39}
\xymatrix{ {\hlhcpy\grso} \ar[r]^-{\id} & {\hlhcpy\grso} }
\eee
and the second one being
\ee
\label{ey.40}
\xymatrix{
\Ctwwv{\arnCcir\grsmo}{\yFm}{\hlhcpy\gszto}
\ar[rr]^-{\stphor} && {\arnCcir\grsmo}.
}
\eee

The tensor product in the \cnvl\ splits
\ee
\label{ey.41}
\arnCcir = \lrbc{\hlhcpy\otimes 1} \oplus
\lrbc{\hlhcpy\otimes\JrWpmz},
%
\eee
and the homomorphism $\yFm$ acts as identity on the first term in this
sum (the term $i=0$ in the sum of \ex{ea2.82x}).
Hence the cone of\rx{ey.40} can be presented in the
form of a `double cone', the inner one being the contractible cone
of the identity homomorphism:
\ee
\label{em.1}
\xymatrix{
{\hlhcpy\otimes\JrWpmz}
\ar[r]
\ar@{}@<0.8cm>[r]^{}="a1"
\ar@{}@<-0.8cm>[r]^{}="a2"
&
{\hlhcpy\otimes 1} \ar[r]^-{\id}="b"
&
{\hlhcpy}
\save
[]+<0.8cm,0cm>*{}="c",
"1,3"."1,2"."b"*[F]\frm{}
\restore
\save
"1,1"."1,3"."b"."a1"."a2"."c"*[F]\frm{}
\restore
}
\eee
(here we omitted the grading shifts).
We \excse\ the contractible cone in accordance with case 2C of Lemma\rw{l.css} and the
original cone becomes
\ee
\label{ey.42x1}
\xymatrix{\Ctwwv{\arnCcir\grsmo}{\yFm}{\hlhcpy\gszto}}
\cchq
\hlhcpy\otimes\JrWpmz\grsmo,
\eee
while the whole complex\rx{ey.40} splits into a direct sum of
complexes
%
%
\ee
\label{ey.43}
\hspace*{-20pt}
\Bigg[
\lrbc{
\hlhcpy\otimes1\;\grsemo
}
\oplus
\lrbc{ \xymatrix{{\hlhcpy\otimes\JrWpmz} \ar[r]^-{\id} &
{\hlhcpy\otimes\JrWpmz}} }
\Bigg]\grsmo.
\eee
The second complex in this sum is contractible and the
first complex consists of only one chain module
\ee
\label{ey.44}
\hlhcpy\otimes 1\;\grsemo\grsmo=\hlhcpy\grsemo\grsv{-2\xN-1}
\eee
(recall that $\dgq 1 = - 2\xN$ in $\JrWp$).
If we substitute it for the complex in the brackets of \ex{ey.34}, then we get
\ex{ea2.86a}.
\qed

In the process of proving \ex{ea2.86a} we also proved
the following formula for the \smcl\ 4-vertex:
\ee
\label{ey.44a}
\hlvery = \bopiztNmo \hlhcpy \grsv{-2\xN+1+2i} \oplus
\hlhcpy\grsmo \oplus\hlhcpy\grso.
\eee
%


\subsection{Second Reidemeister move}
Similar to the first Reidemeister move case, Theorem\rw{th.lr2}
follows from the next
\begin{lemma}
\label{l.r2}
The following objects are \hteqt\ in the category $\Kmplv{\MFWf}$:
\ee
\label{ey.52}
\xymatrixnocompile{ \hbrttw } \cchq
\xymatrixnocompile{ \hbrtutw }\;\;.
\eee
\end{lemma}
\proof
We have already simplified the complex of \Pss, which underlies
the complex $\hrttw$. The result is exhibited in \fg{fg.el.2}, in which the middle
\Ps\ is given by \fg{fg.ey.56}. In order to obtain $\hrttw$, we apply the
\cnvl\ functor to the diagram of \fg{fg.el.2}, thus replacing the dotted
frame boxes with solid ones.

Now it remains to \excse\ \cntrcns\
following the lemmas of Subsection\rw{ss.exc}.
We will refer to the cones in (the \cnvl\ of) the complex in \fg{fg.el.2} and to the constituent
\mfs\ within the cones by pairs of indices $(i,j)$, as if they
were entries of a matrix.

The \cnvls\ (1,2) and (2,3) split, and the splitting \mfs\ $\hlhor$ are
connected by \outr\ identity morphisms with the same \mfs\ at the
corners (1,1) and (3,3) of the complex. These pairs of $\hlhor$ connected by the identity morphisms
form contractible
cones in the category $\KmplMFWf$, and according to
Lemma\rw{l.excc}, they can be excised from the complex. The result
is the complex in \fg{fg.el.50}.
\def\xarru{ \xdar@/^20pt/[rr] }
\def\xaddr{ \xdar@/^20pt/[dd] }
\begin{figure}[hbt]
\ee
\nonumber
\vcenter{
\xymatrix@C=1.8cm@R=2cm{
&
{\hlhrcir}
\ar[r]^-{\yFm}
&
{\hlhor}
\ar[d(0.72)]^-{\ychind}
&
\ar[r]^-{\ychoutu}
\save
[].[ll]*[F]\frm{}
\restore
&
{\hlhrcir}
\\
{\hlpar}
\ar[d]_-{\yF}
\xaddr^-{\yX}="lma"
&
{\hlhor}
\ar[r]
\xarru^-{}="mmao"
\xdar[dr]
\xdar[ddr]
&
{\hlvir}
\ar[r] \ar[d]
\xaddr
&
{\hlpar}
\ar[d]
\xaddr^-{}="mmat"
&
*i{\lrbs{\hlpar}}
\\
{\hlvir}
\ar[d]_-{\yG}
&
{\hlhor}
\ar[r] \ar[d]
\xarru
\xdar[rrd]
\xdar[rd]
&
{\hlpar}
\ar[r] \ar[d]
\xdar[rd]
&
{\hlvir}
\ar[d]
&
{\hlhor}
\ar[d]^-{\yHD}
\\
{\hlhor}
\ar[d]^-{-\ychoutd}
\save
{[].[uu]."lma"}!C="lm"*[F]\frm{}
\restore
\save
\restore
&
{\hlhrcir}
\ar[r]
&
{\hlhor}
\ar[d(0.62)]^-{\ychoutd}
&
{\hlhor}
\save
[].[ll]."mmat"."mmao"*[F]\frm{}
\restore
\save
{"2,2"."4,2"}!C="xo"
\ar "lm";"xo"^-{\ychinu}
\restore
&
{\hlhrcir}
\save
{[].[uu]}!C="rm"*[F]\frm{}
\restore
\save
{"2,4"."4,4"."mmat"}!C="xt"
\ar"xt";"rm"^-{\ychoutu}
\ar"1,5";"rm"^-{-\ychind}
\ar"5,1";"5,2"^-{\ychinu}
\restore
\\
{\hlpar}
&
{\hlhor}
\ar[r]^-{\yF}
\xarru^-{\yX}="mda"
&
{\hlvir}
\ar[r]^-{\yG}
&
{\hlpar}
\save
[].[ll]."mda"*[F]\frm{}
\restore
&
}
}
\eee
\caption{A simplified complex related to the second Reidemeister move}
\label{fg.el.50}
\end{figure}
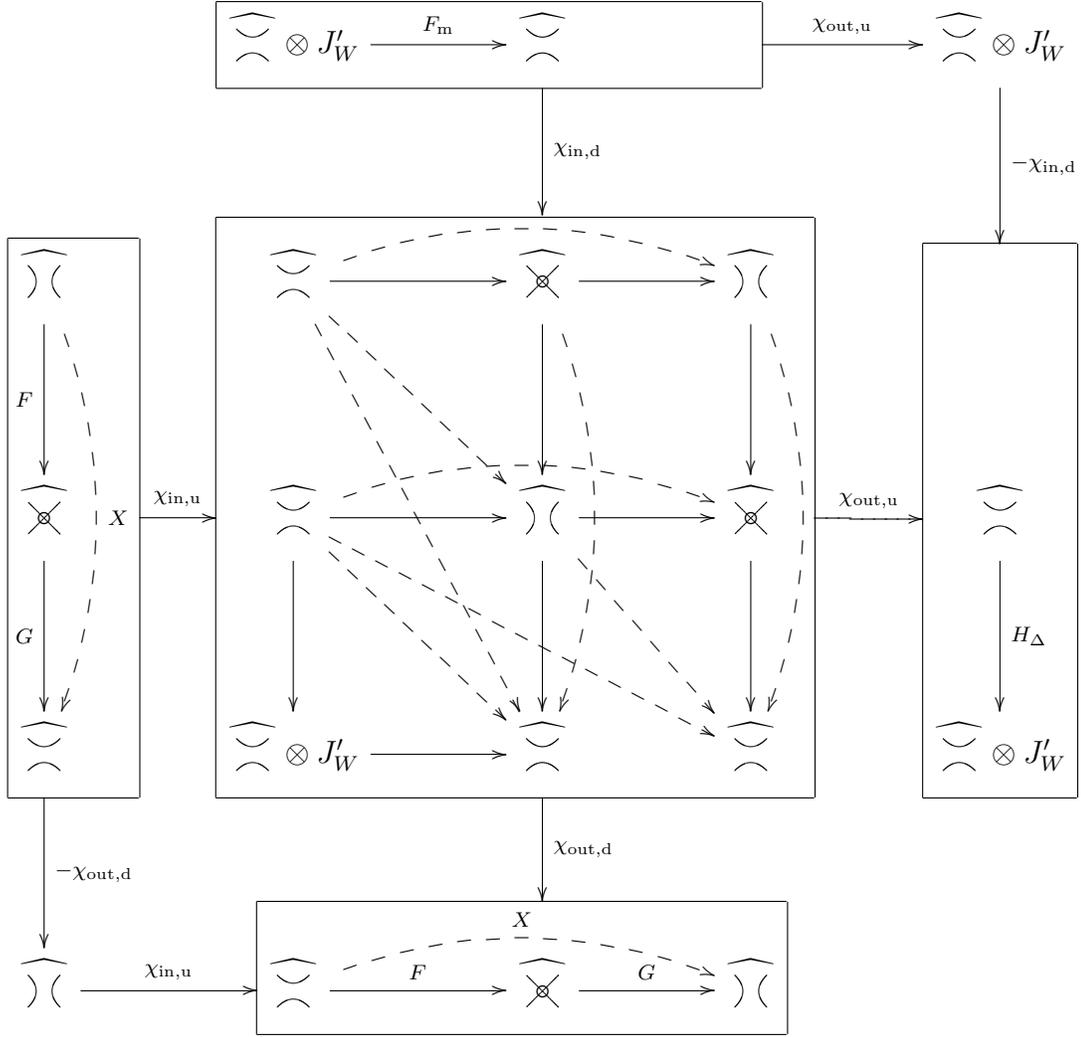
The top line in this complex is similar to the complex\rx{ey.40},
except that $\hlhcpy$ is replaced by $\hlhor$, and we apply to it
a similar simplification procedure. Namely, we split the first \mf\ of the \cnvl\
at (1,2) as
\ee
\label{el.51}
\xymatrix@C=2cm{
{\lrbc{\hlhor\otimes 1}\oplus \lrbc{\hlhor\otimes\JrWpmz}}
\ar[r]^-{\smpm{\id & \xunk}} & {\hlhor}
\save
[].[l]*[F]\frm{}
\restore
}
\eee
The homomorphism $\id$ forms a contractible cone. After its excision,
the top line of the complex becomes
\ee
\label{el.52}
\xymatrix@C=3cm{
{\hlhor\otimes \JrWpmz} \ar[r]^-{\smpm{0\\ \id}}
\ar[d]^-{\xunk}
&
{\lrbc{\hlhor\otimes 1}\oplus\lrbc{\hlhor\otimes\JrWpmz}}
\ar[d]^-{-\ychind}
\\
\cdots
&
\cdots
}
\eee
Now the homomorphism $\id$ forms a contractible cone in the
category $\KmplMF$, and we excise this cone, leaving the \mf\
$\hlhor\otimes 1$ at the position (1,3) of the complex in \fg{fg.el.50}.

Next, we simplify the middle \cnvl\ of the
complex in \fg{fg.el.50}.
 Its detailed structure is
given by (the \cnvl\ of) the diagram of \fg{fg.ey.56}. The \cnst\ \mf\ at the (3,1)
entry of \fg{fg.ey.56} splits:
\ee
\label{el.20}
\hlhrcir =
\lrbc{\hlhor\otimes 1}\oplus
\lrbc{\hlhor\otimes\JrWpmtN}\oplus\lrbc{\hlhor\otimes\zxutN}.
\eee
The term $-\ttNo\,\id\otimes\zxuv{2\xN}$, which appears at $i=0$
in the sum of the expression\rx{ea2.83x} for the homomorphism
$\yHD$, produces a contractible cone
\ee
\label{ey.57}
\xymatrix@C=3cm{
{\hlhor} \ar[r]^-{-\ttNo\,\id} & {\hlhor\otimes\zxuv{2\xN}}
\save [].[l]*[F]\frm{}\restore}
\eee
connecting the \mf\ at the entry (2,1) in the \cnvl\ of \fg{fg.ey.56}
and the third term in the
sum\rx{el.20}. The homomorphisms of the diagram in \fg{fg.ey.56} related to
the cone\rx{ey.57} are directed outwards, so the whole
\cnvl\ of that diagram is a cone of type\rx{el.3x}, where $\xaA$ is
the contractible cone\rx{ey.57}. Its excision modifies the
homomorphism $\ychoutu$ according to case 1A of Lemma\rw{l.css}.

The identity homomorphism $\hmx^0=\id$, which appears at $i=0$ in the sum
of the expression\rx{ea2.82x} for the homomorphism $\yFm$, creates
another contractible cone
\ee
\label{el.54}
\xymatrix@C=2cm{
{\hlhor\otimes 1} \ar[r]^-{\id} & {\hlhor\otimes 1}
\save
[].[l]*[F]\frm{}
\restore
}
\eee
by connecting the first term in the sum\rx{ef.20} with the \mf\ at
the position (3,2)
inside the middle \cnvl\ simplified by the excision of the
cone\rx{ey.57}. All arrows related to this cone are directed
inwards, so the whole simplified \cnvl\ is a cone of the
type\rx{el.3x}, where $\xaB$ is the contractible cone\rx{el.54}.
Its excision modifies the homomorphism $\ychoutu$ according to the
case 2C of Lemma\rw{l.css}. Hence after the excision the whole
complex takes the form of \fg{fg.el.55},
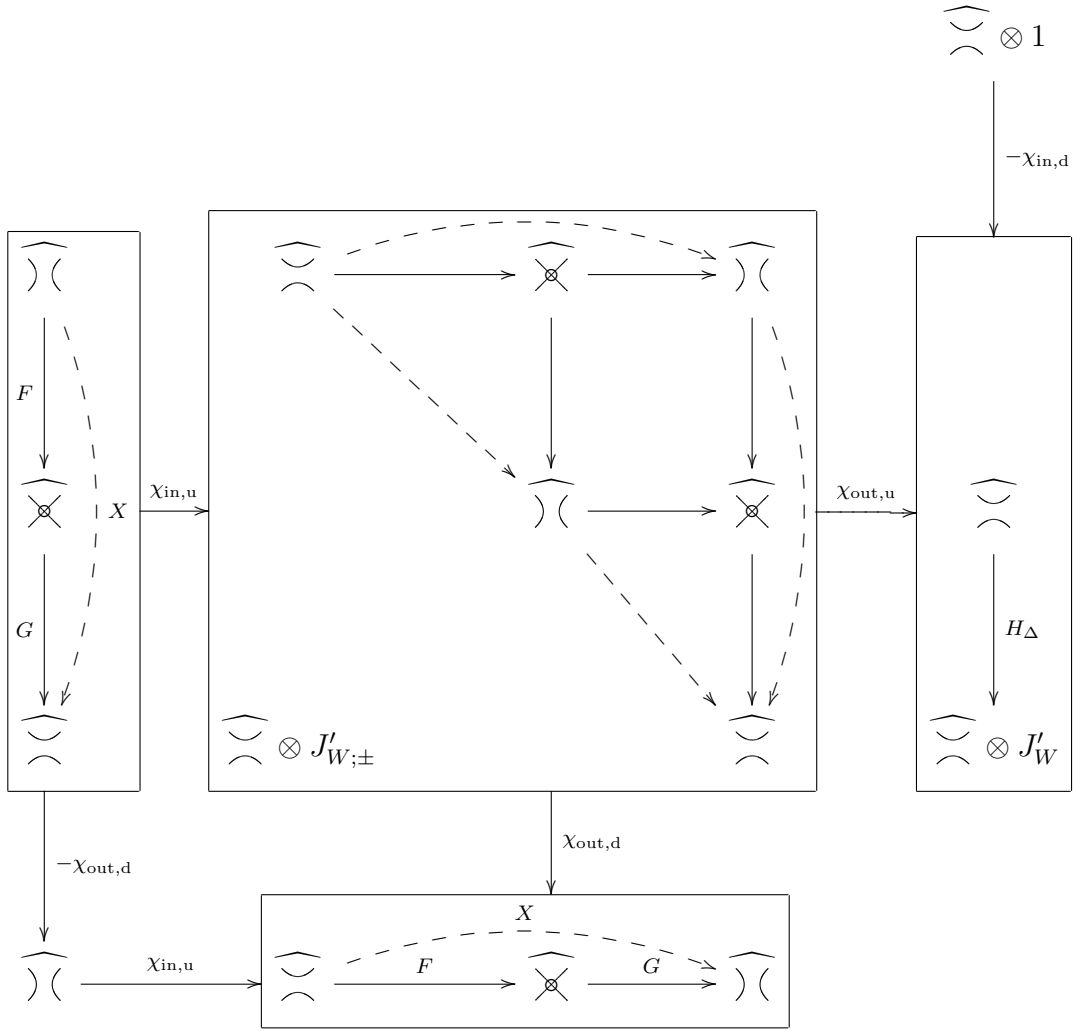
\begin{figure}[hbt]
\ee
\nonumber
\vcenter{
\xymatrix@C=1.7cm@R=2cm{
&
&
&
&
{\hlhor\otimes 1}
\\
{\hlpar}
\ar[d]_-{\yF}
\xaddr^-{\yX}="lma"
&
{\hlhor}
\ar[r]
\xarru^-{}="mmao"
\xdar[dr]
&
{\hlvir}
\ar[r] \ar[d]
&
{\hlpar}
\ar[d]
\xaddr^-{}="mmat"
&
*i{\lrbs{\hlpar}}
\\
{\hlvir}
\ar[d]_-{\yG}
&
&
{\hlpar}
\ar[r]
\xdar[rd]
&
{\hlvir}
\ar[d]
&
{\hlhor}
\ar[d]^-{\yHD}
\\
{\hlhor}
\ar[d]^-{-\ychoutd}
\save
{[].[uu]."lma"}!C="lm"*[F]\frm{}
\restore
\save
\restore
&
{\hlhor\otimes\JrWpmztN}
&
*i{\hlhor}
\save
[]+<0cm,-0.58cm>*{}\ar[d(0.535)]^-{\ychoutd}
\restore
&
{\hlhor}
\save
[].[ll]."mmat"."mmao"*[F]\frm{}
\restore
\save
{"2,2"."4,2"}!C="xo"
\ar "lm";"xo"^-{\ychinu}
\restore
&
{\hlhrcir}
\save
{[].[uu]}!C="rm"*[F]\frm{}
\restore
\save
{"2,4"."4,4"."mmat"}!C="xt"
\ar"xt";"rm"^-{\ychoutu}
\ar"1,5";"rm"^-{-\ychind}
\ar"5,1";"5,2"^-{\ychinu}
\restore
\\
{\hlpar}
&
{\hlhor}
\ar[r]^-{\yF}
\xarru^-{\yX}="mda"
&
{\hlvir}
\ar[r]^-{\yG}
&
{\hlpar}
\save
[].[ll]."mda"*[F]\frm{}
\restore
&
}
}
\eee
\caption{A simplified complex related to the second Reidemeister move}
\label{fg.el.55}
\end{figure}
where the homomorphism $\ychoutu$ acts on the \cnst\ \mf\
$\hlhor\otimes\JrWpmztN$ in the same way as it acted on it in
the diagram of \fg{fg.el.50}.

Now we turn to the \cnvl\ at the position (2,3) in the complex.
After we split its bottom \cnst\ \mf\ as
\ee
\label{el.56}
\hlhor\otimes\JrWp = \lrbc{\hlhor\otimes 1}
\oplus \lrbc{\hlhor\otimes\JrWpmztN}
\oplus \lrbc{\hlhor\otimes \zxutN}
\eee
we see that the term $-\ttNo\;\id$ at $i=0$ in the expression\rx{ea2.83x}
for $\yHD$ forms a contractible cone
\ee
\label{el.57}
\xymatrix@1@C=3cm{ {\hlhor} \ar[r]^-{-\ttNo\id}="a" &
{\hlhor\otimes\zxutN}
\save
[].[l]."a"*[F]\frm{}
\restore
}
\eee
within this convolution. Hence it has the form\rx{el.3x} with
contractible $\xaA$. Its excision transforms the homomorphisms
$-\ychind$ and $\ychoutu$ according to case 1C of Lemma\rw{l.css}.
Hence after the excision they form contractible cones
\ee
\label{el.58}
\xymatrix@C=2cm{{\hlhor\otimes 1} \ar[r]^-{-\id} & {\hlhor\otimes 1}},\qquad
\xymatrix@C=2cm{{\hlhor\otimes\JrWpmztN} \ar[r]^-{\id}
&{\hlhor\otimes\JrWpmztN}}
\eee
within the category $\KmplMF$ with the remaining components of the sum\rx{el.56}.
After we excise them, the complex of \fg{fg.el.55} takes the form of
\fg{fg.el.59}.
\begin{figure}[hbt]
\ee
\nonumber
\vcenter{\xymatrix@C=3cm@R=2.5cm{
{\hlpar}
\ar[d]_-{\yF}
\xaddr^-{\yX}="lma"
&
{\hlhor}
\ar[r]_-{\yF}
\xarru^-{\yX}="mmao"
\xdar[dr]_-{-\yX}
\save
[]+<-0.5cm,0cm>*{}="frb"
\restore
&
{\hlvir}
\ar[r]_-{\yG}
\ar[d]_-{-\yG}
&
{\hlpar}
\ar[d]_-{\yF}
\xaddr^-{\yX}="mmat"
\\
{\hlvir}
\ar[d]_-{\yG}
&
&
{\hlpar}
\ar[r]_-{\yF}
\xdar[rd]_-{\yX}
&
{\hlvir}
\ar[d]_-{\yG}
\\
{\hlhor}
\ar[d]^-{-\ychoutd}
\save
{[].[uu]."lma"}!C="lm"*[F]\frm{}
\restore
\save
\restore
&
*i{\hlhor}
&
*i{\hlhor}
\save
[]+<0cm,-0.58cm>*{}\ar[d(0.6)]^-{\ychoutd}
\restore
&
{\hlhor}
\save
[].[ll]."mmat"."mmao"."frb"*[F]\frm{}
\restore
\save
{"1,2"."3,2"}!C="xo"
\ar "lm";"xo"^-{\ychinu}
\restore
\save
{"1,4"."3,4"."mmat"}!C="xt"
\ar"4,1";"4,2"^-{\ychinu}
\restore
\\
{\hlpar}
&
{\hlhor}
\ar[r]^-{\yF}
\xarru^-{\yX}="mda"
&
{\hlvir}
\ar[r]^-{\yG}
&
{\hlpar}
\save
[].[ll]."mda"*[F]\frm{}
\restore
}
}
\eee
\caption{A simplified complex related to the second Reidemeister move}
\label{fg.el.59}
\end{figure}
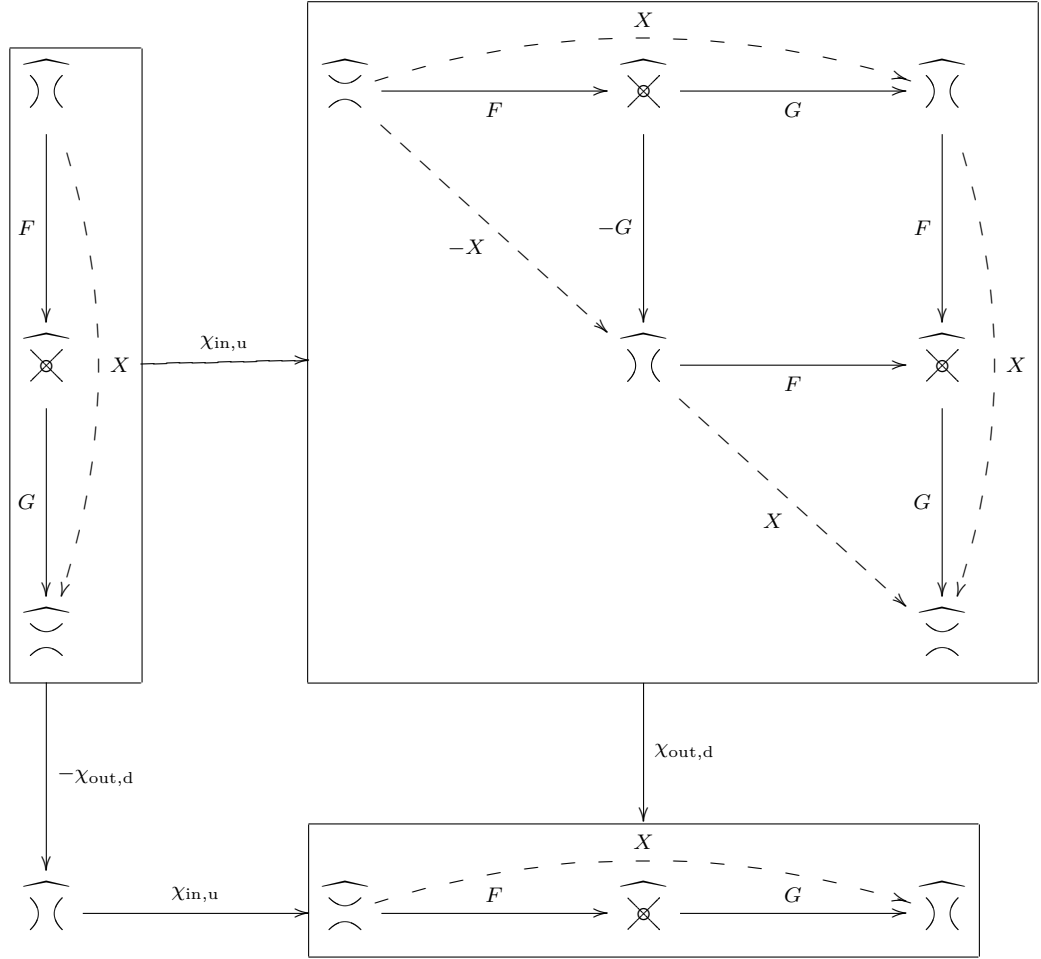

If we apply the homomorphism $\smpm{1 & 1\\ 1 & -1}$
to the sum of the \cnst\ \mfs\ at the positions (1,3) and (2,2) in
the \cnvl\ (1,2) of the diagram of \fg{fg.el.59}, then this \cnvl\
splits into a sum of two \cnvls\ of the
form\rx{ea2.65} and\rx{ea2.65x}. Note that the \cnvls\ at the positions (1,1) and
(2,2) of this complex have the same form, hence they are
all isomorphic to $\hlver$.
Thus the complex of \fg{fg.el.59} takes
the form
\ee
\label{el.30}
\vcenter{\xymatrix@C=2cm@R=2cm{
{\hlver}
\ar[r]^-{\smpm{\ast \\ \id}}
\ar[d]^-{\ychoutd}
&
{\hlver\oplus\hlver}
\ar[d]^-{\smpm{\id & \ast}}
\\
{\hlpar}
\ar[r]^-{\ychinu}
&
{\hlver}
}
}
\eee
where $*$ denotes unspecified homomorphisms, whose precise form is
not important. The components $\id$ of the top and right morphisms
in this diagram form contractible cones, which can be excised in
any order. Thus the complex $\hrttw$ is \hteqt\ to a single \mf\
$\hlpar$ at the position (2,1) in the complex\rx{el.30}. Let us
restore its degree shifts. This \mf\ originates
from the \mf\ $\hparpar$ at (3,1) in the diagram of \fg{fg.ey.53}. The
latter is a tensor product of two \mfs\ $\hlpar$, the top one
being the first \mf\ in the second line of the \rhs of \ex{ey.27} and the
bottom one being the last \mf\ in the second line of the \rhs of \ex{ey.26}.
Hence this \mf\ has no degree shifts, and we proved
\ex{ey.52}.\qed

Note that in process of establishing the \hteq\rx{ey.52}, we
proved the following relation between \mfs\ (just follow the
transformations of the
middle \mf\ in the convolution of the \Ps\ in \fg{fg.ey.53y}):
\ee
\label{el.31z}
\hverver\cchq\lrbc{\hlhor\otimes\JrWpmztN} \oplus
\hlver\grsmo\oplus\hlver\grso.
\eee

\subsection{Virtual and \smvrt\ Reidemeister moves}
\label{ss.vadj}

A diagram of a virtual \grtng\ $\gt$ is a planar graph with 1-valent vertices called \lgs\ and 4-valent
vertices which are of three types: $\lncr$, $\lxver$ and $\lxvir$.
Two edges of $\gt$ are called \vadj, if they are incident to the
same virtual vertex and they are attached to it at opposite sides.
\begin{lemma}
\label{l.vadj}
Suppose that the \lgs\ $i$ and $j$ of a virtual \grtng\ $\gt$ are
connected by a sequence of \vadj\ edges. Then
\ee
\label{vadj.1}
\hgt\cchq\hgtp \otimes \hloarvv{i}{j}\;,
\eee
where $\gtp$ is a virtual \grtng\ constructed by removing the
sequence of \vadj\ edges from $\gt$ and `dissolving' the virtual
vertices, which are incident to these edges.
\end{lemma}
\proof
The lemma follows easily from the definition\rx{eb2.16o} of the
virtual vertex categorification and from the property\rx{ex.13} of
\oarc\ \mfs, which allows us to reduce a \mf\ of a sequence of adjacent
edges into a \mf\ of a single edge connecting the \lgs.\qed

\pr{Theorem}{th.vadj}
This theorem is a simple corollary of Lemma\rw{l.vadj}.

The homotopy invariance of a virtual link complex under the first and second virtual
Reidemeister moves is a particular case of \ex{vadj.1}

The invariance under the third virtual Reidemeister move follows
from the fact that the \mfs\ of both virtual \grtngs\ are homotopy
equivalent to the tensor product of three \oarc\ \mfs.

Finally, the invariance under the \smvrt\ Reidemeister move
follows from the fact that the complexes of \mfs\ corresponding to
both virtual \grtngs\ are equivalent to the tensor product
\ee
\label{vadj.2}
\hlncr\otimes\hloarvv{}{}.
\eee
\qed

\appendix

\section{Virtual crossings, convolutions
and the categorification of the \homflyp\ and of its \sun\ specialization}
\label{ap.1}

A \tarc\ \mf\ corresponding to an \eltr\ virtual crossing can also
be introduced in the categorification\cx{KR1},\cx{KR2},\cx{Ra2} of the
\tvar\ \homflyp\cx{HOMFLY},\cx{PT} and
its \sun\ specialization. Thus both categorifications can be
extended to virtual links. Moreover,
an analog of the relation\rx{e2.x9} between the \mf\ of an \eltr\
graph and the \mfs\ of \tarc\ graphs
also exists.
The \tvar\ and \sun\ cases are similar. First, we present the
\sun\ case in details and then sketch the \tvar\ case.

\subsection{The \sun\ case}
\label{ss.a1}
The basic
algebra of the categorification is $\IQ[x]$ and the basic
polynomial is
\ee
\label{ae.1}
\nWhx = \zx^{\xN+1}.
\eee
The categorification admits a \qgrd\ with $\dgq\zx=2$.
%
Since $\dgq\nW=2(\xN+1)$, we have to set
\ee
\label{ae.3}
\dgq\xD = \xN+1
\eee
on \qdgr s of all \tdiffs\ $\xD$.

Similar to \ex{ea2.13a}, we introduce the finite differences of
$\nW(\zx)$:
\ee
\nWitxot = \frac{\nW(\zxo)-\nW(\zxt)}{\zxo-\zxt},\qquad
\nWithv{\zxo}{\zxt}{\zxth} =
\frac{\nWitv{\zxo}{\zxt}-\nWitv{\zxo}{\zxth}}{\zxt-\zxth}.
\eee
The \oarc\ \mf\ is a 1-row \Kmf\
\ee
\label{ae.4}
\widehat{\loarhvv{1}{2}} = \kmfbv{\zxt-\zxo}{\nWitxot}.
\eee
%
\subsubsection{Virtual crossing and \tarc\ \mfs}
For two legs oriented inwards and two legs oriented outwards there
are two \tarc\ graphs connecting them:
\ee
\label{ae.5}
\xybox{0;/r4pc/:
\xunoverv=<,
(0.8,0.2)*{3},
(2.2,0.2)*{4},
(0.8,-1.2)*{1},
(2.2,-1.2)*{2}
}
\qquad\qquad
\xybox{0;/r4pc/:
(0.5,-0.5)*{\zbendv@(0)=>},
(2,-1)*{}="a";
(1,0)*{}="b";
{\ar"a";"b"},
(1.5,-0.5)*{\vcirc},
(0.8,0.2)*{3},
(2.2,0.2)*{4},
(0.8,-1.2)*{1},
(2.2,-1.2)*{2}
}
\eee
The corresponding \mfs\ of the polynomial
\ee
\label{ae.6}
\zWfx = \nW(\zxth) + \nW(\zxf) - \nW(\zxo) - \nW(\zxt),\qquad\zbx
= \ylst{\zx}{4}
\eee
are 2-row \Kmfs, which are the tensor products of \mfs,
corresponding to individual arcs:
\ee
\label{ae.7}
\hbpar =
\kmftwv{\zxth-\zxo}{\zxf-\zxt}{\nWitv{\zxo}{\zxth}}{\nWitv{\zxt}{\zxf}}
,\qquad
\hbvir =
\kmftwv{\zxth-\zxt}{\zxf-\zxo}{\nWitv{\zxt}{\zxth}}{\nWitv{\zxo}{\zxf}}.
\eee
Following\cx{KR1}, we introduce the third \mf\ of $\zWfx$, which
corresponds to the \Igr:
\ee
\label{ae.7a}
\hbver = \kmftwv{\zzxfm}{\zxth\zxf-\zxo\zxt}{\xunk}{\xunk}\grsmo,
\eee
where $\xunk$ denotes the entries whose exact form is not
important for us here.

The polynomial $\zWfx$ is invariant under the action of the
symmetric group $\St$, permuting the \lgs\ 1 and 2 and their
variables $\zxo$ and $\zxt$; hence this group acts on the category
$\MFWf$ by \enfns. In particular,
\ee
\label{ae.7a1}
\hXot\lrbc{\hbpar} = \hbvir,\qquad
\hXot\lrbc{\hbvir} = \hbpar,\qquad
\hXot\lrbc{\hbver} = \hbver.
\eee
We also need the permutation $\sXoth$, switching
the variables $\zxo$ and $\zxth$, and the corresponding functor
$\hXoth$:
\begin{gather}
\label{ae.7a2}
\sXoth\big(\zWfx\big) = \nW(\zxo) + \nW(\zxf) - \nW(\zxt) - \nW(\zxth),
\\
\hXoth\lrbc{\hbpar} = \hbpro ,\qquad \hXoth\lrbc{\hbvir} = \hbhro.
\end{gather}

The
\mfs\
$\hspar$, $\hsvir$ and $\hsver$
can be presented as the tensor products of \cmn\ and \prpr\ parts. Indeed,
let us introduce the \prpr\ algebra
\ee
\label{ae.8}
\yRprp = \IQ[\zp,\zr,\zC],
\eee
related to the \tarc\ algebra $\yR=\IQ[\zbx]$ by the homomorphism
$\xymatrixnocompile@1{\yRprp\ar[r]^-{\zhomprp} & \yR}$, defined by the
formulas
\ee
\label{ae.9}
\zhomprp(\zp) = \zxf-\zxt,\qquad\zhomprp(\zr) = \zxf-\zxo,\qquad
\zhomprp(\zC) = \ztC(\bfx),
\eee
where
\ee
\label{ae.10}
\ztC(\bfx) = \nWithv{\zxo}{\zxt}{\zxth}
+\nWithv{\zxo}{\zxt}{\zxf}.
\eee
%

We introduce \prpr\ \mfs\
\ee
\label{ae.11}
\hbparprp =\kmfv{\zp}{\zr\zC},\qquad
\hbvirprp =\kmfv{\zr}{\zp\zC},\qquad
\hbverprp =\kmfv{\zp\zr}{\zC}\grsmo
\eee
of the \prpr\ polynomial
\ee
\label{ae.12}
\zWfp = \zp\,\zr\,\zC
\eee
and use the notation
\ee
\label{ea.13}
\ggtprp = \zhhomtprp(\ggprp),
\eee
where $\xg$ is
a graph $\lspar$, $\lsvir$ or $\lsver$, while
$\zhhomtprp$ is the
functor\rx{e2.d3a2}, corresponding to the homomorphism $\zhomprp$.
Finally, we set the \cmn\ \mf\ to be
\ee
\label{ae.14}
\zKcmn = \kmfbv{\zxth+\zxf-\zxo-\zxt}{\zA(\zbx)},
\eee
where
\ee
\label{ae.15}
\zA(\zbx) =
\nWitv{\zxo}{\zxt} + (\zxth-\zxf)\nWithv{\zxo}{\zxt}{\zxth}.
\eee
Note that $\sXot\big(\zA(\zbx)\big) = \zA(\zbx)$, so
\ee
\label{ae.15a}
\hXot(\zKcmn) = \zKcmn.
\eee

The symmetric group $\St$ acts on the \prpr\ algebra $\yRprp$ by
permuting $\zp$ and $\zr$, while leaving $\zC$ intact:
\ee
\label{ae.15a1}
\sXot(\zp)=\zr,\qquad \sXot(\zr)=\zp,\qquad\sXot(\zC) = \zC.
\eee
The homomorphism $\zhhomtprp$ is equivariant with respect to the
simultaneous action of $\St$ on $\zbx$ and on $\zp,\zr,\zC$.
The \prpr\ polynomial $\zWfp$ is invariant under the action of $\St$, so
this group acts on its \mfs\ by \enfns. In particular
\ee
\label{ae.15a2}
\hXot\lrbc{\hbparprp} = \hbvirprp,\qquad
\hXot\lrbc{\hbvirprp}= \hbparprp,\qquad
\hXot\lrbc{\hbverprp} = \hbverprp.
\eee
%
%

\begin{proposition}
\label{pr.suf}
The \tarc\ \mfs\rx{ae.11} and the \Igr\ \mf\rx{ae.7a} factor into the tensor product of the
\cmn\ and \prpr\ \mfs:
\ee
\label{ae.16}
\gg\cong\zKcmn\otR\ggtprp
\eee
\end{proposition}
\proof
The application of the transformation $\rtrv{2}{1}{1}$ (see \ex{ea2.12z2})
to \tarc\
\mfs\ $\hspar$ and $\hsvir$ of \ex{ae.7} and the application of the
composition of transformation $\rtov{2}{-1}\rtrv{1}{2}{-\zxf}$ to the
\Igr\ \mf\ $\hsver$ of \ex{ae.7a} makes their left columns equal
to those of the corresponding \mfs\ in the \rhs of \ex{ae.16}. Since the polynomials
in the left columns form regular sequences, the isomorphisms\rx{ae.16}
follow from Theorem\rw{th.equiv}.\qed

Note that the \tarc\ \mfs\ can be transformed explicitly
into the tensor product\rx{ae.16} by a composition of
two transformations\rx{ea2.12z2} and\rx{ea2.12z3}:
\ee
\label{ae.17}
\xymatrixnocompile@C=4cm{\gg \ar@{|->}[r]^-{\rtrpv{1}{2}{\nWithv{\zxo}{\zxt}{\zxth}}\rtrv{2}{1}{1}}
& \zKcmn\otR\ggtprp}.
\eee

\subsubsection{The \sdlmrp\ and its cone}
We define the \prpr\ \sdlmrp\
\ee
\label{ae.19}
\yFprp\in\ExtmfoqNmo\lrbc{\hsparprp,\hsvirprp}
\eee
by the following commutative diagram
\ee
\label{ae.20}
\xymatrix@C=1.5cm{
{\hsparprp}
\ar[d]^-{\yFprp}
&
{\yRomNo}
\ar[r]^-{\zp}
\ar[rd]^-{1}
&
{\yRz}
\ar[r]^-{\zr\zC}
\ar[dr]^-{-\zC}
&
{\yRomNo}
\\
{\hsvirprp}
&
{\yRomNo}
\ar[r]^-{\zr}
&
{\yRz}
\ar[r]^-{\zp\zC}
&
{\yRomNo}
}
\eee
The corresponding \sdlmrp\ $\yF\in\Extmfo\lrbc{\hspar,\hsvir}$ is
defined with the help of \ex{ae.16} as
\ee
\label{ae.20a1}
\yF = \idcmn\otimes\yFprp.
\eee
We consider two other \sdlmrps, produced by the action of the
functors $\hXot$ and $\hXoth$:
\ee
\label{ae.20a2}
\yG = \hXot(\yF)\in\Extmfo\lrbc{\hsvir,\hspar},\qquad
\yH = \hXoth(\yF)\in\Extmfo\lrbc{\hspro,\hshro}.
\eee

\begin{theorem}
\label{th.a1}
The \sdlmrp\ $\yH$ is equal (up to a possible sign factor) to
the \sdlmrp\ $\eta$ defined in Section~9 of\cx{KR1}.
\end{theorem}
We need the following
\begin{lemma}
\label{lm.a1}
The composition of two \sdlmrps\ is
\ee
\label{ae.15b1}
\yG\yF\mphq -(\xN+1)\sum_{i=0}^{\xN-1}\hxo^i \hxt^{\xN-i-1}.
\eee
\end{lemma}
\proof
The composition of \prpr\ \sdlmrps\ is presented by the diagram
\ee
\label{ae.15b2}
\vcenter{\xymatrix@C=1.5cm{
{\hsparprp}
\ar[d]^-{\yFprp}
&
{\yRomNo}
\ar[r]^-{\zp}
\ar[dr]^-{1}
&
{\yRz}
\ar[r]^-{\zr\zC}
\ar[dr]^-{-\zC}
&
{\yRomNo}
\\
{\hsvirprp}
\ar[d]^-{\yGprp}
&
{\yRomNo}
\ar[r]^-{\zr}
\ar[dr]^-{1}
&
{\yRz}
\ar[r]^-{\zp\zC}
\ar[dr]^-{-\zC}
&
{\yRomNo}
\\
{\hsparprp}
&
{\yRomNo}
\ar[r]^-{\zp}
&
{\yRz}
\ar[r]^-{\zr\zC}
&
{\yRomNo}
}
}
\eee
Obviously, $\yGprp\yFprp=-\zC\,\id$, so in view of \ex{ae.10},
\ee
\label{ae.15b3}
\yGprp\yFprp = -\nWithv{\hxo}{\hxt}{\hxth}
-\nWithv{\hxo}{\hxt}{\hxf}.
\eee
Now \ex{ae.15b1} follows from the relations
\ee
\label{ae.15b4}
\hxth\mphq\hxo,\qquad\hxf\mphq\hxt\qquad\text{in $\EndMF\lrbc{\hspar}$}
\eee
and from the formula
\ee
\label{ae.15b5}
\nWithv{\zxo}{\zxo}{\zxt}
+\nWithv{\zxo}{\zxt}{\zxt} =
(\xN+1)\sum_{i=0}^{\xN-1}\zxo^i \zxt^{\xN-i-1}.
\eee
\qed

\pr{theorem}{th.a1}
In our conventions
\ee
\label{ae.20a3}
\dgq \yH = \dgq \eta = \xN-1.
\eee
Since
\ee
\label{ae.20a4}
\dim\Extmfo_{\qgrsv{\xN-1}}\lrbc{\hspro,\hshro} = 1,
\eee
$\yH$ is proportional to $\eta$:
\ee
\label{ae.15b6}
\yH\mphq a\eta,\qquad a\in\IQ.
\eee
Applying the functor $\hXoth$ to the relation\rx{ae.15b1} we find
that the composition of morphisms
\ee
\label{ae.15b7}
\xymatrix{{\hspro} \ar[r]^-{\yH}& {\hshro}\ar[r]^-{\yHp} &{\hspro}},
\eee
(where $\yHp=\hXoth(\yH)$)
is
\ee
\label{ae.15b8}
\yHp\yH\mphq - (\xN+1)\sum_{i=0}^{\xN-1}\hxo^i \hxt^{\xN-i-1}.
\eee
Since $\eta$ satisfies the same relation, we conclude that
$a^2=1$.\qed

We will consider the cones of \sdlmrps\ $\yF$ and $\yG$. We use
the \cnvl\ notation\rx{ed.42}, because it reveals the internal
structure of the cone.
\begin{theorem}
\label{th.a2}
After an appropriate \qdgr\ shift, the cones of \sdlmrps\ $\yF$ and $\yG$ are \hteqt\
to the \mf\ $\hsver$:
\ee
\label{ae.20a5}
\begin{split}
\hbver\;&\cchq\;\; \xymatrix{\Cbdv{\hbpar\grsmo}{\yF}{\hbvir\grso}}
\\
&\cchq\;\;
\xymatrix{\Cbdv{\hbvir\grsmo}{\yG}{\hbpar\grso}}\;\;.
\end{split}
\eee
\end{theorem}
\proof
It is sufficient to prove the first equivalence of\rx{ae.20a5},
since the second one would follow from the application of the
functor $\hXot$ and from the third equivalence of\rx{ae.7a1}.

We are going to prove the \hteq\ between \prpr\ \mfs
\ee
\label{ae.22}
\xymatrix{\Cbdv{\hbparprp\grsmo}{\yFprp}{\hbvirprp\grso}}\;\;\cchq\;
\hbverprp,
\eee
and the first equivalence of\rx{ae.20a5} follows, if we tensor multiply
both sides by the \cmn\ \mf\ $\zKcmn$.

The following diagram establishes the isomorphism between the
direct sum of \mfs\
$\hsverprp\oplus\kmfbv{1}{\zp\zr\zC}\grso$
and the \cnvl\ in the \lhs of \ex{ae.22}:
\ee
\label{ae.23}
\vcenter{\xymatrix@C=2cm@R=2cm{
{\hsverprp\oplus\kmfbv{1}{\zp\zr\zC}\grso}
\ar[d]^-{\cong}
&
{\yRto}
\ar[r]^-{\smpm{\zp\zr & 0\\0 & 1}}
\ar[d]^-{\smpm{\zr & 1 \\ -1 & 0}}
&
{\yRtz}
\ar[r]^-{\smpm{\zC & 0 \\ 0 & \zp\zr\zC}}
\ar[d]^-{\smpm{1 & p \\ 0 & 1}}
&
{\yRto}
\ar[d]^-{\smpm{\zr & 1 \\ -1 & 0}}
\\
{\hsparprp\grsmo\xrightarrow{\wFp}\hsvirprp\grso}
\save[]*[F]\frm{} \restore
&
{\yRto}
\ar[r]^-{\smpm{\zp & 0 \\ 1 & \zr}}
&
{\yRtz}
\ar[r]^-{\smpm{\zr\zC & 0 \\ -\zC & \zp\zC}}
&
{\yRto}
}
}
\eee
The \Kmf\ $\kmfbv{1}{\zp\zr\zC}$ is contractible, so
$\hsverprp\oplus\kmfbv{1}{\zp\zr\zC}\grso\cchq\hsverprp$.\qed
\subsubsection{Natural morphisms and crossing complexes}
The cones of \ex{ae.20a5} are related by natural morphisms to
their \cnst\ \mfs. We consider two such morphisms:
\ee
\label{ae.23a1}
\xymatrix@C=1.5cm{
{\hbpar\grso}\ar[r]^-{\ychin} &
\Cbdv{\hbvir\grsmo}{\yG}{\hbpar\grso}
}\cchq\hbver
\eee
and
\ee
\label{ae.23a2}
\hbver\cchq
\xymatrix@C=1.5cm{
\Cbdv{\hbpar\grsmo}{\yF}{\hbvir\grso}
\ar[r]^-{\ychout}
&
{\hbpar\grsmo}
}.
\eee

\begin{theorem}
\label{th.aa2}
The morphisms $\ychin$ and $\ychout$ are \hteqt\ up to non-zero
constant factors to the morphisms $\ychiz$ and $\ychio$ of\cx{KR1}.
\end{theorem}
We need the following
\begin{lemma}
\label{lm.a2}
The composition of the natural morphisms is
\ee
\label{ae.23c1}
\ychout\ychin\mphq \hxo-\hxt.
\eee
\end{lemma}

\proof
Obviously,
\ee
\label{ae.23c2}
\ychin=\idcmn\otR\ychinprp,\qquad\ychout=\idcmn\otR\ychoutprp,
\eee
where $\ychinprp$ and $\ychoutprp$ are \prpr\ natural morphisms
\ee
\label{ae.23c3}
\xymatrix@C=1.5cm{
{\hbparprp\grso}\ar[r]^-{\ychinprp} &
\Cbdv{\hbvirprp\grsmo}{\yGprp}{\hbparprp\grso}
}
\eee
and
\ee
\label{ae.23c4}
\xymatrix@C=1.5cm{
\Cbdv{\hbparprp\grsmo}{\yFprp}{\hbvirprp\grso}
\ar[r]^-{\ychoutprp}
&
{\hbparprp\grsmo}
}.
\eee
If we replace the cones by the \hteqt\ \mf\ $\hsverprp$, then the
morphisms $\ychinprp$ and $\ychoutprp$ take the form
\ee
\label{ae.23c5}
\cxy{
{\hsparprp\grso}
\ar[d]^-{\ychinprpp}
&
{\yRotmN}
\ar[r]^-{\zp}
\ar[d]^-{-1}
&
{\yRzso}
\ar[r]^-{\zr\zC}
\ar[d]^-{-\zr}
&
{\yRotmN}
\ar[d]^-{-1}
\\
{\hsverprp}
\ar[d]^-{\ychoutprp}
&
{\yRotmN}
\ar[r]^-{\zp\zr}
\ar[d]^-{\zr}
&
{\yRzmo}
\ar[r]^-{\zC}
\ar[d]^-{1}
&
{\yRotmN}
\ar[d]^-{\zr}
\\
{\hspar\grsmo}
&
{\yRomN}
\ar[r]^-{\zp}
&
{\yRzmo}
\ar[r]^-{\zr\zC}
&
{\yRomN}
}
\eee
Now it is obvious that
\ee
\label{ae.23c6}
\ychoutprp\ychinprp \mphq -\zr,
\eee
and relation\rx{ae.23c1} follows, since according to \ex{ae.9},
$\zhomprp(\zr) = \zxf-\zxo$, while $\zxf\mphq\zxt$ in
$\HomMF\lrbc{\hspar}$. \qed

\pr{theorem}{th.aa2}
In our conventions,
\ee
\label{ae.15b9}
\dgq\ychin=\dgq\ychout=\dgq\ychiz=\dgq\ychio=0,
\eee
but we have shown in\cx{KR1} that
\ee
\label{ae.15b10}
\dim\Extmfzqv{0}\lrbc{\hspar,\hsver}
=\dim\Extmfzqv{0}\lrbc{\hsver,\hspar}=1.
\eee
Therefore
\ee
\label{ae.15b11}
\ychin=a_0\ychiz,\qquad\ychout=a_1\ychio,\qquad a_0,a_1\in\IQ.
\eee
Since $\hxt-\hxo\not\mphq 0$ in
$\HomMF\lrbc{\hspar}$, \ex{ae.23c1} implies that $a_0,a_1\neq
0$.\qed

Now, in view of the equivalences established in Theorems\rw{th.a2}
and\rw{th.aa2}, we can present the categorification formulas of\cx{KR1} for
the \eltr\ crossings as complexes\rx{ae.23c3} and\rx{ae.23c4}:
\begin{align}
\label{ae.d1}
\hbcrn & =
\lrbc{
\xymatrix@C=1.5cm{
{\hbpar\grso}\ar[r]^-{\ychin} &
\Cbdv{\hbvir\grsmo}{\yG}{\hbpar\grso}
}
\;\;}\grsmN\grsemo,
\\\nonumber
\\
\label{ae.d2}
\hbcrp & =
\lrbc{\;\;
\xymatrix@C=1.5cm{
\Cbdv{\hbpar\grsmo}{\yF}{\hbvir\grso}
\ar[r]^-{\ychout}
&
{\hbpar\grsmo}
}\;\;
}\grsN.
\end{align}
Here we assume that the differentials $\ychin$ and $\ychout$ are cohomological (that is,
they have the homological degree $1$) and that the right terms in
the complexes in parentheses have homological degree $0$.

\subsubsection{A categorification complex for a virtual link}
The second formula of\rx{ae.7} allows us to define a
categorification complex for a virtual link in the same way as we
did it for ordinary links in\cx{KR1}.
Let $\xL$ be a virtual link diagram. We split it into a disjoint union of \eltr\ real and virtual
crossings by cutting across all edges. To real crossings we associate complexes of
\mfs\rx{ae.d1} and\rx{ae.d2}. To virtual crossings we associate
the second \mf\ of \ex{ae.7}. Finally, we take the tensor
product of \mfs\ and their complexes, thus forming a chain
complex
$\yCNbL$ of $\ZZ_2\times\ZZ$-graded $\IQ$-vector spaces.
\begin{theorem}
If two diagrams $\xL$ and $\xLp$ represent the same virtual link,
then the corresponding complexes are \hteqt:
\ee
\label{ev.1}
\yCNbLp\cchq \yCNbL.
\eee
\end{theorem}
\proof
If the diagrams $\xL$ and $\xLp$ are related by a Reidemeister
move, then the relation\rx{ev.1} is proved in\cx{KR1}. If the
diagrams $\xL$ and $\xLp$ are related by a virtual Reidemeister
move or by a mixed move\rx{ee.vrm}, then the relation\rx{ev.1} is
obvious.
\qed

The graded Euler characteristic of the complex $\yCNbL$ may serve
as a definition of the \sun\ \homflyp\ for virtual links.

\subsection{The \tvar\ \homflyp\ case}

Consider the categorification\cx{KR2} for the \homflyp, in which we set $a=0$
(see also\cx{Ra2}). Then the basic \qgrd\ algebra is the same as in
the \sun\ case, but the basic polynomial is zero: $\nW=0$. As a
result, the \mfs\ of the \sun\ categorification are replaced by
`\innr'
complexes (we call them \innr\ in order to distinguish them from the
`\outr' complexes associated with the categorification of the
crossings). The \Ztgrdng\ of \mfs\ lifts to homological
\Zgrdng\ of \innr\ complexes. We refer to it as \tgrd, because it generates
powers of $t$ in the formula for the \grdd\ Euler characteristic
of the link categorification complex
\ee
\label{ae.d2a}
\chqt\big(\nCdL\big) = \sum_{i,j,n\in\ZZ}(-1)^{j+n} t^{2j}q^i\rankv\xCnijL
\eee
(\cf \ex{1.6}).

Another important distinction between the categorification of the
\tvar\ \homflyp\ and its \sun\ specialization is that in the
former case the link diagram $\xL$ should be the result of a
(circular) braid closure. Suppose that $\xL$ is a closure of an
$\xnb$-braid $\xbt$, and a \clgr\ $\xGr$ which originates from the
$\xbr$-resolution of the crossings of $\xL$, is the closure of the
$2\xnb$-legged \opgr\ $\xgbr$, which in turn results from the same
resolution of the crossings of $\xbt$. Let $\zhhomn$ be the
functor\rx{ex.9} corresponding to the joining of $\xnb$ incoming
and $\xnb$ outgoing legs of $\xgbr$, which would turn it into $\xGr$.
Then in contrast to the \sun\
categorification formula $\hGr=\zhhomn(\ggbr)$, we have to
introduce an extra degree shift in the \tvar\ case:
\ee
\label{ae.d3}
\hGr=\zhhomn(\ggbr)\grsv{\xnb}\gsztv{-\xnb/2}[\xnb/2].
\eee

Since $\nW=0$, the condition $\dgq\xD=\frac{1}{2}\dgq\nW$ is no
longer imposed on the differentials $\xD$ of \innr\ complexes.
However,
invariance under the first Reidemeister move imposes a requirement
%
\ee
\label{ae.d4}
\dgq\xD=2.
\eee

The \Kmfs\ of Subsection\rw{ss.a1} are replaced by \Kcxs. Let
$\yR$ be a \qgrdd\ polynomial algebra and let $\yp\in\yR$,
$\dgq\yp = k$. Then the \Kcx\ $\kcyp$ is the complex
\ee
\label{ae.d5}
\xymatrixnocompile{\yRo\grsv{k-2}\ar[r]^-{\yp} & \yRz},
\eee
where the index $i$ of $\yRi$ indicates its \tgrd. The \qdgr\
shift $\grsv{k-2}$ is required to satisfy the condition\rx{ae.d3}.
For the column $\ybp$ of polynomials $\ypi\in\yR$, $1\leq i\leq
n$, we define
\ee
\label{ae.d6}
\kcbfp = \bigotimes_{i=1}^n \kcv{\ypi}.
\eee
Now the removal
Removing second columns in all \Kmfs\ of subsection\rw{ss.a1} turns them into \Kcxs\ and thus adapts them
to the \tvar\ \homfly\ case. The only tricky part is the grading shifts, which
have to be adjusted in accordance with \eex{ae.d3} and\rx{ae.d4}.

The \Kcx\ of the \oarc\ graph is
\ee
\label{ae.d7}
\widehat{\loarhvv{1}{2}} = \kcv{\zxt-\zxo}.
\eee
The \mfs\ of the \tarc\ graphs $\hspar$, $\hsvir$ and of the
\eltr\ \fvtx\ graph $\hsver$ have the form\rx{ae.16}, where this
time
\ee
\label{ae.d8}
\zKcmn = \kcv{\zxth+\zxf-\zxo-\zxt},
\eee
while the \prpr\ algebra is $\yRprp=\IQ[\yp,\yr]$ (this time it does not include
$\zC$), the \prpr\ homomorphism is
\ee
\label{ae.d9}
\zhomprp(\zp) = \zxf-\zxt,\qquad\zhomprp(\zr) = \zxf-\zxo
\eee
and the \prpr\ complexes are
\ee
\label{ae.d10}
\hbparprp =\kcv{\zp},\qquad
\hbvirprp =\kcv{\zr},\qquad
\hbverprp =\kcv{\zp\zr}\grsmo.
\eee
The proofs are the same as in subsection\rw{ss.a1}, except that we
ignore the right columns of \Kmfs.

The \sdlmrp\rx{ae.19} is defined by the reduced version of the
diagram\rx{ae.20}
\ee
\label{ae.d11}
\xymatrix@C=1.5cm{
{\hsparprp}
\ar[d]^-{\yFprp}
&
{\yRo}
\ar[r]^-{\zp}
\ar[rd]^-{1}
&
{\yRz}
\\
{\hsvirprp}
&
{\yRo}
\ar[r]^-{\zr}
&
{\yRz}
}
\eee
and it is easy to see that after the appropriate degree shifts its
cone is \hteqt\ to the complex $\hsver$, that is,
isomorphism\rx{ae.22} still
holds. If we tensor multiply both sides by $\zKcmn$ and apply the
$\sXot$ argument of subsection\rw{ss.a1}, then we obtain
isomorphism\rx{ae.20a5}. The natural morphisms\rx{ae.23c3} and\rx{ae.23c4}
associated with the cones\rx{ae.20a5} can be cast in the form
\ee
\label{ae.d11a}
\cxy{
{\hsparprp\grso}
\ar[d]^-{\ychinprpp}
&
{\yRoso}
\ar[r]^-{\zp}
\ar[d]^-{-1}
&
{\yRzso}
\ar[d]^-{-\zr}
\\
{\hsverprp}
\ar[d]^-{\ychoutprp}
&
{\yRoso}
\ar[r]^-{\zp\zr}
\ar[d]^-{\zr}
&
{\yRzmo}
\ar[d]^-{1}
\\
{\hspar\grsmo}
&
{\yRomo}
\ar[r]^-{\zp}
&
{\yRzmo}
}
\eee
and it is easy to verify that they coincide (up to constant
factors) with the morphisms $\ychiz$ and $\ychio$ of\cx{KR2}.
Thus we can present the \eltr\ crossing categorification complexes
of\cx{KR2}
in the form
\begin{align}
\label{ae.d12}
\hbcrn & =
\lrbc{
\xymatrix@C=1.5cm{
{\hbpar\grso}\ar[r]^-{\ychin} &
\Cbdv{\hbvir\grsmo}{\yG}{\hbpar\grso}
}
\;\;}\lrbl{\hlf}\lrbs{-\hlf},
\\\nonumber
\\
\label{ae.d13}
\hbcrp & =
\lrbc{\;\;
\xymatrix@C=1.5cm{
\Cbdv{\hbpar\grsmo}{\yF}{\hbvir\grso}
\ar[r]^-{\ychout}
&
{\hbpar\grsmo}
}\;\;
}\lrbl{-\hlf}\lrbs{-\hlf}.
\end{align}
(fractional shifts were introduced by H.~Wu\cx{Wu}).
Here we assume again that the differentials $\ychin$ and $\ychout$ are cohomological (that is,
they have the homological degree $1$) and that the right terms in
the complexes in parentheses have homological degree $0$.

\section{The isomorphism between the \sutt\ and \sof\
categorifications}
\subsection{The isomorphism theorem}
A simple relation between the groups
\ee
\label{xbe.1}
\msof=\msutt/\ZZ_2
\eee
implies a relation between the corresponding link polynomials: for
a framed oriented link $\xL$
\ee
\label{xbe.2}
\nPoLq = \lrbc{\nPuLq}^2 q^{-3\wrthL}.
\eee
Here $\nP_{\msof}$ is the $\msof$ specialization of the \Kfp,
$\nP_{\msut}$ is the Jones polynomial considered as
the $\msut$ specialization of the \homflyp\ and $\wrth(\xL)$ is
the writhe of $\xL$ defined as
\ee
\label{xbe.2a}
\wrthL = \sum_{i,j=1}^{\ncmpv\xL} l_{ij},
\eee
where $\ncmpv\xL$ is the number of components of $\xL$ and
$l_{ij}$ are the linking numbers between the components, the self-linking numbers
$l_{ii}$ being determined by their framing. If a link diagram
$\xL$ has the blackboard framing, then
\ee
\label{xbe.2a1}
\wrthL = \xnp(\xL)-\xnn(\xL),
\eee
where $\xnp$ and $\xnn$ are the numbers of positive ($\lscrp$) and
negative $(\lscrn)$ crossings of $\xL$. The factor
$q^{-\wrth(\xL)}$ in \ex{xbe.2} reflects the fact that the \sof\
polynomial, as defined in this paper, is covariant (\ex{1.0a1}) rather than invariant with respect
to the first Reidemeister move.

The relation\rx{xbe.2} can be
categorified.
\begin{theorem}
\label{th.suof}
There is a \hteq\ between the categorification complexes
\ee
\label{xbe.3}
\nCobL\cchq\lrbc{\nCubL\otimes\nCubL}\grsv{-3\wrthL}\gsztv{\wrthL}\grsev{-\wrthL)},
\eee
where $\nCubL$ is the \sut\ categorification complex of\cx{Kh} as constructed in\cx{KR1} and
\\
$\nCobL$ is the \sotn\ categorification complex\rx{1.3a} for
$\xN=1$.
\end{theorem}
\subsection{The \sut\ categorification for unoriented link diagrams}
In proving the relation\rx{xbe.3} it is convenient to use the
\sut\ complex construction which differs slightly from the
prescription of\cx{KR1}. Namely, we are going to construct an \sut\
complex $\tCubL$, which is related to the standard one by a degree
shift
\ee
\label{xbe.3a}
\tCubL =
\nCubL\grsv{-\frac{3}{2}}\gsztv{\hlf\,\wrthL}\grsev{-\hlf\,\wrthL}.
\eee
As \ex{xbe.3a} suggests, the complex $\tCubL$ is invariant under the second and third
Reidemeister moves and changes under the first Reidemeister move
according to the formula
\ee
\label{xbe.4}
\tCubs{\akinksk}  &\cchq
\tCubs{\avlinsk}\grsv{\frac{3}{2}}\gsztph\grseph.
\eee

The main feature of the combinatorial construction of $\tCubL$ is
that in contrast to the construction of\cx{KR1} it does not
require the orientation of link components and in this respect it
is similar to our construction of $\nCobL$.

Let us review the \mf\ construction of $\nCubL$ and introduce the
changes that will transform it to $\tCubL$.
%
The basic algebra of \sut\ categorification is $\IQ[\zu]$, its \qgrd\ defined by the condition $\dgq\zu=2$. The
basic polynomial is
\ee
\label{xbe.5}
\nWu(\zu) = \zu^3.
\eee
It has odd degree, hence $\nWu(-\zu)=-\nWu(\zu)$ and in the \sut\ case (as well as in the
\sun\ case with odd $\xN$) we can adopt the same leg orientation
convention, as in the \sotn\ case. Namely, we assume that all graph and \tngl\ legs
are oriented outwards, and if for the purpose of leg joining we
need to switch the orientation of an $i$-th leg, then we change its \mf\ by the functor
\rx{e2.d3a2} associated with the endomorphism of $\IQ[\zu_i]$, which
switches the sign of $\zu_i$.

We are going to review and slightly modify the categorification
formulas of Appendix\rw{ap.1} in view of our new leg orientation
convention. First of all, the \sut\ \oarc\ \mf\ is
\ee
\label{xbe.6}
\zarcxyoo
= \kmfv{\zut+\zuo}{\zut^2 - \zut\zuo + \zuo^2}.
\eee

Next we turn to four \flggd\ graphs $\lxpar$, $\lxvir$, $\lsvru$
and $\lxhor$, the latter \tarc\ graph being included, because the
new leg orientation convention permits it. The \mfs\ of \flggd\
graphs factor into \cmn\ and \prpr\ parts according to \ex{ae.16},
but this time we define both parts slightly differently. Namely,
we choose the polynomial $\zA$ not according to the
expression\rx{ae.15} at $\xN=2$, but rather as
\ee
\label{xbe.7}
\zA(\zbu) = \zuo^2 + \zut^2 + \zuth^2 -
\zuo\zut-\zuo\zuth-\zut\zuth - \zuf(\zuo+\zut+\zuth)-2\zuf^2.
\eee
The \prpr\ algebra this time is
\ee
\label{xbe.8}
\yRprp = \IQ[\zp,\zq,\zr],
\eee
and the \prpr\ homomorphism $\xymatrix@1{\yRprp\ar[r]^-{\zhomprp} &
\yR}$ is defined by the formulas
\ee
\label{xbe.9}
\zhomprp(\zp) = \zut+\zuf,\qquad\zhomprp(\zq)=\zut+\zuth,\qquad
\zhomprp(\zr)=\zuo+\zuf,
\eee
which are similar to \ex{ey.3a}. The \prpr\ \mfs\ of \flggd\
graphs are
\ee
\label{xbe.10}
\hbpruprp =\kmfv{\zp}{3\zq\zr},\qquad
\hbviuprp =\kmfv{\zr}{3\zp\zq},\qquad
\hbhruprp =\kmfv{\zq}{3\zp\zr}
\eee
and
\ee
\label{xbe.11}
\hbvruprp =\kmfv{\zp\zr}{3\zq}\grsmo.
\eee
The proof of the factorization formula\rx{ae.16} in
Proposition\rw{pr.suf} is repeated \emph{verbatim}.

The ordinary and \prpr\ \mfs\ of \tarc\ graphs
have the
same symmetry $\Sf$ and $\Sth$ as in the \sotn\ categorification
case (see subsection\rw{sss.endf}), so it suffices to define
the \sdlmrp\ between one pair of \tarc\ graphs and the morphisms
for other pairs will be dictated by the symmetry. Thus the
\sdlmrp\ is the tensor product\rx{ae.20a1}, where the \prpr\
\sdlmrp\ $\yFprp\in\Extmfoqo\lrbc{\hlparp,\hlvirp}$ is defined by
the diagram
\ee
\label{xbe.12}
\vcenter{\xymatrix@C=1.5cm{
{\hlparp} \ar
@<-0.55ex>
[d]^-{\yFprp}
&
{\yRomo}
\ar[r]^-{\zp}
\ar[dr]^-{1}
&
{\yRz}
\ar[r]^-{3\zq\zr}
\ar[dr]^-{-3\zq}
&
{\yRomo}
\\
{\hlvirp}
&
{\yRomo}
\ar[r]^-{\zr}
&
{\yRz}
\ar[r]^-{3\zp\zq}
&
{\yRomo}
}
}
\eee
This \sdlmrp\ is equal to the \sun\
\sdlmrp\ coming from\rx{ae.20}.

The application of the equivalence transformation $\rtov{1}{3}$
(see\rx{ea2.12z3x}) to the \mf\ $\hsvruprp$ turns it into
$\hshruprp\gszto$, hence
\ee
\label{xbe.13}
\hbvru\cong\hbhru\gszto,
\eee
and we can replace $\hsvru$ by $\hshru$ in the cone
relations\rx{ae.20a5} and in categorification formulas\rx{ae.d1}
and\rx{ae.d2}. The cone relations\rx{ae.20a5} become
\ee
\label{xbe.14}
\begin{split}
\hbhru\gszto\;&\cchq\;\; \xymatrix{\Cbdv{\hbpru\grsmo}{\yF}{\hbviu\grso}}
\\
&\cchq\;\;
\xymatrix{\Cbdv{\hbviu\grsmo}{\yG}{\hbpru\grso}}\;\;.
\end{split}
\eee
and the natural morphisms $\ychin$ and $\ychout$ relating the
cones to constituent \mfs\ are equal (up to a non-zero factor) to
\sdlmrps. Thus if we introduce a non-oriented \eltr\ crossing
complex
\ee
\label{xbe.15}
\hbcnw = \lrbc{
\xymatrix@C=1.5cm{
{\hbpru\grsph\gsztph}
\ar[r]^-{\yF}
&
{\hbhru\grsmh\gsztmh}
}
}{\textstyle \grsemh},
\eee
based on the \sdlmrp\ $\yF$ and define the complex $\tCubL$ for an
unoriented link diagram $\xL$ by joining the \eltr\ complexes\rx{xbe.15}
for each crossing of $\xL$.
The standard \sut\ categorification complex $\nCubL$ comes from
the crossing complexes\rx{ae.d1} and\rx{ae.d2}. These complexes
differ from the unoriented crossing complex\rx{xbe.15} by degree
shifts
\ee
\label{xbe.16}
\hbcrn = \hbcnw\grsv{\hlf-\xN}\gsztmh\grsemh,\qquad
\hbcrp = \hbcpw\grsv{\xN-\hlf}\gsztph\grseph,
\eee
and these shifts result in the relation\rx{xbe.3a} between the
link diagram complexes.

\subsection{Proof of the isomorphism of complexes}

It will be convenient to use the basic \sof\ polynomial, which
has an extra factor $2$ relative to the definition\rx{ea2.4}:
\ee
\label{xbe.17}
\nWoxy = 2(\zx\zy^2 + \zy^3),\qquad\dgq\zx=\dgq\zy=2.
\eee
Consider the algebra homomorphism
\ee
\label{xbe.18}
\xymatrix{\IQ[\zu,\zv]\ar[r]^-{\zxhom} &\IQ[\zx,\zy]},\qquad\zxhom(\zu)
= \zx+\zy,\quad\zxhom(\zv) = \zx-\zy.
\eee
It turns the sum of two \sut\ polynomials into the \sof\ polynomial:
\ee
\label{xbe.19}
\zxhom\lrbc{
\nWu(\zu) + \nWu(\zv)
} = \nWoxy.
\eee
Let $\xg$ and $\xgp$ be two $n$-legged \opgrs, appearing in the unoriented version of the
\sut\ categorification (that is, $\xg$ and $\xgp$ are disjoint
unions of arcs and circles). Suppose that $\xg$ and $\xgp$ have
the same number of legs. We index them and assign variables
$\zu_i$ ($1\leq i\leq n$) to the legs of $\xg$ and $\zv_i$ to the legs of $\xgp$.
Denote by $\zxhom_i$ the homomorphism\rx{xbe.18} applied to
$\IQ[\zu_i,\zv_i]$ and let $\zhhom$ be the composition of all
functors $\zhhom_i$. Then $\zhhom(\gg\otimes\ggp)$ is a \mf\ of
the sum of \sof\ polynomials
\ee
\label{xbe.20}
\sum_{i=1}^n \nWo(\zx_i,\zy_i).
\eee

Introduce a shortcut notation for the \sut\ \mfs:
\ee
\label{xbe.20a}
\ggud = \ggu\otimes\ggu.
\eee
The following lemma leads to a quick proof of Theorem\rw{th.suof}.
\begin{lemma}
\label{lm.a2.1}
The functor $\zhhom$ maps a tensor product of \sut\
crossing \mfs\ into the \sof\ crossing \mf:
\ee
\label{xbe.21}
\zhhom\lrbc{\hlncrud} \cchq \hlncro.
\eee
\end{lemma}
\pr{Theorem}{th.suof}
Since the homomorphism\rx{xbe.18} is an algebra isomorphism,
the relation\rx{xbe.21} implies the \hteq\ of complexes
\ee
\label{xbe.22}
\tCubL\otimes\tCubL \cchq \nCobL.
\eee
The \hteq\rx{xbe.3} follows in view of the
relation\rx{xbe.3a}.\qed

\pr{Lemma}{lm.a2.1}
The functor $\zhhom$ maps \oarc\ graphs into \oarc\ graphs:
\ee
\label{xbe.23}
\zhhom\lrbc{\zarcxyooud} \cong \zarcxyooo.
\eee
Indeed, the \lhs of this equation is a \Kmf\
\ee
\label{xbe.24}
\xymatrixnocompile@C=3cm{
{\kmftwv{\zxo+\zyo+\zxt+\zyt}{\zxo-\zyo+\zxt-\zyt}{\ast}{\ast}}
\ar@{|->}[r]^-{\rtov{1}{2}{1}\rtov{2}{-1}\rtov{1}{1/2}\rtrv{2}{1}{1}}&
{\kmftwv{\zxo+\zxt}{\zyo+\zyt}{\ast}{\ast}}
}
\eee
The right \mf\ in this diagram is isomorphic to the \rhs of
\ex{xbe.23} in view of Theorem\rw{th.equiv}.

The functor $\zhhom$ also maps \sdlmrps\ to \sdlmrps:
\def\xardot{ \ar@{}[d]|{\otimes} }
\ee
\label{xbe.25}
\zhhom\lrbc{
\vcenter{\xymatrix@C=1.5cm{
{\hlparud}
\ar[r]^-{\yFsud}
&
{\hlvirud}
}
}
}
=
\lrbc{
\vcenter{\xymatrix@C=1.5cm{
{\hlparo}
\ar[r]^-{\yFso}
&
{\hlviro}
}
}
}
\eee
This follows from the fact that the \rhs and \lhs morphisms have the same
\qdgr\ and from the uniqueness of the \sdlmrp\ (up to a constant
factor).

Since the functor $\zhhom$ maps \oarc\ \mfs\ to \oarc\ \mfs\ and
\sdlmrps\ to \sdlmrps, then the definition\rx{ey.26} of the \sof\ crossing complex $\hlncro$
indicates that it has an \sutt\ counterpart:
%
\ee
\label{xbe.26}
\zhhom\lrbc{
\xymatrix{
{\hlhorud\grso} \ar[r]^-{
\ychin
}
&
{\hlverudd}
\ar[r]^-{
\ychout
}
&
{\hlparud\grsmo}
}
}\gszto=\hlncro,
\eee
where we used a shortcut notation
\ee
\label{xbe.27}
\hlverudd = \xymatrix@C=1.5cm{
\Cbtwshuv{\hlparud}{\yFsud}{\hlvirud}{\yGsud}{\hlhorud}{\yY}
}
\eee
and $\yY$ is a \scndh\ (its choice does not impact the
homotopy class of the \cnvl, because the condition\rx{ea2.24a9b}
holds). Hence the \hteq\rx{xbe.21} follows from the next
lemma.\qed

From now on we will be dealing only with \mfs\ appearing in the
\sut\ categorification,
hence we drop the indices \sut\ at the graphs and morphisms.
\begin{lemma}
\label{lm.a2.2}
The following complexes of \mfs\ are \hteqt:
\ee
\label{xbe.28}
\hlncrsd\gszto \cchq
\lrbc{
\xymatrix{
{\hlhorsd\grso} \ar[r]^-{
\ychin
}
&
{\hlverudd}
\ar[r]^-{
\ychout
}
&
{\hlparsd\grsmo}
}
}
\eee
\end{lemma}

\proof
According to \ex{xbe.15} (rotated by $\ntydeg$), the \lhs of
\ex{xbe.28} can be presented as a complex
\ee
\label{xbe.29}
\xymatrix{
{\hlhorsd\grso}
\ar[r]^-{\ytchin}
&
{\lrbc{\hlhor\otimes\hlpar}\gszto
\oplus
\lrbc{\hlpar\otimes\hlhor}\gszto}
\ar[r]^-{\ytchout}
&
{\hlparsd\grsmo}
},
\eee
where
\ee
\label{xbe.30}
\ytchin =
\begin{pmatrix}
\yF\otimes\id \\ -\id\otimes\yF
\end{pmatrix}
,
\qquad
\ytchout =
\begin{pmatrix}
\id\otimes\yF &\yF\otimes\id
\end{pmatrix}
\eee
and the middle \mf\ carries the zero homological degree.
Hence we will establish the relation\rx{xbe.28} by proving the
\hteq\ between the central \mfs\ of\rx{xbe.28} and\rx{xbe.29}
\ee
\label{xbe.31}
\hlverudd\cchq\lrbc{\hlhor\otimes\hlpar}\gszto
\oplus
\lrbc{\hlpar\otimes\hlhor}\gszto
\eee
and showing that the homomorphisms $\ychin$ and $\ychout$ are
equivalent to $\ytchin$ and $\ytchout$ respectively.

Let us replace the virtual crossing \mfs\ $\hlvir$ in the
formula\rx{xbe.27} for $\hlverudd$ by the cone expression
\ee
\label{xbe.32}
\hlvir \cchq
\xymatrix{
{\hlhor\grsmo}
\ar[r]^-{\yF}
&
{\hlpar\grso}
\save
[].[l]*[F]\frm{}
\restore
}\gszto
\eee
obtained by
permuting legs 2 and 4 in the second equality of\rx{xbe.14}. Since
this substitution turns the \sdlmrps\ of\rx{xbe.27} to natural
morphisms relating this cone to its \cnst\ \mfs,
the \cnvl\ in the \rhs of \ex{xbe.27} becomes
\begin{multline}
\label{xbe.33}
\hlverudd\cchq
\\
\vcenter{\xymatrix@C=1.5cm@R=0.2cm{
&
{\hlhorsd\grsv{-2}\gszto}
\ar@/^5pc/[rrd]^-{\id}="u"
\ar[r]^-{\yF\otimes\id}
\ar[dd]^-{\id\otimes\yF}
&
{\hlpar\otimes\hlhor\gszto}
\ar[dd]^-{-\id\otimes\yF}
\\
{\hlparsd\grsmo}
\ar@/_5pc/[rrd]_-{\id}="d"
&
&
&
{\hlhorsd\grso}
\\
&
{\hlhor\otimes\hlpar\gszto}
\ar[r]^-{\yF\otimes\id}
&
{\hlparsd\grsv{2}\gszto}
\save
"2,1"."2,4"."u"."d"*[F]\frm{}
\restore
}
}
\end{multline}
The square in the center of this diagram represents the \mf\
$\hlvirsd$ after the substitution\rx{xbe.32}. Note that this
substitution allowed us to set the \scndh\ equal to zero. Let us
rearrange the convolution\rx{xbe.33} by assembling some \cnst\
\mfs\ into `subcones':
\begin{multline}
\label{xbe.34}
\hlverudd\cchq
\\
\!\!\!\!\!\!\!\!\!
\vcenter{
\xymatrix@C=-0.1cm{
&&&&&
{\hlhor\otimes\hlpar\gszto}
\\
{\hlhorsd\grsv{-2}\gszto}\ar[rrrr]^-{\id}
&&&&
{\hlhorsd\grso}
\save
{"2,1"."2,5"}="t"*[F]\frm{},
\ar"t";"1,6"^-{\smpm{-\id\otimes\yF & 0}}
\ar"t";"3,6"_-{\smpm{\yF\otimes\id & 0}}
\restore
&&
{\hlparsd\grsmo}\ar[rrrr]^-{\id}
&&&&
{\hlparsd\grsv{2}\gszto}
\save
{[].[llll]}="b"*[F]\frm{},
\ar"1,6";"b"^-{\smpm{0\\\yF\otimes\id}}
\ar"3,6";"b"_-{\smpm{0\\\id\otimes\yF}}
\restore
\\
&&&&&
{\hlpar\otimes\hlhor\gszto}
\save
"2,1"+<-1.8cm,0cm>*{}="xa",
"2,11"+<1.6cm,0cm>*{}="xb",
"1,6"."3,6"."2,1"."2,11"."xa"."xb"*[F]\frm{}
\restore
}
}
\end{multline}
The cones of identity morphism are contractible. According to Lemma\rw{l.excc},
their excision establishes the \hteq\rx{xbe.31}. It remains to
establish the equivalence of morphisms $\ychin,\ytchin$ and
$\ychout,\ytchout$. We will prove the equivalence for the first
pair, since the second pair can be treated similarly. The morphism
$\ychin$ maps $\hlparsd$ identically to the second \mf\ in the
left subcone of the \cnvl\rx{xbe.34}. Hence its equivalence
to $\ytchin$ follows from the next Lemma, in which
$\xaA$ stands for $\hlparsd$ and $\xaB$
stands for the rest of the
\cnvl\rx{xbe.34}.\qed

\begin{lemma}
For a morphism between two \mfs
\ee
\label{xbe.35}
\xymatrixnocompile{\xaA\ar[r]^-{\yf} &\xaB}
\eee
consider another morphism between $\xaA$ and a cone
\ee
\label{xbe.36}
\xymatrixnocompile@C=2cm{\xaA & \xaA\gszto \ar[r]^-{\smpm{\id&\yf}}="a" & \xaA\oplus\xaB
\save
{[].[l]."a"}="b"*[F]\frm{},
\ar"1,1";"b"^-{\smpm{0&\id&0}}
\restore
}\;\;.
\eee
The target of this morphism is \hteq\ to $\xaB$ and the morphism
itself is equivalent to $\yf$.
\end{lemma}
\proof
Consider a commutative diagram
\ee
\label{xbe.37}
\vcenter{\xymatrix@R=3cm{{\xaA} && {\xaA\gszto} \ar[rr]^-{\smpm{\id&\yf}}="a" &{*i{\hlpar}}
& {\xaA\oplus\xaB}
\save
{[].[ll]."a"}="b"*[F]\frm{},
\ar"1,1";"b"^-{\smpm{0&\id&0}}
\restore
\\
&&
{\xaA\gszto}\ar[r]^-{\id}="c"
&
{\xaA}
\ar@{}[r]|{\oplus}
\save
{[].[l]."c"}="d"*[F]\frm{}
\restore
&
{\xaB}
\save
"2,3"+<-0.8cm,0cm>*{}="s1",
{"2,3"."d"}!C="e",{"e"."s1"}*[F(]\frm{(},
{"2,5"."d"}!C="e"*[F)]\frm{)}
\ar"1,1";"e"_-{\smpm{0&-\id&\yf}}
\ar"1,4";"e"^-{-\smpm{\id & 0 & 0\\0 & \id & 0\\0&-\yf&\id}}
\restore
}
}\;\;.
\eee
The vertical arrow establishes an isomorphism between the upper-right cone
and a sum of a contractible cone and $\xaB$ in the lower-right corner. After we remove the
contractible cone, the morphism between the upper-left and lower-right
\mfs\ will reduce to $\yf$.\qed

\end{document}